\documentclass[11pt]{article}
\usepackage[margin=1in]{geometry}
\usepackage{amsmath,amssymb,amsthm,mathtools}
\usepackage{graphicx,booktabs,longtable,array,multirow}
\usepackage[dvipsnames]{xcolor}
\usepackage{float,algpseudocode}
\floatstyle{ruled}
\newfloat{algorithm}{tbp}{loa}
\floatname{algorithm}{Algorithm}
\usepackage[numbers,sort&compress]{natbib}
\usepackage[colorlinks=true,linkcolor=blue!50!black,citecolor=blue!50!black,urlcolor=blue!50!black,hyperfootnotes=false]{hyperref}
\usepackage{doi}
\usepackage[nameinlink,capitalize]{cleveref}
\usepackage{pgfplots}
\pgfplotsset{compat=1.18}
\usetikzlibrary{patterns}
\makeatletter
\g@addto@macro\UrlBreaks{\do\-}
\makeatother
\usepackage{microtype}
\microtypesetup{expansion=false}
\usepackage{enumitem}
\graphicspath{{./}}
\usepackage[T1]{fontenc}
\usepackage{lmodern}
\hypersetup{hidelinks,hypertexnames=false}
\usepackage[bottom]{footmisc}
\AtBeginDocument{%
  \setlength{\abovedisplayskip}{7pt plus 2pt minus 3pt}%
  \setlength{\belowdisplayskip}{7pt plus 2pt minus 3pt}%
  \setlength{\abovedisplayshortskip}{4pt plus 2pt minus 2pt}%
  \setlength{\belowdisplayshortskip}{4pt plus 2pt minus 2pt}%
}

\usepackage[compact]{titlesec}

\titlespacing*{\section}
  {0pt}{1.55ex plus 0.35ex minus 0.25ex}{0.75ex plus 0.15ex}

\titlespacing*{\subsection}
  {0pt}{1.25ex plus 0.30ex minus 0.20ex}{0.55ex plus 0.12ex}

\titlespacing*{\subsubsection}
  {0pt}{1.00ex plus 0.25ex minus 0.15ex}{0.45ex plus 0.10ex}

\titleformat{\section}
  {\normalfont\large\bfseries}{\thesection.}{0.5em}{}

\titleformat{\subsection}
  {\normalfont\normalsize\bfseries}{\thesubsection.}{0.5em}{}

\titleformat{\subsubsection}
  {\normalfont\normalsize\itshape}{\thesubsubsection.}{0.5em}{}

\setlist{
  topsep=2pt plus 1pt minus 1pt,
  itemsep=1pt plus 0.5pt,
  parsep=0pt,
  partopsep=0pt,
  leftmargin=*
}

\makeatletter
\def\thm@space@setup{%
  \thm@preskip=6pt plus 2pt minus 2pt
  \thm@postskip=6pt plus 2pt minus 2pt
}
\makeatother

\makeatletter
\renewcommand{\maketitle}{%
  \begin{center}
    {\Large\bfseries \@title\par}
    \vskip 0.45em
    {\normalsize
      \lineskip .5em
      \begin{tabular}[t]{c}
        \@author
      \end{tabular}\par}
    \vskip 0.35em
    {\small \@date\par}
  \end{center}
  \par\vskip 0.35em
  \@thanks
  \vskip 0.7em
}
\makeatother

\renewenvironment{abstract}
  {\small\begin{center}\bfseries Abstract\end{center}\vspace{-0.6em}\small}
  {\par\vspace{0.6em}}

\newtheorem{theorem}{Theorem}[section]
\newtheorem{lemma}[theorem]{Lemma}
\newtheorem{proposition}[theorem]{Proposition}
\newtheorem{corollary}[theorem]{Corollary}
\theoremstyle{definition}
\newtheorem{definition}[theorem]{Definition}

\theoremstyle{remark}
\newtheorem{remark}[theorem]{Remark}

\newcommand{\R}{\mathbb R}
\newcommand{\C}{\mathbb C}
\newcommand{\N}{\mathbb N}

\newcommand{\cF}{\mathcal F}
\newcommand{\cM}{\mathcal M}
\newcommand{\cS}{\mathcal N}\newcommand{\Szero}{\mathcal S_0}
\newcommand{\cT}{\mathcal T}
\newcommand{\cA}{\mathcal A}
\newcommand{\cB}{\mathcal B}
\newcommand{\sfP}{\mathsf P}
\newcommand{\ch}{\cosh}
\newcommand{\sh}{\sinh}
\newcommand{\arsinh}{\operatorname{arsinh}}
\newcommand{\Rea}{\operatorname{Re}}
\newcommand{\Ima}{\operatorname{Im}}
\newcommand{\ip}[2]{\langle #1,#2\rangle}
\newcommand{\norm}[1]{\lVert #1\rVert}
\newcommand{\rhoF}{\rho_{\cF}}
\newcommand{\rhoz}{\rho_{\cF}^{0}}
\newcommand{\cQ}{c_{Q}}
\newcommand{\cFF}{c_{\cF}}
\newcommand{\gt}{\tilde g}
\newcommand{\st}{\tilde s}
\newcommand{\clag}{c_{\rm lag}}\newcommand{\clayer}{c_{\rm layer}}\newcommand{\cdisc}{c_{\rm disc}}
\newcommand{\cpi}{c_{\pi}}
\newcommand{\ccorner}{c_{\angle}}
\newcommand{\teff}{\tau_{\rm eff}}
\newcommand{\cFm}{\cF_m}

\crefname{equation}{}{}
\Crefname{equation}{}{}
\crefname{figure}{Figure}{Figures}

\title{Exact (Non-)Accelerations of Heavy-Ball Discretizations}
\author{Dmitry Pasechnyuk-Vilensky\thanks{Mohamed bin Zayed University of Artificial Intelligence, Abu Dhabi, UAE. \texttt{dmivilensky1@gmail.com}}, Martin Tak\'{a}\v{c}\thanks{Mohamed bin Zayed University of Artificial Intelligence, Abu Dhabi, UAE. \texttt{Takac.MT@gmail.com}}}
\date{}
\begin{document}
\maketitle

\begin{abstract}
Momentum methods with the same slow dynamics can have different worst-case
rates on smooth strongly convex functions.  We identify a half-step phase
barrier in the gradient query that governs this difference.
For momentum with extrapolation \(\tau_\kappa\to\tau\),
\(\beta=1-\theta/\sqrt\kappa+o(\kappa^{-1/2})\) and \(hL\to y\), where
\(\kappa=L/\mu\) and \(\theta,y>0\), the exact uniform distance factor is
\(1-c/\sqrt\kappa+o(\kappa^{-1/2})\), with
\(c=\min\{\theta/3,c_Q(\theta,y)\}\), when
\(\tau>1/2\) and \(y(1+2\tau)<4\); \(c_Q\) is the slowest quadratic decay exponent.
Quadratically accelerated families with \(\limsup\tau_\kappa<1/2\) have
divergent trajectories on the full function class.
Upper bounds use finite-lag certificates: nonnegative combinations of
inequalities for values and gradients at finitely separated query points,
together with comparisons to the minimizer.
For any prescribed accuracy, spiral mixtures or quadratic trajectories of one
finite-dimensional objective give matching lower bounds at all iteration horizons.
We resolve the parameter layer near \(\tau=1/2\) and classify all tunings of this family attaining
the optimal leading constant \(1\).  The flow also bounds the leading constant
of every second-order stationary linear family with the same slow parameters
\((\theta,y)\) and fixed finite memory.
\end{abstract}

\section{Introduction}
\label{sec:intro}

Momentum methods can have the same continuous-time limit and different worst-case
behavior.  Polyak's heavy ball is optimal on strongly convex quadratics~\cite{Polyak1964}.
On smooth strongly convex functions it can fail to converge~\cite{LRP16,GTD}.
Nesterov's method accelerates on this whole class~\cite{Nesterov2018}, and the
triple momentum method (TMM) attains the optimal asymptotic distance factor
\(1-\kappa^{-1/2}\), where \(\kappa=L/\mu\)~\cite{VFL18,DT22}.
These methods belong to the family \(\cM\),
\begin{equation}\label{eq:M}
\tilde x_k=x_k+\tau(x_k-x_{k-1}),\qquad
x_{k+1}=x_k+\beta(x_k-x_{k-1})-h\nabla f(\tilde x_k):
\end{equation}
heavy ball has \(\tau=0\), Nesterov's method has \(\tau=\beta\), and TMM has
\(\tau\to1/2\) as \(\kappa\to\infty\).
We determine how extrapolation controls acceleration on the full class.

\paragraph{Slow dynamics and the half-step rule.}
Accelerated tunings have \(1-\beta\asymp\kappa^{-1/2}\).  For families with limits
\begin{equation}\label{eq:scaling-intro}
\beta=1-\theta\kappa^{-1/2}+o(\kappa^{-1/2}),\qquad hL\to y,
\end{equation}
we measure performance by the constant \(c\) in the uniform distance factor
\(\rhoF=1-c\kappa^{-1/2}+o(\kappa^{-1/2})\).
For squared distances the leading constant is \(2c\).
The slow equation is the heavy-ball flow
\(\ddot x+\theta\dot x+\nabla F(x)=0\) with a \(y\)-strongly convex
potential.  We call \((\theta,y)\) its \emph{slow data}: damping \(\theta\) and
strong convexity \(y\) of the scaled flow.  Over all finite smoothness bounds,
its worst distance decay exponent is
\[
\cFF(\theta,y)=\min\{\theta/3,\cQ(\theta,y)\},
\]
where \(\cQ\) is the slowest quadratic exponent.
This exponent bounds the leading constant of every second-order stationary
linear family with the same slow data and fixed finite memory
(\Cref{def:second,thm:ceiling}).  The family \(\cM\) attains it when
\(\tau>1/2\) and \(y(1+2\tau)<4\) (\Cref{thm:B}); acceleration on the
full class requires \(\liminf_{\kappa\to\infty}\tau_\kappa\ge1/2\)
(\Cref{thm:wall}).  We call \(\tau=1/2\) the \emph{half-step wall}.

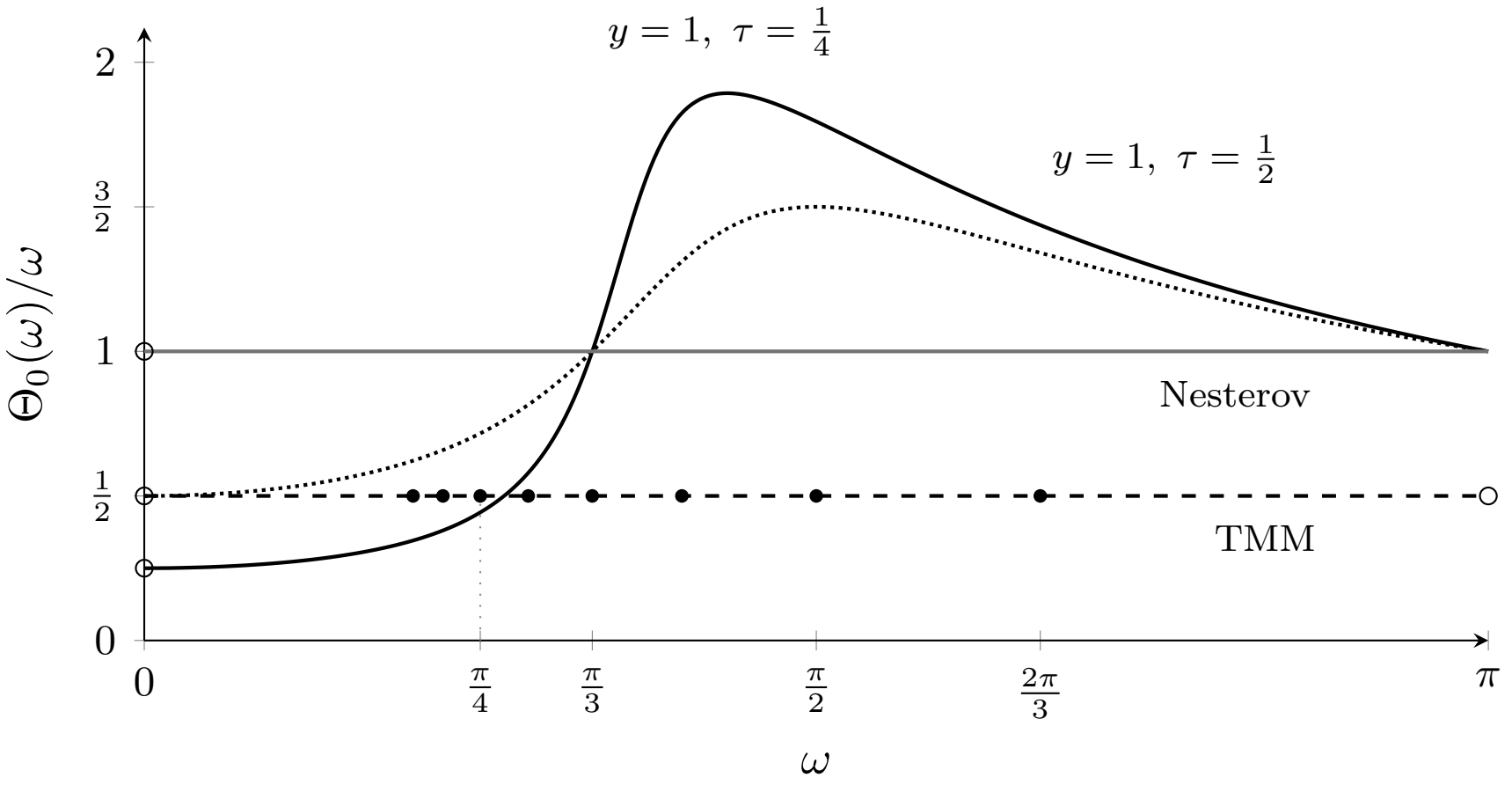
\begin{figure}[t]
\centering
\begin{tikzpicture}
\pgfmathsetmacro{\phasePi}{pi}
\begin{axis}[width=0.68\textwidth,height=0.36\textwidth,trig format plots=rad,
  xmin=0,xmax=3.1415927,ymin=0,ymax=2.12,
  xlabel={\(\omega\)},ylabel={\(\Theta_0(\omega)/\omega\)},
  xtick={0,0.7853982,1.0471976,1.5707963,2.0943951,3.1415927},
  xticklabels={\(0\),\(\frac\pi4\),\(\frac\pi3\),\(\frac\pi2\),\(\frac{2\pi}3\),\(\pi\)},
  ytick={0,0.5,1,1.5,2},yticklabels={\(0\),\(\frac12\),\(1\),\(\frac32\),\(2\)},
  tick label style={font=\footnotesize},label style={font=\small},
  axis lines=left,clip=false]
\addplot[thick,dashed,domain=0:3.1415927,samples=2] {0.5};
\addplot[thick,domain=0:3.1415,samples=320]
  {x==0 ? 0.25 : atan2(0.25*sin(x),1-1.75*(1-cos(x)))/(x==0 ? 1 : x)};
\addplot[thick,densely dotted,domain=0:3.1415,samples=320]
  {x==0 ? 0.5 : atan2(0.5*sin(x),1-1.5*(1-cos(x)))/(x==0 ? 1 : x)};
\addplot[thick,black!55,domain=0:3.1415927,samples=2] {1};
\addplot[only marks,mark=*,mark size=1.25pt,domain=3:10,samples=8]
  ({2*pi/x},{0.5});
\addplot[thin,dotted,black!45] coordinates {(0.7853982,0) (0.7853982,0.5)};
\addplot[only marks,mark=o,mark size=1.8pt,mark options={solid,fill=white}]
  coordinates {(0,0.25) (0,0.5) (0,1) (\phasePi,0.5)};
\node[font=\scriptsize,anchor=south] at (axis cs:1.35,1.96) {\(y=1,\ \tau=\frac14\)};
\node[font=\scriptsize,anchor=south west] at (axis cs:2.08,1.52) {\(y=1,\ \tau=\frac12\)};
\node[font=\scriptsize,anchor=north] at (axis cs:2.55,0.96) {Nesterov};
\node[font=\scriptsize,anchor=north] at (axis cs:2.62,0.46) {TMM};
\end{axis}
\end{tikzpicture}
\caption{Spectral phase lead of limit methods of \(\cM\), in steps.
The horizontal half-step line is TMM; the lag-one symbol is positive above it.
Dots at \(\omega=2\pi/K\), \(3\le K\le10\), mark the one-atom thresholds.
For the solid curve, strict passage below a threshold occurs exactly for
\(K\ge8\), starting at \(\omega=\pi/4\).  Open markers denote endpoint limits.}
\label{fig:halfstep}
\end{figure}

The half step has a spectral explanation.  On a Fourier mode with \(\omega\ne0\),
\[
1-e^{i\omega}
=\frac{2\sin(\omega/2)}{\omega}\,e^{i\omega/2}(-i\omega).
\]
Replacing a derivative by a difference adds half a step of phase.
For limit methods of \(\cM\), the low-frequency lead that compensates for this
shift is \(\tau\).
For larger memory, the full frequency-dependent phase enters the rate
(\Cref{sec:faithful,fig:halfstep}).  Upper bounds use finite-lag certificates:
nonnegative combinations of inequalities for values and gradients at finitely
separated queries and comparisons to the minimizer (\Cref{def:cert}).
We match them with trajectories generated by one objective.

\paragraph{Contributions.}
We extend the spectral characterizations for heavy ball~\cite{TA} and standard
Nesterov~\cite{PS26} to general stationary linear methods.
Finite-dimensional \(C^{1,1}\) witnesses establish distance sharpness for the flow
estimates of~\cite{ADR22} over the union of finite smoothness classes.
A trace identity complements certificate-based method synthesis~\cite{LS20}: quadratic
stability and a nonnegative limiting lag-one symbol force \(y\le2\), with
equality precisely for TMM's limiting transfer function.

The standard Nesterov tuning has \(c=2/3\); its quadratic-optimal tuning has
full-class constant \(4/(3\sqrt3)\).  The best full-class constant within the
Nesterov family is \(\sqrt{2/3}\), attained by
\(h=4/(3L)\), \(\tau=\beta=1-\sqrt{6/\kappa}\)
(\Cref{cor:B,cor:nag-asymptotic-optimal}).
The optimum \(c=1\) over \(\cM\) requires
\(\bigl(\sqrt\kappa(1-\beta),hL,\tau\bigr)\to(3,2,1/2)\).
We call \((hL,\tau)=(2,1/2)\) the \emph{TMM corner}.  In its natural scale,
\(1/2-\tau=a\kappa^{-1/2}+o(\kappa^{-1/2})\),
\(2-hL=b\kappa^{-1/2}+o(\kappa^{-1/2})\), it is attained exactly when
\(\theta=3\), \(a\le5/4\) and \(2a+b\ge7/2\): a wedge with TMM at its
vertex (\Cref{fig:corner}).  We also classify approaches at arbitrary rates.
For some suboptimal corner tunings, integer rotation periods in the lower-bound construction
leave a narrow arithmetic gap between the available bounds
(\Cref{thm:corner,sec:conclusions}).

\paragraph{Related work: flows and exact rates.}
Continuous-time momentum models go back to Polyak~\cite{Polyak1964};
for Nesterov's convex method see~\cite{SBC16}.  High-resolution equations,
Hessian-driven damping and semi-implicit discretizations expose the gradient
corrections separating the methods~\cite{SDJS22,ACFR22,MJ19}; Lyapunov analyses
connect these flows to accelerated schemes~\cite{SZ21,WRJ21}.
For heavy-ball flow, Aujol, Dossal and Rondepierre proved function-value
bounds on quasi-strongly convex functions and sharpness using non-Lipschitz
gradients~\cite{ADR22}; see~\cite{ADR23} for the \L ojasiewicz setting.
\Cref{thm:flow} gives the sharp distance exponent using \(C^{1,1}\) witnesses.
For Nesterov's standard tuning, explicit integral quadratic constraint (IQC)
bounds appear in~\cite{SJFB18};
Gannot obtained the constant \(2/3\) using the Jury--Lee multiplier~\cite{Gan22},
and Pokutta and Stellato proved its sharpness by mixtures of planar
spirals~\cite{PS26}.  Our formula recovers it and identifies the best
asymptotic tuning within the Nesterov family.
For heavy ball, planar cycles establish non-acceleration~\cite{GTD,GDT23}.
The one-dimensional question of~\cite{GTD25} was settled in~\cite{FHDA26};
local acceleration with controlled initialization and a gradient-descent
warm-up is studied in~\cite{LWS26}.  Spectral realizations characterize the
exact rate, separation and one-exponent principle in~\cite{TA}.
\Cref{thm:A} extends this framework to general stationary linear methods;
\Cref{lem:finite,lem:strict,lem:count} give finite-support and atom-count
results under explicit nondegeneracy conditions.
Hernanz~\cite{Her26} identifies the optimal nonaccelerated asymptotic
constant \(8\) in dimension at least two.

\paragraph{Related work: certificates, synthesis and optimal methods.}
The family~\eqref{eq:M} includes robust momentum, which interpolates between TMM
and gradient descent~\cite{CHVL18}, and is a standard synthesis template
\cite{LRP16,SE21,MSE21}.
Rates can be certified by IQCs and dissipativity~\cite{MR97,LRP16,HL17},
Popov and Jury--Lee multipliers~\cite{Pop62,JL64}, and O'Shea--Zames--Falb
multipliers~\cite{OY67,ZF68,WB68,CTH16,ZCH22,SSCK23}.
Exactness in dimension-free lifted settings and limitations for a single
plant are studied in~\cite{KS24,Kha26}; periodic counterexamples and
nonlinear rates beyond linearizations appear in~\cite{SC21,ZSC22}.
Performance estimation and automated Lyapunov analysis compute tight bounds
at fixed parameters~\cite{DT14,THG17,TVSL18,UBTG25}; lifting and interpolation
connect these descriptions~\cite{VSL22}.
Lessard and Seiler derive memory-independent synthesis bounds for a fixed IQC,
with a certificate using consecutive-query data recovering TMM's rate~\cite{LS20}.  Our trace identity
adds an explicit step formula and rigidity of its equality case.
TMM and the information-theoretic exact method attain the optimal asymptotic
factor~\cite{VFL18,TD23,DT22}.
Individual asymptotic root factors under additional local regularity are
studied in~\cite{VSL25,CX26}.  We use the uniform factor, with the supremum
over functions taken before the asymptotic limit (\Cref{rem:uniform-smoothness}).
For broader context on acceleration see~\cite{dST21}.

\paragraph{Organization.}
\Cref{sec:prelim} states the main rate for \(\cM\).
\Cref{sec:flow,sec:closure} develop the flow witnesses and spectral closure;
\Cref{sec:beyond} gives the ceiling, half-step rule and trace formula for
general memory.  \Cref{sec:constant,sec:TMM} resolve the wall layer and the
TMM corner, and \Cref{sec:finite} treats finite condition numbers.
Proofs and detailed algebra are collected in the appendices.

\paragraph{Formal verification.}
The mathematical results of this article have formal proofs in Lean. The package is available at
\href{https://github.com/dmivilensky/PolynomialRealizationsProgram/tree/main/nesterov}{\texttt{PolynomialRealizationsProgram/}\allowbreak\texttt{nesterov}}\footnote{\url{https://github.com/dmivilensky/PolynomialRealizationsProgram/tree/main/nesterov}}.

\section{The Momentum Family and Its Main Rate}
\label{sec:prelim}

\paragraph{Notation and scope.}
Let \(H\) be a real Hilbert space and let \(\cF_{\mu,L}(H)\), \(0<\mu<L\),
consist of differentiable \(\mu\)-strongly convex functions with
\(L\)-Lipschitz gradient; thus \emph{smooth} means \(C^{1,1}\).  We write \(x^\star\) for the minimizer,
\(\kappa=L/\mu\), and \(\N=\{1,2,\ldots\}\).
We identify a real plane with \(\C\), using \(\ip ab=\Rea(a\bar b)\).
Frequencies are \(\omega\in[0,\pi]\), and \(\xi=1-\cos\omega\in[0,2]\).
Bars over complex quantities mean conjugation; bars over real spectral
quantities mean averages with respect to the specified probability measure.

\subsection{Quadratic stability and the rate convention}

For~\eqref{eq:M} set \(y_\mu=h\mu\), \(y=hL\),
\(N_\tau(z)=(1+\tau)z-\tau\), and
\[
\chi_\lambda(z)=(z-1)(z-\beta)+h\lambda N_\tau(z).
\]
This is the characteristic polynomial on \(\lambda\norm x^2/2\).
A \emph{tuning} is a triple \((\beta,\tau,h)\), with
\(\beta\in[0,1)\), \(\tau\ge0\), \(h>0\).
A tuning is \emph{stable} when
\(\rho_q:=\max_{\lambda\in[\mu,L]}\max\{|z|:\chi_\lambda(z)=0\}<1\),
so stability refers to quadratics.

\begin{lemma}\label{lem:qstab}
A tuning of \(\cM\) is stable if and only if \(y(1+2\tau)<2(1+\beta)\).
\end{lemma}

For \(\cM\) the state is \(S_k=(x_{k-1}-x^\star,x_k-x^\star)\).
The following convention also applies to the states of stationary methods
defined in \Cref{def:slm}.
\begin{definition}\label{def:rate}
For \(n\ge0\), let \(\Gamma(n)=\sup\norm{S_{k+n}}/\norm{S_k}\), where the
supremum ranges over real Hilbert spaces \(H\), functions \(f\in\cF_{\mu,L}(H)\),
trajectories and indices \(k\ge0\) with \(S_k\ne0\).  The \emph{worst-case rate} is
\(\rhoF=\lim_{n\to\infty}\Gamma(n)^{1/n}\).  For \(\cM\) we also use the
\emph{rate from rest},
\[
\rhoz=\limsup_{N\to\infty}
\left(\sup\frac{\norm{x_N-x^\star}}{\norm{x_0-x^\star}}\right)^{1/N},
\]
where the supremum ranges over \(d\ge1\), \(f\in\cF_{\mu,L}(\R^d)\) and
\(x_0\ne x^\star\), with \(x_{-1}=x_0\).
\end{definition}

\begin{lemma}\label{lem:envelope}
The envelope \(\Gamma\) is finite and submultiplicative, and
\[
\rhoF=\lim_{n\to\infty}\Gamma(n)^{1/n}
=\inf_{n\ge1}\Gamma(n)^{1/n}\ge\rho_q.
\]
For \(\cM\), \(\rhoz\le\rhoF\).
\end{lemma}

\begin{remark}[Uniform rates and smooth subclasses]\label{rem:uniform-smoothness}
The supremum in \(\Gamma(n)\) is taken \emph{before} the limit in \(n\).
For any fixed stationary linear method and finite \(n\), restricting that supremum to
finite-dimensional \(C^\infty\) objectives in \(\cF_{\mu,L}\) preserves
its value: smooth convolution preserves \(\mu,L\) and approximates every
finite trajectory (\Cref{app:prelim}).
Individual asymptotic root factors, as in~\cite{VSL25,CX26}, place the
limit inside the supremum and describe a different rate convention.  Distance and squared-distance factors have respective
first-order constants \(c\) and \(2c\).
\end{remark}

\subsection{Accelerated scaling and the main theorem}
\label{sec:quadratic}
In the accelerated scaling we normalize \(L=1\) (replace \(f\) by
\(f/L\) and \(h\) by \(hL\)), write \(\mu=\varepsilon^2=1/\kappa\), and
consider families
\begin{equation}\label{eq:scaling}
\beta=1-\theta_\varepsilon\varepsilon,\quad h=y_\varepsilon,\quad
\tau=\tau_\varepsilon,\qquad
(\theta_\varepsilon,y_\varepsilon,\tau_\varepsilon)\to(\theta,y,\tau),\quad
\theta,y>0,\ \tau\ge0,
\end{equation}
stable for all small \(\varepsilon\).  Let
\(\zeta_\pm\) be the roots of \(\zeta^2-\theta\zeta+y\), ordered as
\(\zeta_-\le\zeta_+\) when real, and set
\begin{equation}\label{eq:cQ}
\cQ(\theta,y)=\begin{cases}
\theta/2,&\theta^2\le4y,\\
(\theta-\sqrt{\theta^2-4y})/2,&\theta^2>4y.
\end{cases}
\end{equation}

\begin{lemma}[Quadratic constant]\label{lem:cq}
Under~\eqref{eq:scaling}, \(\limsup_{\kappa\to\infty}\sqrt\kappa\,(1-\rho_q)\le\cQ(\theta,y)\), and if
\(y(1+2\tau)<4\), then \(\sqrt\kappa\,(1-\rho_q)\to\cQ(\theta,y)\).
\end{lemma}

The function \(\cQ\) is continuous, \(\cQ\le\theta/2\) and \(\cQ\le\sqrt y\),
with equality in both at critical damping \(\theta=2\sqrt y\).

A family is \emph{accelerated} on a function class if its worst-case rate
is at most \(1-c\varepsilon\) for some \(c>0\) and all small \(\varepsilon\).
Acceleration on \(\cF_{\mu,L}\) implies acceleration on its quadratic subclass.

In asymptotic statements, \((\theta,y,\tau)\) denotes the limit of
\((\theta_\varepsilon,y_\varepsilon,\tau_\varepsilon)\).
Finite-tuning identities use \(\beta,\tau,h\) and \(y=hL\).  Put
\begin{equation}\label{eq:mainconstant}
\cFF(\theta,y)=\min\{\theta/3,\cQ(\theta,y)\}.
\end{equation}
Spiral mixtures and quadratic modes supply the two branches of this bound.

\begin{theorem}[Exact constant]\label{thm:B}
Consider a family~\eqref{eq:scaling} of stable tunings of \(\cM\).
\begin{enumerate}[label=\textup{(\alph*)}]
\item \textup{(Lower bound, every \(\tau\ge0\).)}  For every \(\delta>0\) there
are \(J=J(\theta,y,\delta)\in\N\) and \(\varepsilon_0>0\) such that for every
\(\varepsilon\in(0,\varepsilon_0)\) there exist
\(f\in\cF_{\varepsilon^2,1}(\R^{2J+1})\) and \(x_0\ne x^\star\) with, from rest,
\(\norm{x_N-x^\star}\ge c_0\,e^{-(\cFF+\delta)\varepsilon N}\norm{x_0-x^\star}\)
for all \(N\ge0\), where \(c_0=c_0(\varepsilon,\delta)>0\) is independent of \(N\).
\item \textup{(Upper bound beyond the wall.)}  Let \(\tau>\frac12\) and \(y(1+2\tau)<4\).  Then for every \(c\in(0,\cFF(\theta,y))\)
there is \(\varepsilon_0>0\) such that \(\rhoF<1-c\varepsilon\) for
\(\varepsilon\in(0,\varepsilon_0)\); this bound is certified by the lag-one
criterion (\Cref{prop:lagone}).
\item Consequently, under the hypotheses of~(b),
\[
\lim_{\kappa\to\infty}\sqrt\kappa\,(1-\rhoF)=\lim_{\kappa\to\infty}\sqrt\kappa\,(1-\rhoz)=\cFF(\theta,y)=\min\bigl(\tfrac\theta3,\ \cQ(\theta,y)\bigr).
\]
\end{enumerate}
\end{theorem}

Every accelerated family satisfies \(\liminf\tau_\varepsilon\ge1/2\)
(\Cref{thm:wall}); TMM reaches this limit from below.
The boundary layer and the stability boundary have their own scales
(\Cref{sec:layer,sec:TMM}).

\paragraph{Designing a target constant.}
For \(c>0\), the condition \(\cFF(\theta,y)\ge c\) is equivalent to
\[
\theta\ge3c,\qquad y\ge c\theta-c^2.
\]
Indeed, under \(\theta\ge3c\), the quadratic condition \(\cQ\ge c\) is
equivalent to \(c^2-\theta c+y\ge0\).
For \(\tau>1/2\) and \(y(1+2\tau)<4\), these are also the level sets of
the method's constant, by \Cref{thm:B}(c).
The smallest step compatible with a target constant is \(y=2c^2\), at
\(\theta=3c\).  Combining this with \(y(1+2\tau)\le4\) explains the
fixed-extrapolation ceiling \(\sqrt{2/(1+2\tau)}\)
(\Cref{cor:frontier}).

\subsection{The quadratic price of extrapolation}
\begin{theorem}[Price of extrapolation]\label{thm:price}
Every tuning of \(\cM\) with \(\rho_q<1\) satisfies
\[
\rho_q\ \ge\ \frac{\sqrt{\kappa_\tau}-1}{\sqrt{\kappa_\tau}+1},\qquad
\kappa_\tau:=(1+2\tau)\kappa-2\tau .
\]
Equality holds exactly at Polyak's tuning,
\[
\tau=0,\qquad
\beta=\left(\frac{\sqrt\kappa-1}{\sqrt\kappa+1}\right)^2,\qquad
h=\frac4{(\sqrt L+\sqrt\mu)^2}.
\]
\end{theorem}

Extrapolation therefore costs at least the replacement of \(\kappa\) by
\(\kappa_\tau\), asymptotically a factor \(1+2\tau\).
Queries also admit a gradient-correction realization, including
quasi-hyperbolic momentum; \Cref{rem:twostep} gives the precise parameter
translation and initialization convention.

\section{The Heavy-Ball Flow}
\label{sec:flow}

In the accelerated scaling~\eqref{eq:scaling}, the slow time
\(t=\varepsilon k\) gives the heavy-ball equation
\[
\ddot x+\theta\dot x+\nabla F(x)=0,\qquad
F=\tfrac{y}{\varepsilon^2}f\in\cF_{y,y\kappa}.
\]
Its smoothness bound grows with \(\kappa\).  We determine the worst decay
exponent over all finite smoothness bounds and transfer sampled spiral
witnesses to discrete time by interpolation (\Cref{lem:flowspirals}).

Let \(H\) be a real Hilbert space and \(m>0\).  Let \(\cFm(H)\) consist of
differentiable \(m\)-strongly convex functions whose gradient is Lipschitz on
bounded sets; thus \(\bigcup_{L>m}\cF_{m,L}(H)\subseteq\cFm(H)\).
Each has a unique minimizer, which we normalize to \(0\), with \(f(0)=0\).
The convex remainder \(f_m(x)=f(x)-\frac m2\norm x^2\) satisfies
\(f_m(0)=0\) and \(\nabla f_m(0)=0\), so
\begin{equation}\label{eq:hineq}
f_m(x)\ge0,\qquad \ip{\nabla f_m(x)}x-f_m(x)\ge0\qquad(x\in H),
\end{equation}
by convexity between \(x\) and \(0\).  For \(\theta>0\) consider
\begin{equation}\label{eq:flow}
\ddot x+\theta\dot x+\nabla f(x)=0 .
\end{equation}
Every initial position and velocity determine a unique solution for all
\(t\ge0\): the vector field is Lipschitz on bounded sets, and the energy
\(f(x)+\frac12\norm{\dot x}^2\) keeps the orbit bounded.
For a nonzero solution define the state decay exponent by
\(\ell(x)=\liminf_{t\to\infty}-\frac1t\log\bigl(\norm{x(t)}+\norm{\dot x(t)}\bigr)\).
With \(\cQ(\theta,m)\) as in~\eqref{eq:cQ} (with \(y\) replaced by \(m\)), put
\begin{equation}\label{eq:cF}
\Delta_c=c^2-\theta c+m,\qquad \cFF(\theta,m)=\min\bigl(\tfrac\theta3,\ \cQ(\theta,m)\bigr).
\end{equation}
For \(f(x)=\frac m2\norm x^2\), the slowest state decay exponent is \(\cQ(\theta,m)\).

\subsection{Exact decay}

\begin{theorem}[Exact decay of the heavy-ball flow]\label{thm:flow}
Let \(\theta,m>0\).
\begin{enumerate}[label=\textup{(\alph*)}]
\item \textup{(Lyapunov identity.)}  For \(c\in\R\) put
\[
V_c=f(x)-\tfrac{c(\theta-c)}2\norm x^2+\tfrac12\norm{\dot x+(\theta-c)x}^2
=f_m(x)+\tfrac{\Delta_c}2\norm x^2+\tfrac12\norm{\dot x+(\theta-c)x}^2 .
\]
Along every solution of~\eqref{eq:flow},
\begin{equation}\label{eq:identity}
\frac{d}{dt}\bigl(e^{2ct}V_c\bigr)=-e^{2ct}\Bigl[(\theta-2c)\Delta_c\norm x^2+(\theta-c)\bigl(\ip{\nabla f_m(x)}x-f_m(x)\bigr)+(\theta-3c)\,f_m(x)\Bigr].
\end{equation}
If \(0<c\le\theta/3\) and \(\Delta_c>0\), the right-hand side is nonpositive, and
\(\norm{x(t)}^2+\norm{\dot x(t)}^2\le K_c\,e^{-2ct}V_c(0)\) with \(K_c\) depending only on \((\theta,m,c)\).
\item \textup{(Quadratic witness.)}  For \(f(x)=\frac m2x^2\) on \(\R\), some solution has \(\ell=\cQ(\theta,m)\).
\item \textup{(Spiral witnesses.)}  For every \(c\in(\theta/3,\cQ(\theta,m))\) there are \(J\in\N\), \(L<\infty\),
\(f\in\cF_{m,L}(\R^{2J})\) and a solution of~\eqref{eq:flow} with \(\norm{x(t)}=e^{-ct}\norm{x(0)}\) and
\(\norm{\dot x(t)}=e^{-ct}\norm{\dot x(0)}\) for all \(t\ge0\).
\item \textup{(Exact rate.)}  Over all real Hilbert spaces \(H\), all
\(f\in\cFm(H)\) and all nonzero solutions,
\[
\inf\ell(x)=\cFF(\theta,m).
\]
The same infimum holds with \(f\in\bigcup_{L>m}\cF_{m,L}(H)\).
On both objective classes, the energy \(f(x)+\frac12\norm{\dot x}^2\)
and the gap \(f(x)\) have logarithmic decay exponents
\(\liminf_{t\to\infty}-t^{-1}\log(\cdot)\) with infimum \(2\cFF(\theta,m)\),
using \(-\log0=+\infty\).
If \(\theta<3\sqrt{m/2}\), every solution satisfies
\(\norm{x(t)}+\norm{\dot x(t)}=O(e^{-\theta t/3})\).
\end{enumerate}
\end{theorem}

The function-value bounds in~(d) are due to Aujol, Dossal and
Rondepierre~\cite{ADR22}, who proved them for quasi-strongly convex \(C^1\)
functions.  Their energy at \(\theta\le3\sqrt{m/2}\) is
\(V_{\theta/3}=f+\frac12\norm{\dot x+\frac{2\theta}3x}^2-\frac{\theta^2}9\norm x^2\).
For \(\theta<3\sqrt{m/2}\), they proved sharpness of the gap exponent
\(2\theta/3\) using one-dimensional examples with non-Lipschitz gradients,
and improved the bound for each fixed Lipschitz constant \(L\).
\Cref{thm:flow}(c) supplies finite-dimensional \(C^{1,1}\) witnesses that are
sharp for distances.  Their dimension may depend on the exponent, and
\(L\) grows as \(c\downarrow\theta/3\).

The coefficient \(\theta/3\) comes from the two convexity inequalities
in~\eqref{eq:hineq}.  After the chain rule cancels the gradient--velocity
term, a homogeneous remainder with \(\ip{\nabla f_m(x)}x=r f_m(x)\)
has coefficient \(r\theta-(r+2)c\) inside the brackets in~\eqref{eq:identity}.
This coefficient vanishes at \(c=r\theta/(r+2)\): at \(\theta/2\) for
quadratic profiles and \(\theta/3\) for conical profiles, whose degree \(1\)
is the minimum allowed by convexity.  The smooth witnesses in~(c) have
nearly conical profiles along their orbits (\Cref{sec:spirals}).

\begin{corollary}[The price of nonlinearity]\label{cor:damping}
For \(m>0\),
\[
\max_{\theta>0}\cFF(\theta,m)=\sqrt{m/2},\qquad
\max_{\theta>0}\cQ(\theta,m)=\sqrt m.
\]
The maxima are attained uniquely at \(\theta=3\sqrt{m/2}\) and
\(\theta=2\sqrt m\), respectively; at the first, \(\cQ=\theta/3\).
Every pair \((\cQ,\cFF)(\theta,m)\) is dominated componentwise by a point of
\(\{(\cQ(\vartheta,m),\vartheta/3):2\sqrt m\le\vartheta\le3\sqrt{m/2}\}\).
\end{corollary}

\begin{figure}[t]
\centering
\begin{tikzpicture}
\pgfmathsetmacro{\flowOptDamping}{3/sqrt(2)}
\pgfmathsetmacro{\flowOptRate}{1/sqrt(2)}
\begin{axis}[width=0.66\textwidth,height=0.35\textwidth,
  xmin=0,xmax=4.2,ymin=0,ymax=1.12,
  xlabel={\(\theta/\sqrt m\)},ylabel={\(c/\sqrt m\)},
  xtick={0,1,2,3,4},ytick={0,0.5,1},
  tick label style={font=\footnotesize},label style={font=\small},
  axis lines=left,clip=false,
  legend style={at={(0.5,1.04)},anchor=south,font=\scriptsize,
    draw=none,legend columns=2},legend cell align=left]
\addplot[thick,dashed,black!65,domain=0:2,samples=2] {x/2};
\addlegendentry{quadratics: \(c_Q\)}
\addplot[thick,domain=0:\flowOptDamping,samples=2] {x/3};
\addlegendentry{full class: \(c_{\cF}\)}
\addplot[thick,dashed,black!65,domain=2:4.2,samples=120,forget plot]
  {(x-sqrt(max(x^2-4,0)))/2};
\addplot[thick,domain=\flowOptDamping:4.2,samples=100,forget plot]
  {(x-sqrt(max(x^2-4,0)))/2};
\addplot[thin,dotted,black!45,forget plot] coordinates {(2,0) (2,1)};
\addplot[thin,dotted,black!45,forget plot]
  coordinates {(\flowOptDamping,0) (\flowOptDamping,\flowOptRate)};
\addplot[only marks,mark=*,mark size=1.5pt,forget plot]
  coordinates {(2,1) (2,0.6666667) (\flowOptDamping,\flowOptRate)};
\node[font=\scriptsize,anchor=south west] at (axis cs:2.04,1.01) {\((2,1)\)};
\node[font=\scriptsize,anchor=east] at (axis cs:1.92,0.64) {\(2/3\)};
\node[font=\scriptsize,anchor=west] at (axis cs:2.6,0.82)
  {\((3/\sqrt2,1/\sqrt2)\)};
\draw[thin,black!60] (axis cs:2.56,0.8)--(axis cs:\flowOptDamping,\flowOptRate);
\end{axis}
\end{tikzpicture}
\caption{Normalized flow exponents (\Cref{thm:flow,cor:damping}).
The marked maxima occur at \(\theta/\sqrt m=2\) for quadratics and
\(3/\sqrt2\) for the full class.  At critical damping, the full-class value
is \(2/3\), corresponding to standard Nesterov (\Cref{cor:B}(i)).}
\label{fig:price}
\end{figure}
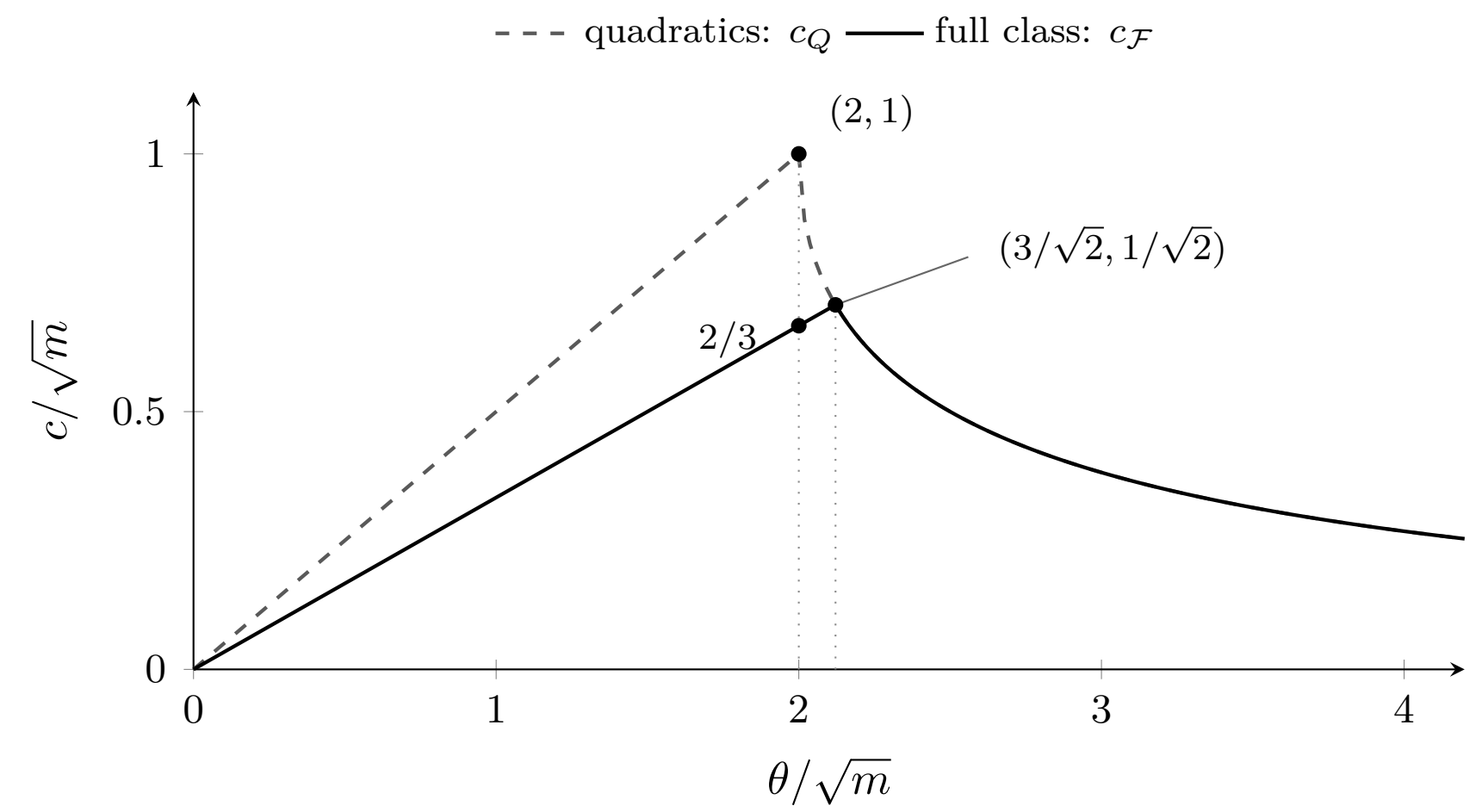

The full-class optimum is \(1/\sqrt2\) times the quadratic optimum and
uses \(3/(2\sqrt2)\approx1.06\) times as much friction.
At critical damping the full-class exponent is \(2\sqrt m/3\).
The arc in \Cref{cor:damping} records the trade-off between the two
exponents (\Cref{fig:price}); discrete tunings also face stability constraints.

\subsection{Why the worst case mixes spirals}
\label{sec:spirals}

The nonquadratic witnesses are orthogonal sums of spirals,
\[
x(t)=\bigl(\sqrt{\nu_j}\,e^{(-c+i\gamma_j)t}\bigr)_{j=1}^J,
\qquad \nu_j\ge0,\quad\sum_j\nu_j=1.
\]
The flow fixes the gradient along each orbit.  Existence of a strongly
convex potential realizing these data is characterized by the interpolation
program \(\mathcal P_c\) in \Cref{lem:creduce}.  Its minimizer constraint is
\[
(3c-\theta)\sum_j\nu_j\gamma_j^2\ge c\Delta_c.
\]
For \(c<\cQ\), \(\Delta_c>0\), so it forces \(c>\theta/3\), matching the coefficient
\(\theta-3c\) in the Lyapunov identity.  The other constraints test convexity between every pair
of orbit points, in both time directions.

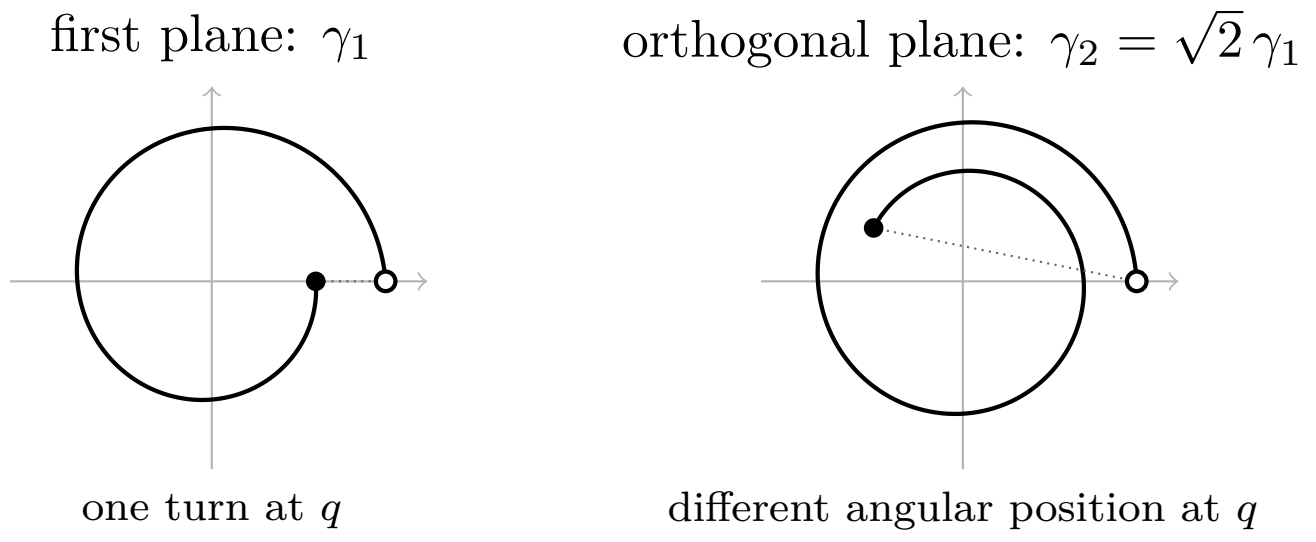
\begin{figure}[t]
\centering
\begin{tikzpicture}[scale=0.92]
\begin{scope}
\draw[->,black!30] (-1.45,0)--(1.55,0);
\draw[->,black!30] (0,-1.35)--(0,1.4);
\draw[thick,domain=0:360,samples=160,variable=\a]
 plot ({1.25*exp(-\a/700)*cos(\a)},{1.25*exp(-\a/700)*sin(\a)});
\draw[densely dotted,black!60] (1.25,0)--({1.25*exp(-360/700)},0);
\draw[fill=white,thick] (1.25,0) circle (2pt);
\fill ({1.25*exp(-360/700)},0) circle (2pt);
\node[font=\small] at (0,1.75) {first plane: \(\gamma_1\)};
\node[font=\scriptsize] at (0,-1.65) {one turn at \(q\)};
\end{scope}
\begin{scope}[xshift=5.4cm]
\draw[->,black!30] (-1.45,0)--(1.55,0);
\draw[->,black!30] (0,-1.35)--(0,1.4);
\draw[thick,domain=0:360,samples=200,variable=\a]
 plot ({1.25*exp(-\a/700)*cos(sqrt(2)*\a)},{1.25*exp(-\a/700)*sin(sqrt(2)*\a)});
\draw[densely dotted,black!60] (1.25,0)--
 ({1.25*exp(-360/700)*cos(sqrt(2)*360)},{1.25*exp(-360/700)*sin(sqrt(2)*360)});
\draw[fill=white,thick] (1.25,0) circle (2pt);
\fill ({1.25*exp(-360/700)*cos(sqrt(2)*360)},{1.25*exp(-360/700)*sin(sqrt(2)*360)}) circle (2pt);
\node[font=\small] at (0,1.75) {orthogonal plane: \(\gamma_2=\sqrt2\,\gamma_1\)};
\node[font=\scriptsize] at (0,-1.65) {different angular position at \(q\)};
\end{scope}
\end{tikzpicture}
\caption{Schematic desynchronization in two orthogonal projections.
At the common lag \(q=2\pi/\gamma_1\), the first returns to its starting ray
and the second has a different angular position.  Open and filled markers
show times \(0\) and \(q\); dotted chords join them.}
\label{fig:spiral-mechanism}
\end{figure}

For \(\theta/3<c<\cQ\), the equally weighted frequencies
\(\gamma_j=\gamma_0(1+(j-\tfrac12)/J)\) satisfy all constraints with a uniform
margin when \(\gamma_0\) and then \(J\) are large (\Cref{lem:spread});
\Cref{lem:extend} supplies a \(C^{1,1}\) potential.
As the frequencies grow, \(\ip{\nabla f_m(x)}x/f_m(x)\to2c/(\theta-c)\).
This limit tends to \(1\) as \(c\downarrow\theta/3\), giving a nearly conical
profile.  The frequency spread separates the components' first full turns
(\Cref{fig:spiral-mechanism});
\Cref{prop:one} shows why one spiral fails near the limiting exponent.

\begin{proposition}[Single-spiral obstruction]\label{prop:one}
Let \(\theta/3<\cQ(\theta,m)\).  There is \(c_\sharp\in(\theta/3,\cQ)\) such that for \(\theta/3<c\le c_\sharp\)
no one-atom measure satisfies~\eqref{eq:P12}--\eqref{eq:P3}.  Consequently, for such \(c\), no
\(f\in\cFm(H)\) has a solution of~\eqref{eq:flow} of the form \(x(t)=e^{-ct}(a\cos\gamma t+b\sin\gamma t)\) with
\(a\perp b\), \(\norm a=\norm b>0\) and \(\gamma\ge0\).
\end{proposition}

Near \(c=\theta/3\), the minimizer constraint forces a single spiral to wind
rapidly, violating convexity at its first return to the starting ray (\Cref{rem:cone}).

We now characterize discrete self-similar trajectories and the certificates
that exclude larger growth factors.

\section{Spectral Closure for Stationary Linear Methods}
\label{sec:closure}

We characterize orbits generated by one smooth strongly convex function.
Interpolation reduces this question to pairwise inequalities.  For orthogonal
rotation--dilation modes, these inequalities depend only on the lag; squared
amplitudes define a spectral measure, and one scalar encodes the function values.
The resulting spectral problem and its finite-lag certificates determine the
same exact rate (\Cref{thm:A}).

\subsection{Normalization and interpolation}

Translate the minimizer to \(0\), set \(f(0)=0\), and define
\begin{equation}\label{eq:psi}
\psi(e)=\tfrac{2}{L-\mu}\bigl(\tfrac{L+\mu}{4}\norm e^2-f(e)\bigr),\qquad
\nabla f(e)=\tfrac{L+\mu}{2}\,e-\tfrac{L-\mu}{2}\,\nabla\psi(e).
\end{equation}
Let \(\cS(H)\) be the set of differentiable \(\psi:H\to\R\) with
\(\psi(0)=0\), \(\nabla\psi(0)=0\) and \(1\)-Lipschitz gradient.

\begin{lemma}[{\cite[Lemma~A.1]{TA}}]\label{lem:norm}
The map \(f\mapsto\psi\) is a bijection from
\(\{f\in\cF_{\mu,L}(H):x^\star=0,\ f(0)=0\}\) onto \(\cS(H)\).
\end{lemma}

The lemma is stated in~\cite{TA} for \(\R^d\); its proof (cocoercivity
of the gradient of a convex \((L-\mu)\)-smooth function~\cite{THG17}) is
dimension-free.  Every \(\psi\in\cS(H)\) satisfies
\(\norm{\nabla\psi(e)}\le\norm e\) and \(|\psi(e)|\le\frac12\norm e^2\).
For data \((\tilde x_i,\sigma_i,\psi_i)_{i\in I}\) in \(H\times H\times\R\) put
\begin{equation}\label{eq:TR}
T_{ij}=\psi_j-\psi_i-\tfrac12\ip{\sigma_i+\sigma_j}{\tilde x_j-\tilde x_i},\qquad
R_{ij}=\tfrac14\bigl(\norm{\tilde x_j-\tilde x_i}^2-\norm{\sigma_j-\sigma_i}^2\bigr).
\end{equation}

\begin{theorem}[{\cite[Theorem~A.2]{TA}}]\label{thm:interp}
Suppose \((\tilde x_\star,\sigma_\star,\psi_\star)=(0,0,0)\) for some
\(\star\in I\).  There is a function \(\psi\in\cS(H)\) interpolating all data,
\(\nabla\psi(\tilde x_i)=\sigma_i\), \(\psi(\tilde x_i)=\psi_i\), if and only if
\(|T_{ij}|\le R_{ij}\) for every \(i,j\in I\).
\end{theorem}

\Cref{thm:interp} is Azagra and Mudarra's extension theorem for convex
\(C^{1,1}\) functions on Hilbert spaces~\cite[Theorem~1.3]{AM17}, expressed
in the normalization \(\psi\).  It applies to arbitrary index sets and hence
to infinite trajectories.  On finite sets it is equivalent to the conditions
of Taylor, Hendrickx and Glineur~\cite{THG17}.

\subsection{Stationary linear methods}

\begin{definition}\label{def:slm}
Let \(D\ge1\), and let \(u(z)=\sum_{j=0}^{D-1}u_jz^j\not\equiv0\) and
\(v(z)=\sum_{j=0}^Dv_jz^j\), \(v_D\neq0\), be real polynomials
with \(u(1)=1\) and \(v(1)=0\).  A \emph{trajectory} of the stationary
linear method \((u,v)\) on \(f\in\cF_{\mu,L}(H)\) is a sequence
\((w_k)_{k\ge0}\subset H\) with
\begin{equation}\label{eq:slm}
\tilde x_k:=\sum_{j=0}^{D-1}u_jw_{k+j},\qquad
\sum_{j=0}^{D}v_jw_{k+j}=\nabla f(\tilde x_k)\qquad(k\ge0).
\end{equation}
We call \(w_k\) the \emph{latent iterates} and \(\tilde x_k\) the \emph{query points},
where gradients are evaluated.  Its \emph{state} is
\(S_k=(w_k-x^\star,\dots,w_{k+D-1}-x^\star)\in H^D\);
\eqref{eq:slm} determines \(S_{k+1}\) from \(S_k\).
\end{definition}

The identities \(u(1)=1\), \(v(1)=0\) preserve the minimizer as a fixed
point; henceforth \(x^\star=0\).  Every recurrence with constant scalar
coefficients
\(x_{k+1}=\sum_{j<D}b_jx_{k-j}-h\nabla f(\sum_{j<D}a_jx_{k-j})\) with
\(\sum_ja_j=\sum_jb_j=1\) fits this definition with \(w_k=x_{k-D+1}\).  Common
roots of \(u\) and \(v\) represent modes that evolve independently of \(f\).  In particular,
\eqref{eq:M} is the case \(D=2\), \(w_k=x_{k-1}\),
\begin{equation}\label{eq:Muv}
u(z)=N_\tau(z)=(1+\tau)z-\tau,\qquad v(z)=-\tfrac1h(z-1)(z-\beta),
\end{equation}
and \(u,v\) are coprime unless \(\beta=\tau/(1+\tau)\).  In normalized
form~\eqref{eq:psi} the method reads
\(\sigma_k:=\nabla\psi(\tilde x_k)=\sum_jp_jw_{k+j}\) with
\begin{equation}\label{eq:p}
p(z)=\tfrac{2}{L-\mu}\bigl(\tfrac{L+\mu}{2}\,u(z)-v(z)\bigr),
\end{equation}
a polynomial of degree \(D\), coprime to \(u\) whenever \(u,v\) are coprime.
The \emph{quadratic rate} and the \emph{middle-curvature radius} are
\[
\rho_q=\max_{\lambda\in[\mu,L]}\max\{|z|:v(z)=\lambda u(z)\},\qquad
r_c=\max\{|z|:p(z)=0\}\le\rho_q.
\]
Thus \(r_c\) is the root radius at curvature \((L+\mu)/2\).
Moreover \(\rho_q>0\): otherwise both \(v-\mu u\) and \(v-Lu\) would
be \(v_Dz^D\), forcing \(u\equiv0\).  We set \(\log0=-\infty\).

For \(\cM\), the endpoint polynomials \(\chi_\mu,\chi_L\) give
explicit normalized symbols.  Indeed,~\eqref{eq:Muv} yields
\(u-p=-2\chi_\mu/(h(L-\mu))\) and \(u+p=2\chi_L/(h(L-\mu))\), hence for
every \(z\in\C\)
\begin{equation}\label{eq:factor}
|u|^2-|p|^2-2i\,\Ima(p\bar u)=(u-p)(\bar u+\bar p)
=-\frac{4\,\chi_\mu(z)\overline{\chi_L(z)}}{h^2(L-\mu)^2},
\end{equation}
and \(r_c\) is the largest root modulus of \(\chi_{(L+\mu)/2}\).  On the
unit circle the phase of \(\chi_\mu\overline{\chi_L}\) is controlled by
\begin{equation}\label{eq:imab}
\Ima\bigl(\chi_\mu(e^{i\omega})\overline{\chi_L(e^{i\omega})}\bigr)
=(y-y_\mu)\sin\omega\,\bigl[(1-\beta)+2\tau\xi\bigr],\qquad\xi=1-\cos\omega,
\end{equation}
which is positive on \((0,\pi)\) (\Cref{app:prelim}).

\begin{remark}[Two-step methods]\label{rem:twostep}
Every scalar stationary method with \(D=2\) has the algebraic form~\eqref{eq:M}:
\[
u=(1+\tau)z-\tau,\quad v=-h^{-1}(z-1)(z-\beta),\qquad
\tau=-u(0),\quad h=-1/v_2,\quad\beta=v(0)/v_2.
\]
The tunings studied here satisfy \(\beta\in[0,1)\), \(\tau\ge0\), \(h>0\).
The queries obey the gradient-correction recurrence
\[
\tilde x_{k+1}=\tilde x_k+\beta(\tilde x_k-\tilde x_{k-1})-h_1\nabla f(\tilde x_k)+h_2\nabla f(\tilde x_{k-1}),\qquad h_1=h(1+\tau),\ h_2=h\tau ,
\]
Conversely, \(h_1>h_2\ge0\) gives the query representation~\eqref{eq:M}
with \(h=h_1-h_2\), \(\tau=h_2/(h_1-h_2)\).
For \(\beta\ne z_\tau:=\tau/(1+\tau)\), the two state descriptions have the same
worst-case rate.  Let \(\rho^0_{\tilde x}\) denote the from-rest rate of the
queries: replace \(x_N\) with \(\tilde x_N\) in the numerator of the from-rest formula
in \Cref{def:rate}.
Then
\[
\rho^0_{\tilde x}\le\rhoz\le\max\{\rho^0_{\tilde x},z_\tau\}.
\]
The two rates coincide above \(z_\tau\), including every accelerated regime
(\Cref{app:prelim}).

For example, the quasi-hyperbolic momentum (QHM) method of Ma and Yarats~\cite{MY19},
\(d_k=\beta d_{k-1}+(1-\beta)\nabla f(x_k)\), \(x_{k+1}=x_k-\alpha[(1-\nu_{\rm qh})\nabla f(x_k)+\nu_{\rm qh}d_k]\) with
\(\nu_{\rm qh}\in[0,1]\), has the query representation~\eqref{eq:M} with
\[
h=\alpha(1-\beta),\qquad \tau=\tau_{\rm qh}:=\frac{(1-\nu_{\rm qh})\beta}{1-\beta};
\]
\(\nu_{\rm qh}=1\) gives heavy ball and \(\nu_{\rm qh}=\beta\) gives Nesterov's method.  With the usual
initialization \(d_{-1}=0\), its iterates are exactly the queries of~\eqref{eq:M} started from rest at \(x_0\)
(\Cref{app:prelim}).  The query-rate comparison transfers the results below
to QHM with this initialization.
\end{remark}

\subsection{Spectra and certificates}

Fix a stationary linear method \((u,v)\) with normalized polynomial \(p\)
from~\eqref{eq:p}.  For \(z=e^{t+i\omega}\) the \emph{wave} \(w_k=z^k\)
rotates by \(\omega\) and dilates by \(e^t\) per step; its queries are
\(\tilde x_k=u(z)z^k\), and~\eqref{eq:slm} requires the normalized gradients
\(\sigma_k=p(z)z^k\).  Put
\begin{equation}\label{eq:gs}
g(t,\omega)=|u(z)|^2-|p(z)|^2,\qquad s(t,\omega)=\Ima\bigl(p(z)\overline{u(z)}\bigr),
\end{equation}
\begin{equation}\label{eq:anbn}
a_n=g\,(\ch nt-\cos n\omega),\qquad b_n=2s\sin n\omega\qquad(n\ge1).
\end{equation}
For waves with a common exponent \(t\), normalized squared amplitudes define
a probability measure \(\nu\) on \([0,\pi]\), called the \emph{spectrum}.
Each atom represents one planar wave component.  One scalar \(C\) encodes the
function values; interpolation with the minimizer bounds it by \(|C|\le\bar g\).
This scalar is essential for sharp boundary-layer witnesses.

\begin{definition}\label{def:feas}
A pair \((\nu,C)\) is \emph{feasible at \(t\)} if
\begin{equation}\label{eq:feas}
\begin{gathered}
|C|\le\bar g,\\
\ell_n:=\bar a_n-\bar b_n+C\sh nt\ge0,\qquad
\upsilon_n:=\bar a_n+\bar b_n-C\sh nt\ge0\qquad(n\ge1).
\end{gathered}
\end{equation}
The numbers \(\ell_n,\upsilon_n\) are the lower and upper \emph{lag-\(n\) slacks},
for queries \(n\) steps apart; \(\bar g\pm C\ge0\) are the two minimizer inequalities.
We call \(|C|<\bar g\) \emph{strict minimizer slack}.
A lag \(n\) is \emph{active} when \(\ell_n=0\) or \(\upsilon_n=0\).
Let \(\cT\) be the set of \(t\in\R\) at which some pair is feasible.
\end{definition}

\begin{definition}\label{def:cert}
A \emph{strict certificate at \(t\)} consists of nonnegative, finitely
supported weights \(\alpha_n,\beta_n\), numbers \(\gamma_\pm\ge0\) and
\(\epsilon>0\) such that
\begin{equation}\label{eq:cert}
\begin{gathered}
\Phi_t(\omega):=\sum_n\bigl\{\alpha_n(a_n-b_n)+\beta_n(a_n+b_n)\bigr\}+(\gamma_++\gamma_-)g\le-\epsilon\quad(0\le\omega\le\pi),\\
\sum_n(\alpha_n-\beta_n)\sh nt+\gamma_+-\gamma_-=0.
\end{gathered}
\end{equation}
\end{definition}

The inequalities~\eqref{eq:feas} are the interpolation conditions between
queries \(n\) steps apart and with the minimizer (\Cref{lem:waves}).
A certificate combines them with nonnegative weights that cancel the function
values.  It is a finite-lag Lyapunov argument, or a discrete-time
O'Shea--Zames--Falb multiplier~\cite{OY67,ZF68,WB68,TA}.  For heavy ball, \(u(z)=z\), and
\eqref{eq:gs} reduces, up to the positive factor \(|z|^2\), to the
quantities of~\cite[eq.~(12) and Definition~2.4]{TA}.

At \(t=0\) the pair \((g,s)\) defines a phase: where \(g-2is\ne0\), let \(\Theta(\omega)\in(-\pi,\pi]\) be given by
\begin{equation}\label{eq:Theta}
g(0,\omega)-2is(0,\omega)=\bigl|g(0,\omega)-2is(0,\omega)\bigr|\,e^{-i\Theta(\omega)} .
\end{equation}
For \(\omega>0\), the \emph{spectral phase lead} is \(\Theta(\omega)/\omega\), measured in steps.
For a limit method of \(\cM\), this lead tends to \(\tau\) at low frequencies
(\Cref{sec:limitphase,sec:faithful}).

\subsection{The closure theorem}

\begin{theorem}[Spectral closure]\label{thm:A}
Let \((u,v)\) be a stationary linear method.
\begin{enumerate}[label=\textup{(\alph*)}]
\item \(\cT\) is nonempty, closed and bounded above, and
\(T:=\max\cT\ge\log\rho_q\).
\item \(\rhoF=e^T\).  Every \(t>T\) admits a strict certificate, and
feasibility at \(t\) excludes one.
\item \textup{(One exponent suffices.)}  If a strict certificate exists at
some \(t\) with \(e^t>r_c\), then \(\rhoF<e^t\).  Consequently
\(\cT\cap(\log r_c,\infty)=(\log r_c,T]\).
\item Every pair \((\nu,C)\) feasible at \(t\) is realized: there are
\(\psi\in\cS(L^2(\nu;\C))\) and a trajectory with \(S_{k+1}=e^tUS_k\) for an
orthogonal \(U\), hence \(\norm{S_N}=e^{Nt}\norm{S_0}\) for all \(N\); if
\(\nu\) has \(J\) atoms, then \(L^2(\nu;\C)\cong\R^{2J}\).
\item If \(T=\log\rho_q\), a quadratic on \(\R^2\) attains the worst case.
If \(T>\log\rho_q\) and \(T\ne0\), every pair feasible at \(T\) has
\(|C|<\bar g\) and finitely many active lags.  Fix any such pair and let
\(M\) be its largest active lag.  Every pair feasible at \(T\) is then
carried by at most \((D+M)/2\) points of \((0,\pi)\), and the worst case
is attained on \(\R^{2J}\) with \(J\le(D+M)/2\).
\end{enumerate}
\end{theorem}

Thus finite-lag certificates approach the exact rate from above, and each
feasible spectrum generates one self-similar trajectory.  Part~(e) gives
finite-dimensional attainment at the spectral boundary under its stated
hypotheses.  The interpolation, separation and one-exponent arguments
for quadratically stable heavy ball are those of~\cite[Theorems~2.5, 3.2
and Lemmas~A.3--A.7]{TA}; finite support, strict minimizer slack and
inactive lag tails were proved there for Polyak's tuning
\cite[Lemma~A.10]{TA}.  For standard Nesterov, a linear program over spirals gives the exact
squared-distance factor~\cite[Theorem~4.1]{PS26}.  The latent polynomials
extend the characterization to arbitrary fixed finite memory and give the
atom bound in~(e).  Finite-dimensional attainment at \(T=0>\log\rho_q\) remains open
(\Cref{sec:conclusions}).

\subsection{Consequences for the class \texorpdfstring{\(\cM\)}{M}}
\label{sec:closureM}

For~\eqref{eq:M} we use \Cref{thm:A} with~\eqref{eq:Muv}.  By
\eqref{eq:factor}, \(g\) and \(s\) are explicit in \(\chi_\mu\) and
\(\chi_L\), and on the unit circle \eqref{eq:imab} gives \(s>0\) on
\((0,\pi)\).  Hence, for a stable tuning, the phase~\eqref{eq:Theta} satisfies \(\Theta(\omega)\in(0,\pi]\) for
\(\omega\in(0,\pi]\): it is the angle from \(-\chi_\mu(e^{i\omega})\) to \(\chi_L(e^{i\omega})\) (for heavy ball this
is~\cite[Definition~3.1]{TA}).

\begin{corollary}[Worst cases from rest]\label{cor:rest}
Let a tuning of \(\cM\) be stable, and let \((\nu,C)\) be feasible at \(t\),
with \(J\) atoms, and with \(|C|<\bar g\) if \(t>0\).  There are
\(f\in\cF_{\mu,L}(\R^{2J+1})\) and \(x_0\ne x^\star\) such that the
trajectory with \(x_{-1}=x_0\) satisfies
\(\norm{x_N-x^\star}\ge c_0e^{Nt}\norm{x_0-x^\star}\) for all \(N\ge0\), with
\(c_0>0\).  In particular, if either \(T=\log\rho_q\) or \(T>\log\rho_q\) and \(T\ne0\), then
\(\rhoz=\rhoF\).  The worst case from rest is attained by a quadratic on \(\R\) in the first case,
and by one function in a fixed dimension in the second.
\end{corollary}

\begin{proposition}[Lag-one criterion]\label{prop:lagone}
Let a tuning of \(\cM\) be stable.
\begin{enumerate}[label=\textup{(\alph*)}]
\item If \(\Rea[(1-z)\chi_\mu(z)\overline{\chi_L(z)}]>0\) on \(|z|=r\) for
some \(r\in(r_c,1)\), then \(\rhoF<r\).
\item If \(\Rea[(1-z)\chi_\mu(z)\overline{\chi_L(z)}]>0\) on \(|z|=1\),
\(z\ne1\), then \(\rhoF<1\).
\item More generally, if \(\Rea[\chi_\mu(z)\overline{\chi_L(z)}\,((1-z)+\lambda(1-\bar z^2))]>0\)
on \(|z|=1\), \(z\ne1\), for some \(\lambda\ge0\), then \(\rhoF<1\).
\end{enumerate}
\end{proposition}

We call \(\Rea[(1-z)\chi_\mu\overline{\chi_L}]\) the \emph{lag-one symbol}
and \Cref{prop:lagone} the \emph{lag-one criterion}.  The corresponding
\emph{lag-one certificate} combines \(\ell_1\ge0\) with a minimizer constraint
to cancel \(C\).  Its Jury--Lee multiplier \(1-z\)~\cite{JL64} underlies
rate bounds in~\cite{LRP16,CHVL18,Gan22} and is exact at the convergence
threshold of Polyak's heavy ball~\cite{BS19}.  On each circle its symbol is
a cubic in \(\xi\) (\Cref{lem:cubic}); \Cref{sec:exact,sec:layer} describe
when the criterion is exact.

For \(r\in(r_c,1]\), a strict certificate at \(t=\log r\) yields a Lyapunov function (\Cref{prop:lyap}).
For lag one, it is a quadratic form in two iterates and a gradient plus a
multiple of the function gap at the query, as
in~\cite{VFL18,LRP16,HL17,CHVL18}.  Longer certificates use a finite window
of query values and gradients.  For a stable tuning, persistence in the exponent (\Cref{lem:persist})
turns a strict certificate at \(t=0\) into a contraction with some \(r<1\).
Lower bounds come from the realizations of \Cref{cor:rest} and periodic rotating orbits (\Cref{prop:cycles}).

\begin{proposition}[Phase locking]\label{prop:cycles}
Let a tuning of \(\cM\) be stable, let \(K\ge3\), \(1\le j<K/2\),
\(\gcd(j,K)=1\) and \(\omega_0=2\pi j/K\).
\begin{enumerate}[label=\textup{(\alph*)}]
\item If \(\Theta(\omega_0)\le\pi/K\), there are \(f\in\cF_{\mu,L}(\R^3)\)
and \(x_0\) such that, from rest,
\(\norm{x_N-x^\star}\ge c_0\norm{x_0-x^\star}>0\) for all \(N\); the
trajectory converges to a \(K\)-periodic orbit that rotates by \(\omega_0\)
in a plane.
\item If \(\Theta(\omega_0)<\pi/K\), then \(\rhoz>1\).
\end{enumerate}
\end{proposition}

We call this periodic rotation \emph{phase locking}, and the corresponding
parameter regions \emph{\(K\)-tongues}.  For heavy ball the proposition
is~\cite[Theorem~3.3(ii)]{TA}, and these regions are the cycles of~\cite{GTD}.

The flow supplies slowly rotating witnesses; the discrete symbol identifies
the phase loss that determines whether the flow rate is attained.

\section{The Half-Step Rule for Fixed Finite Memory}
\label{sec:beyond}

We now allow a fixed finite number of momentum terms, past queries and gradient
corrections.  With the slow dynamics fixed, the flow still bounds the rate.
The limiting spectral phase governs attainment, and a trace identity relates
the step to the lag-one symbol.

For a stationary method \((u,v)\), write
\[
\phi_\lambda(z)=\lambda u(z)-v(z).
\]
This is the characteristic polynomial on \(\lambda\norm\cdot^2/2\): a wave
\(w_k=z^k\) is a trajectory exactly when \(\phi_\lambda(z)=0\).
For \(\cM\), \(\phi_\lambda=\chi_\lambda/h\).  By~\eqref{eq:p},
\(u-p=-\frac2{L-\mu}\phi_\mu\) and \(u+p=\frac2{L-\mu}\phi_L\), so that, as in~\eqref{eq:factor},
\begin{equation}\label{eq:factorgen}
g-2is=(u-p)\overline{(u+p)}=-\frac{4\,\phi_\mu(z)\overline{\phi_L(z)}}{(L-\mu)^2}.
\end{equation}
The proof of \Cref{prop:lagone}(a) uses this factorization and the one-exponent
principle in \Cref{thm:A}(c).  It therefore applies to every stationary method
with \(\phi_\lambda\) replacing \(\chi_\lambda\); we use the same terms
\emph{lag-one symbol}, \emph{certificate} and \emph{criterion}.

\subsection{A ceiling for every memory}
\label{sec:ceiling}

We keep the normalization \(L=1\), \(\mu=\varepsilon^2\).

\begin{definition}[Second-order families]\label{def:second}
Let \(D\ge2\) and, for \(\varepsilon\in(0,\varepsilon_0)\), let \((u_\varepsilon,v_\varepsilon)\) be stationary linear methods
with memory \(D\) whose coefficients extend continuously to \(\varepsilon=0\), with \(\deg v_0=D\).
The family is \emph{second-order with slow data} \((\theta,y)\), \(\theta,y>0\), if
\begin{enumerate}[label=\textup{(\roman*)}]
\item \(v_0(z)=(z-1)^2q_0(z)\) with \(q_0(1)=-1/y\);
\item the root \(\beta_\varepsilon\) of \(v_\varepsilon\) near \(1\) other than \(1\) satisfies \((1-\beta_\varepsilon)/\varepsilon\to\theta\).
\end{enumerate}
It is of class \(C^k\) if the coefficients are \(C^k\) functions of \(\varepsilon\in[0,\varepsilon_0)\).
\end{definition}

The limiting quotient \((z-1)q_0\) has a simple root at \(1\), so for small
\(\varepsilon\),
\[
v_\varepsilon=(z-1)(z-\beta_\varepsilon)q_\varepsilon,\qquad
\theta_\varepsilon:=\frac{1-\beta_\varepsilon}{\varepsilon}\to\theta,\qquad
y_\varepsilon:=-\frac1{q_\varepsilon(1)}\to y.
\]
Here \(\beta_\varepsilon\) is real, inherits the coefficient regularity,
and \(q_\varepsilon\to q_0\).  For \(\cM\), \(q=-1/h\), and the families~\eqref{eq:scaling} are
second-order with slow data \((\theta,y)\).  On a slow wave \(z=e^{-\varepsilon\varpi}\) every such method sees the complex curvature
\(v_\varepsilon(z)/u_\varepsilon(z)=\frac{\varepsilon^2}{y_\varepsilon}\varpi(\theta_\varepsilon-\varpi)+O(\varepsilon^3)\) (\Cref{lem:normal}); after the change \(F=\frac y{\varepsilon^2}f\) this is the
curvature \(\lambda'=\varpi(\theta-\varpi)\) at which \(e^{-\varpi t}\) solves \(\ddot x+\theta\dot x+\lambda'x=0\).  So every second-order family
with slow data \((\theta,y)\) discretizes the same slow flow as the tunings~\eqref{eq:scaling} of \(\cM\).

\begin{theorem}[Ceiling]\label{thm:ceiling}
For every second-order family with slow data \((\theta,y)\),
\[
\limsup_{\kappa\to\infty}\sqrt\kappa\,(1-\rhoF)\le\cFF(\theta,y)=\min\bigl(\tfrac\theta3,\cQ(\theta,y)\bigr).
\]
If \(\theta/3<c<\cQ(\theta,y)\), then for all small \(\varepsilon\) a spectrum with \(J\) atoms, \(J\)
independent of \(\varepsilon\), is feasible at \(t=-c\varepsilon\), and the method has a trajectory in
\(\R^{2J}\) with \(S_{k+1}=e^{-c\varepsilon}US_k\), \(U\) orthogonal.
\end{theorem}

Sampled flow spirals give the same ceiling for every family with these slow
data and any fixed finite memory.

\begin{corollary}[Consequences of the ceiling]\label{cor:ceiling}
Consider second-order families with slow data \((\theta,y)\).
\begin{enumerate}[label=\textup{(\roman*)}]
\item \(\limsup\sqrt\kappa(1-\rhoF)\le\sqrt{y/2}\), and the pair formed by \(\limsup\sqrt\kappa(1-\rho_q)\) and
\(\limsup\sqrt\kappa(1-\rhoF)\) is dominated by a point of the arc of \Cref{cor:damping} with \(m=y\).
\item At critical damping, \(\theta=2\sqrt y\), the constant is at most \(\frac23\sqrt y\).
For any fixed finite memory, the slow data \((2,1)\) of standard Nesterov and
\((\frac4{\sqrt3},\frac43)\) of its quadratic-optimal tuning therefore give ceilings
\(\frac23\) and \(\frac4{3\sqrt3}\), respectively (\Cref{cor:B}(i),~(ii)).
\item A family with \(\liminf\sqrt\kappa(1-\rhoF)\ge1\) has \(\theta\ge3\) and \(y\ge\theta-1\ge2\).
\end{enumerate}
\end{corollary}

\subsection{The limiting phase}
\label{sec:limitphase}

For the \emph{limit method} \((u_0,v_0)\), corresponding to
\(\varepsilon=\mu=0\), \(L=1\), write \(\phi_\lambda=\lambda u_0-v_0\).
Thus the subscript in \(\phi_0,\phi_1\) denotes curvature.
Its lag-one symbol and \emph{effective extrapolation} are
\[
\Szero(\omega)=\Rea\Bigl[(1-e^{i\omega})\,\phi_0(e^{i\omega})\,\overline{\phi_1(e^{i\omega})}\Bigr],\qquad
\teff=u_0'(1)-1+y\,q_0'(1),
\]
with \(\phi_0=-v_0\) and \(\phi_1=u_0-v_0\).  The number \(\teff\) depends only on the transfer function
of the limit method:
\[
y\,\frac{v_0(z)}{u_0(z)}=-(z-1)^2\bigl[1-(1+\teff)(z-1)+O((z-1)^2)\bigr].
\]
Thus multiplying \(u,v\) by a common real polynomial \(r\) with \(r(1)=1\)
preserves \(\teff\).
For \(\cM\), \(u_0=N_\tau\) and \(q_0=-1/y\), so \(\teff=\tau\).

\begin{lemma}[Low frequencies of the limit method]\label{lem:low}
As \(\omega\to0\),
\[
\phi_0\overline{\phi_1}\,(e^{i\omega})=-\frac{\omega^2}y+i\,\frac{\teff}y\,\omega^3+O(\omega^4),\qquad
\Szero(\omega)=\frac{2\teff-1}{2y}\,\omega^4+O(\omega^6).
\]
\end{lemma}

\subsection{The half-step loss}
\label{sec:halfstep-mechanism}

On a Fourier mode with \(\omega\ne0\), the discrete difference and continuous derivative satisfy
\begin{equation}\label{eq:halfstep-multiplier}
1-e^{i\omega}=\frac{2\sin(\omega/2)}{\omega}\,e^{i\omega/2}(-i\omega).
\end{equation}
The factor \(e^{i\omega/2}\) adds half a step of phase, which the spectral
lead must compensate.  The flow Lyapunov function in \Cref{thm:flow} is
associated with a Popov-type multiplier~\cite{Pop62};
\(1-z\) is its discrete Jury--Lee counterpart~\cite{JL64}.
\Cref{thm:faithful} turns this phase relation into a uniform rate theorem
under explicit regularity and quadratic-root conditions.

\paragraph{Phase lead.}
Let \(\Theta_0\) be the phase~\eqref{eq:Theta} of the limit method (\(t=0\), \(\mu=0\), \(L=1\)).  By~\eqref{eq:factorgen} it is defined
where \(\phi_0\phi_1\ne0\) and given by \(-\phi_0\overline{\phi_1}=|\phi_0\phi_1|\,e^{-i\Theta_0}\), and
by \Cref{lem:low}, \(\Theta_0(\omega)/\omega\to\teff\) as \(\omega\to0\).
Thus \(\teff\) is the low-frequency spectral phase lead.
For \(\cM\), the phase is
\(\Theta_0(\omega)=\arg\bigl(y-(2-y\tau)\xi+i\,y\tau\sin\omega\bigr)\).
Two interpolation facts make the phase useful.  At \(\omega\in(0,\pi)\), the lag-one symbol is a positive
multiple of \(\sin(\Theta(\omega)-\omega/2)\) (\Cref{lem:angle}).  For stable tunings of \(\cM\),
\(\Theta\in(0,\pi]\), so positivity means \(\Theta>\omega/2\).  For a general stationary method, the sine
criterion is the branch-independent statement.  Second, a one-atom spectrum at \(2\pi/K\) is feasible
at \(t=0\) exactly when \(|\Theta(2\pi/K)|\le\pi/K\); for \(2\pi j/K\), \(j\ge2\), the threshold remains
\(\pi/K\) (\Cref{lem:single}).

The low-frequency wall is therefore \(\teff=1/2\).
If \(\Szero\not\equiv0\), write
\(\Szero(\omega)=s_{2\mathfrak p}\omega^{2\mathfrak p}+O(\omega^{2\mathfrak p+2})\).
By \Cref{lem:low}, \(\mathfrak p\ge2\), and \(\mathfrak p=2\) exactly when \(\teff\ne1/2\).

\subsection{Faithful discretization}
\label{sec:faithful}

A second-order family with slow data \((\theta,y)\) is \emph{faithful to its slow flow} if
\(\lim_{\kappa\to\infty}\sqrt\kappa(1-\rhoF)=\cFF(\theta,y)\).

\begin{samepage}
\begin{theorem}[Faithful discretization and the half-step rule]\label{thm:faithful}
Let a second-order family have slow data \((\theta,y)\), limit method \((u_0,v_0)\), effective
extrapolation \(\teff\), phase \(\Theta_0\) and, if \(\Szero\not\equiv0\), order \(2\mathfrak p\) and leading coefficient \(s_{2\mathfrak p}\) as above.
\begin{enumerate}[label=\textup{(\alph*)}]
\item \textup{(Beyond the wall.)}  If \(\Szero>0\) on \((0,\pi]\), the roots of \(u_0-2v_0\) lie in the open unit
disk, and the family is of class \(C^{2\mathfrak p}\), then
\[
\lim_{\kappa\to\infty}\sqrt\kappa\,(1-\rhoF)=\cFF(\theta,y),
\]
and for every \(c<\cFF(\theta,y)\) the bound \(\rhoF<1-c\varepsilon\), for small \(\varepsilon\), is certified by the
lag-one criterion.
\item \textup{(Below the wall.)}  If \(|\teff|<\frac12\), or if \(\teff=\frac12\) and \(s_{2\mathfrak p}<0\), there is an integer \(K_0\ge3\)
such that \(|\Theta_0(2\pi/K)|<\pi/K\) for every \(K\ge K_0\), and for each such \(K\) and all small \(\varepsilon\) the one-atom
spectrum at frequency \(2\pi/K\) is feasible at some \(t>0\).
\item \textup{(Behind the iterate.)}  If \(\teff<0\), then \(\liminf_{\kappa\to\infty}\rho_q>1\): the method diverges already on a
quadratic on \(\R\).
\item \textup{(Fast resonances.)}  If \(\phi_0\phi_1(e^{i\omega_0})\ne0\) and \(|\Theta_0(\omega_0)|<\pi/K\) at \(\omega_0=2\pi j/K\) for some \(K\ge3\) and \(1\le j<K/2\)
with \(\gcd(j,K)=1\), then for all small \(\varepsilon\) the one-atom spectrum at frequency \(\omega_0\) is feasible at some \(t>0\).
\end{enumerate}
In~(b)--(d), \(\rhoF>1\): the method has growing trajectories on some \(f\in\cF_{\mu,L}(\R^2)\) for all small
\(\varepsilon\).  These parts require coefficient continuity; part~(a) uses the stated higher regularity.
\end{theorem}
\end{samepage}

\begin{corollary}[Necessary lead]\label{cor:lead}
Every second-order family that is stable on quadratics for all small \(\varepsilon\) has \(\teff\ge0\), and every
second-order family with \(\rhoF<1\) for all small \(\varepsilon\) has \(\teff\ge\frac12\).
\end{corollary}

In~(a), \(u_0-2v_0=2\phi_{1/2}\) controls the fast response at middle curvature.
For \(\teff>1/2\), the regularity hypothesis is \(C^4\); at \(\teff=1/2\),
the quartic term of \(\Szero\) vanishes and the first nonzero term fixes the
required order.  Parts~(b),~(c) detect low-frequency divergence under
coefficient continuity, while~(d) detects fixed-frequency resonances.

\begin{remark}[Extremal leads]\label{rem:halfstep}
For the limit of standard Nesterov, \((y,\tau)=(1,1)\),
\(\Theta_0(\omega)=\omega\) for \(0<\omega\le\pi\); for the limit of TMM, \((y,\tau)=(2,\tfrac12)\),
\(\Theta_0(\omega)=\omega/2\) for \(0<\omega<\pi\).
For TMM, \(\phi_0\phi_1(e^{i\pi})=0\) and \(\Szero\equiv0\) on the entire circle.
For every limit method of \(\cM\), \(\Theta_0(\omega)/\omega\to\tau\) as \(\omega\to0\).

The slow part of the lag-one symbol is the flow identity: in the slow variables \(\omega=\varepsilon\gamma\), the leading part
of the symbol on the circle \(|z|=1-c\varepsilon\) is \(\frac{\varepsilon^3}y\bigl(c\Delta_c+(\theta-3c)\gamma^2\bigr)\), with \(\Delta_c=c^2-\theta c+y\)
(\Cref{lem:twoscale}), which is
the slack \((3c-\theta)\gamma^2-c\Delta_c\) of the minimizer constraint~\eqref{eq:P3} for one atom, with the opposite sign.
For \(0<c<\cQ(\theta,y)\), we have \(\Delta_c>0\): the slow symbol is positive
for every \(\gamma\) exactly when \(c\le\theta/3\), which is also exactly when
the minimizer constraint~\eqref{eq:P3} excludes every flow atom.
The certificate is the discrete form of~\eqref{eq:identity}.
\end{remark}

\subsection{A trace formula for the step}
\label{sec:trace}

The slow step \(y=-1/q_0(1)\), equal to \(\lim hL\) for \(\cM\), satisfies
an exact trace formula in the limiting lag-one symbol.  Recall
\(\xi=1-\cos\omega\).

\begin{theorem}[Trace formula]\label{thm:trace}
Let \((u_0,v_0)\) be a stationary linear method with \(v_0=(z-1)^2q_0\) and \(q_0(1)=-1/y<0\), and let
\(\phi_0=-v_0\), \(\phi_1=u_0-v_0\) and \(\Szero\) be as above.
\begin{enumerate}[label=\textup{(\alph*)}]
\item If the roots of \(\phi_1\) lie in the open unit disk, except possibly a simple root at \(-1\), then
\begin{equation}\label{eq:trace}
\frac2y-1=\frac1{2\pi}\int_0^\pi\frac{(2-\xi)\,\Szero(\omega)}{\xi^2\,|\phi_1(e^{i\omega})|^2}\,d\omega ,
\end{equation}
with a bounded integrand, and \(\Szero\equiv0\) if and only if \(v_0/u_0=-(z-1)^2/(3z-1)\).
\item If the roots of \(\lambda u_0-v_0\) lie in the closed unit disk for all \(\lambda\) in an interval \((\lambda_0,1)\), \(0\le\lambda_0<1\), and
\(\Szero\ge0\) on \((0,\pi)\), then \(y\le2\), with equality if and only if \(v_0/u_0=-(z-1)^2/(3z-1)\).
\end{enumerate}
\end{theorem}

\paragraph{Why equality is rigid.}
For the root hypothesis in~(a), set
\[
Q(z)=-\frac{z(1+z)q_0(z)}{\phi_1(z)}.
\]
It is analytic outside the unit disk, with \(Q(\infty)=1\) and \(Q(1)=2/y\).
On the circle, its imaginary part is minus \(\Szero\) times a positive weight.
The Schwarz formula gives~\eqref{eq:trace}.  If \(\Szero\equiv0\), the imaginary part vanishes on the boundary;
analyticity and \(Q(\infty)=1\) force \(Q\equiv1\), hence
\(u_0=(1-3z)q_0\) and \(v_0/u_0=-(z-1)^2/(3z-1)\).
For part~(b), cancel common factors and apply the same formula at curvature \(\lambda\):
\(0\le Q_\lambda(x)-1\le2/(\lambda y)-1\) for every real \(x>1\).
If \(y=2\), letting \(\lambda\uparrow1\) gives \(Q_1(x)=1\) wherever defined,
and hence the same polynomial identity.

For TMM's limit, \((y,\tau)=(2,\frac12)\) in~\eqref{eq:Muv},
\(\Szero\equiv0\) and \(\teff=\frac12\).
The integrand in~\eqref{eq:trace} tends to \(4(2\teff-1)/y\) as \(\omega\to0\).
For a second-order family stable on quadratics, the quadratic-root hypothesis of~(b) holds
(proof of \Cref{cor:yle2}); hence a step \(y>2\) requires \(\Szero<0\) somewhere.
\begin{corollary}[The price of a faithful lead]\label{cor:yle2}
Consider a second-order family with slow data \((\theta,y)\), stable on quadratics for small \(\varepsilon\), whose
limit method has \(\Szero\ge0\) on \((0,\pi)\).
\begin{enumerate}[label=\textup{(\roman*)}]
\item \(y\le2\), hence \(\limsup\sqrt\kappa(1-\rhoF)\le1\) and
\(\limsup\sqrt\kappa(1-\rho_q)\le\sqrt2\).  The quadratic constant is at
most \(1/\sqrt2\) times Polyak's optimal constant \(2\).
\item If \(\liminf\sqrt\kappa(1-\rhoF)\ge1\), then \((\theta,y)=(3,2)\), and the limit method has the transfer
function \(-(z-1)^2/(3z-1)\) of the limit of TMM; in particular \(\teff=\frac12\) and \(\Szero\equiv0\).
\item Under the hypotheses of \Cref{thm:faithful}(a), \(y<2\), and the constant is \(\cFF(\theta,y)\le\sqrt{y/2}<1\).
\end{enumerate}
\end{corollary}

For \(\cM\) the condition \(\Szero\ge0\) is \(\tau\ge\frac12\), \(y(1+2\tau)\le4\) (\Cref{prop:Ltau}(b)), and (i), (ii) reduce to the
bounds behind \Cref{cor:frontier}(ii) and~(v).  Under the nonnegative-symbol hypothesis,
optimality fixes the limiting transfer function for every fixed finite memory.
Extending the nonnegative-symbol condition from low to all frequencies is an open
question (\Cref{sec:conclusions}).
\section{Two-Step Rates: The Wall, Boundary and Tunings}
\label{sec:constant}

For \(\cM\), \(\teff=\tau\), and the lag-one symbol is cubic in
\(\xi=1-\cos\omega\).  This gives explicit rate formulas, including the layer near
\(\tau=1/2\).
\Cref{sec:TMM} treats the joint limit \((y,\tau)=(2,\tfrac12)\), including all tunings with
constant \(1\) and the arithmetic gaps for suboptimal constants.
Throughout this section, \(L=1\), \(\mu=\varepsilon^2\), and \(\Delta_c=c^2-\theta c+y\).

\subsection{The lag-one polynomial and the wall}
\label{sec:Ltau}

\begin{proposition}[Lag-one polynomial]\label{prop:Ltau}
\leavevmode
\begin{enumerate}[label=\textup{(\alph*)}]
\item For a tuning of \(\cM\) with \(0<y_\mu<y\) and \(\omega\in[0,\pi]\),
\[
\Rea\bigl[(1-e^{i\omega})\chi_\mu(e^{i\omega})\overline{\chi_L(e^{i\omega})}\bigr]=\xi\,L_\tau(\xi),\qquad L_\tau(\xi)=\mathsf A_0+\mathsf A_1\xi+\mathsf A_2\xi^2.
\]
The coefficients are the polynomials in \((\beta,\tau,y_\mu,y)\) given in~\eqref{eq:Acoef}; in particular
\(\mathsf A_0=2(1-\beta)(y-y_\mu)+y_\mu y\), and at \(\beta=1\), \(y_\mu=0\), \(L_\tau(\xi)=2y(2\tau-1)\xi+4(1-\tau y)\xi^2\).
For stable tunings \(L_\tau(0)>0\) and \(L_\tau(2)=\chi_\mu(-1)\chi_L(-1)>0\),
and for \(\omega\in(0,\pi]\), \(\Theta(\omega)>\omega/2\) if and only if
\(L_\tau(\xi)>0\), and \(\Theta(\omega)=\omega/2\) if and only if
\(L_\tau(\xi)=0\).  If \(L_\tau>0\) on \((0,2]\), then \(\rhoF<1\).
\item For a family~\eqref{eq:scaling} of \(\cM\), \(\teff=\tau\), and the lag-one symbol of the limit method is
\(\Szero=\frac{\xi^2}{y^2}\bigl(2y(2\tau-1)+4(1-\tau y)\xi\bigr)\).  It is positive on \((0,\pi]\) if and only if
\(\tau\ge\frac12\) and \(y(1+2\tau)<4\), and nonnegative if and only if \(\tau\ge\frac12\) and \(y(1+2\tau)\le4\); if it is positive, the
roots of \(u_0-2v_0\) lie in the open unit disk.
\end{enumerate}
\end{proposition}

At \(t=0\), positivity of this quadratic on \([0,2]\) gives the lag-one criterion;
a negative value at \(\xi=1-\cos(2\pi/K)\) gives phase locking with \(j=1\)
(\Cref{prop:cycles}).  In the notation of \Cref{sec:faithful}, the symbol in~(b) has \(\mathfrak p=2\) for
\(\tau>\frac12\), while at \(\tau=\frac12\) and \(y<2\), \(\Szero=\frac{2(2-y)}{y^2}\xi^3\) has \(\mathfrak p=3\).

\begin{theorem}[The wall \(\tau=\tfrac12\)]\label{thm:wall}
Let stable tunings \((\beta_\kappa,\tau_\kappa,h_\kappa)\) of \(\cM\), \(\tau_\kappa\ge0\), be
accelerated on quadratics, with \(\limsup_{\kappa\to\infty}\tau_\kappa<\frac12\).  Then for
all large \(\kappa\) there is an integer \(K\ge6\) with \(\Theta(2\pi/K)<\pi/K\); the
method has trajectories from rest that approach \(K\)-periodic orbits in
\(\R^3\), and \(\rhoz>1\).  Consequently every family of stable tunings that is
accelerated on \(\cF_{\mu,L}\) has \(\liminf_{\kappa\to\infty}\tau_\kappa\ge\frac12\).
\end{theorem}

At \(\tau=0\), every heavy-ball tuning that accelerates on quadratics has periodic
orbits on \(\cF_{\mu,L}\) for all large \(\kappa\), strengthening the non-acceleration result
of~\cite{GTD}.  TMM approaches the wall from below:
\(\tau_\kappa<\frac12\) and \(\tau_\kappa\to\frac12\) (\Cref{cor:TMM}).
The leading low-frequency coefficient of the limiting lag-one symbol changes sign at
\(\tau=\frac12\) (\Cref{prop:Ltau}(b)), as predicted by the half-step rule;
the resulting resonances occur at intermediate frequencies.

For fixed \(\tau<\tfrac12\), the resonant band extends from frequencies of order
\(\kappa^{-1/4}\) to a fixed frequency; the number and the maximal size of its integer periods are
of order \(\kappa^{1/4}\).  \Cref{prop:periods} gives the precise band and leading constant.

\subsection{The flow constant and the stability boundary}
\label{sec:exact}

The spiral lower witnesses in \Cref{thm:B} are sampled from the flow.  In the slow variables
\(\omega=\varepsilon\gamma\), \(t=-c\varepsilon\), \(q=n\varepsilon\),
the lag constraints have the flow kernels with uniformly convergent coefficients on bounded frequency sets.
Their common strict margin survives sampling.
For every feasible slow spectrum, the value constant is determined up to \(O(\varepsilon^3)\)
(\Cref{lem:reduce,lem:Cforced,lem:flowspirals}).

\begin{proposition}[Single-slow-frequency obstruction]\label{prop:onedisc}
Consider a family~\eqref{eq:scaling} with \(\theta/3<\cQ(\theta,y)\), and let \(c_\sharp\) be as in \Cref{prop:one} with
\(m=y\).  For every \(\gamma_{\max}>0\) there is \(\varepsilon_1>0\) such that for \(\varepsilon\in(0,\varepsilon_1)\) and
\(\theta/3<c\le c_\sharp\), no pair \((\delta_{\varepsilon\gamma},C)\) with \(0\le\gamma\le\gamma_{\max}\) is feasible at \(t=-c\varepsilon\).
\end{proposition}

The frequency range in this proposition is \(\omega=O(\varepsilon)\).

\paragraph{The stability boundary.}
Several natural tunings, among them the tuning of
Nesterov's method that is optimal on quadratics (\Cref{cor:B}(ii)), sit on the stability boundary
\(y(1+2\tau)=4\) in the limit.  For such families measure the stability margin by
\begin{equation}\label{eq:margin}
\varsigma_\varepsilon=\frac{4-y_\varepsilon(1+2\tau_\varepsilon)}{\varepsilon},\qquad\text{so that}\qquad
\chi_L(-1)=2(1+\beta)-y_\varepsilon(1+2\tau_\varepsilon)=(\varsigma_\varepsilon-2\theta_\varepsilon)\,\varepsilon
\end{equation}
exactly.  Thus the tunings are stable for small \(\varepsilon\) if \(\varsigma_\varepsilon\to\varsigma\in(2\theta,\infty]\), and
unstable if \(\varsigma<2\theta\).  For \(\tau\ge\frac12\) and \(\varsigma\in(2\theta,\infty]\), \(\chi_L\) has a real root close to \(-1\): a
sign-alternating mode at curvature \(L\).  Its coefficient \(\cpi\) complements
the slow quadratic coefficient \(\cQ(\theta,y)\) from~\eqref{eq:cQ}:
\begin{equation}\label{eq:cpi}
\cpi=\frac{(\varsigma-2\theta)(1+2\tau)}{4\tau}\in(0,\infty],
\end{equation}
and we put \(\cpi=\infty\) when \(y(1+2\tau)<4\).

\begin{proposition}[On the stability boundary]\label{prop:boundary}
Consider a family~\eqref{eq:scaling} of stable tunings of \(\cM\) with \(\tau>\frac12\), \(y(1+2\tau)=4\) and
\(\varsigma_\varepsilon\to\varsigma\in(2\theta,\infty]\).  Then for every \(c\in(0,\min(\cFF(\theta,y),\cpi))\) there is
\(\varepsilon_0>0\) such that \(\rhoF<1-c\varepsilon\) for \(\varepsilon\in(0,\varepsilon_0)\), certified by the lag-one criterion, and
\[
\lim_{\kappa\to\infty}\sqrt\kappa\,(1-\rhoF)=\lim_{\kappa\to\infty}\sqrt\kappa\,(1-\rhoz)=\min\bigl(\tfrac\theta3,\ \cQ(\theta,y),\ \cpi\bigr).
\]
The term \(\cpi\) is realized from rest by the quadratic \(\frac L2\norm x^2\) on \(\R\).
\end{proposition}

The alternating mode affects the leading constant when the stability margin is of order
\(\kappa^{-1/2}\).  For \(\tau>\frac12\), at \(\varsigma=2\theta\), stable tunings converge by \Cref{thm:tonguefree}
with constant \(0\), since \(\rho_q\ge1-o(\varepsilon)\) (\Cref{lem:altmode}(b)).
At \((\tau,y)=(\tfrac12,2)\), the limiting lag-one symbol vanishes identically.
For \(\tau=\tfrac12\), \(y<2\), its nonzero sixth-order term determines the wall layer.

\subsection{The boundary layer of the wall}
\label{sec:layer}

\Cref{thm:wall,thm:B} leave open the families with \(\tau_\varepsilon\to\frac12\).  The
transition occurs at \(\frac12-\tau_\varepsilon\asymp\kappa^{-1/4}\).  For a
family~\eqref{eq:scaling} with \(\tau=\frac12\), \(y<2\), put
\[
A_\varepsilon=\bigl(\tfrac12-\tau_\varepsilon\bigr)\varepsilon^{-1/2},\qquad\text{i.e.}\quad
\tau_\varepsilon=\tfrac12-A_\varepsilon\kappa^{-1/4},
\]
and, with \(A_+=\max(A,0)\),
\begin{equation}\label{eq:layer}
A_*=\sqrt{\frac{\theta(2-y)}{y}},\qquad
\Lambda(A)=\theta-\frac{yA_+^2}{2-y},\qquad
\clayer(A)=\min\Bigl(\cQ(\theta,y),\frac{\Lambda(A)}3\Bigr).
\end{equation}

\begin{theorem}[The boundary layer]\label{thm:layer}
Consider a family~\eqref{eq:scaling} of stable tunings of \(\cM\) with
\(\tau=\frac12\) and \(y<2\).
\begin{enumerate}[label=\textup{(\alph*)}]
\item \textup{(Upper bound.)}  If \(A_\varepsilon\to A\in[-\infty,A_*)\), then for every
\(c<\clayer(A)\) there is \(\varepsilon_0>0\) such that \(\rhoF<1-c\varepsilon\) for
\(\varepsilon\in(0,\varepsilon_0)\), certified by the lag-one criterion.
\item \textup{(Lower bound.)}  If \(A_\varepsilon\to A\in(0,\infty)\), then for every
\(c>\Lambda(A)_+/3\) and all small \(\varepsilon\) there are \(d\in\N\),
\(f\in\cF_{\varepsilon^2,1}(\R^d)\) and \(x_0\ne x^\star\) with, from rest,
\(\norm{x_N-x^\star}\ge c_0(1-c\varepsilon)^N\norm{x_0-x^\star}\) for all \(N\ge0\),
where \(c_0=c_0(\varepsilon)>0\) is fixed across all horizons \(N\).
\item \textup{(Divergence.)}  If \(\liminf A_\varepsilon>A_*\), then for all small
\(\varepsilon\) there is an integer \(K\ge3\) with \(\Theta(2\pi/K)<\pi/K\): the method has
trajectories from rest that approach \(K\)-periodic orbits in \(\R^3\), and
\(\rhoz>1\).  If \(A_\varepsilon\to A<\infty\), these \(K\) are of order \(\varepsilon^{-1/4}\).
\end{enumerate}
Consequently, if \(A_\varepsilon\to A\in[-\infty,A_*)\), then
\[
\lim_{\kappa\to\infty}\sqrt\kappa\,(1-\rhoF)=\lim_{\kappa\to\infty}\sqrt\kappa\,(1-\rhoz)=\clayer(A),
\]
and if \(A_\varepsilon\to A\in[A_*,\infty]\) the family is not accelerated on \(\cF_{\mu,L}\).  A family with
\(y<2\) and \(\tau_\varepsilon\to\frac12\) is accelerated on \(\cF_{\mu,L}\) if and only if
\(\limsup A_\varepsilon<A_*\).
\end{theorem}

\begin{figure}[t]
\centering
\begin{tikzpicture}
\begin{axis}[name=frontierleft,width=0.48\textwidth,height=0.36\textwidth,
  title={\small (a) Rate trade-off},
  xmin=0.5,xmax=1.75,ymin=0.48,ymax=1.07,xlabel={\(\cQ\)},ylabel={\(\cFF\)},
  xtick={0.5,1,sqrt(2)},xticklabels={\(\frac12\),\(1\),\(\sqrt2\)},
  ytick={0.5,2/3,1},yticklabels={\(\frac12\),\(\frac23\),\(1\)},
  tick label style={font=\footnotesize},label style={font=\small},
  axis lines=left,clip=false]
\addplot[thin,gray,domain=0.5:1.07] {x};
\addplot[thick,variable=t,domain={2*sqrt(2)}:3,samples=80]
  ({(t-sqrt(max(t^2-8,0)))/2},{t/3});
\addplot[thick,dashed,variable=t,domain={4/sqrt(3)}:{sqrt(6)},samples=80]
  ({(t-sqrt(max(t^2-16/3,0)))/2},{t/3});
\addplot[thick,densely dotted,variable=t,domain={4/sqrt(5)}:{3*sqrt(2/5)},samples=80]
  ({(t-sqrt(max(t^2-16/5,0)))/2},{t/3});
\addplot[only marks,mark=*,mark size=1.5pt] coordinates {(1,1) (1,2/3) ({sqrt(2)},{2*sqrt(2)/3})};
\node[font=\scriptsize,anchor=south west] at (axis cs:1,1) {TMM};
\node[font=\scriptsize,anchor=west] at (axis cs:1.02,0.66) {Nesterov};
\node[font=\scriptsize,anchor=north west] at (axis cs:1.40,0.93) {HNAG\(^{++}\)};
\node[font=\scriptsize,anchor=north] at (axis cs:1.23,0.94) {\(\tau\downarrow\frac12\)};
\node[font=\scriptsize,anchor=south west] at (axis cs:1.12,0.79) {\(\tau=1\)};
\node[font=\scriptsize,anchor=north] at (axis cs:0.75,0.585) {\(\tau=2\)};
\end{axis}
\begin{axis}[at={(frontierleft.east)},anchor=west,xshift=0.08\textwidth,
  width=0.48\textwidth,height=0.36\textwidth,
  title={\small (b) Half-step boundary layer},
  xmin=-0.4,xmax=1.75,ymin=0,ymax=0.76,xlabel={\(A\)},ylabel={constant},
  xtick={0,sqrt(2)},xticklabels={\(0\),\(A_*\)},
  ytick={0,1/3,2/3},yticklabels={\(0\),\(\frac13\),\(\frac23\)},
  tick label style={font=\footnotesize},label style={font=\small},
  axis lines=left,clip=true]
\addplot[draw=none,fill=gray!12] coordinates {({sqrt(2)},0) (1.75,0) (1.75,0.76) ({sqrt(2)},0.76)} --cycle;
\addplot[thick,domain=-0.4:{sqrt(2)},samples=100] {(2-max(x,0)^2)/3};
\addplot[thick,dashed,domain=-0.4:{sqrt(2)},samples=100] {min(2/3,(2-max(x,0)^2)/2)};
\addplot[thin,gray,densely dotted] coordinates {({sqrt(2)},0) ({sqrt(2)},0.76)};
\node[font=\scriptsize,anchor=north east] at (axis cs:0.65,0.50) {\(\clayer\)};
\node[font=\scriptsize,anchor=south west] at (axis cs:0.90,0.61) {one atom};
\node[font=\scriptsize,rotate=90] at (axis cs:1.59,0.36) {nonaccelerated};
\end{axis}
\end{tikzpicture}
\caption{Left: the Pareto arcs \(\mathcal C_\tau\) of \Cref{cor:frontier}(iv),
with the diagonal \(\cFF=\cQ\), TMM, standard Nesterov and
HNAG\(^{++}\)~\cite{CX26} (parameters in \Cref{cor:B}(v)).
Right: for \(\theta=2\), \(y=1\), the exact wall-layer constant
\(\clayer(A)=(2-A_+^2)/3\) on \(A<A_*=\sqrt2\) (solid), compared with the one-atom
bound \(\min(\frac23,\Lambda(A)_+/2)\) (dashed).
The shaded region is nonaccelerated (\Cref{thm:layer}).}
\label{fig:frontier}
\end{figure}
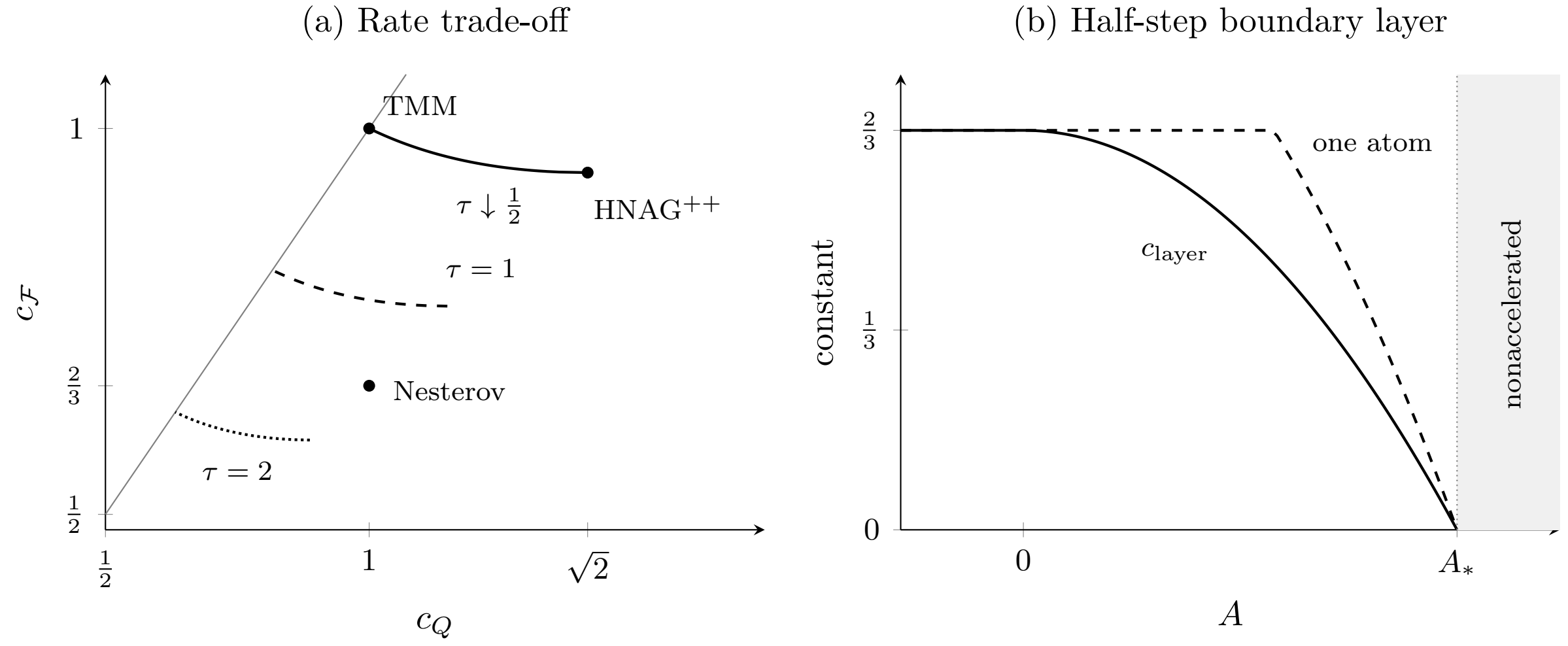

The wall layer has width \(A_*\kappa^{-1/4}\).  Its constant equals the flow constant
\(\cFF(\theta,y)\) whenever \(\Lambda(A)\ge3\cFF(\theta,y)\), in particular for \(A\le0\),
and then decreases to \(0\) at \(A=A_*\) (\Cref{fig:frontier}, right).
The lag-one certificate is exact throughout the accelerated layer.  The mechanism is the lag-one polynomial at
the scale \(\xi\asymp\kappa^{-1/4}\), i.e.\ at frequencies \(\omega\asymp\kappa^{-1/8}\): for the lag-one symbol
\(H_r=\Rea[(1-z)\chi_\mu(z)\overline{\chi_L(z)}]\), \(z=re^{i\omega}\) (\Cref{prop:lagone,lem:cubic}), with
\(\xi=\varepsilon^{1/2}\varphi\) and \(r=1-c\varepsilon\),
\[
\frac{H_r}{\xi}=\varepsilon\bigl[2y(\theta-3c)-4Ay\varphi+2(2-y)\varphi^2\bigr]+o(\varepsilon).
\]
Minimizing the bracket over \(\varphi\ge0\) gives \(3c=\Lambda(A)\): the minimum occurs at a double root
for \(A>0\), and at \(\varphi=0\) for \(A\le0\).  At \(c=0\), this gives \(A=A_*\).
The stabilizing term corresponds to \(s_6\omega^6\) in the limiting symbol,
with \(s_6=(2-y)/(4y^2)\) (\Cref{prop:Ltau}(b)).
The resonant periods \(\kappa^{1/8}\) interpolate between the periods up to
order \(\kappa^{1/4}\) at fixed \(\tau<\frac12\) (\Cref{prop:periods}) and the absence
of resonances beyond the wall.

The witnesses in~(b) spread frequencies over a smooth band around \(2\pi/K\),
\(K\asymp\kappa^{1/8}\), desynchronizing the returns of the planar spirals
(\Cref{lem:spreadhalf}).  At these frequencies the gradient has a half-step lead up to
relative error \(O(\kappa^{-1/4})\).  The spectra asymptotically saturate the lower lag-one
and minimizer constraints, with \(C\) approaching its minimizer bound, and retain strict slack
at the other lags.  This matches the lag-one rate.  Each witness is finite-dimensional;
the number of spirals may grow with \(\kappa\).
The one-atom spectra of \Cref{lem:mixed}, with \(C=0\), give the bound
\(\min(\cFF,\Lambda(A)_+/2)\) (\Cref{rem:onelayer}).

The scales behind the different witnesses are summarized in \Cref{tab:scales}.
For sufficiently smooth families, \Cref{thm:faithful} recovers \Cref{thm:B}(c)
and the layer case \(A=0\);
the direct two-step analysis also covers nonsmooth parameter families and trajectories from rest.

\begin{table}[t]
\centering\small
\begin{tabular}{@{}p{0.23\textwidth}p{0.25\textwidth}p{0.22\textwidth}p{0.20\textwidth}@{}}
\toprule
Regime & Parameter scale & Relevant frequency & Rotation or cycle scale\\
\midrule
Sampled flow & fixed slow data & \(\omega\asymp\kappa^{-1/2}\) & rotation time \(\kappa^{1/2}\)\\
Below the wall & fixed \(\tau<1/2\) & lower edge \(\kappa^{-1/4}\) & maximal \(K\asymp\kappa^{1/4}\)\\
Wall layer, \(y<2\) & \(1/2-\tau_\varepsilon\asymp\kappa^{-1/4}\) & \(\omega\asymp\kappa^{-1/8}\) & \(K\asymp\kappa^{1/8}\)\\
TMM corner & \(1/2-\tau_\varepsilon,\ 2-y_\varepsilon=O(\kappa^{-1/2})\) & fixed \(2\pi/K\) also matter & integer \(K\) leaves a gap\\
\bottomrule
\end{tabular}
\caption{Scales of the witnesses.  In the first row, rotation time is the angular period of a decaying spiral component.}
\label{tab:scales}
\end{table}

\subsection{Known tunings and ceilings}
\label{sec:known}

We apply the rate formulas to standard tunings and optimize at fixed extrapolation.

\begin{corollary}[Known tunings]\label{cor:B}
\leavevmode
\begin{enumerate}[label=\textup{(\roman*)}]
\item \textup{(Nesterov's method.)}  The tuning \(h=1/L\), \(\tau=\beta=\frac{\sqrt\kappa-1}{\sqrt\kappa+1}\)
has \((\theta_\varepsilon,y_\varepsilon,\tau_\varepsilon)\to(2,1,1)\), \(\cQ=1\) and
\[
\rhoF=1-\tfrac23\kappa^{-1/2}+o(\kappa^{-1/2}),
\]
and the same for \(\rhoz\); for squared distances the factor is
\(1-\frac43\kappa^{-1/2}+o(\kappa^{-1/2})\).
\item \textup{(Nesterov's method, quadratic-optimal tuning.)}  The tuning
\(h=4/(3L+\mu)\), \(\tau=\beta=(\sqrt{3\kappa+1}-2)/(\sqrt{3\kappa+1}+2)\), which is optimal for Nesterov's
method on quadratics with rate \(1-2/\sqrt{3\kappa+1}\)~\cite[Proposition~1]{LRP16}, has
\((\theta_\varepsilon,y_\varepsilon,\tau_\varepsilon)\to(\frac4{\sqrt3},\frac43,1)\), on the stability boundary, with margin \(\varsigma=\frac{32}{3\sqrt3}\) and
\(\cpi=\cQ=\frac2{\sqrt3}\).  On \(\cF_{\mu,L}\),
\[
\rhoF=1-\frac4{3\sqrt3}\,\kappa^{-1/2}+o(\kappa^{-1/2}),\qquad \frac4{3\sqrt3}\approx0.770,
\]
and the same for \(\rhoz\).
\item \textup{(Quasi-hyperbolic momentum.)}  For the method of~\cite{MY19} in the form of
\Cref{rem:twostep}, with step \(\alpha>0\), discount \(\beta\in[0,1)\), weight \(\nu_{\rm qh}\in[0,1]\) and \(\tau_{\rm qh}\) as there, rates
refer to its iterates from the usual initialization \(d_{-1}=0\).
\begin{enumerate}[label=\textup{(\alph*)}]
\item If a family of tunings accelerates on quadratics, then \(1-\nu_{\rm qh}=O(\kappa^{-1/2})\); in
particular a fixed \(\nu_{\rm qh}<1\) with \(\beta\to1\) does not accelerate.  If it accelerates on
\(\cF_{\mu,L}\), then moreover \(\liminf\tau_{\rm qh}\ge\frac12\).
\item If \(1-\beta=\theta\kappa^{-1/2}\), \(\alpha L(1-\beta)\to y\) and \(\tau_{\rm qh}\to\tau>\frac12\) with \(y(1+2\tau)<4\), the
constant is \(\min(\theta/3,\cQ(\theta,y))\).
\end{enumerate}
\item \textup{(Over-relaxed heavy ball.)}  The method
\(x_{k+1}=x_k+\beta(x_k-x_{k-1})-\gamma\bigl(2\nabla f(x_k)-\nabla f(x_{k-1})\bigr)\) with \(\gamma=1/(\sqrt L+\sqrt\mu)^2\) and
\(\beta=L\gamma\) has \(\rhoF=1-\frac23\kappa^{-1/2}+o(\kappa^{-1/2})\).
\item \textup{(HNAG\(^{++}\).)}  Put \(\eta=\sqrt{2/\kappa}\).
The gradient-correction form of HNAG\(^{++}\)~\cite{CX26} has
\[
\beta=\frac1{(1+\eta)^2},\qquad hL=\frac{2+\eta}{(1+\eta)^2},\qquad \tau=\frac1{2+\eta}.
\]
For arbitrary initial states, its uniform distance factor satisfies
\(\rhoF=1-\frac{2\sqrt2}{3}\kappa^{-1/2}+o(\kappa^{-1/2})\), and its quadratic factor satisfies
\(\rho_q=1-\sqrt2\kappa^{-1/2}+o(\kappa^{-1/2})\).
\end{enumerate}
\end{corollary}

Part~(i) gives \(\cFF(2,1)=2/3\), recovering Pokutta and Stellato~\cite{PS26} and
Gannot~\cite{Gan22}.  Both Nesterov tunings are critically damped in the slow flow;
their quadratic constant is \(3/2\) times the full-class constant.
In~(ii), the quadratic endpoint modes balance at \(\cQ=\cpi=2/\sqrt3\),
and the spiral term determines the slower full-class rate.
Squared-distance constants are twice the distance constants.

For the default quasi-hyperbolic weight \(\nu_{\rm qh}=0.7\)~\cite{MY19},
\Cref{cor:QHM} gives the finite-condition-number threshold.
The gradient correction in~(iv) is exactly \(\tau=1\) under the same translation, so the over-relaxed
method of Wei and Chen~\cite{WC25} shares the slow data and the constant of standard Nesterov.
For HNAG\(^{++}\), the squared-distance coefficient is therefore \(4\sqrt2/3\), improving the global
coefficient \(\sqrt2\) in~\cite[Theorem~5.1]{CX26}.  Their faster coefficient \(2\sqrt2\) for the auxiliary sequence \(y_k\) is an individual
asymptotic estimate under local asymptotic symmetry: the two first-order Taylor
remainders of \(f\) between \(x\) and \(x^\star\) have ratio tending to \(1\) as
\(x\to x^\star\)~\cite[Theorem~5.2]{CX26};
\Cref{rem:uniform-smoothness} describes the uniform rate convention.
Optimizing the damping for the full class gives the following tuning.

\begin{corollary}[Asymptotically optimal Nesterov tuning]\label{cor:nag-asymptotic-optimal}
For \(\kappa\ge6\), the tuning
\[
h=\frac4{3L},\qquad \tau=\beta=1-\sqrt{\frac6\kappa}
\]
is stable, converges linearly on \(\cF_{\mu,L}\), and satisfies
\[
\rhoF=1-\sqrt{\frac23}\,\kappa^{-1/2}+o(\kappa^{-1/2}),\qquad
\rhoz=1-\sqrt{\frac23}\,\kappa^{-1/2}+o(\kappa^{-1/2}).
\]
The constant \(\sqrt{2/3}\) is maximal among families with \(\tau=\beta\) on \(\cF_{\mu,L}\).
\end{corollary}

\Cref{tab:nag-tunings} compares the three asymptotic Nesterov tunings;
\Cref{fig:levels} shows their geometry.

\begin{table}[t]
\centering
\begin{tabular}{lccc}
\toprule
Nesterov tuning & \(\theta\) & \(y\) & Full-class constant \(c\)\\
\midrule
Standard & \(2\) & \(1\) & \(2/3\)\\
Quadratic-optimal & \(4/\sqrt3\) & \(4/3\) & \(4/(3\sqrt3)\)\\
Class-optimal asymptotically & \(\sqrt6\) & \(4/3\) & \(\sqrt{2/3}\)\\
\bottomrule
\end{tabular}
\caption{Asymptotic distance constants within \(\tau=\beta\), so \(\tau\to1\)
(\Cref{cor:B,cor:nag-asymptotic-optimal}).  Optimizing on quadratics and on the full function class
gives different damping parameters.}
\label{tab:nag-tunings}
\end{table}

The target-rate inequalities in \Cref{sec:quadratic} give the fixed-extrapolation ceiling.
Thus the supremum at fixed \(\tau>1/2\) is \(\sqrt{2/(1+2\tau)}\), reached on the stability boundary
when the sign-alternating mode is no slower.  Every family accelerated on \(\cF_{\mu,L}\)
loses at least a factor \(\sqrt2\) in its quadratic constant relative to Polyak's optimum.
\Cref{cor:frontier} gives the Pareto arcs and the necessary limit
\((\theta,y,\tau)=(3,2,1/2)\) for constant \(1\).

At TMM, \(y_\varepsilon\to2\), so the stabilizing wall-layer coefficient vanishes.
The joint limit \(A_\varepsilon\to0\), \(y_\varepsilon\to2\) motivates the corner scale below.

\section{The Corner of the Triple Momentum Method}
\label{sec:TMM}

The optimal constant \(1\) forces \((\theta_\varepsilon,y_\varepsilon,\tau_\varepsilon)\to(3,2,\tfrac12)\)
(\Cref{cor:frontier}(v)).  We classify all approaches that attain it, first in the
natural corner scale and then at arbitrary rates of approach.  The lag-one
certificate tests all frequencies; the lower-bound construction spreads spectra
around \(2\pi/K\) for integers \(K\ge3\).  For some suboptimal tunings, the resulting
bounds leave an \emph{arithmetic gap}.

We begin with TMM's exact finite-\(\kappa\) certificate.

\subsection{An exact certificate at TMM}
\label{sec:robust}

For TMM, let \(\varrho=1-\kappa^{-1/2}\).  Its endpoint polynomials factor as
\[
\chi_\mu(z)=(z-\varrho)(z-\varrho^2),\qquad \chi_L(z)=z(z+\varrho).
\]
The lag-one symbol vanishes on \(|z|=\varrho\) and is strictly positive on every circle
\(\varrho<r<1\).  Letting \(r\downarrow\varrho\) gives the exact rate.
\Cref{thm:robust} gives the full factorization and the exact rate for the robust momentum family
of~\cite{CHVL18}, including TMM; the rate itself is already known from~\cite{CHVL18,VFL18}.

For robust momentum with design rate \(\varrho=1-c\kappa^{-1/2}\), \(0<c<1\), the limiting data are
\((\theta,y,\tau)=(3c,2c^2,1/(2c^2))\).  They lie on the balance curve \(\theta=3\sqrt{y/2}\),
where the quadratic and spiral bounds agree.

\begin{corollary}[TMM]\label{cor:TMM}
TMM~\cite{VFL18} is the tuning~\eqref{eq:robust} at \(\varrho=1-\kappa^{-1/2}\).  It has \((\theta_\varepsilon,y_\varepsilon,\tau_\varepsilon)\to(3,2,\frac12)\), with
\(2-hL=\kappa^{-1/2}\) and \(\frac12-\tau=\frac54\kappa^{-1/2}(1+O(\kappa^{-1/2}))\), so that \(\tau<\frac12\) for every \(\kappa\), and
\(\rhoF=\rhoz=\rho_q=1-\kappa^{-1/2}\) exactly (the upper bound is~\cite[Theorem~1]{VFL18}).  By~\cite{DT22} this
rate, and the constant \(1=\cFF(3,2)=\cQ(3,2)\), are optimal on \(\cF_{\mu,L}\) among deterministic
first-order methods with value and gradient access.
The HNAG\(^{+}\) tuning of~\cite[Theorem~4.1]{CX26} has the same gradient-correction parameters.
\end{corollary}

\subsection{The corner}
\label{sec:corner}

By \Cref{cor:TMM}, TMM approaches the corner along
\(2-y_\varepsilon=\varepsilon\) and \(\frac12-\tau_\varepsilon=\frac54\varepsilon(1+o(1))\), at the same scale \(\varepsilon=\kappa^{-1/2}\) in both
coordinates.  The natural scaling of the corner is
\begin{equation}\label{eq:corner}
\beta=1-\theta_\varepsilon\varepsilon,\qquad hL=y_\varepsilon=2-b_\varepsilon\varepsilon,\qquad
\tau_\varepsilon=\tfrac12-a_\varepsilon\varepsilon,\qquad
(\theta_\varepsilon,a_\varepsilon,b_\varepsilon)\to(\theta,a,b),
\end{equation}
with \(\theta>0\), \(a,b\in\R\).  Expanding \(y_\varepsilon(1+2\tau_\varepsilon)<2(1+\beta)\) shows that
these tunings are stable for small \(\varepsilon\) if \(2a+b>\theta\) and unstable if
\(2a+b<\theta\); we assume \(2a+b>\theta\).  Put
\begin{equation}\label{eq:cornerN}
\begin{gathered}
\mathsf N(\xi)=4\theta+2(\theta-4a)\,\xi+2(4a+b-2\theta)\,\xi^2,\qquad
\mathsf q(\xi)=\frac{\mathsf N(\xi)}{4(3-\xi)},\qquad
\cpi=2(2a+b-\theta),\\
\clag=\min_{0\le\xi\le2}\mathsf q(\xi),\qquad
\cdisc=\min\bigl\{\mathsf q(\xi):\ \xi\in\{0,2\}\cup\{\xi_K:K\ge3\}\bigr\},\qquad \xi_K=1-\cos\tfrac{2\pi}K .
\end{gathered}
\end{equation}
Here \(\cpi\) is~\eqref{eq:cpi} at \(\tau=\frac12\): the margin~\eqref{eq:margin} of~\eqref{eq:corner} tends to
\(\varsigma=4a+2b\).  We call \(\mathsf q\) the \emph{lag-one ratio}.  Its values at the ends of \([0,2]\) are \(\mathsf q(0)=\theta/3\)
and \(\mathsf q(2)=\cpi\), so \(\clag\le\cdisc\le\min(\theta/3,\cpi)\).

\begin{theorem}[The corner]\label{thm:corner}
Consider a family~\eqref{eq:corner} with \(2a+b>\theta\).
\begin{enumerate}[label=\textup{(\alph*)}]
\item \textup{(Upper bound.)}  For every \(c\) with \(0<c<\min(\cQ(\theta,2),\clag)\) there is \(\varepsilon_0>0\)
such that \(\rhoF<1-c\varepsilon\) for \(\varepsilon\in(0,\varepsilon_0)\), certified by the lag-one
criterion.
\item \textup{(Lower bound.)}  For every \(c>\max\bigl(0,\min(\cQ(\theta,2),\cdisc)\bigr)\) and all small \(\varepsilon\),
\(\rhoz>1-c\varepsilon\).
\item Consequently, if \(\clag>0\),
\[
\min\bigl(\cQ(\theta,2),\clag\bigr)\le\liminf_{\kappa\to\infty}\sqrt\kappa\,(1-\rhoF)\le
\limsup_{\kappa\to\infty}\sqrt\kappa\,(1-\rhoz)\le\min\bigl(\cQ(\theta,2),\cdisc\bigr),
\]
and the constant equals \(\min(\cQ(\theta,2),\clag)\) whenever \(\mathsf q\) attains its minimum over \([0,2]\) at \(0\),
at \(2\) or at a point \(\xi_K\), or \(\cQ(\theta,2)\le\clag\).  In particular, with
\(\ccorner=\min(\theta/3,\cQ(\theta,2),\cpi)\), the constant equals \(\ccorner\) whenever \(\mathsf N(\xi)\ge4\ccorner(3-\xi)\) on \([0,2]\),
for instance if \(4a+b\le2\theta\).
\item \textup{(Optimal tunings.)}  The constant equals \(1\) if and only if \(\theta=3\), \(a\le\frac54\) and
\(2a+b\ge\frac72\).
\end{enumerate}
\end{theorem}

The lower bounds in~(b) have three sources, all realized from rest:
flow or quadratic witnesses for \(\min(\theta/3,\cQ(\theta,2))\) (\Cref{thm:B}(a));
the sign-alternating mode of \(\frac L2\norm x^2\) for \(\cpi\);
and spectra spread around \(2\pi/K\) for \(\mathsf q(\xi_K)\) (\Cref{lem:spreadhalf}).

\begin{remark}[Explicit regions]\label{rem:regions}
Let \(\theta\le3\), so that \(\min(\theta/3,\cQ(\theta,2))=\theta/3\).  By \Cref{thm:corner}(c), the constant
equals \(\theta/3\) if \(a\le\frac{5\theta}{12}\) and \(2a+b\ge\frac{7\theta}6\), and it equals \(\cpi=2(2a+b-\theta)<\theta/3\) if
\(2a+b<\frac{7\theta}6\) and \(5a+2b\le\frac{11\theta}4\) (\Cref{app:corner}).
\end{remark}

\begin{figure}[t]
\centering
\begin{tikzpicture}
\begin{axis}[name=cornerleft,width=0.49\textwidth,height=0.39\textwidth,
  title={\small (a) Optimal corner wedge},
  xmin=-1,xmax=2.6,ymin=-2.4,ymax=5.6,xlabel={\(a\)},ylabel={\(b\)},
  ylabel style={rotate=-90},
  xtick={-1,0,5/4,2},xticklabels={\(-1\),\(0\),\(\frac54\),\(2\)},
  ytick={-2,0,1,3,5},
  tick label style={font=\footnotesize},label style={font=\small},
  axis lines=left,clip=true]
\addplot[draw=none,fill=gray!35] coordinates {(-1,5.5) (1.25,1) (1.25,5.6) (-1,5.6)} --cycle;
\addplot[draw=none,fill=gray!10] coordinates {(1.25,1) (2.6,-1.7) (2.6,5.6) (1.25,5.6)} --cycle;
\addplot[draw=none,pattern=north east lines,pattern color=gray!60] coordinates {(-1,5) (2.6,-2.2) (2.6,-1.7) (-1,5.5)} --cycle;
\addplot[thick,domain=-1:2.6] {3.5-2*x};
\addplot[thick,dashed,domain=-1:2.6] {3-2*x};
\addplot[thick] coordinates {(1.25,1) (1.25,5.6)};
\addplot[only marks,mark=*,mark size=1.8pt] coordinates {(1.25,1)};
\addplot[only marks,mark=o,mark size=2pt] coordinates {(2,1)};
\node[font=\scriptsize,anchor=north east,fill=white,inner sep=1pt] at (axis cs:1.2,0.95) {TMM};
\node[font=\scriptsize] at (axis cs:0.1,4.8) {constant \(1\)};
\node[font=\scriptsize] at (axis cs:1.9,4.1) {suboptimal};
\node[font=\scriptsize,anchor=south] at (axis cs:2,1.15) {example};
\node[font=\scriptsize] at (axis cs:-0.2,-0.8) {unstable};
\end{axis}
\begin{axis}[at={(cornerleft.east)},anchor=west,xshift=0.08\textwidth,
  width=0.49\textwidth,height=0.39\textwidth,
  title={\small (b) Arithmetic gap: \((a,b)=(2,1)\)},
  xmin=0.3,xmax=0.82,ymin=0.8465,ymax=0.8685,
  xlabel={\(\xi\)},ylabel={\(\mathsf q(\xi)\)},
  xtick={0.3765101981,0.5,0.5505102572,0.6909830056},xticklabels={\(\xi_7\),\(\xi_6\),\(\xi_*\),\(\xi_5\)},
  ytick={0.85,0.86},
  tick label style={font=\footnotesize},label style={font=\small},
  axis lines=left,clip=true]
\addplot[draw=none,fill=gray!18] coordinates {(0.3,{3*sqrt(6)-13/2}) (0.82,{3*sqrt(6)-13/2}) (0.82,0.85) (0.3,0.85)} --cycle;
\addplot[thick,domain=0.3:0.82,samples=160] {(6*x^2-10*x+12)/(4*(3-x))};
\addplot[thin,dashed] coordinates {(0.3,0.85) (0.82,0.85)};
\addplot[thin,densely dotted] coordinates {({3-sqrt(6)},0.8465) ({3-sqrt(6)},{3*sqrt(6)-13/2})};
\addplot[only marks,mark=*,mark size=1.8pt] coordinates
{({1-cos(360/7)},{(6*(1-cos(360/7))^2-10*(1-cos(360/7))+12)/(4*(3-(1-cos(360/7))))})
 (0.5,0.85)
 ({(5-sqrt(5))/4},{(6*((5-sqrt(5))/4)^2-10*((5-sqrt(5))/4)+12)/(4*(3-(5-sqrt(5))/4))})};
\addplot[only marks,mark=diamond*,mark size=2.4pt] coordinates {({3-sqrt(6)},{3*sqrt(6)-13/2})};
\node[font=\scriptsize,anchor=south west] at (axis cs:0.56,0.8503) {\(\cdisc\)};
\end{axis}
\end{tikzpicture}
\caption{The corner at \(\theta=3\).
Left: the dark wedge \(a\le5/4\), \(2a+b\ge7/2\) has constant \(1\),
with TMM at its vertex; the dashed line is the stability boundary \(2a+b=3\),
and the hatched strip has \(\cpi<1\).
Right: for the marked example, the continuous minimum
\(\clag=3\sqrt6-13/2\) at \(\xi_*=3-\sqrt6\) (diamond) lies between \(\xi_6\) and \(\xi_5\);
the resonant minimum is \(\cdisc=17/20\) (filled point at \(\xi_6\)).
The shaded band is the gap between the certificate and witness bounds.}
\label{fig:corner}
\end{figure}
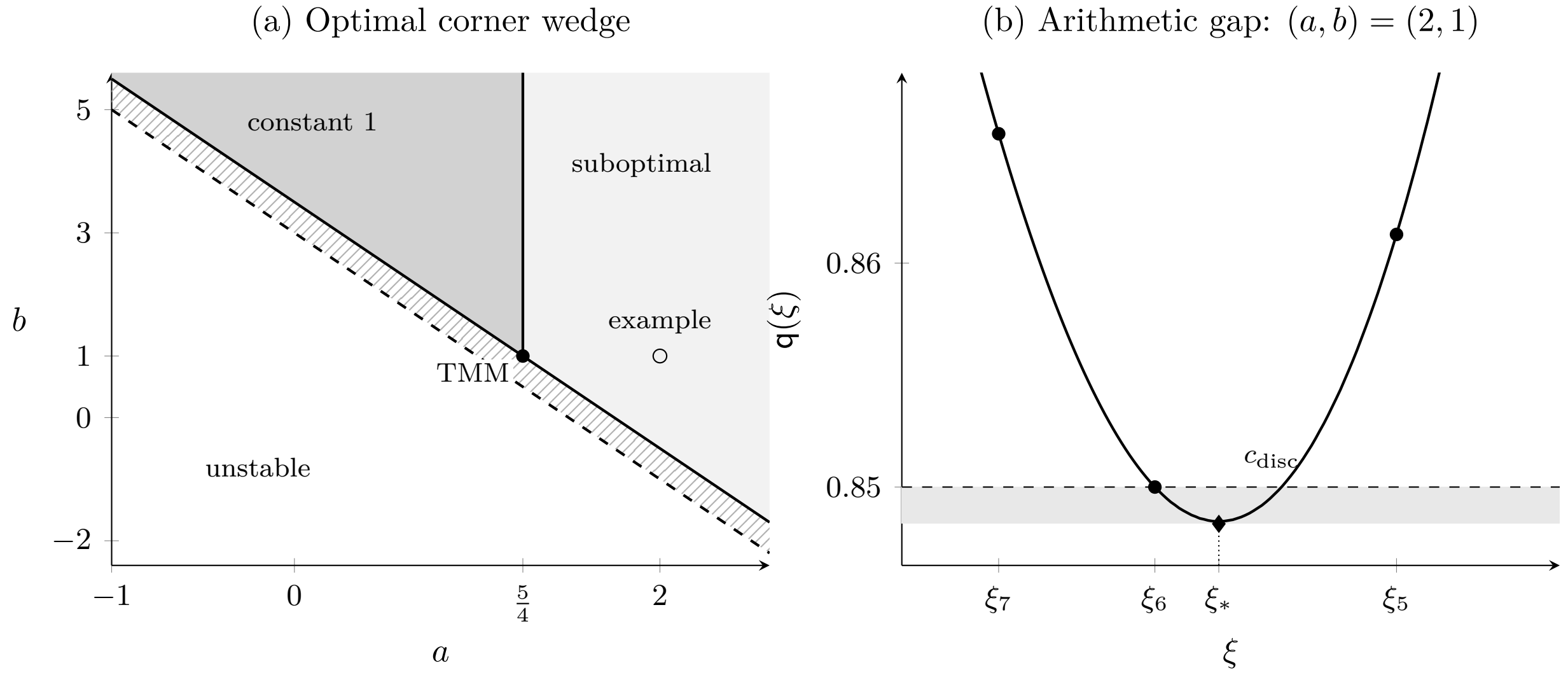

TMM is the vertex \((\theta,a,b)=(3,\frac54,1)\) of the optimal wedge
\(\{a\le\frac54,\ 2a+b\ge\frac72\}\) in~(d) (\Cref{fig:corner}).
It has the largest admissible deficit \(a\) and, at that deficit, the largest leading step
\(y_\varepsilon=2-b\varepsilon+o(\varepsilon)\).
The flow and both quadratic endpoint constants coincide:
\(\theta/3=\cQ(3,2)=\cpi=1\).  Also,
\[
\mathsf N(\xi)=12-4\xi=4(3-\xi),
\]
so \(\mathsf q\equiv1\): to first order in \(\varepsilon\) the lag-one
certificate at rate \(1-\varepsilon\) is tight at every frequency simultaneously.
For \(\theta=3\),
\[
\mathsf N(\xi)-4(3-\xi)=2\xi\bigl[(5-4a)+(4a+b-6)\xi\bigr],
\]
an affine bracket whose values at \(\xi=0\) and \(\xi=2\) are \(5-4a\) and \(2(2a+b)-7\): inside the wedge
\(\mathsf q\ge1\) on \([0,2]\), and beyond its edge \(a=\frac54\) the ratio dips below \(1\) at low frequencies, where the points \(\xi_K\) accumulate.  The limit method of TMM, \((y,\tau)=(2,\frac12)\),
leads the iterate by half a step for \(0<\omega<\pi\) (\Cref{rem:halfstep}).
\Cref{thm:corner} recovers TMM's constant \(1\)~\cite{VFL18}, and \Cref{thm:robust}
gives its rate \(1-\kappa^{-1/2}\) at every \(\kappa\), where the design-circle symbol vanishes identically.  The optimal set also contains families with \(\tau_\varepsilon\ge\frac12\):
these include \(a<0\), \(b\ge\frac72-2a\), and exact-wall families with
\(a=0\), \(b\ge\frac72\).
\Cref{rem:regions} gives the corresponding wedge and alternating-mode region for \(\theta\le3\).

The spread spectra for \(\mathsf q(\xi_K)\) use the boundary-layer construction (\Cref{sec:layer});
the limit method has a half-step lead on \((0,\pi)\).
A gap between the bounds in~(c) requires a minimum of \(\mathsf q\) strictly between consecutive
\(\xi_{K+1}<\xi_K\), below the quadratic contribution.  This occurs for some tunings with \(\theta=3\) and
\(a>\frac54\), where~(b) bounds every limiting leading coefficient strictly below \(1\).  For instance, for \(\theta=3\) and \((a,b)=(2,1)\),
\(\mathsf q(\xi)=(6\xi^2-10\xi+12)/(4(3-\xi))\) is minimal at \(\xi=3-\sqrt6\), strictly between \(\xi_6=\frac12\) and \(\xi_5\), and every limiting leading coefficient lies in
\([3\sqrt6-\frac{13}2,\frac{17}{20}]\approx[0.8485,0.85]\)
(\Cref{fig:corner}, right; see also \Cref{sec:conclusions}).

\begin{corollary}[The arc at \(\tau=\frac12\)]\label{cor:arc}
\leavevmode
\begin{enumerate}[label=\textup{(\roman*)}]
\item The tunings \(\tau=\frac12\), \(\beta=1-2\sqrt2\,\kappa^{-1/2}\),
\(hL=2-4\kappa^{-1/2}\) have constant \(\frac{2\sqrt2}3\) on \(\cF_{\mu,L}\) and \(\sqrt2\) on quadratics; so the bound \(\sqrt2\) of
\Cref{cor:frontier}(ii) is attained.
\item Every point of the arc \(\mathcal C_{1/2}\) of \Cref{cor:frontier}(iv), from \((\sqrt2,\frac{2\sqrt2}3)\) to TMM's point \((1,1)\), is
attained by a family~\eqref{eq:corner} with \(a=0\) and \(b=2\theta\).
\end{enumerate}
\end{corollary}

\subsection{Optimal tunings}
\label{sec:optimal}

\Cref{thm:corner}(d) describes optimal families~\eqref{eq:corner}, with
\(\frac12-\tau_\varepsilon\) and \(2-y_\varepsilon\) both \(O(\kappa^{-1/2})\).
  By \Cref{cor:frontier}(v), every family with constant \(1\)
approaches the corner, possibly at other rates.  For a family~\eqref{eq:scaling} with
\((\theta_\varepsilon,y_\varepsilon,\tau_\varepsilon)\to(\theta,2,\frac12)\) put
\begin{equation}\label{eq:alphavs}
\alpha_\varepsilon=\frac{(\frac12-\tau_\varepsilon)\,y_\varepsilon}{2\varepsilon},\qquad
\varsigma_\varepsilon=\frac{4-y_\varepsilon(1+2\tau_\varepsilon)}{\varepsilon},\qquad
\mathsf N_\varepsilon(\xi)=4\theta_\varepsilon+2(\theta_\varepsilon-4\alpha_\varepsilon)\xi+(\varsigma_\varepsilon+4\alpha_\varepsilon-4\theta_\varepsilon)\xi^2 .
\end{equation}
The tunings are stable if and only if \(\varsigma_\varepsilon>2\theta_\varepsilon\), and for~\eqref{eq:corner}, \(\alpha_\varepsilon\to a\),
\(\varsigma_\varepsilon\to4a+2b\) and \(\mathsf N_\varepsilon\to\mathsf N\).
The criterion below allows unbounded \(\alpha_\varepsilon\) and \(\varsigma_\varepsilon\).

\begin{theorem}[Optimal tunings of \(\cM\)]\label{thm:optimal}
A family of stable tunings of \(\cM\) satisfies \(\lim_{\kappa\to\infty}\sqrt\kappa(1-\rhoF)=1\) if and only if, in the
notation of~\eqref{eq:scaling} and~\eqref{eq:alphavs},
\[
(\theta_\varepsilon,y_\varepsilon,\tau_\varepsilon)\to\bigl(3,2,\tfrac12\bigr),\qquad
\min_{0\le\xi\le2}\frac{\mathsf N_\varepsilon(\xi)}{4(3-\xi)}\to1 .
\]
Under the first condition, put \(d_\varepsilon=(4\alpha_\varepsilon-5)_+\).
The second condition is equivalent to
\[
\liminf_{\varepsilon\to0}\varsigma_\varepsilon\ge7,\qquad
\frac{d_\varepsilon^2}{\varsigma_\varepsilon+5+d_\varepsilon}\to0.
\]
The denominator is positive by stability.  In this case also \(\sqrt\kappa(1-\rhoz)\to1\).
\end{theorem}

For the families~\eqref{eq:corner}, the criterion reduces to
\(\theta=3\), \(4a+2b\ge7\), \(a\le\frac54\), recovering \Cref{thm:corner}(d).  If \(b_\varepsilon=(2-y_\varepsilon)/\varepsilon\to\infty\)
and \(\alpha_\varepsilon\) stays bounded, then the exact identity
\(\varsigma_\varepsilon=2b_\varepsilon+4\alpha_\varepsilon\) gives
\(\varsigma_\varepsilon\to\infty\), and the ratio condition holds: any fixed multiple of
\(\kappa^{-1/2}\) for the deficit \(\frac12-\tau_\varepsilon\) is compatible with constant \(1\).
The edge \(a=\frac54\) and TMM's vertex therefore belong to the scale
\(2-y_\varepsilon=O(\kappa^{-1/2})\).  The proof combines the expansion of the lag-one symbol near the
corner, valid for arbitrary \(\alpha_\varepsilon\) and \(\varsigma_\varepsilon\), with three lower-bound mechanisms: the
sign-alternating mode when \(\liminf\varsigma_\varepsilon<7\), spread spectra (\Cref{lem:spreadhalf}) at frequencies
where the lag-one ratio dips below \(1\), and, when this dip is deep, planar periodic orbits
(\Cref{prop:cycles}).

\section{Finite Condition Numbers}
\label{sec:finite}

At fixed \(\kappa\), the polynomial \(L_\tau\) gives explicit convergence regions.
Call a stable tuning \emph{lag-one convergent} if \(L_\tau>0\) on \((0,2]\).
The periodic tongues of \Cref{prop:cycles,rem:hill} lie outside this region.
The key facts are that \(L_\tau\) is affine in \(y_\mu\) and equals
\(|\chi_L(e^{i\omega})|^2\) at \(y_\mu=y\).

\begin{samepage}
\begin{theorem}[Tongue-free extrapolations]\label{thm:tonguefree}
Let \(\beta\in[0,1)\), \(\tau\ge0\), and
\[
\tau_3(\beta)=\frac{\beta^2+3\beta-1}6,\qquad \tau_*(\beta)=\frac{\beta^2+4\beta-1}8 .
\]
\begin{enumerate}[label=\textup{(\alph*)}]
\item If \(\tau\ge\tau_3(\beta)\), every stable tuning \((\beta,\tau,h)\), for every
\(0<\mu<L\), converges: \(\rhoF<1\).  It is certified by the lag-two certificate with multiplier
\((1-z)+\lambda(1-\bar z^2)\), with \(\lambda=0\) (the lag-one certificate) if
\(\tau\ge\tau_*(\beta)\) and
\(\lambda=\frac{\beta-2\tau}{1+2\beta-\beta^2+2\tau(1+3\beta)}\) if
\(\tau_3(\beta)\le\tau<\tau_*(\beta)\).
\item If \(\tau<\tau_3(\beta)\), there are \(0<\mu<L\) and \(h>0\) such that
\((\beta,\tau,h)\) is stable and lies in the \(3\)-tongue, \(\Theta(2\pi/3)<\pi/3\):
the method has trajectories from rest approaching \(3\)-periodic orbits, and
\(\rhoz>1\).
\item If \(\tau<\tau_*(\beta)\), there are \(0<\mu<L\) and \(h>0\) such that
\((\beta,\tau,h)\) is stable and not lag-one convergent.
\end{enumerate}
\end{theorem}
\end{samepage}

The \(3\)-tongue sets the exact threshold for convergence of every stable tuning.
All stable tunings admit lag-one certificates precisely when
\(\tau\ge\tau_*(\beta)=\tau_3(\beta)+(1-\beta^2)/24\).
Between the thresholds, negative lag-one values occur only in \((2\pi/3,\pi)\);
adding the upper lag-two constraint gives a strict certificate.
Since \(\tau_3<\tau_*<\frac12\) and \(\tau_*(\beta)\le\beta\) on \([0,1)\), every stable
tuning with \(\tau\ge\frac12\), and every stable tuning of the Nesterov family
\(\tau=\beta\), converges linearly on \(\cF_{\mu,L}\), whatever \(h\), \(\mu\) and \(L\); this agrees with the automated
search of~\cite{GDT23}, which found cycles of Nesterov's method on smooth convex functions only on the
boundary of its stability region.
For heavy ball (\(\tau=0\)), every stable tuning converges at every condition
number exactly when \(\beta\le(\sqrt{13}-3)/2\approx0.303\); the lag-one threshold is
\(\beta\le\sqrt5-2\approx0.236\).  The lag-two certificate covers the interval between them.
Earlier sufficient conditions depend on \(h,\mu,L\)~\cite{GFJ15}.
Heavy-ball \(3\)-cycles are described by~\cite[Theorem~3.5]{GTD}, and each rotation-number tongue
by~\cite[Corollary~D.11(b)]{TA}; see also~\cite{GDT23}.
For \(\cM\), all such tongues are absent exactly when \(\tau\ge\tau_3(\beta)\).
\Cref{rem:twostep} transfers this threshold to QHM.

\begin{corollary}[Quasi-hyperbolic momentum]\label{cor:QHM}
For the method of~\cite{MY19} in the form of \Cref{rem:twostep}, with fixed discount \(\beta\in[0,1)\) and weight
\(\nu_{\rm qh}\in[0,1]\), and with the usual initialization \(d_{-1}=0\), every stable choice of the step \(\alpha\) converges linearly on
\(\cF_{\mu,L}\) for all \(0<\mu<L\) if and only if \(\tau_{\rm qh}\ge(\beta^2+3\beta-1)/6\).
\end{corollary}

\begin{figure}[t]
\centering
\begin{tikzpicture}
\begin{axis}[width=0.70\textwidth,height=0.35\textwidth,
  xmin=0,xmax=1,ymin=0,ymax=0.62,xlabel={\(\beta\)},ylabel={\(\tau\)},
  ylabel style={rotate=-90},
  xtick={0,0.5,1},xticklabels={\(0\),\(\frac12\),\(1\)},
  ytick={0,0.25,0.5},yticklabels={\(0\),\(\frac14\),\(\frac12\)},
  tick label style={font=\footnotesize},label style={font=\small},
  axis lines=left,clip=false]
\addplot[draw=none,fill=gray!12,domain={sqrt(5)-2}:1,samples=120] {(x^2+4*x-1)/8} \closedcycle;
\addplot[draw=none,fill=gray!35,domain={(sqrt(13)-3)/2}:1,samples=120] {(x^2+3*x-1)/6} \closedcycle;
\addplot[thick,domain={(sqrt(13)-3)/2}:1,samples=120] {(x^2+3*x-1)/6};
\addplot[thick,dashed,domain={sqrt(5)-2}:1,samples=120] {(x^2+4*x-1)/8};
\addplot[thin,densely dotted,domain=0:1] {0.5};
\addplot[thin,domain=0:0.62] {x};
\addplot[only marks,mark=o,mark size=2pt] coordinates {(1,0.5)};
\node[font=\scriptsize,anchor=west] at (axis cs:1.01,0.5) {TMM limit};
\node[font=\scriptsize,anchor=west] at (axis cs:0.625,0.6) {Nesterov \(\tau=\beta\)};
\node[font=\scriptsize,anchor=south] at (axis cs:0.24,0.5) {half-step wall};
\node[font=\scriptsize] at (axis cs:0.75,0.11) {\(3\)-cycle witnesses};
\node[font=\scriptsize] at (axis cs:0.66,0.43) {lag one};
\node[font=\scriptsize,anchor=south east] at (axis cs:0.86,0.40) {\(\tau_*\)};
\node[font=\scriptsize,anchor=north west] at (axis cs:0.80,0.33) {\(\tau_3\)};
\end{axis}
\end{tikzpicture}
\caption{Regions for convergence of every stable tuning at fixed \((\beta,\tau)\)
(\Cref{thm:tonguefree}).  Lag one suffices for \(\tau\ge\tau_*\) (dashed),
and lag two covers \(\tau_3\le\tau<\tau_*\), with \(\tau_3\) solid.
For \(\tau<\tau_3\), some stable choices of \(h,\mu,L\) admit \(3\)-cycle witnesses.
The open endpoint \((1,\frac12)\) is the limit of TMM's parameters.}
\label{fig:plane}
\end{figure}
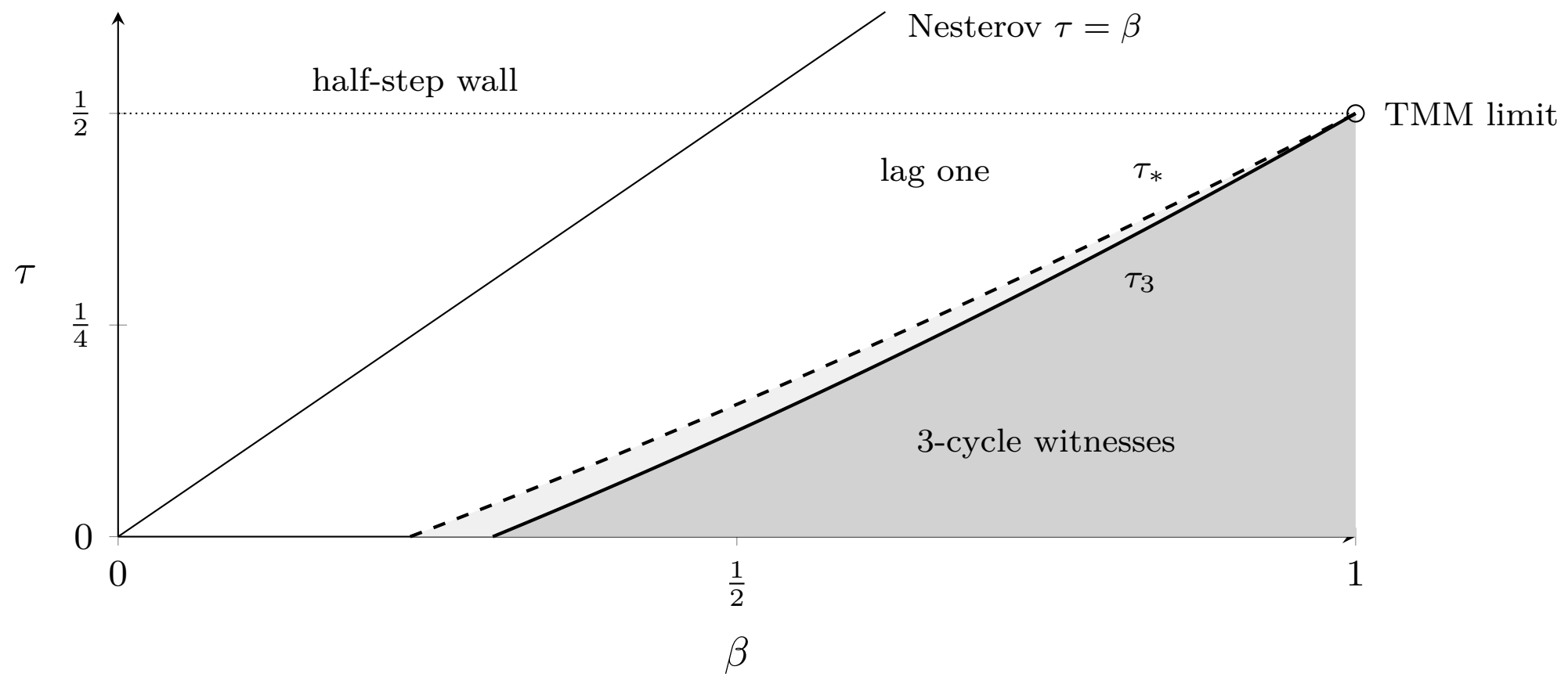

\begin{samepage}
\begin{proposition}[Monotonicity in the condition number]\label{prop:mono}
Fix a tuning \((\beta,\tau,h)\) of \(\cM\) and \(L>0\) with \(hL(1+2\tau)<2(1+\beta)\).  The set
of \(\mu\in(0,L)\) for which \(\rhoF<1\) is an interval \((\mu_\partial,L)\), possibly
\((0,L)\), and it contains the interval \((\mu_*,L)\) on which the tuning is
lag-one convergent.
\end{proposition}
\end{samepage}

Thus uniform linear convergence holds below a critical condition number.
For every fixed stationary method, class inclusion
\(\cF_{\mu',L}\subseteq\cF_{\mu,L}\), \(\mu'\ge\mu\), makes the rate nonincreasing in \(\mu\).
Strict-certificate persistence makes the convergence interval open.
For heavy ball, the interval statement also appears in~\cite[Theorem~3.12]{TA}.

Both thresholds approach the half-step wall as \(\beta\to1\), with
\(\tau_3(1-\theta\varepsilon)=\tfrac12-\tfrac{5\theta}{6}\varepsilon+O(\varepsilon^2)\).
In the corner scale, stable margins reaching the \(3\)-tongue exist for
\(a>5\theta/6\) and are excluded for \(a<5\theta/6\).
At \(\theta=3\), the threshold \(a=5/2\) is twice TMM's deficit.

\section{Conclusions}
\label{sec:conclusions}

The flow ceiling \(\min(\theta/3,\cQ)\) organizes the discrete rates.
For \(\cM\), we determine the constants for \(\tau>1/2\), \(y(1+2\tau)<4\)
and in the wall layer, and classify every family attaining the optimal leading constant \(1\).
For fixed finite memory, positivity of the limiting lag-one symbol gives
faithful discretization under the stated regularity and root conditions.
Quadratic stability and nonnegativity of that symbol force \(y\le2\),
with equality fixing TMM's limiting transfer function.
Spectral closure connects these rates to trajectories generated by one
objective and to finite-lag certificates.
Four directions remain open.
\begin{itemize}
\item \emph{The arithmetic gap at the corner.}  The bounds
\(\min(\cQ,\clag)\) and \(\min(\cQ,\cdisc)\) can differ when the lag-one ratio
\(\mathsf q\) is minimized between consecutive points \(\xi_{K+1}<\xi_K\).
Spectra concentrated in the corresponding frequency gap violate the leading-order
lag-\(K+1\) constraint.  After radial reweighting, the limiting interpolation constraints at the corner
define the cone of finite nonnegative measures \(\nu\) on \((0,\pi]\) satisfying
\(\int\sin\frac{n\omega}2\sin\frac{(n\mp1)\omega}2\,d\nu\ge0\) for every \(n\ge1\);
its Dirac measures are precisely those at \(2\pi/K\), \(K\ge2\).
Does linear-fractional optimization over this cone give the exact constant in
the gaps, and how does it depend on the ordering of rational frequencies with bounded denominators?
\item \emph{Memory \(D\ge3\).}  Does acceleration of a second-order family
force \(\Szero\ge0\) at every frequency?  The unresolved case has a negative
symbol somewhere and \(|\Theta_0(2\pi j/K)|\ge\pi/K\) at every rational
frequency where the phase is defined.  An affirmative answer would extend \(y\le2\) and TMM rigidity
in \Cref{cor:yle2} to every accelerated second-order family.
How does the constant depend on the approach to a zero-symbol limit,
\(\Szero\equiv0\)?
\item \emph{The size of the worst case.}  For \(\theta/3<\cQ\), must the
number of spirals in \Cref{thm:B}(a) grow as \(\delta\to0\)?  Can the
dimension in \Cref{thm:A}(e) be bounded by the memory alone?
At \(T=0>\log\rho_q\), can the attained worst case in \(L^2(\nu;\C)\)
always be realized by a finite spectrum?
\item \emph{Periodic and time-varying methods.}  Which parts of the flow
ceiling, phase rule and optimal-constant classification extend to these methods?
\end{itemize}

\section*{Use of AI tools}
The authors used Claude Opus 5.5 and ChatGPT-6.1 Sol
to assist with exposition, formal proof verification and literature checks.
The authors verified the statements and take responsibility for the results
and accompanying code.

\appendix

\section{Elementary Rate and Quadratic Facts}
\label{app:prelim}

\begin{proof}[Proof of \Cref{rem:uniform-smoothness}]
Fix a stationary method \((u,v)\) and a finite horizon \(n\ge1\).
Translate the minimizer to \(0\) and restart at the state used in the ratio defining \(\Gamma(n)\).
For a trajectory in an arbitrary real Hilbert space, let
\[
V=\operatorname{span}\{w_0,\ldots,w_{D-1},
\nabla f(\tilde x_0),\ldots,\nabla f(\tilde x_{n-1})\}.
\]
This space is finite-dimensional.  Since \(v_D\ne0\), the recursion~\eqref{eq:slm} shows inductively
that \(w_0,\ldots,w_{n+D-1}\) and all queries used to compute them lie in \(V\).
The restriction \(f|_V\) belongs to \(\cF_{\mu,L}(V)\), has minimizer \(0\), and its gradient at each
of these queries equals the original gradient.  Thus the restriction preserves the finite trajectory
and its ratio of state norms.

Identify \(V\) with \(\R^d\).  Let \(\eta\ge0\) be a smooth compactly supported probability density
on the unit ball, and set \(f_\delta(x)=\int f(x-\delta z)\eta(z)\,dz\).
Integrating the strong-convexity inequality and the gradient Lipschitz bound gives
\(f_\delta\in C^\infty\cap\cF_{\mu,L}(\R^d)\), with
\(\sup_x\norm{\nabla f_\delta(x)-\nabla f(x)}\le L\delta\).
If \(a_\delta\) is its minimizer, strong monotonicity gives
\(\norm{a_\delta}\le\norm{\nabla f_\delta(0)}/\mu\le L\delta/\mu\).
Hence the recentered function
\(g_\delta(x)=f_\delta(x+a_\delta)-f_\delta(a_\delta)\) has minimizer \(0\), belongs to the same smooth
class, and satisfies
\[
\sup_x\norm{\nabla g_\delta(x)-\nabla f(x)}
\le L\delta\bigl(1+L/\mu\bigr)\longrightarrow0.
\]
Run the method on \(g_\delta\) from the same initial coordinates
\(w_0,\ldots,w_{D-1}\) of~\eqref{eq:slm}.
Induction in the explicit recursion~\eqref{eq:slm}, using \(v_D\ne0\) and the common Lipschitz bound,
shows convergence of every iterate up to \(w_{n+D-1}\).
The initial state is unchanged and nonzero, so the state-norm ratio converges to the original ratio.
Every ratio in \(\Gamma(n)\) is therefore approximated by finite-dimensional \(C^\infty\) objectives;
the reverse inequality follows by class inclusion.  The case \(n=0\) is immediate.
Equality at every finite horizon gives equality of the corresponding uniform rates.
For \(\cM\), the same approximation preserves the initialization \(x_{-1}=x_0\), so it also applies
to each finite-horizon supremum defining \(\rhoz\).
\end{proof}

\begin{proof}[Proof of \Cref{lem:envelope}]
Center at \(x^\star=0\).  The explicit recursion~\eqref{eq:slm} and
\(\norm{\nabla f(x)}\le L\norm x\) give
\[
\norm{w_{k+D}}\le
\frac{\sum_{j<D}|v_j|+L\sum_{j<D}|u_j|}{|v_D|}\norm{S_k}.
\]
Thus \(\Gamma(n)\le C_\Gamma^n<\infty\), with \(C_\Gamma\) independent of the objective and dimension.
Restarting at \(S_m\) gives
\(\norm{S_{m+n}}\le\Gamma(n)\norm{S_m}\le\Gamma(n)\Gamma(m)\norm{S_0}\),
including zero states, which are absorbing.
Taking suprema proves submultiplicativity; the multiplicative Fekete lemma gives the root limit and infimum formula.

For a root \(v(z)=\lambda u(z)\), \(\lambda\in[\mu,L]\), the real planar quadratic
\(f(w)=\lambda|w|^2/2\) has the trajectory \(w_k=z^k\), hence
\(\norm{S_n}=|z|^n\norm{S_0}\).  Taking the largest root gives \(\Gamma(n)\ge\rho_q^n\).
For \(\cM\) started from rest, \(\norm{S_0}=\sqrt2\norm{x_0}\), so
\(\norm{x_N}\le\sqrt2\Gamma(N)\norm{x_0}\); taking root limits gives \(\rhoz\le\rhoF\).
\end{proof}

\begin{proof}[Proof of \Cref{lem:qstab}]
A real monic quadratic \(z^2-Pz+Q\) has both roots in the open unit disk if and
only if \(|Q|<1\), \(1+Q-P>0\) and \(1+Q+P>0\) (the Schur--Cohn, or Jury,
conditions).  For \(\chi_\lambda\) we have \(P=1+\beta-h\lambda(1+\tau)\) and
\(Q=\beta-h\lambda\tau\le\beta<1\).  Further \(1+Q-P=\chi_\lambda(1)=h\lambda>0\), and
\(1+Q+P=\chi_\lambda(-1)=2(1+\beta)-h\lambda(1+2\tau)\), which decreases in \(\lambda\); it is
positive for all \(\lambda\in[\mu,L]\) if and only if \(y(1+2\tau)<2(1+\beta)\).  Under this
condition also \(h\lambda\tau\le y\tau<1+\beta\), i.e.\ \(Q>-1\).
\end{proof}

\begin{proof}[Proof of identities~\eqref{eq:factor} and~\eqref{eq:imab}]
Substitution of~\eqref{eq:Muv} in~\eqref{eq:p} gives
\[
u-p=-\frac{2\chi_\mu}{h(L-\mu)},\qquad
u+p=\frac{2\chi_L}{h(L-\mu)}.
\]
Their product \((u-p)(\bar u+\bar p)=|u|^2-|p|^2-2i\Ima(p\bar u)\) gives~\eqref{eq:factor};
the roots of \(p\) are those of \((\chi_\mu+\chi_L)/2=\chi_{(L+\mu)/2}\).
For \(\Omega=(z-1)(z-\beta)\) and \(N=N_\tau\),
\(\Ima(\chi_\mu\overline{\chi_L})=(y-y_\mu)\Ima(\Omega\bar N)\).
On \(|z|=1\), expansion with \(z\bar z=1\) gives
\(\Ima(\Omega\bar N)=\sin\omega[(1-\beta)+2\tau\xi]\), proving~\eqref{eq:imab}.
\end{proof}

\begin{proof}[Proof of \Cref{thm:price}]
Let \(\rho=\rho_q\), and let \(w_1,w_2\) and \(v_1,v_2\) be the roots of \(\chi_\mu\) and
\(\chi_L\).  Then \(h\mu=\chi_\mu(1)=(1-w_1)(1-w_2)\ge(1-\rho)^2\): for real roots both
factors are at least \(1-\rho\), for complex ones \(|1-w_1|^2\ge(1-|w_1|)^2\).
Similarly \(\chi_L(-1)=(1+v_1)(1+v_2)\ge(1-\rho)^2\), and \(\beta-h\mu\tau=w_1w_2\le\rho^2\).
Since \(\chi_L(-1)=2(1+\beta)-hL(1+2\tau)\),
\(hL(1+2\tau)\le2(1+\beta)-(1-\rho)^2\le(1+\rho)^2+2h\mu\tau\).  Dividing by
\(h\mu\ge(1-\rho)^2\) gives \(\kappa(1+2\tau)\le(\frac{1+\rho}{1-\rho})^2+2\tau\), i.e.\ the
claim.  Equality forces equality in each step: \(w_1=w_2=\rho\) (for complex roots
\(|1-w|>1-|w|\)), \(v_1=v_2=-\rho\), so \(\beta-h\mu\tau=\rho^2=\beta-hL\tau\) and \(\tau=0\);
then \(\chi_\mu=(z-\rho)^2\) and \(\chi_L=(z+\rho)^2\) give \(\beta=\rho^2\),
\(h\mu=(1-\rho)^2\), \(hL=(1+\rho)^2\), which is Polyak's tuning.
\end{proof}
\begin{proof}[Proof of \Cref{lem:cq}]
In this proof write \(\beta=\beta_\varepsilon\), \(\tau=\tau_\varepsilon\), and let
\(P_\ell(z)=(z-1)(z-\beta)+\ell N_\tau(z)\), so
\(\chi_\lambda=P_{h\lambda}\) and \(h\lambda\in[y_\varepsilon\varepsilon^2,y_\varepsilon]\).
(i) At \(\lambda=\mu\), \(\chi_\mu(1-\varepsilon\zeta)=\varepsilon^2(\zeta^2-\theta'\zeta+y_\varepsilon)\) exactly,
with \(\theta'=\theta_\varepsilon+y_\varepsilon(1+\tau_\varepsilon)\varepsilon\), because
\((z-1)(z-\beta)=\varepsilon^2\zeta(\zeta-\theta_\varepsilon)\) and
\(N_\tau(1-\varepsilon\zeta)=1-(1+\tau)\varepsilon\zeta\).  So the roots at \(\mu\) are \(1-\varepsilon\zeta\), \(\zeta\) a
root of \(\zeta^2-\theta'\zeta+y_\varepsilon\); if real, the largest is
\(1-\varepsilon\zeta_-(\theta',y_\varepsilon)\), and if complex both have modulus
\(\sqrt{\beta-y_\mu\tau}=1-\frac{\theta_\varepsilon}2\varepsilon+O(\varepsilon^2)\).  In both cases the largest root
modulus at \(\mu\) is \(1-\varepsilon\cQ(\theta',y_\varepsilon)+O(\varepsilon^2)\) (if the roots are real,
\(\zeta_-(\theta',y_\varepsilon)\le\theta'/2\)), so \(\rho_q\ge1-\varepsilon\cQ(\theta',y_\varepsilon)+O(\varepsilon^2)\); as \(\cQ\) is continuous,
\(\theta'\to\theta\) and \(y_\varepsilon\to y\), this proves the first claim.  For the second claim, assume the strict inequality in the statement.
(ii) Complex roots at any \(\ell\) have modulus \(\sqrt{\beta-\ell\tau}\le\sqrt\beta=1-\frac{\theta_\varepsilon}2\varepsilon+O(\varepsilon^2)\).
(iii) The discriminant of \(P_\ell\) is a convex quadratic in \(\ell\); its smaller
zero is \(\ell_a=O(\varepsilon^2)\) and its larger zero is \(4/(1+\tau_\varepsilon)^2+o(1)\).  On
\([y_\varepsilon\varepsilon^2,\ell_a]\) (empty if the roots at \(\mu\) are complex) the roots are real and positive, the larger root \(z_+\) exceeds
\(\tau/(1+\tau)\) (it is close to \(1\)), so \(N_\tau(z_+)>0\) and, in the interior,
\(dz_+/d\ell=-N_\tau(z_+)/P_\ell'(z_+)<0\); by continuity the maximum over this interval is at
\(\ell=y_\varepsilon\varepsilon^2\).
(iv) Fix \(\ell_*>0\) below the larger zero of the discriminant.  On \([\ell_a,\ell_*]\)
the roots are complex.  On \([\ell_*,y_\varepsilon]\), real roots \(w\) satisfy \(|w|<1\) and
\(P_\ell(1)=\ell\ge\ell_*\), \(P_\ell(-1)\ge2(1+\beta)-y_\varepsilon(1+2\tau)\ge\delta_0>0\) for small
\(\varepsilon\) (by this strict inequality); since \(P_\ell(\pm1)=(1\mp w)(1\mp w')\le2(1\mp w)\) for the
other root \(w'\in[-1,1]\), real roots lie in \([-1+\delta_1,1-\delta_1]\),
\(\delta_1=\frac12\min(\ell_*,\delta_0)\).
Combining: the root moduli on \([\ell_*,y_\varepsilon]\) are at most \(1-\delta_1\) or
\(\sqrt\beta\); on \([\ell_a,\ell_*]\) they are at most \(\sqrt\beta\); and on \([y_\varepsilon\varepsilon^2,\ell_a]\) they are at
most \(1-\varepsilon\zeta_-(\theta',y_\varepsilon)\), attained at \(\ell=y_\varepsilon\varepsilon^2\).  If the roots at \(\mu\) are complex,
\(\rho_q=\sqrt{\beta-y_\mu\tau}=1-\frac{\theta_\varepsilon}2\varepsilon+O(\varepsilon^2)\); if they are real,
\(\zeta_-(\theta',y_\varepsilon)\le\theta'/2\) and \(\rho_q=1-\varepsilon\zeta_-(\theta',y_\varepsilon)+O(\varepsilon^2)\).  In both cases
\(1-\rho_q=\varepsilon\cQ(\theta',y_\varepsilon)+O(\varepsilon^2)\), and \(\cQ\) is continuous, with \(\theta'\to\theta\) and
\(y_\varepsilon\to y\).
\end{proof}

\begin{proof}[Proof of the claims in \Cref{rem:twostep}]
Applying the operator \(N_\tau\) to~\eqref{eq:M} shows that the queries satisfy the gradient-correction
recursion.  Conversely, if \(h_1>h_2\ge0\), \(h=h_1-h_2\), \(\tau=h_2/(h_1-h_2)\) and \(\beta\ne\tau/(1+\tau)\), each trajectory of the
gradient-correction method is the query sequence of exactly one trajectory of~\eqref{eq:M}.
The two descriptions have the same rates.  Indeed, \(\bigl(\tilde x_{k-1},\,\tilde x_k+(1+\tau)h\nabla f(\tilde x_{k-1})\bigr)\) is the image of \((x_{k-2},x_{k-1})\)
under a fixed linear map with determinant \((1+\tau)\beta-\tau\); as \(\norm{\nabla f(x)}\le L\norm x\), the
states \((x_{k-2},x_{k-1})\) and \((\tilde x_{k-1},\tilde x_k)\) bound each other up to factors independent of \(f\)
when \(\beta\ne\tau/(1+\tau)\), and the worst-case
rate of the queries, defined as in \Cref{def:rate} with the state \((\tilde x_{k-1},\tilde x_k)\), equals
\(\rhoF\).  From rest, \(x_{-1}=x_0\), the queries start with \(\tilde x_0=x_0\) and
\(\tilde x_1=\tilde x_0-h_1\nabla f(\tilde x_0)\), and with \(z_\tau=\tau/(1+\tau)<1\), inverting \(\tilde x_k=(1+\tau)x_k-\tau x_{k-1}\) gives
\[
x_k=\frac1{1+\tau}\sum_{j=0}^{k}z_\tau^{\,j}\,\tilde x_{k-j}+z_\tau^{\,k+1}x_{-1};
\]
hence \(\rho^0_{\tilde x}\le\rhoz\le\max(\rho^0_{\tilde x},z_\tau)\).  The two rates are equal whenever they exceed \(z_\tau\), in
particular in every accelerated regime, where the rates tend to \(1\) while \(\tau\) stays bounded: by \Cref{thm:price},
\(c\kappa^{-1/2}\le1-\rho_q\le2\kappa_\tau^{-1/2}\), and \(\kappa_\tau\ge(1+2\tau)(\kappa-1)\).

Quasi-hyperbolic momentum satisfies the gradient-correction recursion for \(k\ge1\) with
\(h_1=\alpha(1-\nu_{\rm qh}\beta)\) and \(h_2=\alpha(1-\nu_{\rm qh})\beta\).
It therefore has the query representation~\eqref{eq:M} with
\(h=\alpha(1-\beta)\) and \(\tau=\tau_{\rm qh}\).  With \(d_{-1}=0\) its first step is \(x_1=x_0-h_1\nabla f(x_0)\), so its iterates are
exactly the queries of~\eqref{eq:M} started from rest at \(x_0\).
\end{proof}
\section{Proofs for Section~\ref{sec:flow}}
\label{app:flow}

We use the following level-set and scaling identities for \(\cQ\) and \(\cFF\).

\begin{lemma}[The quadratic constant]\label{lem:cQprop}
Let \(\theta,m>0\) and \(\Delta_c=c^2-\theta c+m\).
\begin{enumerate}[label=\textup{(\alph*)}]
\item For \(c\in\R\), \(c<\cQ(\theta,m)\) if and only if \(c<\theta/2\) and \(\Delta_c>0\).
\item For \(\alpha>0\), \(\cQ(\alpha\theta,\alpha^2m)=\alpha\,\cQ(\theta,m)\) and \(\cFF(\alpha\theta,\alpha^2m)=\alpha\,\cFF(\theta,m)\).
\end{enumerate}
\end{lemma}

\begin{proof}
(a) If \(\theta^2<4m\), then \(\Delta_c>0\) for every \(c\) and
\(\cQ=\theta/2\).  Otherwise \(\Delta_c\) has roots
\(\zeta_-\le\theta/2\le\zeta_+\), so \(\Delta_c>0\) together with
\(c<\theta/2\) is equivalent to \(c<\zeta_-=\cQ\).
(b) Rescaling \((\theta,m,c)\) to \((\alpha\theta,\alpha^2m,\alpha c)\)
preserves the two inequalities in (a), and therefore rescales \(\cQ\)
by \(\alpha\).  The same holds for \(\cFF=\min(\theta/3,\cQ)\).
\end{proof}

The witnesses of \Cref{thm:flow}(c) are orthogonal sums of self-similar spirals.  For
\(c>0\) and \(\gamma\ge0\) put
\begin{equation}\label{eq:AB}
\cA_\gamma(q)=\ch cq-\cos\gamma q\ge0,\qquad
\cB_\gamma(q)=\frac{\gamma\sh cq}{c}-\sin\gamma q\ge0\qquad(q>0).
\end{equation}
Let \(\nu=\sum_{j=1}^J\nu_j\delta_{\gamma_j}\) be a probability measure on \([0,\infty)\), let
\(z_j=-c+i\gamma_j\), and in \(\C^J\cong\R^{2J}\) let
\begin{equation}\label{eq:cspiral}
x(t)=\bigl(\sqrt{\nu_j}\,e^{z_jt}\bigr)_j,\qquad G(t)=\bigl(w_j\sqrt{\nu_j}\,e^{z_jt}\bigr)_j,\qquad
w_j=-z_j^2-\theta z_j-m,
\end{equation}
\begin{equation}\label{eq:cvalue}
\eta(t)=e^{-2ct}\sum_j\nu_j\eta_j,\qquad \eta_j=\frac{(\theta-c)\gamma_j^2-c\Delta_c}{2c}.
\end{equation}
Then \(\ddot x+\theta\dot x+mx+G=0\) and \(\dot\eta=\ip G{\dot x}\): the triples \((x(t),G(t),\eta(t))\) are the
candidate values of \((x,\nabla f_m(x),f_m(x))\) along an orbit of~\eqref{eq:flow}.  Consider the
program \(\mathcal P_c\) for such measures:
\begin{gather}
\textstyle\sum_j\nu_j\bigl[(\gamma_j^2-\Delta_c)\cA_{\gamma_j}(q)\mp(\theta-2c)\gamma_j\cB_{\gamma_j}(q)\bigr]\ge0\qquad(q>0),\label{eq:P12}\\
\textstyle(3c-\theta)\sum_j\nu_j\gamma_j^2\ge c\,\Delta_c .\label{eq:P3}
\end{gather}

\begin{lemma}[Spirals of the flow]\label{lem:creduce}
For \(t\ge0\) and \(q>0\),
\begin{align*}
\eta(t+q)-\eta(t)-\ip{G(t)}{x(t+q)-x(t)}&=e^{-2ct-cq}\textstyle\sum_j\nu_j\bigl[(\gamma_j^2-\Delta_c)\cA_{\gamma_j}(q)-(\theta-2c)\gamma_j\cB_{\gamma_j}(q)\bigr],\\
\eta(t)-\eta(t+q)-\ip{G(t+q)}{x(t)-x(t+q)}&=e^{-2ct-cq}\textstyle\sum_j\nu_j\bigl[(\gamma_j^2-\Delta_c)\cA_{\gamma_j}(q)+(\theta-2c)\gamma_j\cB_{\gamma_j}(q)\bigr],\\
\norm{x(t+q)-x(t)}^2&=2e^{-2ct-cq}\textstyle\sum_j\nu_j\cA_{\gamma_j}(q),
\end{align*}
and \(\ip{G(t)}{x(t)}-\eta(t)=e^{-2ct}\sum_j\nu_j\bigl((3c-\theta)\gamma_j^2-c\Delta_c\bigr)/(2c)\).  Hence, if
\(c\le\theta/2\), there is a convex \(f_m\) on \(\C^J\) with \(f_m(0)=0\), \(0\in\partial f_m(0)\), \(f_m(x(t))=\eta(t)\) and
\(G(t)\in\partial f_m(x(t))\) for all \(t\ge0\) if and only if \(\nu\) satisfies~\eqref{eq:P12}--\eqref{eq:P3}.
\end{lemma}

\begin{proof}[Proof of \Cref{lem:creduce}]
Orthogonality makes every quantity linear in the weights \(\nu_j\),
so it suffices to take one atom.  Write
\(z=-c+i\gamma\) and
\(w=-z^2-\theta z-m=(\gamma^2-\Delta_c)+i(2c-\theta)\gamma\).
Then
\[
\ip{G(t)}{\dot x(t)}
=e^{-2ct}\Rea(w\bar z)
=e^{-2ct}\bigl(c\Delta_c+(c-\theta)\gamma^2\bigr)=\dot\eta(t).
\]
For the forward lag, multiply
\(\eta(t+q)-\eta(t)-\ip{G(t)}{x(t+q)-x(t)}\)
by \(e^{2ct+cq}\).  The result is
\[
-2\eta_1\sh cq+(\gamma^2-\Delta_c)(e^{cq}-\cos\gamma q)
-(2c-\theta)\gamma\sin\gamma q
=(\gamma^2-\Delta_c)\cA_\gamma(q)
-(\theta-2c)\gamma\cB_\gamma(q),
\]
using \(2c\eta_1=(\theta-c)\gamma^2-c\Delta_c\).
Reversing the pair changes the sign of the \(\cB_\gamma\)-term.
The distance identity follows from
\(|e^{zq}-1|^2=2e^{-cq}(\ch cq-\cos\gamma q)\);
subtracting \(\eta(t)\) from
\(\ip{G(t)}{x(t)}=e^{-2ct}(\gamma^2-\Delta_c)\)
gives the last identity.

Convex interpolation of values and subgradients is equivalent to
\(\eta_b\ge\eta_a+\ip{G_a}{x_b-x_a}\) for every pair.
Sufficiency follows by taking the supremum of these affine supports;
here the orbit data are bounded, so the supremum is finite everywhere.
Along the orbit, these inequalities are the two lag inequalities just
computed; with the minimizer datum they add
\(\eta\ge0\) and \(\ip Gx-\eta\ge0\).
The latter is~\eqref{eq:P3}, and for \(c\le\theta/2\) it implies
the former because
\[
\eta(t)-(\ip{G(t)}{x(t)}-\eta(t))
=e^{-2ct}\frac{\theta-2c}{c}\sum_j\nu_j\gamma_j^2\ge0.
\]
\end{proof}

\begin{lemma}[Smooth extension with a margin]\label{lem:extend}
Let \(\{(x_a,G_a,\eta_a)\}_{a\in I}\) be data in \(H\times H\times\R\) containing \((0,0,0)\), and suppose that
for all \(a,b\in I\), with \(\Delta x=x_b-x_a\),
\[
\eta_b-\eta_a-\ip{G_a}{\Delta x}\ge\epsilon_0\norm{\Delta x}^2,\qquad \norm{G_b-G_a}\le K\norm{\Delta x},
\]
for some \(\epsilon_0>0\), \(K<\infty\).  Then there is a convex \(\tilde f_m\in C^{1,1}(H)\) with
\(\mathrm{Lip}(\nabla\tilde f_m)\le K^2/(2\epsilon_0)\), \(\tilde f_m(x_a)=\eta_a\) and \(\nabla\tilde f_m(x_a)=G_a\) for all \(a\), and, for every \(m>0\),
\(f=\tilde f_m+\frac m2\norm\cdot^2\) is \(m\)-strongly convex with
\(\mathrm{Lip}(\nabla f)\le m+K^2/(2\epsilon_0)\).
In particular it belongs to \(\cF_{m,L}(H)\) for every \(L>m\) at least this large.
\end{lemma}

\begin{proof}
If \(K\le0\), the gradient data are constant, and adding the two margin inequalities for a pair forces
\(x_a=x_b\); since the data contain \((0,0,0)\), the zero function is an
interpolant.  Suppose \(K>0\).  By~\cite[Theorem~1.3]{AM17}, a convex \(C^{1,1}\) interpolant with \(\mathrm{Lip}(\nabla\tilde f_m)\le M\) exists
if and only if \(\eta_b\ge\eta_a+\ip{G_a}{x_b-x_a}+\frac1{2M}\norm{G_b-G_a}^2\) for all \(a,b\).  With
\(M=K^2/(2\epsilon_0)\) the last term is at most \(\epsilon_0\norm{\Delta x}^2\).
\end{proof}

\begin{lemma}[Spread spectra]\label{lem:spread}
Let \(\theta,m>0\), \(\theta/3<c<\cQ(\theta,m)\), \(\Delta_c=c^2-\theta c+m\), and \(W_\gamma=(\gamma^2-\Delta_c)\cA_\gamma-(\theta-2c)\gamma\cB_\gamma\).
There are \(\gamma_0\ge1\), \(J\in\N\) and \(\epsilon_*\in(0,1]\) such that, for
\(\gamma_i=\gamma_0\lambda_i\), \(\lambda_i=1+(i-\frac12)/J\), \(\nu_i=1/J\) \((1\le i\le J)\), and all \(q>0\),
\begin{align}
\textstyle\sum_i\nu_iW_{\gamma_i}(q)&\ge\epsilon_*\textstyle\sum_i\nu_i\bigl(\cA_{\gamma_i}(q)+\cB_{\gamma_i}(q)\bigr),\label{eq:S1}\\
\textstyle\sum_i\nu_i\bigl[(\gamma_i^2-\Delta_c)\cA_{\gamma_i}+(\theta-2c)\gamma_i\cB_{\gamma_i}\bigr]&\ge\epsilon_*\textstyle\sum_i\nu_i(\cA_{\gamma_i}+\cB_{\gamma_i}),\label{eq:S2}\\
(3c-\theta)\textstyle\sum_i\nu_i\gamma_i^2-c\Delta_c&\ge1.\label{eq:S3}
\end{align}
\end{lemma}

\begin{proof}
We separate short lags, a compact band where the spread suppresses
returns, and long lags where hyperbolic growth dominates.  The frequency
scale \(\gamma_0\) is chosen first, and the quadrature size \(J\) last.
Write \(c=\theta/3+\delta\), \(\delta>0\).  By \Cref{lem:cQprop}(a), \(\Delta_c>0\) and \(c<\theta/2\), so
\(\theta-2c=c-3\delta>0\).  Put \(a=3\delta/c\), \(\mathcal D(q)=e^{-cq}+a\sh cq\) and
\(\mathcal D_\delta=\min(e^{-1},a\sh1)\le\inf \mathcal D\).  The identity
\begin{equation}\label{eq:Lid}
W_\gamma(q)=\gamma^2\bigl(e^{-cq}-\cos\gamma q\bigr)+a\gamma^2\sh cq-\Delta_c\cA_\gamma(q)+(c-3\delta)\gamma\sin\gamma q
\end{equation}
follows from \(\ch x-\frac{c-3\delta}c\sh x=e^{-x}+a\sh x\).  All atoms satisfy
\(\gamma_0\le\gamma\le2\gamma_0\).  Let \(\mathcal U(x)=\int_1^2\lambda^2(1-\cos\lambda x)\,d\lambda\),
\(\mathcal V(x)=\int_1^2\lambda^2\cos\lambda x\,d\lambda\), and \(\mathcal U_J\), \(\mathcal V_J\) their midpoint sums
\(\frac1J\sum_i\lambda_i^2(1-\cos\lambda_ix)\), \(\frac1J\sum_i\lambda_i^2\cos\lambda_ix\); then
\(|\mathcal U-\mathcal U_J|,|\mathcal V-\mathcal V_J|\le(1+x)^2/(6J^2)\) (the second derivatives in \(\lambda\) are bounded
by \(4(1+x)^2\)) and, integrating by parts once, \(|\mathcal V(x)|\le8/x\).  Fix
\(q_I=\frac12\), \(q_{II}=\max(1,32/\mathcal D_\delta)\), \(\mathcal U_{\min}=\min_{[q_I,q_{II}]}\mathcal U>0\) and
\(q_{III}=c^{-1}\arsinh(2/a)\).

\emph{Regime I: \(q\le q_I/\gamma_0\).}  Assume \(\gamma_0\ge\max(2c,\sqrt{12\Delta_c})\).  Then
\(\gamma q\le1\), \(cq\le\frac14\), and \(1-\cos x\ge\frac{x^2}2(1-\frac{x^2}{12})\),
\(\ch x-1\ge\frac{x^2}2\), \(\sh x-x\le\frac{x^3}6\ch x\) give
\(\cA_\gamma\ge\frac{11}{24}(c^2+\gamma^2)q^2\) and
\(\cB_\gamma\le\frac13\gamma(c^2+\gamma^2)q^3\), while \(\cA_\gamma,\cB_\gamma\le(c^2+\gamma^2)q^2\).  Using
\(\gamma^2-\Delta_c\ge\frac{11}{12}\gamma^2\) and \((c-3\delta)q/3\le c/(6\gamma_0)\le\frac1{12}\),
\[
W_\gamma(q)\ge(c^2+\gamma^2)q^2\bigl[\tfrac{11}{12}\cdot\tfrac{11}{24}-\tfrac1{12}\bigr]\gamma^2
\ge\tfrac13\gamma^2(c^2+\gamma^2)q^2\ge\tfrac{\gamma_0^2}6(\cA_\gamma+\cB_\gamma)
\]
atom by atom.

\emph{Regime II: \(q_I/\gamma_0\le q\le q_{II}/\gamma_0\), \(x=\gamma_0 q\).}  Assume \(\gamma_0\ge cq_{II}\).
By~\eqref{eq:Lid}, using \(\gamma^2(1-e^{-cq})\le4\gamma_0^2cq\), \(\cA\le3\) and
\(|\gamma\sin\gamma q|\le2\gamma_0\),
\(\sum_i\nu_iW_{\gamma_i}\ge \gamma_0^2\mathcal U_J(x)-4c\gamma_0x-3\Delta_c-2c\gamma_0\).  With \(J\) so large
that \((1+q_{II})^2/(6J^2)\le \mathcal U_{\min}/2\) and \(\gamma_0\) so large that
\(\gamma_0^2\mathcal U_{\min}/4\ge(4cq_{II}+2c)\gamma_0+3\Delta_c\), this is at least \(\gamma_0^2\mathcal U_{\min}/4\), while
\(\sum_i\nu_i(\cA+\cB)\le4+4q_{II}\) (as \(\gamma\sh(cq)/c\le2\gamma_0 q\ch1\le4x\)).

\emph{Regime III: \(q\ge q_{II}/\gamma_0\).}  By~\eqref{eq:Lid},
\(\sum_i\nu_iW_{\gamma_i}=\gamma_0^2[m_2\mathcal D(q)-\mathcal V_J(\gamma_0 q)]-\Delta_c\sum_i\nu_i\cA_{\gamma_i}+(c-3\delta)\sum_i\nu_i\gamma_i\sin\gamma_iq\),
\(m_2=\frac1J\sum_i\lambda_i^2\in[1,4]\).  If \(q\ge q_{III}\), then \(\mathcal D\ge2\) and \(|\mathcal V_J|\le m_2\),
so \(m_2\mathcal D-\mathcal V_J\ge m_2\mathcal D/2\ge \mathcal D/2\).  If \(q_{II}/\gamma_0\le q\le q_{III}\), take \(J\) with
\((1+\gamma_0 q_{III})^2/(6J^2)\le \mathcal D_\delta/4\); then
\(|\mathcal V_J(\gamma_0 q)|\le8/q_{II}+\mathcal D_\delta/4\le \mathcal D_\delta/2\le \mathcal D/2\), and again \(m_2\mathcal D-\mathcal V_J\ge \mathcal D/2\).
Since \(\ch x\le1+\sh x\), \(\mathcal D\ge \mathcal D_\delta\) and \(\mathcal D\ge a\sh cq\),
\(\Delta_c\sum_i\nu_i\cA_{\gamma_i}+2c\gamma_0\le(2\Delta_c+2c\gamma_0)+\Delta_c\sh cq
\le[(2\Delta_c+2c\gamma_0)/\mathcal D_\delta+\Delta_c/a]\mathcal D\le\frac{\gamma_0^2}4\mathcal D\) for \(\gamma_0\) large; so
\(\sum_i\nu_iW_{\gamma_i}\ge\frac{\gamma_0^2}4\mathcal D\), while
\(\sum_i\nu_i(\cA+\cB)\le3+(1+2\gamma_0/c)\sh cq\le[3/\mathcal D_\delta+(1+2\gamma_0/c)/a]\mathcal D\).

Put
\[
K_3=\frac3{\mathcal D_\delta}+\frac{1+2\gamma_0/c}{a},\qquad
\epsilon_*=
\min\left\{1,\frac{\gamma_0^2}6,
\frac{\gamma_0^2\mathcal U_{\min}}{16(1+q_{II})},
\frac{\gamma_0^2}{4K_3},(\theta-2c)\gamma_0\right\}>0.
\]
All lower bounds on \(\gamma_0\) above are independent of \(J\).
Impose also
\(\gamma_0^2\ge\max\{\Delta_c+1,(c\Delta_c+1)/(3\delta)\}\),
then choose \(J\) to satisfy the two midpoint-error bounds.
The three regimes prove~\eqref{eq:S1} with this common \(\epsilon_*\).
For~\eqref{eq:S2}, each atom has \(\gamma^2-\Delta_c\ge1\) and
\((\theta-2c)\gamma\ge(\theta-2c)\gamma_0\), while
\(\cA,\cB\ge0\).  Finally,
\((3c-\theta)\sum_i\nu_i\gamma_i^2-c\Delta_c
=3\delta\gamma_0^2m_2-c\Delta_c\ge1\), proving~\eqref{eq:S3}.
\end{proof}

\begin{proof}[Proof of \Cref{thm:flow}]
(a) Put \(v=\dot x\), \(G=\nabla f_m(x)\), and
\(a=\theta-c\).  Since \(\dot v=-\theta v-mx-G\),
\[
\dot V_c=\ip Gv+\Delta_c\ip xv+\ip{v+ax}{-cv-mx-G}.
\]
The gradient--velocity terms cancel.  Using
\(\Delta_c=m-ca\) and adding \(2cV_c\) yields
\[
\dot V_c+2cV_c
=-(\theta-2c)\Delta_c\norm x^2
-(\theta-c)(\ip Gx-f_m)-(\theta-3c)f_m,
\]
which is~\eqref{eq:identity}.  Under the stated conditions all three
terms on the right are nonpositive by~\eqref{eq:hineq}.  Moreover
\[
V_c\ge\tfrac{\Delta_c}{2}\norm x^2
+\tfrac12\norm{v+(\theta-c)x}^2
\ge K_c^{-1}(\norm x^2+\norm v^2),
\]
with \(K_c\) depending only on \((\theta,m,c)\), proving (a).

(b) The solutions on \(\frac m2x^2\) are combinations of \(e^{st}\), \(s^2+\theta s+m=0\), with \(te^{st}\) at a
double root.  The root with the larger real part has real part \(-\cQ(\theta,m)\), and the
solution with this mode alone (or its real part) has \(\ell=\cQ(\theta,m)\).

(c) Let \(\theta/3<c<\cQ(\theta,m)\).  \Cref{lem:spread} gives \(\gamma_0\ge1\), \(J\), \(\epsilon_*\in(0,1]\) and atoms
\(\gamma_j\in[\gamma_0,2\gamma_0]\) with equal weights such that, for all \(q>0\),
\[
\textstyle\sum_j\nu_j\bigl[(\gamma_j^2-\Delta_c)\cA_{\gamma_j}\mp(\theta-2c)\gamma_j\cB_{\gamma_j}\bigr](q)\ge\epsilon_*\sum_j\nu_j\cA_{\gamma_j}(q),\qquad
(3c-\theta)\textstyle\sum_j\nu_j\gamma_j^2-c\Delta_c\ge1 ,
\]
where we dropped the terms \(\epsilon_*\cB_{\gamma_j}\ge0\).  Use these atoms in~\eqref{eq:cspiral}--\eqref{eq:cvalue}.  By
\Cref{lem:creduce}, for two times \(t<t+q\) both lag quantities are at least
\(e^{-2ct-cq}\epsilon_*\sum_j\nu_j\cA_{\gamma_j}(q)=\frac{\epsilon_*}2\norm{\Delta x}^2\).  For a time and the minimizer,
\(\ip{G(t)}{x(t)}-\eta(t)\ge\frac1{2c}e^{-2ct}=\frac1{2c}\norm{x(t)}^2\), and \(\eta(t)\ge\ip{G(t)}{x(t)}-\eta(t)\), as in
the proof of \Cref{lem:creduce}, since \(c<\cQ\le\theta/2\).  Moreover
\(\norm{G(t')-G(t)}\le K\norm{x(t')-x(t)}\) and \(\norm{G(t)}\le K\norm{x(t)}\) with \(K=\max_j|w_j|\), since \(G\)
acts on each coordinate as multiplication by \(w_j\).  \Cref{lem:extend} with
\(\epsilon_0=\min(\frac{\epsilon_*}2,\frac1{2c})\) gives \(f\in\cF_{m,L}(\R^{2J})\) with \(\nabla f(x(t))=mx(t)+G(t)\) for
\(t\ge0\).  Since \(\ddot x+\theta\dot x+mx+G=0\), the curve \(x\) solves~\eqref{eq:flow}, and
\(\norm{x(t)}=e^{-ct}\), \(\norm{\dot x(t)}^2=e^{-2ct}\sum_j\nu_j|z_j|^2\).

(d) Let \(0<c<\cFF(\theta,m)\).  Then \(c<\theta/3\), and \(\Delta_c>0\) by \Cref{lem:cQprop}(a), as \(c<\cQ\).  So (a) gives \(\ell\ge c\) for every solution, and \(\inf\ell\ge\cFF\).
Conversely (b) gives \(\ell=\cQ\), and (c) gives \(\ell=c\) for every \(c\in(\theta/3,\cQ)\); hence \(\inf\ell\le\cFF\),
and the witnesses lie in \(\bigcup_L\cF_{m,L}\).  For the energy and the gap: by (a),
\(f(x)=f_m+\frac m2\norm x^2\le(1+\frac m{\Delta_c})V_c\), and \(\frac12\norm{\dot x}^2\le K_ce^{-2ct}V_c(0)\); along the spirals
\(f(x(t))=e^{-2ct}(\sum_j\nu_j\eta_j+\frac m2)\) and the energy is a fixed multiple of \(e^{-2ct}\); on the quadratic
\(\frac m2x^2\) the slowest mode gives \(2\cQ\).  If \(\theta<3\sqrt{m/2}\), then \(\Delta_{\theta/3}=m-\frac29\theta^2>0\), and (a)
applies with \(c=\theta/3\).
\end{proof}

\begin{proof}[Proof of \Cref{cor:damping}]
If \(\theta\le2\sqrt m\), then \(\cQ=\theta/2\le\sqrt m\) and \(\cFF=\theta/3\le\frac23\sqrt m<\sqrt{m/2}\).  On \(\theta\ge2\sqrt m\),
\(\cQ=\zeta_-=2m/(\theta+\sqrt{\theta^2-4m})\) decreases from \(\sqrt m\) and \(\theta/3\) increases; they are equal if and
only if \(\frac{\theta^2}9-\frac{\theta^2}3+m=0\), i.e.\ \(\theta=3\sqrt{m/2}\), with common value \(\sqrt{m/2}\).  This gives both
maxima.  For the domination: if \(\theta\le2\sqrt m\), \((\cQ,\cFF)=(\theta/2,\theta/3)\) is dominated by
\((\sqrt m,\frac23\sqrt m)\), the endpoint \(\theta=2\sqrt m\) of the arc; for \(2\sqrt m\le\theta\le3\sqrt{m/2}\) the point
\((\zeta_-,\theta/3)\) lies on the arc; for \(\theta\ge3\sqrt{m/2}\), \(\cQ=\cFF=\zeta_-\le\sqrt{m/2}\), dominated by the other
endpoint \((\sqrt{m/2},\sqrt{m/2})\).
\end{proof}

\begin{proof}[Proof of \Cref{prop:one}]
Take \(c_\sharp\in(\theta/3,\cQ)\); then every \(c\in(\theta/3,c_\sharp]\) has \(c<\cQ\), hence
\(\Delta_c>0\) and \(\theta-2c>0\) by \Cref{lem:cQprop}(a).  The atom \(\gamma=0\) violates~\eqref{eq:P3} because \(\Delta_c>0\).  For \(\gamma>0\),
\eqref{eq:P3} gives \(\gamma^2\ge c\Delta_c/(3c-\theta)\).  At \(q=2\pi/\gamma\) the
inequality with the minus sign in~\eqref{eq:P12} reads
\((\gamma^2-\Delta_c)(\ch v-1)\ge(\theta-2c)\gamma^2\sh v/c\), \(v=2\pi c/\gamma\); it fails if
\(\gamma^2\le\Delta_c\), and otherwise, since \((\ch v-1)/\sh v=\tanh(v/2)\le v/2\), it
forces \(\pi c/\gamma\ge(\theta-2c)/c\), i.e.\ \(\gamma\le\pi c^2/(\theta-2c)\).  As
\(c\downarrow\theta/3\), \(c\Delta_c/(3c-\theta)\to\infty\) (because
\(\Delta_{\theta/3}>0\) when \(\theta/3<\cQ\), by \Cref{lem:cQprop}(a)) while \(\pi c^2/(\theta-2c)\to\pi\theta/3\); choose
\(c_\sharp\) so that the first bound exceeds the square of the second on
\((\theta/3,c_\sharp]\).

For the consequence, suppose such a solution exists.  Replacing \(f\) by \(f(r\,\cdot)/r^2\in\cFm(H)\)
and \(x\) by \(x/r\), \(r=\norm a\), we may assume \(\norm a=\norm b=1\), and we identify the plane \(P\) spanned
by \(a,b\) with \(\C\) so that \(x(t)=e^{zt}\), \(z=-c+i\gamma\).  Then
\(\nabla f_m(x(t))=-\ddot x-\theta\dot x-mx=G(t)\) with \(\nu=\delta_\gamma\) in~\eqref{eq:cspiral}; it lies in \(P\), so it is also the
gradient of the restriction of \(f_m\) to \(P\), a convex function with \(f_m(0)=0\) and \(\nabla f_m(0)=0\).
Since \(\frac d{dt}f_m(x(t))=\ip{G}{\dot x}=\dot\eta\) and both \(f_m(x(t))\) and \(\eta(t)\) tend to \(0\), \(f_m(x(t))=\eta(t)\).  As \(c\le\theta/2\),
\Cref{lem:creduce} shows that \(\delta_\gamma\) satisfies~\eqref{eq:P12}--\eqref{eq:P3}, a contradiction.
\end{proof}

\begin{remark}[The single-spiral cone]\label{rem:cone}
Let \(0<c<\theta/2\).  For one atom \(\gamma>0\) with \(\gamma^2>\Delta_c\), put
\(\hat\gamma=\gamma/c\) and \(\tan\phi=(\theta-2c)\gamma/(\gamma^2-\Delta_c)\), the (unsigned) angle between \(G\) and \(x\).
The inequality with the minus sign in~\eqref{eq:P12} and~\eqref{eq:P3} read
\[
\tan\phi\le\Upsilon(\hat\gamma):=\inf_{x>0}\frac{\ch x-\cos\hat\gamma x}{\hat\gamma\sh x-\sin\hat\gamma x},\qquad
\hat\gamma\tan\phi\le1 .
\]
A spiral with gradient twist \(\phi\) and winding ratio \(\hat\gamma\) is
admissible if and only if both hold.  When \(\theta/3<\cQ\), approaching
\(c=\theta/3\) from above requires
\(\hat\gamma\to\infty\) by the second inequality, while the first fails at
the first return \(q=2\pi/\gamma\) of a fast winding.  Indeed, for one atom the inequality with the minus sign
in~\eqref{eq:P12}, divided by \((\gamma^2-\Delta_c)\cB_\gamma(q)\), with \(x=cq\), reads
\((\ch x-\cos\hat\gamma x)/(\hat\gamma\sh x-\sin\hat\gamma x)\ge\tan\phi\) for all \(x>0\), and
\eqref{eq:P3} reads \((\theta-2c)\gamma^2\le c(\gamma^2-\Delta_c)\), i.e.\
\(\hat\gamma\tan\phi\le1\); the inequality with the plus sign holds because \(\cA_\gamma,\cB_\gamma\ge0\).
\end{remark}

\section{Proofs for Section~\ref{sec:closure}}
\label{app:closure}

Throughout this appendix, \(x^\star=0\).  We compute the wave interpolation
conditions, prove duality and the one-exponent principle, and construct the
realizations.  The first computation extends~\cite[Lemma~A.3]{TA} to arbitrary
polynomial amplitudes \(u,p\).

\begin{lemma}[Exact spectral interpolation]\label{lem:waves}
Let \(\nu\) be a probability measure on \([0,\pi]\), \(C\in\R\), \(t\in\R\), and
let \(H_\nu=L^2([0,\pi],\nu;\C)\) with \(\ip v{v'}=\Rea\int v\,\overline{v'}\,d\nu\).
For \(z=e^{t+i\omega}\) and complex polynomials (or constants) \(u,p\), put
\(w_k(\omega)=z^k\), \(\tilde x_k=u(z)z^k\), \(\sigma_k=p(z)z^k\) and
\(\psi_k=\frac14\bigl(C+2\int\Rea(p\bar u)\,d\nu\bigr)e^{2kt}\).  Then, with
\(g,s,a_n,b_n\) as in~\eqref{eq:gs}--\eqref{eq:anbn}, for \(j\ge0\), \(n\ge1\),
\begin{equation}\label{eq:slacks}
R_{j,j+n}\pm T_{j,j+n}=\tfrac12e^{(2j+n)t}\bigl(\bar a_n\mp\bar b_n\pm C\sh nt\bigr),\qquad
R_{j\star}\pm T_{\star j}=\tfrac14e^{2jt}(\bar g\pm C).
\end{equation}
Hence the data \((\tilde x_k,\sigma_k,\psi_k)_{k\ge0}\) together with \((0,0,0)\) are
interpolated by some \(\psi\in\cS(H_\nu)\) if and only if \((\nu,C)\) is
feasible at \(t\).
\end{lemma}

\begin{proof}
All expressions are \(\nu\)-averages of the one-frequency expressions,
so first take \(\nu=\delta_\omega\), and put \(m=j+n\).
The identities
\[
|z^n-1|^2=2e^{nt}(\ch nt-\cos n\omega),\qquad
(1+z^n)(\bar z^n-1)=(e^{2nt}-1)-2ie^{nt}\sin n\omega
\]
give
\[
R_{jm}=\tfrac12e^{(2j+n)t}a_n,\qquad
\ip{\sigma_j+\sigma_m}{\tilde x_m-\tilde x_j}
=e^{2jt}\bigl[\Rea(p\bar u)(e^{2nt}-1)
+2s e^{nt}\sin n\omega\bigr].
\]
Subtracting half of the second expression from
\(\psi_m-\psi_j=\tfrac14(C+2\Rea(p\bar u))
e^{2jt}(e^{2nt}-1)\) cancels \(\Rea(p\bar u)\) and yields
\[
T_{jm}=\tfrac12e^{(2j+n)t}(C\sh nt-b_n).
\]
For the minimizer pair,
\(R_{j\star}=\tfrac14ge^{2jt}\) and
\(T_{\star j}=\psi_j-\tfrac12\ip{\sigma_j}{\tilde x_j}
=\tfrac14Ce^{2jt}\).
Averaging proves~\eqref{eq:slacks}.  Since \(R\) is symmetric and
\(T\) antisymmetric, nonnegativity of these slacks is exactly the
pairwise interpolation condition of \Cref{thm:interp}, including the
minimizer pairs.  This proves the equivalence with feasibility.
\end{proof}

\begin{lemma}[Phase of the lag-one symbol]\label{lem:angle}
Let \((u,v)\) be a stationary linear method, \(t\in\R\) and \(\omega\in(0,\pi]\), and let \((g,s)=(g,s)(t,\omega)\ne(0,0)\),
\(g-2is=Re^{-i\Theta}\) with \(R>0\).  Then
\begin{equation}\label{eq:defect}
\Rea\bigl[(1-e^{i\omega})(g-2is)\bigr]=2R\sin\tfrac\omega2\,\sin\Bigl(\frac{\omega-2\Theta}2\Bigr),\qquad
\Rea\bigl[(1-e^{-i\omega})(g-2is)\bigr]=2R\sin\tfrac\omega2\,\sin\Bigl(\frac{\omega+2\Theta}2\Bigr).
\end{equation}
In particular, at \(t=0\) the lag-one symbol \(\Rea[(1-e^{i\omega})\phi_\mu\overline{\phi_L}]\) is a positive multiple of
\(\sin(\Theta(\omega)-\frac\omega2)\).
\end{lemma}

\begin{proof}
Since \(1-e^{\pm i\omega}=2\sin\frac\omega2\,e^{\pm i(\omega-\pi)/2}\), the left-hand sides are \(2R\sin\frac\omega2\cos(\pm\frac{\omega-\pi}2-\Theta)\),
and \(\cos(\frac{\omega-\pi}2-\Theta)=\sin(\frac\omega2-\Theta)\), \(\cos(\frac{\pi-\omega}2-\Theta)=\sin(\frac\omega2+\Theta)\).  At \(t=0\),
\(g-2is=-\frac4{(L-\mu)^2}\phi_\mu\overline{\phi_L}\) by~\eqref{eq:factorgen} (for \(\cM\), \eqref{eq:factor}), so the lag-one symbol is
\(\frac{(L-\mu)^2}2R\sin\frac\omega2\sin(\Theta-\frac\omega2)\).
\end{proof}

\begin{lemma}[Compactness]\label{lem:compact}
For every compact interval \(I\subset\R\), the set of \((t,\nu,C)\) with \(t\in I\) and
\((\nu,C)\) feasible at \(t\) is compact (weak topology on \(\nu\)).
\end{lemma}

\begin{proof}
Put \(G=\max_{I\times[0,\pi]}|g|\).  Feasibility gives \(|C|\le G\), so the set lies in the
compact product \(I\times\mathcal P([0,\pi])\times[-G,G]\), where
\(\mathcal P([0,\pi])\) denotes probability measures with the weak topology.
Every constraint is closed: its integrand is jointly continuous and hence varies uniformly with \(t\),
while integration is continuous under weak convergence of probability measures.
\end{proof}

\begin{lemma}[Duality]\label{lem:duality}
\eqref{eq:feas} is infeasible at \(t\) if and only if a strict certificate
at \(t\) exists.
\end{lemma}

\begin{proof}
Apply the compact-moment separation argument of~\cite[Theorem~2.5(a)]{TA}.
By \Cref{lem:compact}, infeasibility is witnessed by finitely many constraints, including the two minimizer
constraints.  Finite-dimensional strict separation of their moment vectors gives nonnegative weights
whose coefficient of \(C\) is zero and whose combined integrand is uniformly negative.
These are exactly the two conditions in~\eqref{eq:cert}.
Conversely, integrating a strict certificate against a feasible measure would give
\(0\le\int\Phi_t\,d\nu<0\).
\end{proof}

\begin{lemma}[Strict certificates persist]\label{lem:persist}
Let \((\alpha_n,\beta_n,\gamma_\pm)\) be a strict certificate at \(t\), with
margin \(\epsilon\).  There is \(\delta>0\) such that, for any method and
curvature interval with symbols \(g',s'\) at \(t'\),
\[
|t'-t|\le\delta,\qquad
\sup_{\omega\in[0,\pi]}\max\{|g'(t',\omega)-g(t,\omega)|,
|s'(t',\omega)-s(t,\omega)|\}\le\delta
\]
imply a strict certificate with the same lag weights and margin \(\epsilon/2\).
Put \(e=\sum_n(\alpha_n-\beta_n)(\sh nt'-\sh nt)\).
Increase \(\gamma_-\) by \(e\) if \(e\ge0\), and \(\gamma_+\) by \(-e\) otherwise.
For \(t'=t=0\), both minimizer weights are unchanged.
\end{lemma}

\begin{proof}
With finitely many fixed lag weights, the symbol is continuous in
\((t,g,s)\), uniformly in \(\omega\).  The residual \(e\) defined in the statement
tends to zero as \(t'\to t\).  The prescribed correction restores cancellation
and adds \(|e|g'\) to the symbol, also tending uniformly to zero.
For sufficiently small \(\delta\), the strict margin is therefore at least
\(\epsilon/2\).  At \(t'=t=0\), \(e=0\).
\end{proof}

\begin{lemma}[Certificate IQC]\label{lem:iqc}
Let a strict certificate at \(t\) be given, put \(r=e^t\), and define
\[
\begin{gathered}
\mathsf A(\omega)=\sum_n(\alpha_n+\beta_n)(\ch nt-\cos n\omega)+\gamma_++\gamma_-,
\\
\mathsf B(\omega)=\sum_n(\alpha_n-\beta_n)\sin n\omega,\qquad
\Pi(\omega)=\begin{pmatrix}\mathsf A&i\mathsf B\\-i\mathsf B&-\mathsf A\end{pmatrix}.
\end{gathered}
\]
\textup{(a)} For every \(\psi\in\cS(H)\) and sequence with
\(\check x_k=r^{-k}\tilde x_k\in\ell^2\), put
\(\check\sigma_k=r^{-k}\nabla\psi(\tilde x_k)\), and extend both sequences
by zero to negative times.  Then
\[
\mathcal J(\check x,\check\sigma)
:=\frac{1}{2\pi}\int_{-\pi}^{\pi}
[\hat x;\hat\sigma]^*\Pi(\omega)[\hat x;\hat\sigma]\,d\omega\ge0,
\qquad \hat x(\omega)=\sum_k\check x_ke^{-ik\omega},
\]
The Fourier series are understood in \(L^2\), with the same convention for \(\hat\sigma\).
\textup{(b)} For any \(\mathsf u,\mathsf p\in\C\),
\([\mathsf u;\mathsf p]^*\Pi(\omega)[\mathsf u;\mathsf p]\)
equals the certificate symbol computed with
\(g=|\mathsf u|^2-|\mathsf p|^2\),
\(s=\Ima(\mathsf p\bar{\mathsf u})\) and signed frequency \(\omega\).
In particular, for the real method polynomials,
\[
[u(re^{i\omega});p(re^{i\omega})]^*\Pi(\omega)
[u(re^{i\omega});p(re^{i\omega})]=\Phi_t(|\omega|)\qquad(|\omega|\le\pi).
\]
\end{lemma}

\begin{proof}
Use the interpolation-to-IQC construction of~\cite[Lemma~D.4]{TA}.
With \(M=\max(\{n:\alpha_n+\beta_n>0\}\cup\{0\})\), extend
\(\tilde x_k,\sigma_k,\psi_k\) by zero for \(k<0\).  The interpolation inequalities give
\[
\mathcal Q=\sum_{k\ge-M}r^{-2k}\Bigl[
\sum_{n\le M}2r^{-n}\{\alpha_n(R+T)_{k,k+n}+\beta_n(R-T)_{k,k+n}\}
+4\gamma_+(R+T)_{\star k}+4\gamma_-(R-T)_{\star k}\Bigr]\ge0.
\]
This sum converges absolutely by \(\norm{\sigma_k}\le\norm{\tilde x_k}\),
\(|\psi_k|\le\frac12\norm{\tilde x_k}^2\) and \(\check x\in\ell^2\).
The coefficient of each \(\psi_m\), including the zero-extension terms, is
\[
4r^{-2m}\Bigl[\sum_n(\alpha_n-\beta_n)\sh nt+\gamma_+-\gamma_-\Bigr]=0.
\]
On harmonics \((\check x_k,\check\sigma_k)=(\mathsf u,\mathsf p)e^{ik\omega}\),
expansion of the remaining finite-lag form gives
\(\mathsf A(|\mathsf u|^2-|\mathsf p|^2)-2\mathsf B\Ima(\mathsf p\bar{\mathsf u})
=[\mathsf u;\mathsf p]^*\Pi[\mathsf u;\mathsf p]\).
Parseval gives \(\mathcal Q=\mathcal J\), proving both parts.
For real \(u,p\), the lag terms are even in \(\omega\), giving the last identity.
\end{proof}

\begin{lemma}[Homotopy]\label{lem:homotopy}
Let \(\mathcal G\) be convolution with a real impulse response \(\mathsf h=(\mathsf h_n)_{n\ge0}\),
\(\sum_{n\ge0}|\mathsf h_n|<\infty\), \(\mathsf h_0=0\), acting coordinatewise on
\(\R^d\)-valued sequences, and let \(G\) be its transfer function.  This includes every real rational transfer function whose poles lie in the open unit disk and whose numerator has lower degree than its denominator (stable and strictly proper).  Let \(\Delta\) be a causal map on all
one-sided sequences satisfying
\(\norm{P_N\Delta v}\le\norm{P_Nv}\), where \(P_N\) retains the indices
\(0,\ldots,N-1\).  Let \(\Pi\) be a bounded measurable Hermitian \(2\times2\)
matrix function, and suppose
\(\mathcal J(v,\vartheta\Delta v)\ge0\) for every \(v\in\ell^2\),
\(\vartheta\in[0,1]\), and
\[
[G(e^{i\omega});1]^*\Pi(\omega)[G(e^{i\omega});1]
\le-\epsilon'<0\qquad(-\pi\le\omega\le\pi).
\]
Then \(v=\mathcal G\vartheta\Delta v+v^0\) has a unique causal
solution and \(\norm v_{\ell^2}\le c_*\norm{v^0}_{\ell^2}\) for every
\(\vartheta\in[0,1]\), \(v^0\in\ell^2\).  The gain \(c_*\) depends only
on \(G\), \(\sup\norm\Pi\) and \(\epsilon'\), and is independent of
\(d\), \(\Delta\) and \(\vartheta\).
\end{lemma}

\begin{proof}
Use the discrete-time IQC continuation argument~\cite[Lemma~D.5]{TA}, following~\cite[Theorem~1]{MR97}.
The zero direct term gives a unique causal solution.  Parseval and the strict frequency inequality give
an a priori gain \(\norm v\le C\norm{v^0}\) for square-summable solutions, uniformly in \(\vartheta\).
For continuation from \(\vartheta_0\) to \(\vartheta_0+\delta\), truncation and causality give
\[
\norm{P_Nv}\le C\norm{v^0}+C\delta\norm{\mathcal G}\norm{P_N\Delta v}
\le C\norm{v^0}+\tfrac12\norm{P_Nv}
\]
when \(\delta\le(2C\norm{\mathcal G})^{-1}\).
This proves square summability; the a priori estimate restores the same gain \(C\).
Finitely many steps cover \([0,1]\).  The scalar transfer function acts coordinatewise,
so every constant is independent of \(d\).
\end{proof}

\begin{proof}[Proof of \Cref{thm:A}(c)]
Let a strict certificate exist at \(t\), and let \(r=e^t>r_c\).

\emph{Reduction to finite dimension.}  The compression argument in the proof of
\Cref{rem:uniform-smoothness} preserves each finite trajectory and its state norms.
It therefore suffices to obtain a bound in \(\R^d\) uniform in \(d\).

\emph{Feedback form.}  Fix \(H=\R^d\), \(\psi\in\cS(H)\) and a trajectory with initial
state \(S_0\).  Write \(w=w^\sharp+w^0\), where \(w^0\) solves \(\sum_jp_jw^0_{k+j}=0\) with
initial state \(S_0\), and \(w^\sharp\) solves \(\sum_jp_jw^\sharp_{k+j}=\sigma_k\) with zero initial
state.  Then \((\tilde x_k)=\mathcal G(\sigma_k)+(\tilde x^0_k)\), where \(\tilde x^0_k=\sum_ju_jw^0_{k+j}\) is a free response
and \(\mathcal G\) is the causal convolution with transfer function \(u(z)/p(z)\).
After the scaling \(\check x_k=r^{-k}\tilde x_k\), \(\check\sigma_k=r^{-k}\sigma_k\), the transfer function
becomes \(G_t(z)=u(r z)/p(r z)\).  It is stable, because the roots of \(p\) have
modulus at most \(r_c<r\), and strictly proper, because \(\deg u<\deg p\).  The
free response \(\check x^0\) is a combination of \(r^{-k}\)-scaled modes of \(p\), so
\(\norm{\check x^0}_{\ell^2}\le c_1\norm{S_0}\).  The map
\(\Delta:\check x\mapsto(r^{-k}\nabla\psi(r^k\check x_k))_k\) is causal, memoryless, \(1\)-Lipschitz in
each step, and \(\Delta0=0\); moreover \(\vartheta\Delta\) is the same map for
\(\vartheta\psi\in\cS(H)\).  Hence \Cref{lem:iqc}(a) gives \(\mathcal J(\check x,\vartheta\Delta\check x)\ge0\) for
all \(\check x\in\ell^2\) and \(\vartheta\in[0,1]\).  By \Cref{lem:iqc}(b) with \(\mathsf u=u(r e^{i\omega})\),
\(\mathsf p=p(r e^{i\omega})\),
\[
[G_t;1]^*\Pi[G_t;1]=\frac{\Phi_t(|\omega|)}{|p(r e^{i\omega})|^2}\le-\frac{\epsilon}{\max_{|z|=r}|p|^2}.
\]

\emph{An \(\ell^2\) bound.}  \Cref{lem:homotopy} establishes square summability of the
causal feedback solution.
With \(\vartheta=1\) and \(v^0=\check x^0\), it yields
\(\norm{\check x}\le c_*\norm{\check x^0}\le c_*c_1\norm{S_0}\).  Then \(\norm{\check\sigma}\le\norm{\check x}\), and since
\(w^\sharp\) is the output of the stable filter \(1/p(r z)\) applied to \(\check\sigma\),
\begin{equation}\label{eq:l2}
\sum_{k\ge0}r^{-2k}\norm{w_k}^2\le C_0\norm{S_0}^2,
\end{equation}
with \(C_0\) independent of \(H\), \(\psi\) and the trajectory.

\emph{From \(\ell^2\) to a rate.}  Put \(\mathsf E_k=\sum_{j\ge k}r^{-2j}\norm{w_j}^2\).
Apply~\eqref{eq:l2} after restarting at \(S_k\):
\[
\mathsf E_k\le C_0r^{-2k}\norm{S_k}^2\le C'(\mathsf E_k-\mathsf E_{k+D}),\qquad C'=C_0\max(1,r^{2(D-1)}),
\]
because \(r^{-2k}\norm{S_k}^2=\sum_{j<D}r^{-2k}\norm{w_{k+j}}^2\le\max(1,r^{2(D-1)})(\mathsf E_k-\mathsf E_{k+D})\).
Increasing \(C_0\) if necessary, assume \(C'>1\).  Then
\(\mathsf E_{k+D}\le \rho_D\mathsf E_k\) with \(\rho_D=1-1/C'<1\).
Hence \(\mathsf E_{mD}\le \rho_D^m\mathsf E_0\le \rho_D^mC_0\norm{S_0}^2\), and
\(\norm{S_N}^2\le\max(1,r^{2(D-1)})r^{2N}\mathsf E_N\le C''r^{2N}\rho_D^{\lfloor N/D\rfloor}\norm{S_0}^2\).  Thus
\(\Gamma(N)\le(C'')^{1/2}r^N\rho_D^{(N/D-1)/2}\) and \(\rhoF\le r \rho_D^{1/(2D)}<r\).
The second statement of (c) is proved after (b).
\end{proof}

\begin{proof}[Proof of \Cref{thm:A}(d)]
By \Cref{lem:waves} there is \(\psi\in\cS(H_\nu)\) with \(\nabla\psi(\tilde x_k)=\sigma_k\) for all \(k\ge0\).
Since \(\sigma_k=p(z)z^k=\sum_jp_jw_{k+j}\), this is the normalized form
of~\eqref{eq:slm} (proof of \Cref{lem:envelope}), and by \Cref{lem:norm}, \((w_k)\) is a
trajectory on the corresponding \(f\in\cF_{\mu,L}(H_\nu)\).  Multiplication by
\(e^{i\omega}\) is an orthogonal map \(U\) of \(H_\nu^D\), and \(S_{k+1}=e^tUS_k\), so
\(\norm{S_N}=e^{Nt}\norm{S_0}\).  If \(\nu\) has \(J\) atoms, \(H_\nu\cong\C^J\cong\R^{2J}\).
\end{proof}

\begin{proof}[Proof of \Cref{thm:A}(a), (b) and the end of (c)]
(a) Let \(|z|=\rho_q\) with \(v(z)=\lambda u(z)\), \(\lambda\in[\mu,L]\); replacing \(z\) by \(\bar z\) we
may assume \(\omega=\arg z\in[0,\pi]\).  If \(u(z)\ne0\), then \(p(z)=m\,u(z)\) with
\(m=\frac{L+\mu-2\lambda}{L-\mu}\in[-1,1]\); if \(u(z)=0\), then also \(v(z)=p(z)=0\).  In both cases \(g\ge0\) and \(s=0\), so for
\(\nu=\delta_\omega\), \(C=0\), one has
\(a_n\ge0\), \(b_n=0\), and the pair is feasible at \(\log\rho_q\).  Thus \(\log\rho_q\in\cT\).  If \(t\in\cT\), part (d) gives
\(\Gamma(n)\ge e^{nt}\), so \(t\le\log C_\Gamma\) by the proof of \Cref{lem:envelope}.  By
\Cref{lem:compact}, \(\cT\cap I\) is compact for every compact \(I\); hence \(\cT\) is closed,
and it is bounded above.
(b) By (d), \(\rhoF\ge e^T\).  Let \(t>T\).  Then \(t\notin\cT\), and \Cref{lem:duality} gives a
strict certificate at \(t\); since \(e^t>e^T\ge\rho_q\ge r_c\), part (c) gives \(\rhoF<e^t\).
Letting \(t\downarrow T\), \(\rhoF\le e^T\).  At \(t\in\cT\) a strict certificate is excluded by
\Cref{lem:duality}.
(c), second statement: if \(\log r_c<t\le T\) and \(t\notin\cT\), a strict certificate at \(t\)
would give \(\rhoF<e^t\le e^T\), contradicting (b).  Hence
\((\log r_c,T]\subset\cT\), and the reverse inclusion holds since \(T=\max\cT\).
\end{proof}

\paragraph{Normalized lag functionals.}
For part (e), normalize the lag functionals to obtain a limit as \(n\to\infty\).  For \(t\ne0\) let \(N_n(t)=\sh(n|t|)/\sh|t|\) and
\(\mathsf L_n=(a_n-b_n)/N_n\), \(\mathsf U_n=(a_n+b_n)/N_n\).  Then
\[
\ell_n/N_n=\bar{\mathsf L}_n+\operatorname{sgn}(t)\sh|t|\,C,\qquad
\upsilon_n/N_n=\bar{\mathsf U}_n-\operatorname{sgn}(t)\sh|t|\,C,
\]
so the coefficient of \(C\) is independent of \(n\).  As \(n\to\infty\), uniformly for
\(\omega\in[0,\pi]\) and \(|t|\) in compact subsets of \((0,\infty)\),
\(\ch(nt)/N_n\to\sh|t|\) while \(\cos n\omega/N_n,\sin n\omega/N_n\to0\); hence
\(\mathsf L_n,\mathsf U_n\to g\sh|t|\), and the limit constraints are the minimizer
inequalities \(\bar g\pm C\ge0\) multiplied by \(\sh|t|\).
Let \(\mathcal I=\{\mathsf L,\mathsf U\}\times(\N\cup\{\infty\})\),
where \(\N\cup\{\infty\}\) is the one-point compactification.
For \(i=(\mathsf L,n)\), put
\(\phi_i=\mathsf L_n\), \(\varpi_i=\operatorname{sgn}(t)\sh|t|\);
for \(i=(\mathsf U,n)\), put
\(\phi_i=\mathsf U_n\), \(\varpi_i=-\operatorname{sgn}(t)\sh|t|\).
Thus constraint \(i\) is
\(\int\phi_i(t,\omega)\,d\nu+\varpi_i(t)C\ge0\), including the two
minimizer constraints at infinity.  The two distinct labels have the same limiting function,
\(\phi_{(\mathsf L,\infty)}=\phi_{(\mathsf U,\infty)}\).
The functionals and coefficients are jointly continuous for \(t\ne0\).

\begin{lemma}[Finite support]\label{lem:finite}
Let \(T=\max\cT\ne0\) and \(T>\log\rho_q\).  There are a probability measure \(\Lambda\) on
\(\mathcal I\) and a function \(\Phi_T=\int\phi_i(T,\cdot)\,d\Lambda(i)\le0\) on
\([0,\pi]\), analytic in a complex neighbourhood of \([0,\pi]\), with
\(\Phi_T(0)<0\), \(\Phi_T(\pi)<0\), such that every pair feasible at \(T\) is
carried by \(\{\Phi_T=0\}\), a finite subset of \((0,\pi)\), and \(\Lambda\) is
carried by constraints active for every pair feasible at \(T\).
\end{lemma}

\begin{proof}
\emph{A limiting certificate.}  Choose \(t_k\downarrow T\) with
\(0<t_k-T<|T|/2\).  Part (b) gives strict certificates
at \(t_k\).  Written in the normalized functionals and divided by their total
weight, they are probability measures \(\Lambda_k\) on \(\mathcal I\) whose coefficient of
\(C\) vanishes and with \(\int\phi_i(t_k,\omega)\,d\Lambda_k(i)<0\) for all \(\omega\).  A weak limit point
\(\Lambda\) exists by compactness of \(\mathcal I\); by joint continuity,
\(\Phi_T(\omega)=\int\phi_i(T,\omega)\,d\Lambda(i)\le0\), and the coefficient of \(C\) still vanishes.

\emph{Support.}  For \((\nu,C)\) feasible at \(T\), integrate against \(\nu\)
and use cancellation of the \(C\)-coefficient:
\(\int\Phi_T\,d\nu=\int\mathrm{slack}_i\,d\Lambda(i)\ge0\), where \(\mathrm{slack}_i\ge0\) is the
normalized slack of constraint \(i\).  Since \(\Phi_T\le0\), both sides vanish: \(\nu\) is
carried by \(\{\Phi_T=0\}\), and \(\Lambda\) by the constraints with zero slack.

\emph{Analyticity.}  The functions \(g(T,\cdot)\), \(s(T,\cdot)\) are trigonometric
polynomials.  In the strip \(|\Ima\omega|<|T|/2\), \(|\cos n\omega|,|\sin n\omega|\le\ch(n|T|/2)\), and
\(\ch(n|T|/2)/N_n(T)\to0\) geometrically, so the normalized functionals are analytic
and uniformly bounded there and converge uniformly on compact subsets.
Hence \(\Phi_T\) is analytic in the strip.

\emph{Endpoints.}  At \(\omega=0\), \(z=e^T\) is real, \(s=0\) and \(g(T,0)=u(e^T)^2-p(e^T)^2\).  If
\(u(e^T)\ne0\) and \(g(T,0)\ge0\), then \(p(e^T)=m\,u(e^T)\) with \(m\in[-1,1]\), i.e.\
\(v(e^T)=\lambda'u(e^T)\) with \(\lambda'=\frac{L+\mu}2-\frac{L-\mu}2m\in[\mu,L]\), so \(e^T\le\rho_q\), a
contradiction.  If \(u(e^T)=0=p(e^T)\), then also \(v(e^T)=0\), making \(e^T\) a root of
\(v-\lambda u\) for every \(\lambda\in[\mu,L]\), again contrary to \(e^T>\rho_q\).  Thus
\(u(e^T)=0\) implies \(g(T,0)=-p(e^T)^2<0\).  In all cases \(g(T,0)<0\); every normalized functional is a positive multiple of \(g\) at \(\omega=0\)
(there \(b_n=0\) and \(\ch nT-1>0\)), and \(\Phi_T(0)<0\).  The same argument at
\(z=-e^T\) gives \(\Phi_T(\pi)<0\).  A nonzero analytic function has finitely many
zeros in \([0,\pi]\).
\end{proof}

\begin{lemma}[Strict minimizer slack]\label{lem:strict}
Under the hypotheses of \Cref{lem:finite}, every pair feasible at \(T\) has
\(|C|<\bar g\).
\end{lemma}

\begin{proof}
By \Cref{lem:finite}, \(\nu=\sum_{j\le J}\nu_j\delta_{\omega_j}\) with distinct \(\omega_j\in(0,\pi)\).  Put
\(F_\pm(n)=\sum_j\nu_j(g_j\cos n\omega_j\pm2s_j\sin n\omega_j)\).  Then
\(\ell_n=\bar g\ch nT-F_+(n)+C\sh nT\) and \(\upsilon_n=\bar g\ch nT-F_-(n)-C\sh nT\).  Suppose
\(|C|=\bar g\).  If \(C=-\operatorname{sgn}(T)\bar g\), then \(\ell_n\ge0\) reads \(F_+(n)\le\bar ge^{-n|T|}\);
if \(C=\operatorname{sgn}(T)\bar g\), then \(\upsilon_n\ge0\) reads \(F_-(n)\le\bar ge^{-n|T|}\).  In either case
write \(F=F_\pm=\Rea\sum_jc_je^{in\omega_j}\) with \(c_j=\nu_j(g_j\mp2is_j)\).

Each \(c_j\ne0\): \(g_j=s_j=0\) means \(|u|=|p|\) and \(\Ima(p\bar u)=0\) at \(z_j\), i.e.\ \(p\bar u=\pm|u|^2\).
If \(u(z_j)\ne0\), this gives \(p/u=\pm1\); if \(u(z_j)=0\), then \(p(z_j)=v(z_j)=0\).
In either case \(v(z_j)=\mu u(z_j)\) or \(v(z_j)=Lu(z_j)\), a quadratic mode with
\(|z_j|=e^T>\rho_q\), which is impossible.

Distinct \(\omega_j\in(0,\pi)\) make every nonconstant harmonic in \(F\) and \(F^2\)
nonzero modulo \(2\pi\).  Thus the Ces\`aro means of \(F\) tend to \(0\), and
those of \(F^2\) tend to \(\frac12\sum_j|c_j|^2>0\).  On the other hand,
\(F\le\bar ge^{-n|T|}\) makes the Ces\`aro means of the positive part \(F^+\) tend to \(0\);
then so do those of \(F^-=F^+-F\) and of \(|F|=F^++F^-\), and since \(|F|\le K\), the
means of \(F^2\le K|F|\) tend to \(0\).  This contradiction proves \(|C|<\bar g\).
\end{proof}

\begin{lemma}[Atom count]\label{lem:count}
Under the hypotheses of \Cref{lem:finite}, every pair feasible at \(T\) has
finitely many active lags; if \(M\) is the largest active lag of one of
them, every pair feasible at \(T\) is carried by at most \(\frac12(D+M)\)
points.
\end{lemma}

\begin{proof}
With \(\eta=\bar g-|C|>0\) (\Cref{lem:strict}) and \(K=\max|g|+2\max|s|\),
\(\ell_n,\upsilon_n\ge\bar g\ch nT-|C|\sh n|T|-K\ge\eta\ch nT-K>0\) for large \(n\), so only finitely
many lags are active.
The probability measure \(\Lambda\) of \Cref{lem:finite} is carried
by active constraints, and strict minimizer slack excludes its two
infinity labels.  Hence every feasible pair has at least one active
finite lag.  Fix one such pair and let \(M\) be its largest active lag.
Then \(\Lambda\) is carried by
\(\{(\mathsf L,n),(\mathsf U,n):n\le M\}\).
Since \(\deg u<D\) and
\(\deg p=D\), \(g(T,\cdot)\) is an even and \(s(T,\cdot)\) an odd trigonometric polynomial of
degree at most \(D\).  So \(a_n\) and \(b_n\) are even trigonometric polynomials of
degree at most \(D+n\), and \(\Phi_T=P(\cos\omega)\) with a polynomial \(P\) of degree at
most \(D+M\).  Now \(P\le0\) on \([-1,1]\), \(P\not\equiv0\) since \(\Phi_T(0)<0\), and every zero of
\(P\) in \((-1,1)\) is a local maximum, hence of even multiplicity.  So \(P\) has at
most \(\frac12(D+M)\) zeros in \((-1,1)\), and by \Cref{lem:finite} every feasible pair is
carried by them.
\end{proof}

\begin{proof}[Proof of \Cref{thm:A}(e)]
If \(T=\log\rho_q\), choose a maximizing root
\(v(z)=\lambda u(z)\), \(\lambda\in[\mu,L]\).
The wave \(w_k=z^k\) is a trajectory on the quadratic
\(\frac\lambda2\norm\cdot^2\) on \(\C\cong\R^2\), and its state
norm grows exactly by \(e^T\) at each step.  If \(T>\log\rho_q\) and \(T\ne0\), \Cref{lem:finite,lem:strict,lem:count}
show that every pair feasible at \(T\) has \(|C|<\bar g\), finitely many active lags, and
at most \(\frac12(D+M)\) atoms in \((0,\pi)\); part (d) realizes it on \(\R^{2J}\), \(J\le\frac12(D+M)\).
\end{proof}

\begin{lemma}[Start from rest in \(\cM\)]\label{lem:startup}
Let a tuning of \(\cM\) be stable, and let \((\nu,C)\) be feasible at \(t\), with
\(\nu=\sum_{j=1}^J\nu_j\delta_{\omega_j}\), \(\nu_j>0\), and \(|C|<\bar g\) if \(t>0\).
Then there are \(f\in\cF_{\mu,L}(\C^J\oplus\R)\) and \(x_0\ne0\) such that the
trajectory of~\eqref{eq:M} with \(x_{-1}=x_0\) satisfies
\(\norm{x_N}\ge c_0e^{Nt}\norm{x_0}\) for all \(N\ge0\), with \(c_0>0\).  The
trajectory \((x_k)_{k\ge0}\) is the orthogonal sum of the waves
\(x_k^{(j)}=\sqrt{\nu_j}z_j^{k+1}\), \(z_j=e^{t+i\omega_j}\), and a sequence tending
to \(0\) in the last coordinate.
\end{lemma}

\begin{proof}
Following~\cite[Lemma~A.8]{TA}, prescribe the trajectory and its query data in
latent coordinates.  We verify their interpolation by a normalized objective;
\Cref{lem:norm} then gives \(f\).

\emph{Waves.}  Let \(x^w_k=(\sqrt{\nu_j}z_j^{k+1})_j\in\C^J\) for \(k\ge-1\), \(z_j=e^{t+i\omega_j}\).
The latent sequence \(w_k=x^w_{k-1}\) is the wave \(z^k\) of \Cref{lem:waves} in
\(H_\nu\cong\C^J\) (with the weights \(\nu_j\) absorbed into the amplitudes).  Hence the
queries \(\tilde x^w_k=(1+\tau)x^w_k-\tau x^w_{k-1}=(\sqrt{\nu_j}u(z_j)z_j^k)_j\), the normalized
gradients \(\sigma^w_k=(\sqrt{\nu_j}p(z_j)z_j^k)_j\) and the values \(\psi^w_k\) of \Cref{lem:waves},
\(k\ge0\), together with \((0,0,0)\), satisfy \(|T|\le R\) for all pairs, and the
waves satisfy~\eqref{eq:M} with these gradients.

\emph{Extra coordinate.}  For \(\lambda\in(\mu,L)\) let \(P_{-1}=P_0=1\),
\(Q_k=(1+\tau)P_k-\tau P_{k-1}\) and \(P_{k+1}=(1+\beta)P_k-\beta P_{k-1}-h\lambda Q_k\): this is the run
of~\eqref{eq:M} from rest on the scalar quadratic \(\frac\lambda2v^2\), with iterates \(P_k\)
and queries \(Q_k\).  By induction \(Q_k\) is a polynomial in \(\lambda\) of degree \(k\) with
leading coefficient \((-h)^k(1+\tau)^k\ne0\), so \(Q_k(\lambda)=1\) has finitely many roots
for each \(k\ge1\); choose \(\lambda\) outside this countable set.  By stability \(Q_k\to0\), so
\((1-Q_k)^2\ge\frac14\) for large \(k\), and \(\eta_0=\inf_{k\ge1}(1-Q_k)^2>0\).  In normalized
form the quadratic \(\frac\lambda2v^2\) is \(\psi(v)=\frac{m_0}2v^2\) with
\(m_0=\frac{L+\mu-2\lambda}{L-\mu}\in(-1,1)\).  For data of a quadratic, \(T_{ik}=0\) for every
pair, and \(R_{ik}=\frac{1-m_0^2}4(v_i-v_k)^2\ge0\).

\emph{Trajectory.}  Put \(x_{-1}=x_0=(x^w_0,Y)\) with \(Y>0\) to be chosen.  Then
\(\tilde x_0=x_0\).  At \(\tilde x_0\) prescribe the gradient \(((x^w_0-x^w_1)/h,\ \lambda Y)\) and the
normalized value \(\psi_0=\frac12\ip{\sigma_0}{\tilde x_0}\), where \(\sigma_0\) is the normalized gradient.
Then \(x_1=x_0-h\nabla f(\tilde x_0)=(x^w_1,(1-h\lambda)Y)=(x^w_1,P_1Y)\).  For \(k\ge1\), \(x_{k-1}\) and \(x_k\)
have wave parts \(x^w_{k-1},x^w_k\), so \(\tilde x_k=(\tilde x^w_k,Q_kY)\); prescribe the gradient
(wave gradient at \(\tilde x^w_k\), \(\lambda Q_kY\)) and the normalized value
\(\psi_k=\psi^w_k+\frac{m_0}2Q_k^2Y^2\).  By induction \(x_k=(x^w_k,P_kY)\) for all \(k\ge0\).

\emph{Interpolation.}  The data split orthogonally into the wave coordinates
and the extra coordinate, and \(R\) and \(T\) are sums of the two contributions.
In the extra coordinate, \(T=0\) and \(R\ge0\) for every pair, since
\(Q_0=1\) and the extra part of \(\psi_0\) is \(m_0Y^2/2\).  In the
wave coordinates, pairs indexed in \(\{\star,1,2,\ldots\}\) satisfy \(|T|\le R\) by the
first step.  For the pairs \((0,k)\), \(k\ge1\), and \((\star,0)\), the separately
prescribed index-\(0\) datum may cause a wave deficit \(\max(0,|T|-R)\).  This
deficit is bounded uniformly in \(k\).  For \(t\le0\) all data are bounded.  For
\(t>0\), expand \(R_{0k}\) and \(T_{0k}\) around \(R_{k\star}\) and \(T_{\star k}\):
\(\norm{\tilde x_k-\tilde x_0}^2=\norm{\tilde x_k}^2-2\ip{\tilde x_k}{\tilde x_0}+\norm{\tilde x_0}^2\) and similarly for \(\sigma\), and
\(T_{0k}=T_{\star k}-T_{\star0}-\frac12\ip{\sigma_0}{\tilde x_k}+\frac12\ip{\sigma_k}{\tilde x_0}\); the cross terms are
\(O(e^{kt})\) and the rest \(O(1)\).  By~\eqref{eq:slacks},
\(R_{k\star}-|T_{\star k}|=\frac14(\bar g-|C|)e^{2kt}\), so
\(R_{0k}-|T_{0k}|\ge\frac14(\bar g-|C|)e^{2kt}-O(e^{kt})-O(1)\), which is bounded below because
\(\bar g-|C|>0\).  The extra coordinate adds \(\frac{1-m_0^2}4(1-Q_k)^2Y^2\ge\frac{1-m_0^2}4\eta_0Y^2\)
to \(R_{0k}\) and \(\frac{1-m_0^2}4Y^2\) to \(R_{0\star}\), and nothing to \(T\).  One large \(Y\) makes all
pairs interpolable, and \Cref{thm:interp,lem:norm} give \(\psi\) and \(f\) whose
gradients at the queries are the prescribed ones; the trajectory of~\eqref{eq:M}
on \(f\) from \(x_{-1}=x_0\) is the constructed one.

Finally \(\norm{x_N}\ge\norm{x^w_N}=(\sum_j\nu_j|z_j|^{2(N+1)})^{1/2}=e^{(N+1)t}\) and
\(\norm{x_0}=(e^{2t}+Y^2)^{1/2}\), so \(c_0=e^t(e^{2t}+Y^2)^{-1/2}\).
\end{proof}

\begin{lemma}[From rest on a quadratic]\label{lem:restquad}
Let \((\beta,\tau,h)\) be a tuning of \(\cM\), \(\lambda>0\), and let \(\rho(\chi_\lambda)\) be the largest modulus of a root of
\(\chi_\lambda\).  The run of~\eqref{eq:M} from rest on \(\frac\lambda2x^2\) on \(\R\), with \(x_{-1}=x_0\ne0\), satisfies
\(\limsup_{N\to\infty}|x_N|^{1/N}=\rho(\chi_\lambda)\).  In particular \(\rhoz\ge\rho(\chi_\lambda)\) for every \(\lambda\in[\mu,L]\).  If
\(\beta\ne\tau/(1+\tau)\), the queries \(\tilde x_N=(1+\tau)x_N-\tau x_{N-1}\) of this run also satisfy \(\limsup_{N\to\infty}|\tilde x_N|^{1/N}=\rho(\chi_\lambda)\).
\end{lemma}

\begin{proof}
Normalize \(x_0=1\), and write
\(\chi_\lambda(z)=z^2-Sz+Q=(z-w_1)(z-w_2)\).
The recurrence and \(x_{-1}=x_0=1\) give the generating function
\[
X(\zeta):=\sum_{N\ge0}x_N\zeta^N
=\frac{1-Q\zeta}{1-S\zeta+Q\zeta^2}.
\]
A pole at \(\zeta=1/w_j\), \(w_j\ne0\), cannot cancel:
the numerator there is \(1-w_{3-j}\ne0\), since
\(\chi_\lambda(1)=h\lambda>0\).  The Cauchy--Hadamard formula
therefore gives \(\limsup_N|x_N|^{1/N}=\max_j|w_j|\).
If both roots vanish, the series is a polynomial and the same conclusion
holds with value zero.  This proves the lower bound on \(\rhoz\).

For the queries,
\[
\sum_{N\ge0}\tilde x_N\zeta^N
=((1+\tau)-\tau\zeta)X(\zeta)-\tau.
\]
At a nonzero mode \(w_j\), its multiplier is
\(N_\tau(w_j)/w_j\).  This vanishes only if
\(w_j=\tau/(1+\tau)\); substituting in \(\chi_\lambda\) then gives
\(\beta=\tau/(1+\tau)\).  Under the stated exclusion no pole
cancels, proving the query claim as well.
\end{proof}

\begin{proof}[Proof of \Cref{cor:rest}]
The first statement is \Cref{lem:startup}.  If \(T>\log\rho_q\) and \(T\ne0\), \Cref{thm:A}(e) provides a pair
feasible at \(T\) with finitely many atoms and \(|C|<\bar g\).  \Cref{lem:startup} realizes
it from rest, so \(\rhoz\ge e^T=\rhoF\ge\rhoz\) by \Cref{lem:envelope}.  If \(T=\log\rho_q\), let \(\rho_q\) be
attained by a root of \(\chi_\lambda\), \(\lambda\in[\mu,L]\); \Cref{lem:restquad} gives \(\rhoz\ge\rho_q=e^T=\rhoF\ge\rhoz\), the worst case
from rest being attained by \(\frac\lambda2x^2\) on \(\R\).
\end{proof}

\begin{proof}[Proof of \Cref{prop:lagone}]
Let \(r=e^t\le1\).  Take \(\alpha_1=r\), all other \(\alpha_n\) and all \(\beta_n\) equal to \(0\),
\(\gamma_-=0\) and \(\gamma_+=-r\sh t\ge0\).  The second condition of~\eqref{eq:cert} holds:
\(\alpha_1\sh t+\gamma_+=0\).  Since \(r(\ch t-\sh t)=re^{-t}=1\),
\[
\Phi_t=r\bigl[g(\ch t-\cos\omega)-2s\sin\omega\bigr]-r\sh t\,g=g(1-r\cos\omega)-2rs\sin\omega
=\Rea\bigl[(g-2is)(1-z)\bigr],
\]
with \(z=re^{i\omega}\), and by~\eqref{eq:factor}
\[
\Phi_t=-\frac{4\Rea[(1-z)\chi_\mu(z)\overline{\chi_L(z)}]}{h^2(L-\mu)^2}.
\]
(a) The right-hand side is continuous and negative on the circle \(|z|=r\), so
\(\Phi_t\le-\epsilon\) for some \(\epsilon>0\), and \Cref{thm:A}(c) applies since \(r>r_c\).
(b) is (c) with \(\lambda=0\).
(c) At \(t=0\) take in addition \(\beta_2=\lambda\); the second condition of~\eqref{eq:cert} still
holds since \(\sh0=0\).  At \(t=0\), \(a_2+b_2=g(1-\cos2\omega)+2s\sin2\omega=\Rea[(g-2is)(1-e^{-2i\omega})]\),
so the symbol becomes
\[
\Phi_0=-\frac{4\Rea[\chi_\mu\overline{\chi_L}\,((1-z)+\lambda(1-\bar z^2))]}{h^2(L-\mu)^2},\qquad z=e^{i\omega}.
\]
By hypothesis \(\Phi_0<0\) on \((0,\pi]\) and \(\Phi_0(0)=0\).
To obtain strict negativity on the full circle, add \(\eta>0\) to both \(\gamma_+\) and \(\gamma_-\); this keeps the second condition and
adds \(2\eta g\) to the symbol.  At \(z=1\), \(u(1)=1\) and \(p(1)=\frac{L+\mu}{L-\mu}>1\), so
\(g(0,0)<0\), and \(g\le-g_1<0\) on some \([0,\omega_1]\), where \(\Phi_0\le0\); on \([\omega_1,\pi]\),
\(\Phi_0\le-\epsilon_1<0\).  For \(\eta\le\epsilon_1/(4\max|g|)\) the new symbol is at most
\(-\min(2\eta g_1,\epsilon_1/2)<0\).  Stability gives \(r_c<1\), and \Cref{thm:A}(c) applies at
\(t=0\).
\end{proof}

\begin{lemma}[Separation of quadratic modes]\label{lem:separate}
Let a tuning of \(\cM\) and \(r\in(0,1)\) satisfy \(\Rea[(1-z)\chi_\mu(z)\overline{\chi_L(z)}]>0\) on \(|z|=r\).  Then no \(\chi_\lambda\),
\(\lambda\in[\mu,L]\), has a root on \(|z|=r\).  If moreover the roots of \(\chi_\mu\) lie in \(|z|<r\), then \(\rho_q<r\); in particular
the tuning is stable, \(r_c<r\), and \(\rhoF<r\).
\end{lemma}

\begin{proof}
As \(\chi_\lambda\) is affine in \(\lambda\), \(\chi_\lambda=(1-\vartheta)\chi_\mu+\vartheta\chi_L\) with \(\vartheta=(\lambda-\mu)/(L-\mu)\in[0,1]\).  Suppose
\(\chi_\lambda(z)=0\) with \(|z|=r\).  If \(\vartheta=0\) or \(\vartheta=1\), then \(\chi_\mu(z)\overline{\chi_L(z)}=0\).  Otherwise
\(\chi_\mu(z)=-\frac\vartheta{1-\vartheta}\chi_L(z)\), so \(\chi_\mu\overline{\chi_L}=-\frac\vartheta{1-\vartheta}|\chi_L|^2\le0\) at \(z\), and
\(\Rea[(1-z)\chi_\mu\overline{\chi_L}]\le0\) because \(\Rea(1-z)\ge1-r>0\).  Both contradict the hypothesis.  The roots of the monic
quadratics \(\chi_\lambda\) depend continuously on \(\lambda\in[\mu,L]\) and never meet \(|z|=r\); if those of \(\chi_\mu\) lie in \(|z|<r\),
so do those of every \(\chi_\lambda\).  Hence \(\rho_q<r<1\), and \(r_c\le\rho_q\) since \(r_c\) is the largest root modulus of
\(\chi_{(L+\mu)/2}\).  \Cref{prop:lagone}(a) gives \(\rhoF<r\).
\end{proof}

\begin{lemma}[One-atom spectra]\label{lem:single}
Let \((u,v)\) be a stationary linear method, \(K\ge3\), \(1\le j<K/2\), \(\gcd(j,K)=1\) and \(\omega_0=2\pi j/K\).  For \(t\in\R\) with
\((g,s)(t,\omega_0)\ne(0,0)\) write \(g-2is=|g-2is|e^{-i\Theta^{(t)}}\) at \((t,\omega_0)\), \(\Theta^{(t)}\in(-\pi,\pi]\).
\begin{enumerate}[label=\textup{(\alph*)}]
\item If \(|\Theta^{(t)}|\le\pi/K\), then \(g(t,\omega_0)>0\) and \((\delta_{\omega_0},0)\) is feasible at \(t\).  At \(t=0\) the converse holds.
\item If \((g,s)(0,\omega_0)\ne(0,0)\) and \(|\Theta(\omega_0)|<\pi/K\), then \((\delta_{\omega_0},0)\) is feasible at some \(t>0\), with \(0=|C|<\bar g\).
Consequently \(\rhoF\ge e^t>1\), attained by a trajectory in \(\C\cong\R^2\) with \(S_{k+1}=e^tUS_k\) (\Cref{thm:A}(d)), and for a stable
tuning of \(\cM\) also \(\rhoz\ge e^t>1\) (\Cref{lem:startup}).
\end{enumerate}
\end{lemma}

\begin{proof}
At \(t=0\), the criterion is~\cite[Proposition~D.9]{TA}.
(a) For \(\nu=\delta_{\omega_0}\) and \(C=0\), \(\bar g=g\),
\(g=|g-2is|\cos\Theta^{(t)}\), and \(2s=|g-2is|\sin\Theta^{(t)}\), so
\[
\begin{aligned}
\ell_n&=g(\ch nt-1)+|g-2is|\bigl(\cos\Theta^{(t)}-\cos(n\omega_0-\Theta^{(t)})\bigr),\\
\upsilon_n&=g(\ch nt-1)+|g-2is|\bigl(\cos\Theta^{(t)}-\cos(n\omega_0+\Theta^{(t)})\bigr),
\end{aligned}
\]
Since \(\gcd(j,K)=1\), the residues \(m=nj\bmod K\) run through \(0,\ldots,K-1\).
The cosine differences above equal
\(2\sin(\pi m/K)\sin(\pi m/K\mp\Theta^{(t)})\).
For \(|\Theta^{(t)}|\le\pi/K\), these are nonnegative, \(g>0=|C|\), and \(\ch nt\ge1\);
hence the pair is feasible.  At \(t=0\), the terms \(g(\ch nt-1)\) vanish,
and \(|\Theta|>\pi/K\) makes the lower or upper slack at residue \(m=1\) negative.
This proves the converse.
(b) The values \(g,s\) at \((t,\omega_0)\) are continuous in \(t\), and \(g-2is\ne0\) at \(t=0\) with \(\Theta(\omega_0)\in(-\pi/K,\pi/K)\), away from the
branch point \(\pm\pi\); so \(\Theta^{(t)}\) is continuous near \(t=0\), and \(|\Theta^{(t)}|<\pi/K\) for small \(t>0\).  By (a), \((\delta_{\omega_0},0)\) is
feasible at such \(t\), with \(\bar g=g>0=|C|\).  \Cref{thm:A}(d) realizes it by a trajectory with \(S_{k+1}=e^tUS_k\) in \(\C\cong\R^2\), so
\(\Gamma(n)\ge e^{nt}\) and \(\rhoF\ge e^t\); for a stable tuning of \(\cM\), \Cref{lem:startup} applies, as \(|C|<\bar g\), and gives \(\rhoz\ge e^t\).
\end{proof}

\begin{proof}[Proof of \Cref{prop:cycles}]
On the unit circle \(u=N_\tau\ne0\) (its root \(\tau/(1+\tau)\) lies inside the disk), and,
by~\eqref{eq:factor}--\eqref{eq:imab}, \(g-2is=|g-2is|e^{-i\Theta(\omega)}\) with \(\Theta\in(0,\pi)\) on
\((0,\pi)\).
(a) By \Cref{lem:single}(a), \((\delta_{\omega_0},0)\) is feasible at \(t=0\), and \Cref{lem:startup} with
one atom gives \(f\in\cF_{\mu,L}(\C\oplus\R)=\cF_{\mu,L}(\R^3)\) and \(x_0\).  The trajectory is
\(x_k=(z^{k+1},P_kY)\), \(k\ge0\), with \(z=e^{i\omega_0}\) and \(z^K=1\): its first component rotates by
\(\omega_0\) in the plane \(\C\) and is \(K\)-periodic, and \(P_k\to0\).  Hence the trajectory
approaches the \(K\)-periodic orbit \((z^{k+1},0)\), and
\(\norm{x_N}\ge1=(1+Y^2)^{-1/2}\norm{x_0}\).
(b) This is \Cref{lem:single}(b), since \(\Theta(\omega_0)\in(0,\pi/K)\).
\end{proof}

\begin{samepage}
\begin{proposition}[Lyapunov form]\label{prop:lyap}
Let a stable tuning of \(\cM\) admit a strict certificate at \(t=\log r\),
\(r\in(r_c,1]\), with lags at most \(M\ge1\).  For any trajectory on
\(f\in\cF_{\mu,L}(H)\), put
\[
\sigma_k=\nabla\psi(\tilde x_k-x^\star),\qquad
\mathbf X_k=(x_{k-1}-x^\star,x_k-x^\star,\sigma_k,\ldots,\sigma_{k+M-1}).
\]
There is a function \(V_k\), quadratic in \(\mathbf X_k\) minus a linear
combination of \(\psi(\tilde x_{k+j}-x^\star)\), \(0\le j<M\), such that
\[
c_-\norm{\mathbf X_k}^2\le V_k\le c_+\norm{\mathbf X_k}^2,\qquad
V_{k+1}\le r^2V_k-\eta\norm{\mathbf X_k}^2.
\]
Its coefficients and the constants \(c_-,c_+,\eta>0\) depend only on the tuning,
\(\mu,L\), \(M\), and certificate, and work for every \(H,f\) and trajectory.  From rest,
\(\norm{x_N-x^\star}\le K_Vr^N\norm{x_0-x^\star}\), with the same dependence
for \(K_V\).  For the lag-one certificate of \Cref{prop:lagone}(a), \(M=1\),
and \(V_k\) is quadratic in \(\mathbf X_k\) plus
\[
\frac4{L-\mu}\left(f(\tilde x_k)-f(x^\star)
-\frac{L+\mu}4\norm{\tilde x_k-x^\star}^2\right).
\]
\end{proposition}
\end{samepage}

\begin{proof}[Proof of \Cref{prop:lyap}]
\emph{State space.}  Write \(X_k=x_k-x^\star\) and \(e_k=\tilde x_k-x^\star\), so that \(\sigma_k=\nabla\psi(e_k)\) and
\(\mathbf X_k=(X_{k-1},X_k,\sigma_k,\dots,\sigma_{k+M-1})\).  Let \(\lambda_+=\frac{L+\mu}2\) and \(\lambda_-=\frac{L-\mu}2\), so that
\(e_k=(1+\tau)X_k-\tau X_{k-1}\) and, by~\eqref{eq:psi}, \eqref{eq:M} becomes
\(X_{k+1}=S_+X_k-Q_+X_{k-1}+h\lambda_-\sigma_k\) with \(S_+=1+\beta-h\lambda_+(1+\tau)\) and
\(Q_+=\beta-h\lambda_+\tau\).  With the input \(d_k=\sigma_{k+M}\),
\(\mathbf X_{k+1}=\tilde A\mathbf X_k+\tilde Bd_k\), where \(\tilde A\) updates \((X_{k-1},X_k)\) by this recursion
and shifts the gradients, and \(\tilde B\) feeds \(d_k\) into the last slot.  The
matrix \(\tilde A\) is block triangular; its eigenvalues are \(0\) (from the shift) and
the roots of \(z^2-S_+z+Q_+=\chi_{\lambda_+}(z)\), of modulus at most \(r_c<r\).  The
pair \((\tilde A,\tilde B)\) is controllable: \(\tilde A^j\tilde B\), \(0\le j<M\), are the unit vectors
of the gradient slots, \(\tilde A^M\tilde B=(0,h\lambda_-,0,\dots)\) and
\(\tilde A^{M+1}\tilde B=(h\lambda_-,S_+h\lambda_-,0,\dots)\), and these \(M+2\) vectors are linearly
independent.  Each \(e_{k+j}\), \(0\le j\le M\), is a linear function of \(\mathbf X_k\): \(X_{k+j}\)
depends on \(X_{k-1},X_k\) and \(\sigma_k,\dots,\sigma_{k+j-1}\).

\emph{Supply and telescoping.}  Let \(\mathcal Q_k\) be the combination
\[
\mathcal Q_k=\sum_n2r^{-n}\bigl(\alpha_n(R+T)_{k,k+n}+\beta_n(R-T)_{k,k+n}\bigr)+4\gamma_+(R+T)_{\star k}+4\gamma_-(R-T)_{\star k}\ \ge0
\]
of interpolation inequalities in the window \(\{\star,k,\dots,k+M\}\).  Its function values
are \(\sum_j\mathfrak a_j\psi(e_{k+j})\) with \(\mathfrak a_j=2r^{-j}(\alpha_j-\beta_j)\) for \(1\le j\le M\) and
\(\mathfrak a_0=4(\gamma_+-\gamma_-)-\sum_{j\ge1}\mathfrak a_j\), and
\(\sum_j\mathfrak a_jr^{2j}=4[\sum_n(\alpha_n-\beta_n)\sh nt+\gamma_+-\gamma_-]=0\), as in the proof of
\Cref{lem:iqc}.  So \(\sum_j\mathfrak a_j\zeta^j\) is divisible by \(\zeta-r^2\); define real numbers \(\mathfrak b_0,\dots,\mathfrak b_{M-1}\) by
\(\sum_j\mathfrak a_j\zeta^j=(\zeta-r^2)\sum_{j<M}\mathfrak b_j\zeta^j\), and comparing coefficients
(\(\mathfrak a_0=-r^2\mathfrak b_0\), \(\mathfrak a_j=\mathfrak b_{j-1}-r^2\mathfrak b_j\), \(\mathfrak a_M=\mathfrak b_{M-1}\)) gives
\(\sum_j\mathfrak a_j\psi(e_{k+j})=\Psi_{k+1}-r^2\Psi_k\) with \(\Psi_k=\sum_{j<M}\mathfrak b_j\psi(e_{k+j})\).  Let \(\Sigma_r\) be
the symmetric matrix of the remaining quadratic part of \(\mathcal Q_k\) as a form in
\(\mathbf Z_k=(\mathbf X_k,d_k)\), so that \(\mathcal Q_k=\ip{\mathbf Z_k}{\Sigma_r\mathbf Z_k}+\Psi_{k+1}-r^2\Psi_k\).

\emph{Frequency inequality.}  Let
\(\mathcal L(\mathsf W)=[\tilde A\ \tilde B]^\top\mathsf W[\tilde A\ \tilde B]-r^2E^\top\mathsf W E+\Sigma_r\), \(E=[I\ 0]\).  On a complex
solution \(\mathbf Z_k=z^k\mathbf Z\), \(|z|=r\), of the state equation, the vector \(\mathbf Z\) is determined
by \(d\) (as \(z\) is not an eigenvalue of \(\tilde A\)), and since \(p(z)\ne0\), the latent
sequence is a wave \(cz^k\), \(c\in\C\), with \(\norm{\mathbf Z}\le K|c|\).  The computation in the
proof of \Cref{lem:waves}, with the function values removed, gives
\(\mathbf Z^*\Sigma_r\mathbf Z=\Phi_t(|\omega|)|c|^2\le-\epsilon|c|^2\); the terms in \(\Rea(p\bar u)\) cancel by the
second condition of~\eqref{eq:cert}.  Hence the strict frequency inequality
holds on \(|z|=r\), and it persists with \(\Sigma_r\) replaced by \(\Sigma_r+\eta I\) for
\(\eta<\epsilon/K^2\).  The discrete-time Kalman--Yakubovich--Popov lemma on the circle
\(|z|=r\)~\cite[Lemma~2.4]{Gan22} (see also~\cite{Ran96}), applicable because \((\tilde A,\tilde B)\) is
controllable, gives a real symmetric \((M+2)\times(M+2)\) matrix \(\mathsf W\) with \(\mathcal L(\mathsf W)\preceq-\eta I\).

\emph{Dissipation.}  With \(V_k=\ip{\mathbf X_k}{\mathsf W\mathbf X_k}-\Psi_k\) (\(\mathsf W\) acting as \(\mathsf W\otimes I\) on \(H^{M+2}\)),
\[
V_{k+1}-r^2V_k=\ip{\mathbf Z_k}{\mathcal L(\mathsf W)\mathbf Z_k}-\mathcal Q_k\le-\eta\norm{\mathbf Z_k}^2,
\]
because \(\ip{\mathbf X_{k+1}}{\mathsf W\mathbf X_{k+1}}-r^2\ip{\mathbf X_k}{\mathsf W\mathbf X_k}=\ip{\mathbf Z_k}{(\mathcal L(\mathsf W)-\Sigma_r)\mathbf Z_k}\) and
\(\mathcal Q_k\ge0\).  For lag one, \(\mathcal Q_k=2(R+T)_{k,k+1}+2(1-r^2)(R+T)_{\star k}\) and \(\Psi_k=2\psi(e_k)\), i.e.\ \(\mathfrak b_0=2\), and
\(-2\psi(e_k)=\frac{4}{L-\mu}\bigl(f(\tilde x_k)-f(x^\star)-\frac{L+\mu}4\norm{e_k}^2\bigr)\) by~\eqref{eq:psi}.

\emph{Bounds on \(V\).}  Upper bound: \(|\psi(e)|\le\frac12\norm e^2\) and \(e_{k+j}\) is linear in
\(\mathbf X_k\), so \(|V_k|\le\eta^+\norm{\mathbf X_k}^2\) with \(\eta^+\) depending only on \(\mathsf W\), \(\mathfrak b\) and the tuning.  Lower bound: by \Cref{thm:A}(c), \(\rhoF<r\), so for some \(r'<r\) every
trajectory satisfies \(\norm{\mathbf X_k}\le c'r'^k\norm{\mathbf X_0}\) (the gradients in \(\mathbf X_k\) are
bounded by the queries, hence by the states \(S_k,\dots,S_{k+M}\)).  Summing
\(r^{-2(k+1)}V_{k+1}\le r^{-2k}V_k-\eta r^{-2(k+1)}\norm{\mathbf X_k}^2\) over \(k<n\),
\(r^{-2n}V_n\le V_0-\eta r^{-2}\norm{\mathbf X_0}^2\), and \(r^{-2n}|V_n|\le\eta^+c'^2(r'/r)^{2n}\norm{\mathbf X_0}^2\to0\).
Hence \(V_0\ge\eta r^{-2}\norm{\mathbf X_0}^2\), and the same holds at every \(k\) by restarting.
From rest, \(\norm{\mathbf X_0}\le K''\norm{X_0}\), because \(\norm{\sigma_j}\le\norm{e_j}\) and
\(\norm{S_j}\le\Gamma(j)\norm{S_0}\); so
\(\frac\eta{r^2}\norm{X_N}^2\le V_N\le r^{2N}V_0\le\eta^+K''^2r^{2N}\norm{X_0}^2\).
\end{proof}

\begin{remark}[Tongues as phase locking of a Hill recursion]\label{rem:hill}
On the planar orbit of \Cref{prop:cycles}(a), the gradient field acts by a
rotating curvature.  With \(\mathsf R\) the rotation by \(\omega_0\), the
query-normalized gradient amplitude \(\hat p=p(z)/u(z)\)
(\(|\hat p|\le1\) by feasibility), and a
symmetric matrix \(\mathsf H_0\) of norm at most \(1\) with first column \((\Rea\hat p,\Ima\hat p)\), for
instance \(\bigl(\begin{smallmatrix}\Rea\hat p&\Ima\hat p\\\Ima\hat p&-\Rea\hat p\end{smallmatrix}\bigr)\), of norm \(|\hat p|\),
in a basis in which \(\tilde x_0\) is real and positive the orbit satisfies
\(\sigma_n=\mathsf H_n\tilde x_n\) with \(\mathsf H_n=\mathsf R^n\mathsf H_0\mathsf R^{-n}\), a symmetric matrix of
norm at most \(1\) whose trace-free part rotates by \(2\omega_0\) per step.  In the
plane of the orbit the method is therefore a linear recursion with
periodic symmetric coefficients (period dividing \(K\)), a discrete Hill
equation.  In the rotating frame \(\mathsf R^{-n}x_n\), the orbit is a fixed point.
Thus a \(K\)-tongue is a resonance region for periodically rotating curvatures
that one function of the class can realize.  One may take \(\mathsf H_0\), and then every
\(\mathsf H_n\), scalar if and only if \(\Ima\hat p=0\), i.e.\ for a quadratic mode.
\end{remark}

\paragraph{Spread spectra at a half step.}
The next lemma constructs spectra for a method with phase lead close to half a
step on a frequency band.  It gives strict margins for both minimizer constraints,
all upper constraints and all lower constraints with lag at least two;
lag one is checked separately.  We use it in \Cref{app:classM,app:corner}.
For \(t<0\), put \(C=-(1-\eta)\bar g\), \(\eta\in(0,1)\), and
\[
\mathcal E_\pm(x)=\ch x\pm(1-\eta)\sh x\qquad(x\ge0),\qquad
\Phi^{(1)}_t(\omega)=\Rea\bigl[(1-e^{t+i\omega})(g-2is)(t,\omega)\bigr].
\]
Since \(\ell_n=\bar g\ch nt-\Rea\int(g-2is)e^{in\omega}d\nu+C\sh nt\) and
\(\upsilon_n=\bar g\ch nt-\Rea\int(g+2is)e^{in\omega}d\nu-C\sh nt\), one has, for
\(C=-(1-\eta)\bar g\),
\begin{equation}\label{eq:Epm}
\ell_n=\bar g\,\mathcal E_+(n|t|)-\Rea\int(g-2is)e^{in\omega}\,d\nu,\qquad
\upsilon_n=\bar g\,\mathcal E_-(n|t|)-\Rea\int(g+2is)e^{in\omega}\,d\nu .
\end{equation}
For \(n=1\), \(\mathcal E_+(|t|)=e^{-t}-\eta\sh|t|\) and \(e^{-t}g-\Rea[(g-2is)e^{i\omega}]=e^{-t}\Phi^{(1)}_t\), so
\begin{equation}\label{eq:l1extreme}
\ell_1=e^{-t}\int\Phi^{(1)}_t\,d\nu-\eta\sh|t|\,\bar g .
\end{equation}
Thus the lower lag-one slack is the average of the certificate symbol
\(\Phi^{(1)}_t\), up to its explicit positive factor and the correction in \(\eta\).  For \(\cM\), \(\Phi^{(1)}_t=-4H_r/(h^2(L-\mu)^2)\),
\(r=e^t\), with \(H_r\) as in \Cref{lem:cubic}.

Choose once and for all an even smooth cutoff \(\mathsf b_0\ge0\), equal to \(1\) on \([-\frac18,\frac18]\) and zero for \(|u|\ge\frac14\), and set \(\mathsf b=\mathsf b_0/\int\mathsf b_0\).  Thus \(\mathsf b\in C^\infty_c((-\frac12,\frac12))\), \(\mathsf b\ge0\) and \(\int\mathsf b=1\).  Fix this choice throughout and put
\(\mathsf k=\mathsf b*\mathsf b\), supported in \([-1,1]\), and \(\hat{\mathsf k}(x)=\int e^{ixu}\mathsf k(u)\,du=\hat{\mathsf b}(x)^2\).

\begin{lemma}[Spread spectra at a half step]\label{lem:spreadhalf}
There are \(\mathsf s_1\in(0,\frac16]\), \(c_{\mathsf k}>0\) and \(C_{\mathsf k}<\infty\), depending only on \(\mathsf b\),
with the following property.  Let \((u,v)\) be a stationary linear method, \(t<0\),
\(K\ge3\), \(\omega_0=2\pi/K\), \(\mathsf s\in(0,\mathsf s_1]\), \(\mathsf w=\mathsf s\omega_0\) and
\(B=[\omega_0-\mathsf w,\omega_0+\mathsf w]\).  Suppose that \((g,s)(t,\omega)\ne(0,0)\) on \(B\), write
\(g-2is=Re^{-i\Theta}\) with \(R>0\) and \(\Theta\) continuous on \(B\), and put
\(\mathsf d=\max_B|\omega-2\Theta(\omega)|\), the phase defect from a half-step lead.
Let \(\eta\in(0,1)\),
\(\mathsf m=\min(\frac14,\frac{\omega_0^2}6,c_{\mathsf k}\mathsf s^2)\), and assume
\begin{equation}\label{eq:spreadhyp}
\mathsf d\le\frac{\mathsf m}2,\qquad C_{\mathsf k}\Bigl(\frac{|t|}{\mathsf m\,\mathsf w}\Bigr)^2+\mathsf d<\frac{\sqrt\eta}8 .
\end{equation}
Let \(\nu\) be the probability measure on \(B\) with density proportional to
\(\mathsf k((\omega-\omega_0)/\mathsf w)/R(\omega)\), and \(C=-(1-\eta)\bar g\).  Then \(\bar g>0\), \(|C|<\bar g\),
\(\upsilon_n>0\) for every \(n\ge1\), and \(\ell_n>0\) for every \(n\ge2\).  If moreover
\begin{equation}\label{eq:spreadl1}
\int_B\Phi^{(1)}_t\,d\nu>\eta\,e^t\sh|t|\,\bar g ,
\end{equation}
then \((\nu,C)\) is feasible at \(t\), and some pair \((\nu',C')\) with \(\nu'\) carried by finitely many
points of \(B\) and \(|C'|<\bar g'\) is feasible at \(t\).
\end{lemma}

\begin{proof}
\emph{The kernel.}  \(\hat{\mathsf b}\) is real and even, \(\hat{\mathsf b}(0)=1\), \(|\hat{\mathsf b}(x)|<1\) for \(x\ne0\) (the
support of \(\mathsf b\) has interior points), \(\hat{\mathsf b}(x)=1-\frac{m_{\mathsf b}}2x^2+O(x^4)\) with
\(m_{\mathsf b}=\int u^2\mathsf b>0\), and \(|\hat{\mathsf b}(x)|\le C_1/|x|\).  Hence \(\hat{\mathsf k}=\hat{\mathsf b}^2\ge0\),
\(\hat{\mathsf k}(x)=1-m_{\mathsf b}x^2+O(x^4)\) and \(\hat{\mathsf k}(x)\le C_1^2/x^2\).  Choose \(x_\circ\in(0,1]\) with
\(1-\frac32m_{\mathsf b}x^2\le\hat{\mathsf k}(x)\le1-\frac12m_{\mathsf b}x^2\) on \([0,x_\circ]\), and let
\(M_0=\max_{x\ge x_\circ}\hat{\mathsf k}(x)<1\).  Put
\(a_0=\min\bigl(\frac{x_\circ}3,((1-M_0)/(3m_{\mathsf b}))^{1/2},(3m_{\mathsf b})^{-1/2}\bigr)\),
\(c_{\mathsf k}=\min(\frac{m_{\mathsf b}}3,\frac{1-M_0}2)\), \(\mathsf s_1=\min(\frac16,\frac{3a_0}\pi)\) and \(C_{\mathsf k}=16C_1^2\).
For \(0<a\le a_0\) and \(x\ge3a\),
\begin{equation}\label{eq:kgap}
\hat{\mathsf k}(a)\ge\tfrac12,\qquad \hat{\mathsf k}(x)\le\hat{\mathsf k}(a)-c_{\mathsf k}\min(x^2,1).
\end{equation}
Indeed \(\hat{\mathsf k}(a)\ge1-\frac32m_{\mathsf b}a^2\ge\frac12\).  If \(x\le x_\circ\), then
\(\hat{\mathsf k}(a)-\hat{\mathsf k}(x)\ge\frac12m_{\mathsf b}x^2-\frac32m_{\mathsf b}a^2\ge\frac13m_{\mathsf b}x^2\), because \(a^2\le x^2/9\).  If
\(x\ge x_\circ\), then \(\hat{\mathsf k}(a)-\hat{\mathsf k}(x)\ge1-\frac32m_{\mathsf b}a^2-M_0\ge\frac12(1-M_0)\).

\emph{Reduction.}  Let \(\varrho(\omega)=\mathsf w^{-1}\mathsf k((\omega-\omega_0)/\mathsf w)\), so that
\(\hat\varrho(\zeta):=\int e^{i\zeta\omega}\varrho(\omega)\,d\omega=e^{i\zeta\omega_0}\hat{\mathsf k}(\zeta \mathsf w)\), and \(d\nu=\varrho\,d\omega/(ZR)\) with
\(Z=\int\varrho/R\).  Since \(g\pm2is=Re^{\pm i\Theta}\), multiplying~\eqref{eq:Epm} by \(Z\) shows that
\(\ell_n>0\) and \(\upsilon_n>0\) are equivalent to
\[
\Rea\int e^{i(n\omega-\Theta)}\varrho\,d\omega<\mathcal E_+(n|t|)\int\cos\Theta\,\varrho\,d\omega,\qquad
\Rea\int e^{i(n\omega+\Theta)}\varrho\,d\omega<\mathcal E_-(n|t|)\int\cos\Theta\,\varrho\,d\omega .
\]
Write \(\Theta=\frac\omega2-\frac\delta2\) with \(|\delta|\le\mathsf d\) on \(B\).  Since \(|e^{ix}-1|\le|x|\),
\[
\Bigl|\int e^{i(n\omega\mp\Theta)}\varrho\,d\omega-\hat\varrho\bigl(n\mp\tfrac12\bigr)\Bigr|\le\frac{\mathsf d}2,\qquad
\int\cos\Theta\,\varrho\,d\omega\ge G_0-\frac{\mathsf d}2,\quad G_0:=\Rea\hat\varrho(\tfrac12)=\cos\tfrac{\omega_0}2\,\hat{\mathsf k}(\tfrac{\mathsf w}2).
\]
As \(\omega_0\le\frac{2\pi}3\) and \(\frac{\mathsf w}2\le\frac{\pi\mathsf s}3\le a_0\), \(G_0\ge\frac14\).
Also \(\mathsf d\le\omega_0^2/12<(5/6)\omega_0\le\omega\) on \(B\),
so \(\Theta>0\).  Finally \(\mathsf d\le\frac{\mathsf m}2\le\frac18\) and
\(\omega_0+\mathsf w\le\frac{7\pi}9\), so \(0<\Theta\le\frac12(\frac{7\pi}9+\frac18)<\frac{5\pi}{12}\) on \(B\); hence \(\cos\Theta>\frac14\),
\(\bar g=\int R\cos\Theta\,d\nu>0\), and \(|C|=(1-\eta)\bar g<\bar g\).

\emph{The margin.}  For \(n\ge1\) let \(\mathfrak m(n)=G_0-\Rea\hat\varrho(n+\frac12)=G_0-\cos((n+\frac12)\omega_0)\,\hat{\mathsf k}((n+\frac12)\mathsf w)\).
We claim \(\mathfrak m(n)\ge\mathsf m\).  Put \(a=\frac{\mathsf w}2\) and \(x=(n+\frac12)\mathsf w\ge3a\).  If
\(\cos((n+\frac12)\omega_0)\le0\), then \(\mathfrak m(n)\ge G_0\ge\frac14\), since \(\hat{\mathsf k}\ge0\).  If \(n\equiv0\) or
\(n\equiv-1\pmod K\), then \(n\ge K-1\), \(\cos((n+\frac12)\omega_0)=\cos\frac{\omega_0}2\), and by~\eqref{eq:kgap}
\(\mathfrak m(n)=\cos\frac{\omega_0}2(\hat{\mathsf k}(a)-\hat{\mathsf k}(x))\ge\frac12c_{\mathsf k}\min(x^2,1)\); here
\(x\ge(K-\frac12)\mathsf w=2\pi\mathsf s(1-\frac1{2K})\ge\frac{5\pi}3\mathsf s\), so \(x^2\ge27\mathsf s^2\), and \(27\mathsf s^2\le1\), whence
\(\mathfrak m(n)\ge c_{\mathsf k}\mathsf s^2\).  Otherwise \((n+\frac12)\omega_0\) is congruent modulo \(2\pi\) to
\((i+\frac12)\omega_0\) with \(1\le i\le K-2\), so \(\cos((n+\frac12)\omega_0)\le\cos\frac{3\omega_0}2\); if this cosine is
positive, then \(\hat{\mathsf k}(x)\le\hat{\mathsf k}(a)\) gives
\(\mathfrak m(n)\ge\hat{\mathsf k}(a)(\cos\frac{\omega_0}2-\cos\frac{3\omega_0}2)\ge\sin\omega_0\sin\frac{\omega_0}2\ge\frac{\omega_0^2}6\), by concavity of
\(\sin\) on \([0,\pi]\) (\(\sin\varphi\ge\frac{3\sqrt3}{4\pi}\varphi\) on \([0,\frac{2\pi}3]\) and \(\sin\varphi\ge\frac{3\sqrt3}{2\pi}\varphi\) on \([0,\frac\pi3]\), and
\(\frac{27}{16\pi^2}>\frac16\)).

\emph{Upper constraints.}  Let \(n\ge1\) and \(x=n|t|\).  If \(x\le\frac{\mathsf m}4\), then \(\mathcal E_-(x)\ge e^{-x}\ge1-x\), so
the right-hand side is at least \((1-x)(G_0-\frac{\mathsf d}2)\ge G_0-\frac{\mathsf d}2-x\), while the left-hand side is
at most \(G_0-\mathfrak m(n)+\frac{\mathsf d}2\le G_0-\mathsf m+\frac{\mathsf d}2\); since \(\mathsf m>\frac{\mathsf m}4+\mathsf d\), \(\upsilon_n>0\).  If
\(x>\frac{\mathsf m}4\), then \((n+\frac12)\mathsf w>\mathsf m \mathsf w/(4|t|)\), so \(|\Rea\hat\varrho(n+\frac12)|\le C_1^2(4|t|/(\mathsf m \mathsf w))^2\), and
\(\mathcal E_-(x)=(1-\frac\eta2)e^{-x}+\frac\eta2e^x\ge\sqrt{\eta(2-\eta)}\ge\sqrt\eta\); the left-hand side is at most
\(C_{\mathsf k}(|t|/(\mathsf m \mathsf w))^2+\frac{\mathsf d}2\), the right-hand side at least \(\sqrt\eta(G_0-\frac{\mathsf d}2)\ge\frac{\sqrt\eta}8\), and
\(\upsilon_n>0\) by~\eqref{eq:spreadhyp}.

\emph{Lower constraints, \(n\ge2\).}  \(\mathcal E_+\ge1\), so the right-hand side is at least \(G_0-\frac{\mathsf d}2\), and
the left-hand side is at most \(\Rea\hat\varrho(n-\frac12)+\frac{\mathsf d}2=G_0-\mathfrak m(n-1)+\frac{\mathsf d}2\le G_0-\mathsf m+\frac{\mathsf d}2\); as
\(\mathsf m>\mathsf d\), \(\ell_n>0\).

\emph{Lag one and feasibility.}  By~\eqref{eq:l1extreme}, \(\ell_1>0\) is~\eqref{eq:spreadl1}.  With \(|C|<\bar g\), all
constraints of~\eqref{eq:feas} hold strictly.

\emph{Finitely many atoms.}  For every probability measure \(\nu'\) on \(B\) and
\(C'=-(1-\eta)\bar g'\), the identities~\eqref{eq:Epm} hold, and
\(|\Rea\int(g\mp2is)e^{in\omega}d\nu'|\le\int R\,d\nu'<4\bar g'\) because \(\cos\Theta>\frac14\) on \(B\).  Since
\(\mathcal E_\pm(x)\ge\frac\eta2e^x\), all constraints with \(n\ge n_2:=\lceil|t|^{-1}\log(8/\eta)\rceil\) hold for every such
\(\nu'\).  The finitely many remaining constraints, \(\ell_n,\upsilon_n\) with \(n<n_2\), are integrals of
continuous functions on \(B\) against \(\nu'\) and hold strictly for \(\nu\); they persist for the
Riemann sums \(\nu'\) of \(\nu\) on a fine enough grid, which converge weakly to \(\nu\).  Finally
\(|C'|=(1-\eta)\bar g'<\bar g'\).
\end{proof}

For \(\cM\) in the accelerated scaling the hypotheses of the lemma are checked once, in \Cref{lem:spreadwit}, where the
defect \(\mathsf d\) is read off the mixed symbol of \Cref{lem:mixed} through~\eqref{eq:defect}.

\section{Proofs for Section~\ref{sec:beyond}}
\label{app:beyond}

Throughout, a second-order family with slow data \((\theta,y)\) is fixed, \(L=1\), \(\mu=\varepsilon^2\), and
\(\theta_\varepsilon=(1-\beta_\varepsilon)/\varepsilon\), \(y_\varepsilon=-1/q_\varepsilon(1)\).  For the member \(\varepsilon\) we write \(\phi^\varepsilon_\lambda=\lambda u_\varepsilon-v_\varepsilon\); at
\(\varepsilon=0\), \(\phi^0_0=\phi_0\) and \(\phi^0_1=\phi_1\) are the polynomials of the limit method in \Cref{sec:faithful}.

\begin{lemma}[Slow normal form]\label{lem:normal}
For \(z=e^{-\varepsilon\varpi}\), uniformly for \(\varpi\) in compact subsets of \(\C\),
\[
\hat\lambda_\varepsilon(z):=\frac{v_\varepsilon(z)}{u_\varepsilon(z)}=\frac{\varepsilon^2}{y_\varepsilon}\,\varpi(\theta_\varepsilon-\varpi)+O(\varepsilon^3).
\]
For small \(\varepsilon\) the polynomial \(\phi^\varepsilon_{\varepsilon^2}=\varepsilon^2u_\varepsilon-v_\varepsilon\) has two roots \(e^{-\varepsilon\varpi_\pm}\) with \(\varpi_\pm\) converging to the roots
of \(\varpi^2-\theta\varpi+y\).  Hence \(\rho_q\ge1-\varepsilon\cQ(\theta,y)+o(\varepsilon)\).
\end{lemma}

\begin{proof}
With \(E(x)=(1-e^{-x})/x\), continuously extended by \(E(0)=1\), \(z-1=-\varepsilon\varpi E(\varepsilon\varpi)\) and \(z-\beta_\varepsilon=\varepsilon(\theta_\varepsilon-\varpi E(\varepsilon\varpi))\).  The
coefficients are bounded and \(|z-1|=O(\varepsilon)\), so \(q_\varepsilon(z)=q_\varepsilon(1)+O(\varepsilon)\) and
\(u_\varepsilon(z)=u_\varepsilon(1)+O(\varepsilon)=1+O(\varepsilon)\).  Multiplying, \(v_\varepsilon(z)=-\varepsilon^2\varpi(\theta_\varepsilon-\varpi)q_\varepsilon(1)+O(\varepsilon^3)\), which gives
the first claim.  Next, \(\varepsilon^{-2}\phi^\varepsilon_{\varepsilon^2}(e^{-\varepsilon\varpi})=u_\varepsilon(z)\bigl(1-\varepsilon^{-2}\hat\lambda_\varepsilon(z)\bigr)\) converges, uniformly
on compact sets, to \(\frac1y(\varpi^2-\theta\varpi+y)\); Rouch\'e's theorem, counting multiplicities, gives two zeros converging
to the roots of \(\varpi^2-\theta\varpi+y\).  The moduli of the corresponding roots of \(\phi^\varepsilon_{\varepsilon^2}\) are \(e^{-\varepsilon\Rea\varpi_\pm}\),
and the smaller of \(\Rea\varpi_\pm\) tends to \(\cQ(\theta,y)\) (to \(\zeta_-\) if the roots are real, to \(\theta/2\)
otherwise).  Since \(\rho_q\) is the largest modulus of a root of \(\phi^\varepsilon_\lambda\), \(\lambda\in[\mu,L]\), the last claim
follows.
\end{proof}

\paragraph{Slow spectra.}
Here a \emph{slow spectrum} is supported in \([0,\varepsilon\gamma_{\max}]\)
for some \(\gamma_{\max}\) independent of \(\varepsilon\).
The lower bounds of \Cref{thm:ceiling,thm:B}(a) use \(t=-c\varepsilon\) and
\(\omega_j=\varepsilon\gamma_j\), with \(\Delta_c=c^2-\theta c+y\).
We choose the value constant \(C=C_*\) of \Cref{lem:split}.
In query normalization, where each wave has query amplitude \(1\), the lag inequalities have exactly the
kernels \(\cA_\gamma(q),\cB_\gamma(q)\) of the flow program at
\(q=n\varepsilon\); their coefficients converge uniformly to the flow
coefficients on bounded frequency sets (\Cref{lem:reduce}).
The uniform margin of \Cref{lem:spread} absorbs these perturbations
for all lags simultaneously (\Cref{lem:flowspirals}).
Conversely, \Cref{lem:Cforced} shows that every feasible slow spectrum
has \(C=C_*+O(\varepsilon^3)\), which will exclude a single slow frequency.

\begin{lemma}[Per-atom splitting]\label{lem:split}
Let \(\nu=\sum_j\nu_j\delta_{\omega_j}\) and \(t<0\).  Set
\(C_*=t^{-1}\sum_j2\nu_js_j\omega_j\).  For \(C=C_*\),
\[
\ell_n,\ \upsilon_n=\textstyle\sum_j\nu_j\bigl[g_j(\ch nt-\cos n\omega_j)\pm2s_j\bigl(\tfrac{\omega_j\sh n|t|}{|t|}-\sin n\omega_j\bigr)\bigr]
\quad(\text{signs }+,-\text{ respectively}),
\]
and \(|C_*|\le\bar g\) if and only if \(|\sum_j\nu_js_j\omega_j|\le\frac{|t|}2\bar g\).
\end{lemma}

\begin{proof}
Substitute
\(C_*\sh nt=\sum_j\nu_j2s_j\omega_j\sh(n|t|)/|t|\)
in~\eqref{eq:feas}.  The two lag identities follow, and
\(|C_*|\le\bar g\) is exactly
\(|\sum_j\nu_j2s_j\omega_j|\le|t|\bar g\).
\end{proof}

For a finitely supported \(\nu\) with \(u(z_j)\ne0\), the constraints~\eqref{eq:feas}
are homogeneous in \((\nu,C)\); replacing \(g,s\) by \(g/|u|^2=1-|\hat p|^2\) and
\(s/|u|^2=\Ima\hat p\), \(\hat p=p/u\), and \(\nu_j\) by \(\nu_j|u(z_j)|^2\) (renormalized, with \(C\)
divided by the same constant) gives an equivalent system.  We call this the
\emph{query normalization}.  The original amplitudes \(u,p\) and spectral
weights in~\eqref{eq:feas} give the \emph{latent normalization}.
By~\eqref{eq:p} with \(L=1\), \(\mu=\varepsilon^2\),
\(\hat p=(1+\varepsilon^2-2\hat\lambda)/(1-\varepsilon^2)\) for every stationary linear method, where
\(\hat\lambda=v(z)/u(z)\) is the complex curvature of the wave (for \(\cM\), \(\hat\lambda=-(z-1)(z-\beta)/(hN_\tau(z))\)).
Define \(\gt,\st\) by
\begin{equation}\label{eq:gtst}
1-|\hat p|^2=\tfrac{4\varepsilon^2}{y_\varepsilon}\gt,\qquad 2\Ima\hat p=\tfrac{4\varepsilon^2}{y_\varepsilon}\st .
\end{equation}

\begin{lemma}[Exact reduction and expansion]\label{lem:reduce}
Let \(c>0\), \(0<\varepsilon<1\) with \(y_\varepsilon>0\), and put \(t=-c\varepsilon\), \(\omega_j=\varepsilon\gamma_j\in[0,\pi]\) with \(\gamma_j\ge0\).  Let \(\nu\) be a finitely supported probability
measure on \(\{\omega_j\}\) in the query normalization.  The pair \((\nu,C_*)\)
of \Cref{lem:split} is feasible at \(t\) if and only if, for all
\(q\in\varepsilon\N\),
\[
\textstyle\sum_j\nu_j\bigl[\gt_j\cA_{\gamma_j}(q)\pm\st_j\cB_{\gamma_j}(q)\bigr]\ge0,\qquad
\bigl|\sum_j\nu_j\gamma_j\st_j\bigr|\le c\sum_j\nu_j\gt_j .
\]
Moreover \(u_\varepsilon(z_j)\to1\), so the query normalization is defined for small
\(\varepsilon\), and \(\gt\to\gamma^2-\Delta_c\), \(\st\to(\theta-2c)\gamma\) as \(\varepsilon\to0\),
uniformly for \((c,\gamma)\) in compact subsets of \((0,\infty)\times[0,\infty)\).
\end{lemma}

\begin{proof}
\emph{Reduction.}  With \(q=n\varepsilon\), \(|t|=c\varepsilon\) and \(\omega_j=\varepsilon\gamma_j\) we have, exactly,
\[
\ch nt-\cos n\omega_j=\ch cq-\cos\gamma_jq=\cA_{\gamma_j}(q),\qquad
\frac{\omega_j\sh n|t|}{|t|}-\sin n\omega_j=\frac{\gamma_j\sh cq}c-\sin\gamma_jq=\cB_{\gamma_j}(q).
\]
In the query normalization \(g_j=1-|\hat p_j|^2=\frac{4\varepsilon^2}{y_\varepsilon}\gt_j\) and \(2s_j=2\Ima\hat p_j=\frac{4\varepsilon^2}{y_\varepsilon}\st_j\).
\Cref{lem:split}, divided by \(4\varepsilon^2/y_\varepsilon\), gives the lag constraints, and its minimizer
condition \(|\sum_j\nu_j2s_j\omega_j|\le c\varepsilon\,\bar g\) becomes \(|\sum_j\nu_j\gamma_j\st_j|\le c\sum_j\nu_j\gt_j\).

\emph{Expansion.}  Put \(\varpi=c-i\gamma\), so that \(z=e^{t+i\omega}=e^{-\varepsilon\varpi}\).  The coefficients of \(u_\varepsilon\) are
bounded and \(|z-1|=O(\varepsilon)\), so \(u_\varepsilon(z)=u_\varepsilon(1)+O(\varepsilon)=1+O(\varepsilon)\); in particular the query normalization is
defined for small \(\varepsilon\).  By \Cref{lem:normal}, the complex curvature is
\(\hat\lambda=\varepsilon^2\varpi(\theta_\varepsilon-\varpi)/y_\varepsilon+O(\varepsilon^3)\).  Writing
\(\hat p=1-2\varepsilon^2\tilde w\), we get \(\tilde w=(\hat\lambda/\varepsilon^2-1)/(1-\varepsilon^2)\to\varpi(\theta-\varpi)/y-1\).  Since
\(1-|\hat p|^2=4\varepsilon^2(\Rea\tilde w-\varepsilon^2|\tilde w|^2)\) and \(2\Ima\hat p=-4\varepsilon^2\Ima\tilde w\),~\eqref{eq:gtst} gives
\(\gt=y_\varepsilon(\Rea\tilde w-\varepsilon^2|\tilde w|^2)\) and \(\st=-y_\varepsilon\Ima\tilde w\).  Finally
\(\varpi(\theta-\varpi)=(\gamma^2-c^2+\theta c)+i(2c-\theta)\gamma\), so \(\gt\to\gamma^2-\Delta_c\) and \(\st\to(\theta-2c)\gamma\).
All remainders are uniform for \((c,\gamma)\) in compact sets, as in \Cref{lem:normal}.
\end{proof}

\begin{lemma}[The value constant is forced]\label{lem:Cforced}
Fix \(\gamma_{\max}>0\), a compact set \(I_c\subset(0,\infty)\), and a
second-order family with slow data \((\theta,y)\).  There are
\(\varepsilon_2,K_*>0\) such that every pair \((\nu,C)\) feasible at
\(t=-c\varepsilon\), \(c\in I_c\), \(0<\varepsilon<\varepsilon_2\), in query
normalization with support in \([0,\varepsilon\gamma_{\max}]\), satisfies
\[
|C-C_*|\le K_*\varepsilon^3,\qquad
C_*:=\frac1t\int2s\omega\,d\nu
=-\frac{4\varepsilon^2}{y_\varepsilon c}\int\gamma\st\,d\nu.
\]
\end{lemma}

\begin{proof}
By \Cref{lem:reduce}, \(|g|,|s|\le K_1\varepsilon^2\) on the support of \(\nu\), uniformly for \(c\in I_c\).
The lag-one constraints \(\ell_1,\upsilon_1\ge0\) read \(\bar a_1-\bar b_1+C\sh t\ge0\) and
\(\bar a_1+\bar b_1-C\sh t\ge0\), i.e.\ \(|\bar b_1-C\sh t|\le\bar a_1\).  Here
\(\bar a_1=\int g(\ch c\varepsilon-\cos\varepsilon\gamma)\,d\nu\le K_2\varepsilon^4\), since the bracket is \(O(\varepsilon^2)\).  Next,
\[
\bar b_1-C_*\sh t=\int2s\Bigl(\sin\varepsilon\gamma-\varepsilon\gamma\frac{\sh c\varepsilon}{c\varepsilon}\Bigr)d\nu ,
\]
and the bracket is \(O(\varepsilon^3)\) uniformly, so \(|\bar b_1-C_*\sh t|\le K_3\varepsilon^5\).  Hence
\(|C-C_*|\,|\sh t|\le|\bar b_1-C\sh t|+|\bar b_1-C_*\sh t|\le K_4\varepsilon^4\), and since \(|\sh t|\ge c\varepsilon\),
\(|C-C_*|\le(K_4/\min I_c)\varepsilon^3\).  The formula for \(C_*\) follows from
\(2s=\frac{4\varepsilon^2}{y_\varepsilon}\st\), \(\omega=\varepsilon\gamma\) and \(t=-c\varepsilon\).
\end{proof}

\begin{lemma}[Sampled flow spirals]\label{lem:flowspirals}
Let \(\theta/3<c<\cQ(\theta,y)\), and let \(J\), \(\gamma_i=\gamma_0\lambda_i\) and \(\nu_i=1/J\) \((1\le i\le J)\) be given by \Cref{lem:spread}
with \(m=y\).  For all small \(\varepsilon\), the spectrum with atoms \(\omega_i=\varepsilon\gamma_i\) and weights \(\nu_i\) in the query
normalization, with \(C=C_*\), is feasible at \(t=-c\varepsilon\) and has \(|C_*|<\bar g\).  Consequently the corresponding pair in the
latent normalization, with weights \(\nu_i|u_\varepsilon(z_i)|^{-2}\) renormalized and \(C\) rescaled accordingly, is feasible at \(t\),
has \(J\) atoms and satisfies \(|C|<\bar g\).
\end{lemma}

\begin{proof}
Let \(\epsilon_*\) be as in \Cref{lem:spread}.
The quantities \(c,\gamma_0,J,\epsilon_*\) are fixed before
\(\varepsilon\) is chosen; take \(\varepsilon\) small enough that also
\(\varepsilon\gamma_i\in[0,\pi]\) for every \(i\).
By \Cref{lem:reduce} choose \(\varepsilon_0\) such that, for \(\varepsilon<\varepsilon_0\) and all \(i\),
\[
|\gt_i-(\gamma_i^2-\Delta_c)|\le\eta',\qquad |\st_i-(\theta-2c)\gamma_i|\le\eta',\qquad
\eta'=\min\bigl(\tfrac{\epsilon_*}2,\tfrac1{2(c+2\gamma_0)}\bigr).
\]
Since \(\cA_\gamma,\cB_\gamma\ge0\), \eqref{eq:S1} and~\eqref{eq:S2} give, for every \(q>0\),
\[
\sum_i\nu_i\bigl[\gt_i\cA_{\gamma_i}(q)\pm\st_i\cB_{\gamma_i}(q)\bigr]\ge(\epsilon_*-\eta')\sum_i\nu_i\bigl(\cA_{\gamma_i}(q)+\cB_{\gamma_i}(q)\bigr)\ge0
\]
(the minus sign uses~\eqref{eq:S1}, the plus sign~\eqref{eq:S2}).  For the minimizer
condition, using \(\gamma_i\le2\gamma_0\) and~\eqref{eq:S3},
\[
\begin{aligned}
c\sum_i\nu_i\gt_i-\Bigl|\sum_i\nu_i\gamma_i\st_i\Bigr|
&\ge c\sum_i\nu_i(\gamma_i^2-\Delta_c)-(\theta-2c)\sum_i\nu_i\gamma_i^2-(c+2\gamma_0)\eta'\\
&=(3c-\theta)\sum_i\nu_i\gamma_i^2-c\Delta_c-(c+2\gamma_0)\eta'\ \ge\ \tfrac12 .
\end{aligned}
\]
\Cref{lem:reduce,lem:split} now give feasibility and \(|C_*|<\bar g\).
Homogeneity preserves both conclusions under the inverse change of spectral
weights from query to latent normalization.
\end{proof}

\begin{proof}[Proof of \Cref{thm:ceiling}]
If \(\cQ\le\theta/3\), \Cref{lem:normal} and \(\rhoF\ge\rho_q\) (\Cref{lem:envelope}) give the claim.  Let
\(\theta/3<c<\cQ\).  \Cref{lem:flowspirals} gives \(J\), independent of \(\varepsilon\), and for small \(\varepsilon\) a pair with \(J\) atoms that is
feasible at \(t=-c\varepsilon\).  By \Cref{thm:A}(d), it is realized by a trajectory on
\(H=L^2(\nu;\C)\cong\R^{2J}\), with \(S_{k+1}=e^{-c\varepsilon}US_k\).
Hence \(\Gamma(n)\ge e^{-c\varepsilon n}\), \(\rhoF\ge e^{-c\varepsilon}\), and
\(\limsup\sqrt\kappa(1-\rhoF)\le c\).
Letting \(c\downarrow\theta/3\) completes the nonlinear branch.
\end{proof}

\begin{proof}[Proof of \Cref{cor:ceiling}]
(i) By \Cref{lem:normal} the quadratic constant is at most \(\cQ(\theta,y)\); apply \Cref{thm:ceiling}
and \Cref{cor:damping}.  (ii) \(\cFF(2\sqrt y,y)=\min(\frac23\sqrt y,\sqrt y)\).  For \((\frac4{\sqrt3},\frac43)\),
\(\theta^2=4y\), so \(\cFF=\min(\frac4{3\sqrt3},\frac2{\sqrt3})\).  (iii) \(\cFF\ge1\) requires \(\theta\ge3\) and \(\cQ(\theta,y)\ge1\).
If \(\theta^2\ge4y\), \(\zeta_-\ge1\) is equivalent to \(\theta-2\ge\sqrt{\theta^2-4y}\), i.e.\ \(y\ge\theta-1\); if \(\theta^2<4y\), then
\(y>\theta^2/4\ge\theta-1\).
\end{proof}

\begin{proof}[Proof of \Cref{lem:low}]
First, the transfer-function form of \(\teff\) in \Cref{sec:faithful} follows from \(v_0=(z-1)^2q_0\), \((q_0/u_0)(1)=-1/y\) and
\((q_0/u_0)'(1)=q_0'(1)+u_0'(1)/y\).  Here \(\phi_0\overline{\phi_1}=-v_0\overline{(u_0-v_0)}\).  Since \(v_0=O(\omega^2)\) on the circle near \(1\),
\(\phi_0\overline{\phi_1}=-(z-1)^2q_0(z)\overline{u_0(z)}+O(\omega^4)\), and \(\overline{u_0(z)}=u_0(1/z)\) for \(|z|=1\).  Put
\(G(z)=q_0(z)u_0(1/z)\): \(G(1)=-1/y\) and \(G'(1)=q_0'(1)-q_0(1)u_0'(1)=q_0'(1)+u_0'(1)/y\).  With
\((z-1)^2=-\omega^2-i\omega^3+O(\omega^4)\) and \(G(e^{i\omega})=G(1)+i\omega G'(1)+O(\omega^2)\),
\[
-(z-1)^2G(z)=\omega^2G(1)+i\omega^3\bigl(G(1)+G'(1)\bigr)+O(\omega^4),
\]
and \(G(1)+G'(1)=\teff/y\).  For \(\Szero\), multiply by \(1-e^{i\omega}=-i\omega+\frac{\omega^2}2+O(\omega^3)\): the real part is
\(\frac{\teff}y\omega^4-\frac{\omega^4}{2y}+O(\omega^5)\), and \(\Szero\) is even in \(\omega\), since the coefficients are real.
\end{proof}

\begin{lemma}[Two-scale expansion]\label{lem:twoscale}
Fix \(c\in\R\), put \(r=1-c\varepsilon\), \(z=re^{i\omega}\), and
\(\mathcal H(\varepsilon,\omega)=\Rea\bigl[(1-z)\phi^\varepsilon_{\varepsilon^2}(z)\overline{\phi^\varepsilon_1(z)}\bigr]\).
For every second-order family, \(\mathcal H\) is even in \(\omega\), and, as \((\varepsilon,\omega)\to(0^+,0)\),
\[
\mathcal H(\varepsilon,\omega)=\frac\varepsilon y\bigl(c\Delta_c\,\varepsilon^2+(\theta-3c)\,\omega^2\bigr)+o\bigl((\varepsilon+|\omega|)^3\bigr).
\]
If the family is of class \(C^3\), then \(\mathcal H\) is also of class \(C^3\) on the half-neighbourhood \(\varepsilon\ge0\).
\end{lemma}

\begin{proof}
Write \(s=\varepsilon+|\omega|\).  The factorization following \Cref{def:second} and coefficient continuity give
\begin{align*}
z-1&=-c\varepsilon+i\omega+O(s^2),&
z-\beta_\varepsilon&=(\theta-c)\varepsilon+i\omega+o(s),\\
q_\varepsilon(z)&=-\frac1y+o(1),& u_\varepsilon(z)&=1+o(1).
\end{align*}
Consequently, with \(a=c\varepsilon\) and \(b=(\theta-c)\varepsilon\),
\[
\phi^\varepsilon_{\varepsilon^2}(z)
=\varepsilon^2+\tfrac1y(i\omega-a)(i\omega+b)+o(s^2),\qquad
\phi^\varepsilon_1(z)=1+o(1).
\]
The other factor is \(1-z=a-i\omega+O(s^2)\), so their product has real part
\[
a\varepsilon^2-\frac{a^2b}{y}+\frac{\omega^2(b-2a)}y+o(s^3)
=\frac\varepsilon y\bigl(c\Delta_c\varepsilon^2+(\theta-3c)\omega^2\bigr)+o(s^3).
\]
This proves the joint expansion for every relative scale of \(\omega\) and \(\varepsilon\).
Replacing \(\omega\) by \(-\omega\) conjugates the bracket, proving evenness.
Finally, \(\mathcal H\) is a polynomial in \(r\cos\omega\), \(r\sin\omega\) and the real coefficients;
hence \(C^3\) coefficient regularity implies \(C^3\) regularity of \(\mathcal H\).
\end{proof}

\begin{proof}[Proof of \Cref{thm:faithful}]
(a) Let \(0<c<\cFF\); then \(\Delta_c>0\) and \(\theta-3c>0\).  Here \(\Szero>0\) near \(0\), so \(s_{2\mathfrak p}>0\).  The family
is of class \(C^{2\mathfrak p}\), so the polynomial dependence on its coefficients makes \(\mathcal H\) of class \(C^{2\mathfrak p}\), and by Taylor's formula at \((0,0)\),
\(\mathcal H=\sum_{i+j\le2\mathfrak p}\mathfrak h_{ij}\varepsilon^i\omega^j+R\) with \(|R|\le\eta(\varepsilon+|\omega|)^{2\mathfrak p}\) for \(\varepsilon+|\omega|\le\rho_\eta\), for every \(\eta>0\).  Here
\(\mathfrak h_{ij}=0\) for odd \(j\), since \(\mathcal H\) is even in \(\omega\); \(\mathfrak h_{0j}\) are the Taylor coefficients of \(\mathcal H(0,\omega)=\Szero(\omega)\), so
\(\mathfrak h_{0j}=0\) for \(j<2\mathfrak p\) and \(\mathfrak h_{0,2\mathfrak p}=s_{2\mathfrak p}\); and by \Cref{lem:twoscale} there are no terms of degree below \(3\), and the
cubic ones are \(\frac\varepsilon y(c\Delta_c\varepsilon^2+(\theta-3c)\omega^2)\).  Every other monomial has \(i\ge1\) and either \(j=0\), \(i\ge4\), or
\(j\ge2\), \(i+j\ge4\); for \(\varepsilon+|\omega|\le1\) it is bounded by \(\varepsilon^3(\varepsilon+|\omega|)\) in the first case and by \(\varepsilon\omega^2(\varepsilon+|\omega|)\)
in the second.  With \(M_0\) the sum of the absolute values of their coefficients, \(\eta=s_{2\mathfrak p}2^{-2\mathfrak p}\) and
\((\varepsilon+|\omega|)^{2\mathfrak p}\le2^{2\mathfrak p-1}(\varepsilon^4+\omega^{2\mathfrak p})\) (as \(2\mathfrak p\ge4\), \(\varepsilon\le1\)), for \(\varepsilon+|\omega|\le\rho\le\min(\rho_\eta,1)\),
\[
\mathcal H\ge\varepsilon^3\Bigl(\frac{c\Delta_c}y-\bigl(M_0+\tfrac{s_{2\mathfrak p}}2\bigr)\rho\Bigr)+\varepsilon\omega^2\Bigl(\frac{\theta-3c}y-M_0\rho\Bigr)+\frac{s_{2\mathfrak p}}2\,\omega^{2\mathfrak p}.
\]
Fix \(\rho\) so small that both brackets are positive; then \(\mathcal H>0\) for \(0<\varepsilon\le\rho/2\), \(|\omega|\le\rho/2\).  On
\(\rho/2\le|\omega|\le\pi\), \(\mathcal H(\varepsilon,\omega)\to \Szero(\omega)\) uniformly, and \(\min_{[\rho/2,\pi]}\Szero>0\); so \(\mathcal H>0\) there for small \(\varepsilon\).
Hence the lag-one symbol is positive on \(|z|=1-c\varepsilon\).  By~\eqref{eq:p} the normalized
polynomial is proportional to \((1+\varepsilon^2)u_\varepsilon-2v_\varepsilon\), whose roots converge to those of \(u_0-2v_0\)
(the leading coefficient is \(-2\) times that of \(v_\varepsilon\), which is nonzero at \(\varepsilon=0\)); so
\(r_c<1-c\varepsilon\) for small \(\varepsilon\), and the lag-one criterion gives \(\rhoF<1-c\varepsilon\).  With
\Cref{thm:ceiling} the limit follows.

\emph{A criterion for (b).}  At \(t=0\) let \(\Xi=\phi^\varepsilon_{\varepsilon^2}\overline{\phi^\varepsilon_1}\) on the unit circle, let
\(\mathcal H=\Rea[(1-e^{i\omega})\Xi]=(1-\cos\omega)\Rea \Xi+\sin\omega\,\Ima \Xi\) be the lag-one symbol, and let
\(\check{\mathcal H}=\Rea[(1-e^{-i\omega})\Xi]=\mathcal H-2\sin\omega\,\Ima \Xi\) be the symbol of the backward multiplier \(1-\bar z\).  By~\eqref{eq:factorgen}, \(g-2is\)
is a positive multiple of \(-\Xi\), so by \Cref{lem:angle}, where \(\Xi\ne0\), \(\mathcal H\) is a positive multiple of \(\sin(\Theta-\frac\omega2)\) and \(\check{\mathcal H}\) a
positive multiple of \(-\sin(\Theta+\frac\omega2)\).  If \(\mathcal H<0\) and \(\check{\mathcal H}<0\) at some \(\omega\in(0,\pi)\), then, as \(\Theta\in(-\pi,\pi]\),
\(\Theta\in(\frac\omega2-\pi,\frac\omega2)\cap(-\frac\omega2,\pi-\frac\omega2)\), i.e.\ \(|\Theta(\omega)|<\omega/2\).

(b) At \(\varepsilon=0\) the symbols of the criterion are those of the limit method,
\(\mathcal H=\Szero\) and \(\check{\mathcal H}=\check{\mathcal H}_0:=\Szero-2\sin\omega\,\Ima\Xi_0\) with \(\Xi_0=[\phi_0\overline{\phi_1}](e^{i\omega})\), and \(-\Xi_0=|\Xi_0|\,e^{-i\Theta_0}\) by the definition
of \(\Theta_0\).  By \Cref{lem:low}, \(\Ima\Xi_0=\frac\teff y\omega^3+O(\omega^5)\) (it is odd in \(\omega\)), so
\[
\Szero=\frac{2\teff-1}{2y}\,\omega^4+O(\omega^6),\qquad \check{\mathcal H}_0=-\frac{1+2\teff}{2y}\,\omega^4+O(\omega^6).
\]
If \(|\teff|<\frac12\), both coefficients are negative.  If \(\teff=\frac12\) and \(s_{2\mathfrak p}<0\), then \(\Szero\le\frac12s_{2\mathfrak p}\omega^{2\mathfrak p}\) and
\(\check{\mathcal H}_0\le-\frac1{2y}\omega^4\) near \(0\).  In both cases there is \(\omega_1\in(0,\frac{2\pi}3]\) with \(\Szero<0\) and \(\check{\mathcal H}_0<0\) on \((0,\omega_1]\), and the
criterion gives \(|\Theta_0(\omega)|<\omega/2\) there.  With \(K_0=\lceil2\pi/\omega_1\rceil\ge3\), every \(K\ge K_0\) has \(2\pi/K\le\omega_1\), so
\(|\Theta_0(2\pi/K)|<\pi/K\); as \(\Szero(2\pi/K)\ne0\), also \(\phi_0\phi_1(e^{2\pi i/K})\ne0\), and (d), with \(j=1\), gives the feasible
spectra.

(c) Let \(\mathcal R=v_0/u_0\), holomorphic near \(z=1\) as \(u_0(1)=1\).  By the transfer-function form of \(\teff\),
\(\mathcal R(1+w)=-\frac{w^2}y\bigl(1-(1+\teff)w+O(w^2)\bigr)\).  The function \(G(\zeta,\delta)=\delta^{-2}\mathcal R(1+\delta\zeta)-1\) is holomorphic near
\((\zeta,\delta)=(i\sqrt y,0)\), with \(G(\zeta,0)=-\zeta^2/y-1\) and \(\partial_\zeta G(i\sqrt y,0)=-2i/\sqrt y\ne0\).  The implicit function theorem
gives a holomorphic \(\zeta(\delta)\) with \(G(\zeta(\delta),\delta)=0\) and \(\zeta(0)=i\sqrt y\), and the terms of order \(\delta\) give
\(\zeta'(0)=\frac{1+\teff}2\zeta(0)^2=-\frac{(1+\teff)y}2\).  So for small \(\lambda=\delta^2>0\), \(z=1+\delta\zeta(\delta)\) is a root of \(\lambda u_0-v_0\) with
\[
|z|^2=1+2\delta\Rea\zeta(\delta)+\delta^2|\zeta(\delta)|^2=1-\teff\,\lambda y+O(\lambda^{3/2}).
\]
If \(\teff<0\), fix such a \(\lambda_1\in(0,1]\) with a root of modulus \(r_1>1\).  For small \(\varepsilon\), \(\lambda_1\in[\varepsilon^2,1]\), and the roots of
\(\lambda_1u_\varepsilon-v_\varepsilon\) converge to those of \(\lambda_1u_0-v_0\), since its leading coefficient \(-v_D\) is nonzero at
\(\varepsilon=0\); so \(\rho_q\ge(1+r_1)/2\) for small \(\varepsilon\), and the root gives growing solutions on \(\frac{\lambda_1}2x^2\) on \(\R\).  By
\Cref{lem:envelope}, \(\rhoF\ge\rho_q>1\).

(d) At \(t=0\) and \(\omega_0=2\pi j/K\) the values \(g,s\) depend continuously on \(\varepsilon\), with a nonzero limit
by~\eqref{eq:factorgen} and \(\phi_0\phi_1(e^{i\omega_0})\ne0\).
Hence, for small \(\varepsilon\), \((g,s)(0,\omega_0)\ne(0,0)\) and \(|\Theta(\omega_0)|<\pi/K\) (as \(\Theta_0(\omega_0)\) lies away from the branch point \(\pm\pi\)), and
\Cref{lem:single}(b) gives a pair feasible at some \(t>0\) and \(\rhoF>1\).
\end{proof}

\begin{proof}[Proof of \Cref{cor:lead}]
If \(\teff<0\), \Cref{thm:faithful}(c) gives \(\rho_q>1\) for small \(\varepsilon\).  If \(\rhoF<1\) for all small \(\varepsilon\), then
\(\teff\ge0\) by \Cref{lem:envelope} and the first claim, and \(\teff\ge\frac12\) by \Cref{thm:faithful}(b), whose conclusion is
\(\rhoF>1\) whenever \(|\teff|<\frac12\).
\end{proof}

\begin{proof}[Proof of \Cref{thm:trace}]
(a) Put
\[
Q(z)=-\frac{z(1+z)\,q_0(z)}{\phi_1(z)}=\frac{z(1+z)\,\phi_0(z)}{(z-1)^2\phi_1(z)} .
\]
Since \(\deg u_0<\deg v_0\), \(\phi_1\) and \(z(1+z)q_0\) have the same degree and opposite leading coefficients, so
\(Q(\infty)=1\); and \(Q(1)=-2q_0(1)/\phi_1(1)=2/y\), as \(\phi_1(1)=u_0(1)=1\).  By the hypothesis, \(Q\) is holomorphic in
\(|z|>1\) and continuous on \(|z|\ge1\), the simple root of \(\phi_1\) at \(-1\), if any, being cancelled by \(1+z\).  On
the unit circle, \(z/(z-1)^2=-1/(4\sin^2\frac\omega2)\) and \((1+z)/(1-z)=i\cot\frac\omega2\), so
\[
Q(e^{i\omega})=-\frac{i\cot\frac\omega2}{4\sin^2\frac\omega2}\cdot\frac{(1-z)\phi_0(z)\overline{\phi_1(z)}}{|\phi_1(z)|^2},\qquad
\Ima Q(e^{i\omega})=-\frac{\cos\frac\omega2\;\Szero(\omega)}{4\sin^3\frac\omega2\;|\phi_1(e^{i\omega})|^2}
\]
for \(0<\omega<\pi\).  Let \(F(\zeta)=Q(1/\zeta)-1\): it is holomorphic in the unit disk, continuous on its closure, has
real Taylor coefficients and \(F(0)=0\).  By the Schwarz formula,
\(F(\zeta)=\frac i{2\pi}\int_{-\pi}^{\pi}\Ima F(e^{it})\,\frac{e^{it}+\zeta}{e^{it}-\zeta}\,dt\), and since \(\Ima F(e^{-it})=-\Ima F(e^{it})\), pairing \(t\) with \(-t\) gives
\[
F(\zeta)=\frac2\pi\int_0^\pi\Ima F(e^{it})\,\frac{\zeta\sin t}{1-2\zeta\cos t+\zeta^2}\,dt\qquad(|\zeta|<1).
\]
Here \(\Ima F(e^{it})=\Ima Q(e^{-it})=-\Ima Q(e^{it})\).  Let \(\zeta=x\uparrow1\).  As \(1-2x\cos t+x^2\ge2x(1-\cos t)\), the kernel is at
most \(\frac1{2}\cot\frac t2\) and tends to it, while \(\Ima F(e^{it})\cot\frac t2=\cos^2\frac t2\,\Szero(t)/(4\sin^4\frac t2\,|\phi_1|^2)\) is bounded: near
\(t=0\) because \(\Szero=O(t^4)\) (\Cref{lem:low}), and near \(t=\pi\) because \(Q\) is continuous at \(-1\).  By dominated
convergence and the continuity of \(F\) at \(1\),
\[
\frac2y-1=F(1)=\frac1\pi\int_0^\pi\frac{\cos^2\frac t2\;\Szero(t)}{4\sin^4\frac t2\;|\phi_1(e^{it})|^2}\,dt ,
\]
which is~\eqref{eq:trace}, since \(\sin^2\frac t2=\frac\xi2\) and \(\cos^2\frac t2=\frac{2-\xi}2\).  If \(\Szero\equiv0\), then \(\Ima F=0\) on the circle, so the
harmonic function \(\Ima F\) vanishes, \(F\) is a real constant, and \(F=F(0)=0\): \(Q\equiv1\).  Then
\(u_0-(z-1)^2q_0=\phi_1=-z(1+z)q_0\), i.e.\ \(u_0=(1-3z)q_0\), and \(v_0/u_0=-(z-1)^2/(3z-1)\).  Conversely, if
\(v_0/u_0=-(z-1)^2/(3z-1)\), then \(u_0=(1-3z)q_0\), \(\phi_0/\phi_1=(v_0/u_0)/(v_0/u_0-1)=(z-1)^2/(z(z+1))\), and on the circle
\((1-z)\phi_0/\phi_1=-(z-1)^3/(z(z+1))=4i\sin^2\frac\omega2\tan\frac\omega2\) is imaginary, so \(\Szero=|\phi_1|^2\Rea[(1-z)\phi_0/\phi_1]\equiv0\);
moreover \(q_0(1)=-\frac12u_0(1)=-\frac12\), so \(y=2\), and \(u_0'(1)=-2q_0'(1)+\frac32\), so \(\teff=\frac12\).

(b) Divide \(u_0\) and \(v_0\) by their greatest common divisor \(r\), normalized by \(r(1)=1\) (possible since \(u_0(1)=1\)).
This preserves \(v_0/u_0\) and \(y\), keeps the roots of \(\lambda u_0-v_0\) among the old ones, and gives a reduced symbol
\(\Szero^{\rm red}\) with \(\Szero=|r|^2\Szero^{\rm red}\).  Thus \(\Szero^{\rm red}\ge0\) away from the finitely many roots of \(r\),
and everywhere by continuity.  We may therefore assume \(u_0\) and \(v_0\) coprime.  Put \(\phi_\lambda=\lambda u_0-v_0=\lambda\phi_1+(1-\lambda)\phi_0\).  On the circle
\begin{equation}\label{eq:lambdasym}
\Rea\bigl[(1-z)\phi_0\overline{\phi_\lambda}\bigr]=\lambda\,\Szero+(1-\lambda)\,\xi\,|\phi_0|^2 .
\end{equation}
Let \(\lambda\in(\lambda_0,1)\).  If \(\phi_\lambda(e^{i\alpha})=0\), then \(\alpha\ne0\), since \(\phi_\lambda(1)=\lambda\), and~\eqref{eq:lambdasym} vanishes at \(\alpha\); both terms
are nonnegative (\(\Szero\) is even and continuous), so \(\phi_0(e^{i\alpha})=0\), and then \(u_0(e^{i\alpha})=0\), against coprimality.
So the roots of \(\phi_\lambda\) lie in the open unit disk, and (a) applies to the method \((u_0,v_0/\lambda)\), whose step is \(\lambda y\)
and whose \(\phi_0,\phi_1\) are \(\phi_0/\lambda\), \(\phi_\lambda/\lambda\):
\begin{equation}\label{eq:tracelambda}
\frac2{\lambda y}-1=\frac1{2\pi}\int_0^\pi\frac{(2-\xi)\bigl(\lambda\Szero+(1-\lambda)\xi|\phi_0|^2\bigr)}{\xi^2\,|\phi_\lambda(e^{i\omega})|^2}\,d\omega\ \ge0 .
\end{equation}
For the same strictly stable curvatures, put
\[
Q_\lambda(z)=-\frac{z(1+z)q_0(z)}{\phi_\lambda(z)},\qquad F_\lambda(\zeta)=Q_\lambda(1/\zeta)-1.
\]
The paired Schwarz representation in~(a) has a nonnegative density by~\eqref{eq:lambdasym}.
For \(0<\zeta<1\), its kernel is nonnegative and at most its limiting value at \(\zeta=1\).
Thus, for every real \(x>1\),
\[
0\le Q_\lambda(x)-1\le F_\lambda(1)=\frac2{\lambda y}-1.
\]
Letting \(\lambda\uparrow1\) gives \(y\le2\).  If \(y=2\), fix \(x>1\) with \(\phi_1(x)\ne0\).
Continuity in \(\lambda\) and the displayed bounds give \(Q_1(x)=1\), hence
\(-x(1+x)q_0(x)=\phi_1(x)\).  Since \(\phi_1(1)=1\), only finitely many real \(x\) are excluded.
The two polynomials agree at infinitely many points, so
\(\phi_1(z)=-z(1+z)q_0(z)\), or \(u_0=(1-3z)q_0\).
Therefore \(v_0/u_0=-(z-1)^2/(3z-1)\).  Conversely, this transfer function has \(y=2\), as computed in~(a).
\end{proof}

\begin{proof}[Proof of \Cref{cor:yle2}]
The limit method satisfies the hypotheses of \Cref{thm:trace}(b): for small \(\varepsilon\) the roots of
\(\lambda u_\varepsilon-v_\varepsilon\) lie in the open unit disk for \(\lambda\in[\varepsilon^2,1]\), and they depend continuously on \(\varepsilon\) (\(\deg v_0=D\)),
so for \(\lambda\in(0,1)\) the roots of \(\lambda u_0-v_0\) lie in the closed unit disk.
(i) \(y\le2\) by \Cref{thm:trace}(b), and \Cref{cor:ceiling}(i) gives \(\limsup\sqrt\kappa(1-\rhoF)\le\sqrt{y/2}\le1\) and
\(\limsup\sqrt\kappa(1-\rho_q)\le\sqrt y\le\sqrt2\), since the first coordinates \(\zeta_-(\theta',y)\), \(\theta'\ge2\sqrt y\), of the arc of
\Cref{cor:damping} are at most \(\zeta_-(2\sqrt y,y)=\sqrt y\).  (ii) By
\Cref{cor:ceiling}(iii), \(\theta\ge3\) and \(y\ge\theta-1\ge2\); with (i), \(y=2\) and \(\theta=3\), and the equality case of
\Cref{thm:trace}(b) gives the transfer function, for which \(\teff=\frac12\) and \(\Szero\equiv0\) by \Cref{thm:trace}(a).  (iii) By
\Cref{thm:faithful}(a), \(\rhoF<1\) for small \(\varepsilon\), so the family is stable on quadratics (\Cref{lem:envelope}), and
\(\Szero>0\) on \((0,\pi]\).  So \(y\le2\), and \(y\ne2\), since \(y=2\) would force \(\Szero\equiv0\).  The constant is \(\cFF(\theta,y)\) by
\Cref{thm:faithful}(a), and \(\cFF(\theta,y)\le\sqrt{y/2}\) by \Cref{cor:damping}.
\end{proof}

\section{Two-Step Rates and Boundary Layers}
\label{app:classM}

Write \(\chi_\mu(z)=z^2-S_\mu z+Q_\mu\) and \(\chi_L(z)=z^2-S_Lz+Q_L\) with
\(S_\mu=1+\beta-y_\mu(1+\tau)\), \(Q_\mu=\beta-y_\mu\tau\),
\(S_L=1+\beta-y(1+\tau)\), \(Q_L=\beta-y\tau\).

\begin{lemma}[Lag-one cubic]\label{lem:cubic}
For \(z=re^{i\omega}\), let \(H_r(\omega)=\Rea[(1-z)\chi_\mu(z)\overline{\chi_L(z)}]\),
\(E_0=r^4+S_\mu S_Lr^2+Q_\mu Q_L\), \(E_1=S_Lr^2+S_\mu Q_L\) and
\(E_2=S_\mu r^2+Q_\mu S_L\).  Then
\[
H_r=E_0+r^2E_2-r(E_0+E_1+E_2+Q_\mu r^2)\cos\omega+r^2(E_1+Q_\mu+Q_L)\cos2\omega-Q_Lr^3\cos3\omega .
\]
In \(\xi=1-\cos\omega\), this becomes
\(H_r=\sfP_0+\sfP_1\xi+\sfP_2\xi^2+\sfP_3\xi^3\), where
\(\sfP_0=(1-r)\chi_\mu(r)\chi_L(r)\), \(\sfP_3=4Q_Lr^3\),
\(\sfP_2=2r^2(E_1+Q_\mu+Q_L)-12Q_Lr^3\) and
\(\sfP_1=r(E_0+E_1+E_2+Q_\mu r^2)-4r^2(E_1+Q_\mu+Q_L)+9Q_Lr^3\).  At \(\omega=\pi\) (\(\xi=2\)),
\(H_r=(1+r)\chi_\mu(-r)\chi_L(-r)\).
\end{lemma}

\begin{proof}
Using \(z\bar z=r^2\), expand
\[
\chi_\mu\overline{\chi_L}=E_0-E_1z-E_2\bar z+Q_Lz^2+Q_\mu\bar z^2.
\]
Multiplication by \(1-z\), followed by
\(\Rea z^n=\Rea\bar z^n=r^n\cos n\omega\), gives the trigonometric polynomial.
Substitute \(\cos\omega=1-\xi\),
\(\cos2\omega=1-4\xi+2\xi^2\) and
\(\cos3\omega=1-9\xi+12\xi^2-4\xi^3\) and collect powers of \(\xi\).
The evaluations at \(z=r\) and \(z=-r\) give \(\sfP_0\) and the endpoint formula.
\end{proof}
\begin{lemma}[Mixed test]\label{lem:mixed}
Let \((\beta,\tau,h)\) be a tuning of \(\cM\), \(r>0\), and for \(z=re^{i\omega}\) put
\(H^{\rm mix}_r(\omega)=\Rea[(1-e^{i\omega})\chi_\mu(z)\overline{\chi_L(z)}]\).  Then
\[
\Ima\bigl(\chi_\mu(z)\overline{\chi_L(z)}\bigr)=(y-y_\mu)\,r\sin\omega\,
\bigl[(1+\tau)(r^2-\beta)+\tau(1+\beta)-2\tau r\cos\omega\bigr],
\]
and \(H^{\rm mix}_r=\xi\,\mathsf M_r(\xi)\), \(\xi=1-\cos\omega\), with the quadratic
\(\mathsf M_r(\xi)=E_0-(E_1+E_2)r(1-\xi)+(Q_\mu+Q_L)r^2(1-4\xi+2\xi^2)
+(2-\xi)(y-y_\mu)r[(1+\tau)(r^2-\beta)+\tau(1+\beta)-2\tau r(1-\xi)]\) (notation of
\Cref{lem:cubic}).  If \(K\ge3\), \(\omega_0=2\pi/K\),
\(\Ima(\chi_\mu\overline{\chi_L})(re^{i\omega_0})>0\) and \(H^{\rm mix}_r(\omega_0)\le0\), then \((\delta_{\omega_0},0)\)
is feasible at \(t=\log r\).
\end{lemma}

\begin{proof}
Writing \(\chi_\lambda=\Omega+y_\lambda N_\tau\) gives
\(\Ima(\chi_\mu\overline{\chi_L})=(y-y_\mu)\Ima(\Omega\overline{N_\tau})\).
Substitution of \(\Omega=(z-1)(z-\beta)\), \(N_\tau=(1+\tau)z-\tau\) and
\(z\bar z=r^2\) gives the stated imaginary part.  The product expansion in \Cref{lem:cubic} gives the real part, so
\[
H^{\rm mix}_r=\xi\Rea(\chi_\mu\overline{\chi_L})
+\sin\omega\,\Ima(\chi_\mu\overline{\chi_L})=\xi\mathsf M_r,
\]
using \(\sin^2\omega=\xi(2-\xi)\).
At \(t=\log r\), \eqref{eq:factor} gives
\(g-2is=-\mathfrak c\chi_\mu\overline{\chi_L}\), \(\mathfrak c>0\).
Thus \(s>0\), and \(g-2is=Re^{-i\Theta_r}\) with \(R>0\), \(\Theta_r\in(0,\pi)\).
By \Cref{lem:angle} and the sign hypothesis,
\[
0\le-\mathfrak c H^{\rm mix}_r(\omega_0)
=2R\sin\frac{\omega_0}2\sin\frac{\omega_0-2\Theta_r}2.
\]
Since \((\omega_0-2\Theta_r)/2\in(-\pi,\pi/2)\), this implies
\(\Theta_r\le\omega_0/2=\pi/K\).
\Cref{lem:single}(a), with \(j=1\), gives the claimed feasible atom.
\end{proof}

\begin{proof}[Proof of \Cref{prop:Ltau}(a)]
The coefficients of \(L_\tau\) are
\begin{equation}\label{eq:Acoef}
\begin{gathered}
\mathsf A_0=2(1-\beta)(y-y_\mu)+y_\mu y,\qquad \mathsf A_2=4(\beta-\tau y),\\
\mathsf A_1=2(1-\beta)^2-2y_\mu(\beta+\beta\tau+\tau)+y(6\tau-2-2\beta\tau)+2\tau(1+\tau)y_\mu y ;
\end{gathered}
\end{equation}
at \(\beta=1\), \(y_\mu=0\) they are \(0\), \(2y(2\tau-1)\) and \(4(1-\tau y)\).  On \(|z|=1\), put
\(\mathcal R_0=\Rea(\chi_\mu\overline{\chi_L})\) and
\(q=(y-y_\mu)[(1-\beta)+2\tau\xi]\).
By~\eqref{eq:imab},
\[
H_1=\xi\mathcal R_0+q\sin^2\omega
=\xi\bigl[\mathcal R_0+(2-\xi)q\bigr],
\qquad L_\tau=\mathcal R_0+(2-\xi)q.
\]
Substituting \(\chi_\lambda=(z-1)(z-\beta)+y_\lambda((1+\tau)z-\tau)\),
\(z\bar z=1\) and \(\Rea z=1-\xi\) and collecting powers of \(\xi\)
gives~\eqref{eq:Acoef}.

Stability gives \(\chi_\mu(-1)\chi_L(-1)>0\) (\Cref{lem:qstab}), hence
\(L_\tau(2)=\chi_\mu(-1)\chi_L(-1)>0\).  Also \(L_\tau(0)=\mathsf A_0>0\).  For the
angle, let \(\omega\in(0,\pi]\) with \(g-2is\ne0\) at \((0,\omega)\) (for \(\omega<\pi\) this holds by~\eqref{eq:imab}, and at \(\omega=\pi\) for stable
tunings).  By~\eqref{eq:factor} and~\eqref{eq:imab}, \(2s=\mathfrak c q\sin\omega\ge0\) with \(\mathfrak c=4/(h^2(L-\mu)^2)\), so \(\Theta\in[0,\pi]\) and
\(\Theta-\frac\omega2\in[-\frac\pi2,\pi)\).  By \Cref{lem:angle}, \(H_1=\xi L_\tau\) is a positive multiple of \(\sin(\Theta-\frac\omega2)\); hence \(L_\tau(\xi)>0\) if
and only if \(\Theta>\omega/2\), and \(L_\tau(\xi)=0\) if and only if \(\Theta=\omega/2\).  The last claim is \Cref{prop:lagone}(b), since \(H_1=\xi L_\tau\).
\end{proof}
\begin{lemma}[Resonance on a band]\label{lem:band}
\leavevmode
\begin{enumerate}[label=\textup{(\alph*)}]
\item Let \(\mathsf s\in[0,1)\) and \(0<\omega_a<\omega_b\le\frac{2\pi}3\) with \(2\pi(1-\mathsf s)/\omega_a-2\pi(1+\mathsf s)/\omega_b\ge1\).  Then there is an
integer \(K\ge2\pi(1+\mathsf s)/\omega_b\ge3\) with \([(1-\mathsf s)\frac{2\pi}K,(1+\mathsf s)\frac{2\pi}K]\subset[\omega_a,\omega_b]\).
\item Let a tuning of \(\cM\) be stable, \(K\ge3\), and \(L_\tau(1-\cos\frac{2\pi}K)<0\).  Then \(\Theta(2\pi/K)<\pi/K\), the method has trajectories
from rest that approach \(K\)-periodic orbits in \(\R^3\), and \(\rhoz>1\).  In particular this holds for some \(K\ge2\pi/\omega_b\) if
\(L_\tau(1-\cos\omega)<0\) for all \(\omega\in[\omega_a,\omega_b]\), with \(\omega_a,\omega_b\) as in~(a) for \(\mathsf s=0\).
\end{enumerate}
\end{lemma}

\begin{proof}
(a) The interval \([2\pi(1+\mathsf s)/\omega_b,\,2\pi(1-\mathsf s)/\omega_a]\) has length at least \(1\), so it contains an integer \(K\); then
\((1+\mathsf s)2\pi/K\le\omega_b\), \((1-\mathsf s)2\pi/K\ge\omega_a\), and \(K\ge2\pi/\omega_b\ge3\).
(b) By \Cref{prop:Ltau}(a), \(\Theta(2\pi/K)<\pi/K\), and \Cref{prop:cycles}, with \(j=1\), gives the orbits and \(\rhoz>1\).  For the
last claim, (a) with \(\mathsf s=0\) gives \(K\) with \(2\pi/K\in[\omega_a,\omega_b]\).
\end{proof}

\begin{lemma}[Accelerated tunings]\label{lem:accel}
Let a stable tuning of \(\cM\) on \([\mu,L]\), \(L=1\), \(\mu=\varepsilon^2\), \(0<\varepsilon<1\), have
\(\rho_q\le1-c\varepsilon\), \(c>0\).
\begin{enumerate}[label=\textup{(\alph*)}]
\item \(y_\mu\ge c^2\varepsilon^2\), \(c^2\le y<4\), \(\tau<2/c^2\) and \(1-\beta<\max(4,8/c)\,\varepsilon\).
\item If such tunings are given along a sequence \(\varepsilon\to0^+\), with
parameters \(\beta=1-\theta_\varepsilon\varepsilon\), \(h=y_\varepsilon\) and \(\tau_\varepsilon\), a subsequence has \((\theta_\varepsilon,y_\varepsilon,\tau_\varepsilon)\to(\theta,y,\tau)\) with \(\theta\ge2c\),
\(y\ge c^2\), \(\tau\ge0\) and \(y(1+2\tau)\le4\); along it the tunings form a family~\eqref{eq:scaling}.
\end{enumerate}
\end{lemma}

\begin{proof}
(a) As in the proof of \Cref{thm:price}, \(y_\mu\ge(1-\rho_q)^2\ge c^2\varepsilon^2\) and \(y=\kappa y_\mu\ge c^2\).
Stability gives \(y(1+2\tau)<2(1+\beta)<4\), so \(\tau<(4-y)/(2y)<2/c^2\).  Put
\(\delta=1-\beta\).  For \(\eta\ge0\),
\(\chi_\mu(1-\eta)=\eta^2-\eta(\delta+y_\mu(1+\tau))+y_\mu\le\eta^2-\eta\delta+y_\mu\).  If
\(y_\mu<\delta^2/4\), then \(\chi_\mu(1-2y_\mu/\delta)<0<\chi_\mu(1)\), so \(\chi_\mu\) has a root
in \((1-2y_\mu/\delta,1)\), and \(c\varepsilon<2y_\mu/\delta\), i.e.\
\(\delta<2y_\mu/(c\varepsilon)<8\varepsilon/c\).  Otherwise \(\delta\le2\sqrt{y_\mu}<4\varepsilon\).
(b) The roots of \(\chi_\mu\) have moduli at most \(\rho_q\le1-c\varepsilon\), and their product is \(Q_\mu=\beta-y_\mu\tau\).  Hence
\(1-\theta_\varepsilon\varepsilon-y_\varepsilon\tau\varepsilon^2=Q_\mu\le(1-c\varepsilon)^2\), and \(y_\varepsilon\tau<8/c^2\) by (a); so \(\theta_\varepsilon\ge2c-(c^2+8/c^2)\varepsilon\).  By (a),
\((\theta_\varepsilon,y_\varepsilon,\tau_\varepsilon)\) stays in a bounded set, so a subsequence converges; the limits satisfy \(\theta\ge2c>0\), \(y\ge c^2>0\), \(\tau\ge0\),
and \(y(1+2\tau)\le4\) by the stability condition \(y_\varepsilon(1+2\tau_\varepsilon)<2(1+\beta)\) (\Cref{lem:qstab}).
\end{proof}

\begin{proof}[Proof of \Cref{thm:wall}]
If the first claim failed, \Cref{lem:accel} would give a subsequence
with \(\beta\to1\), \(y_\mu\to0\),
\(y_\varepsilon\to y>0\), and
\(\tau_\varepsilon\to\tau<1/2\).
By \Cref{prop:Ltau}, its lag-one polynomials converge coefficientwise to
\[
L_0(\xi)=2y(2\tau-1)\xi+4(1-\tau y)\xi^2.
\]
The negative linear coefficient makes \(L_0<0\) on a sufficiently
small interval \((0,\xi_0)\).  Choose one integer \(K\ge6\) with
\(1-\cos(2\pi/K)<\xi_0\).  Coefficientwise convergence gives
\(L_{\tau_\varepsilon}(1-\cos(2\pi/K))<0\) for all sufficiently
small \(\varepsilon\) on the subsequence.  \Cref{prop:Ltau,prop:cycles}
then give the periodic orbits and \(\rhoz>1\), contradicting the
choice of that subsequence.

If an accelerated family on \(\cF_{\mu,L}\) had
\(\liminf\tau_\kappa<1/2\), a subsequence with
\(\limsup\tau_\kappa<1/2\) would satisfy the first claim, contrary
to \(\rhoz\le\rhoF<1\).
\end{proof}

\begin{proposition}[Range of resonant periods]\label{prop:periods}
Let \(\tau\in[0,\frac12)\), \(\theta,y>0\), \(y(1+2\tau)<4\), and consider
\(\beta=1-\theta\varepsilon\), \(h=y\).  Then \(\tau y<1\).  Put
\(\xi_*=\frac{y(1-2\tau)}{2(1-\tau y)}\in(0,2)\) and
\(\xi_{\min}=\frac{\theta\varepsilon}{1-2\tau}\).  For every
\(\delta\in(0,\frac12)\) and all small \(\varepsilon\),
\(L_\tau<0\) on \([(1+\delta)\xi_{\min},(1-\delta)\xi_*]\), while \(L_\tau>0\)
on \((0,(1-\delta)\xi_{\min}]\) and on \([(1+\delta)\xi_*,2]\).
Consequently the method has \(K\)-periodic orbits for every \(K\) with
\(1-\cos(2\pi/K)\in[(1+\delta)\xi_{\min},(1-\delta)\xi_*]\); the largest
such \(K\) is
\(\pi\bigl(2(1-2\tau)/((1+\delta)\theta)\bigr)^{1/2}\varepsilon^{-1/2}(1+o(1))\),
and the number of such \(K\) grows like \(\varepsilon^{-1/2}\).
\end{proposition}

\begin{proof}[Proof of \Cref{prop:periods}]
Stability implies \(\tau y<4\tau/(1+2\tau)<1\), and
\(\xi_*<2\).  The coefficients of \Cref{prop:Ltau} satisfy
\[
\mathsf A_0=2\theta y\varepsilon+O(\varepsilon^2),\quad
\mathsf A_1=-2y(1-2\tau)+O(\varepsilon),\quad
\mathsf A_2=4(1-\tau y)+O(\varepsilon)>0.
\]
The limit polynomial has the two simple roots \(0\) and \(\xi_*\).
Its perturbed roots are therefore
\[
\xi_1=\frac{\theta\varepsilon}{1-2\tau}+O(\varepsilon^2)
=\xi_{\min}+O(\varepsilon^2),\qquad
\xi_2=\xi_*+O(\varepsilon).
\]
For small \(\varepsilon\), the polynomial is negative between these
roots and positive outside; the claimed intervals follow for every
fixed \(\delta>0\).  \Cref{lem:band}(b) gives the periodic orbits.
Finally,
\(1-\cos(2\pi/K)=2\pi^2K^{-2}(1+O(K^{-2}))\)
gives the largest admitted period and the asserted count.
\end{proof}

For the remainder of this appendix, use~\eqref{eq:scaling} with \(L=1\),
\(\mu=\varepsilon^2\), and \(\Delta^\varepsilon_c=c^2-\theta_\varepsilon c+y_\varepsilon\).

\begin{lemma}[Lag-one and mixed coefficients in the accelerated scaling]\label{lem:Pscaled}
Let \(r=1-c\varepsilon\) with \(c\in\R\).  The following expansions hold uniformly for \(c\) and the parameters in compact sets.
\begin{enumerate}[label=\textup{(\alph*)}]
\item The coefficients of \Cref{lem:cubic} satisfy
\(\sfP_0=c\varepsilon\,\chi_\mu(r)\chi_L(r)\) with
\(\chi_\mu(r)=\varepsilon^2\bigl(\Delta_c^\varepsilon-y_\varepsilon(1+\tau_\varepsilon)c\varepsilon\bigr)\) and
\(\chi_L(r)=y_\varepsilon+O(\varepsilon)\), and
\[
\sfP_1=2y_\varepsilon(\theta_\varepsilon-3c)\varepsilon+O(\varepsilon^2),\qquad
\sfP_2=2y_\varepsilon(2\tau_\varepsilon-1)+O(\varepsilon),\qquad
\sfP_3=4(1-\tau_\varepsilon y_\varepsilon)+O(\varepsilon).
\]
\item For \(z=re^{i\omega}\),
\(\Ima(\chi_\mu\overline{\chi_L})=(y_\varepsilon-y_\mu)r\sin\omega\,[\varepsilon(\theta_\varepsilon-2c)+(1+\tau_\varepsilon)c^2\varepsilon^2+2\tau_\varepsilon r\xi]\),
and the quadratic of \Cref{lem:mixed} is \(\mathsf M_r=\mathsf Q_0+\mathsf Q_1\xi+\mathsf Q_2\xi^2\) with
\(\mathsf Q_0=2y_\varepsilon(\theta_\varepsilon-2c)\varepsilon+O(\varepsilon^2)\),
\(\mathsf Q_1=2y_\varepsilon(2\tau_\varepsilon-1)+O(\varepsilon)\) and \(\mathsf Q_2=4(1-\tau_\varepsilon y_\varepsilon)+O(\varepsilon)\).
\item Regard \((\theta',y',\tau')=(\theta_\varepsilon,y_\varepsilon,\tau_\varepsilon)\) as independent parameters.  Then
\[
\begin{gathered}
\sfP_1=\varepsilon \dot{\mathsf P}_1+O(\varepsilon^2),\quad \sfP_2=2y'(2\tau'-1)+\varepsilon \dot{\mathsf P}_2+O(\varepsilon^2),\quad \sfP_3=4(1-\tau'y')+\varepsilon \dot{\mathsf P}_3+O(\varepsilon^2),\\
\mathsf Q_0=\varepsilon \dot{\mathsf Q}_0+O(\varepsilon^2),\quad \mathsf Q_1=2y'(2\tau'-1)+\varepsilon \dot{\mathsf Q}_1+O(\varepsilon^2),\quad \mathsf Q_2=4(1-\tau'y')+\varepsilon \dot{\mathsf Q}_2+O(\varepsilon^2),
\end{gathered}
\]
where \(\dot{\mathsf P}_k,\dot{\mathsf Q}_k\) are polynomials in \((c,\theta',y',\tau')\), and at \((y',\tau')=(2,\frac12)\)
\[
\dot{\mathsf P}_1=4(\theta'-3c),\quad \dot{\mathsf P}_2=2(\theta'+2c),\quad \dot{\mathsf P}_3=-4\theta',\qquad \dot{\mathsf Q}_0=4(\theta'-2c),\quad \dot{\mathsf Q}_1=2(\theta'+c),\quad \dot{\mathsf Q}_2=-4\theta' .
\]
\end{enumerate}
\end{lemma}

\begin{proof}
Treat \(p=(c,\theta',y',\tau')\) as independent parameters and substitute
\(r=1-c\varepsilon\), \(\beta=1-\theta'\varepsilon\) and
\(y_\mu=y'\varepsilon^2\) in the coefficient formulas of
\Cref{lem:cubic,lem:mixed}.  Each of the six coefficients \(F\) is polynomial, so
\[
F(\varepsilon,p)=F(0,p)+\varepsilon\partial_\varepsilon F(0,p)
+\varepsilon^2R_F(\varepsilon,p),
\]
with \(R_F\) polynomial.  Collecting the constant and linear powers of
\(\varepsilon\) gives the displayed coefficients in (a)--(c);
setting \((y',\tau')=(2,1/2)\) in the linear coefficients gives the six displayed derivative values.
For every compact parameter set \(K\), both the first coefficients and
\(R_F\) are bounded on \([-1,1]\times K\), proving all the stated uniform remainders.

The exact characteristic formulas in (a) follow from
\(r-1=-c\varepsilon\), \(r-\beta=(\theta'-c)\varepsilon\) and
\(N_{\tau'}(r)=1-(1+\tau')c\varepsilon\).
For (b), the imaginary-part bracket of \Cref{lem:mixed} reduces exactly to
\(\varepsilon(\theta'-2c)+(1+\tau')c^2\varepsilon^2+2\tau'r\xi\).
The same coefficient extraction then gives the assertions about \(\mathsf M_r\).
\end{proof}

\begin{lemma}[Lag-one certificates in the accelerated scaling]\label{lem:lagscaled}
Let tunings of \(\cM\) with \(L=1\), \(\mu=\varepsilon^2\), \(\beta=1-\theta_\varepsilon\varepsilon\), \(h=y_\varepsilon\), \(\tau=\tau_\varepsilon\) and
\((\theta_\varepsilon,y_\varepsilon,\tau_\varepsilon)\to(\theta,y,\tau)\), \(\theta,y>0\), be given for small \(\varepsilon\) (or along a sequence \(\varepsilon\to0\)).  Let
\(0<c<\cQ(\theta,y)\) and \(r=1-c\varepsilon\).  Then, for small \(\varepsilon\), \(\sfP_0>0\) and the roots of \(\chi_\mu\) lie in \(|z|<r\).  If moreover
\(\sfP_1+\sfP_2\xi+\sfP_3\xi^2>0\) on \([0,2]\) for small \(\varepsilon\), then \(H_r>0\) on \([0,2]\), and \(\rho_q<r\), \(r_c<r\) and \(\rhoF<r\); the last bound is
certified by the lag-one criterion.
\end{lemma}

\begin{proof}
Since \(c<\cQ\), \Cref{lem:cQprop}(a) gives
\(\Delta_c>0\), and \Cref{lem:Pscaled}(a) yields
\(\sfP_0=c\varepsilon^3(\Delta_c+o(1))(y+o(1))>0\).
Exactly,
\[
\chi_\mu(1-\varepsilon\zeta)
=\varepsilon^2(\zeta^2-\theta''_\varepsilon\zeta+y_\varepsilon),
\qquad
\theta''_\varepsilon=
\theta_\varepsilon+y_\varepsilon(1+\tau_\varepsilon)\varepsilon\to\theta.
\]
The two \(\zeta\)-roots are bounded, and their minimum real part tends to \(\cQ>c\).
Hence, for some \(\eta>0\) and \(K<\infty\),
\(|1-\varepsilon\zeta|^2\le1-2(c+\eta)\varepsilon
+K^2\varepsilon^2<(1-c\varepsilon)^2\) for small \(\varepsilon\).
This puts the roots of \(\chi_\mu\) inside \(|z|<r\).
The additional coefficient hypothesis gives
\(H_r=\sfP_0+\xi(\sfP_1+\sfP_2\xi+\sfP_3\xi^2)>0\)
on \([0,2]\); \Cref{lem:separate} then gives all three rate bounds.
\end{proof}

\begin{lemma}[The sign-alternating mode]\label{lem:altmode}
\leavevmode
\begin{enumerate}[label=\textup{(\alph*)}]
\item Let tunings of \(\cM\) with \(L=1\) be given along \(\varepsilon\to0\), with \(\chi_L(-1)>0\), \(\chi_L(-1)\to0\) and \(S_L\to S_\infty\in(-2,0)\).  Then, for small
\(\varepsilon\), \(\chi_L\) has the real roots \(z_L=-1+\delta\) and \(z_2=1+S_L-\delta\), with
\(\delta=\frac{\chi_L(-1)}{2+S_L}\bigl(1+O(\chi_L(-1))\bigr)>0\) and \(z_2\to1+S_\infty\in(-1,1)\); so \(|z_2|<|z_L|=1-\delta\), and \(\rho_q\ge1-\delta\).  From rest with \(x_0=1\)
on \(\frac12x^2\) on \(\R\), \(\limsup_N|x_N|^{1/N}=1-\delta\), so \(\rhoz\ge1-\delta\); the dominant mode \(z_L^k\) alternates in sign.
\item Let a family~\eqref{eq:scaling} have \(\tau\ge\frac12\), \(y(1+2\tau)=4\) and \(\varsigma_\varepsilon\to\varsigma\in[2\theta,\infty]\) (notation~\eqref{eq:margin}), and be stable for
small \(\varepsilon\).  Then (a) applies, with \(S_\infty=-\frac2{1+2\tau}\), \(z_2\to\frac{2\tau-1}{1+2\tau}\in[0,1)\) and
\(\delta/\varepsilon\to(\varsigma-2\theta)\frac{1+2\tau}{4\tau}\), which is \(\cpi\) if \(\varsigma>2\theta\) (\(\infty\) if \(\varsigma=\infty\)) and \(0\) if \(\varsigma=2\theta\).
\end{enumerate}
\end{lemma}

\begin{proof}
(a) Put \(e=\chi_L(-1)\) and \(\sigma=2+S_L\to2+S_\infty\in(0,2)\).  Then \(\chi_L(-1+\delta)=e-\sigma\delta+\delta^2\), whose roots for \(e<\sigma^2/4\) are
\(\delta=\frac{2e}{\sigma+\sqrt{\sigma^2-4e}}=\frac e\sigma(1+O(e))>0\) and \(\sigma-\delta\).  So \(z_L=-1+\delta\) and \(z_2=-1+\sigma-\delta=1+S_L-\delta\to1+S_\infty\in(-1,1)\), and
\(|z_2|<1-\delta=|z_L|\) for small \(\varepsilon\).  Hence \(\rho(\chi_L)=1-\delta\), and \Cref{lem:restquad} at \(\lambda=L\) gives
\(\rhoz\ge1-\delta\), with dominant mode \(z_L^k\).
(b) By~\eqref{eq:margin}, \(\chi_L(-1)=(\varsigma_\varepsilon-2\theta_\varepsilon)\varepsilon\), which is positive by stability (\Cref{lem:qstab}) and tends to \(0\), since
\(\varsigma_\varepsilon\varepsilon=4-y_\varepsilon(1+2\tau_\varepsilon)\to0\).  Further \(S_L=1+\beta-y_\varepsilon(1+\tau_\varepsilon)\to2-\frac{4(1+\tau)}{1+2\tau}=-\frac2{1+2\tau}\in[-1,0)\), so
\(2+S_L\to\frac{4\tau}{1+2\tau}\), \(z_2\to1+S_\infty=\frac{2\tau-1}{1+2\tau}\), and \(\delta/\varepsilon=\frac{\varsigma_\varepsilon-2\theta_\varepsilon}{2+S_L}(1+o(1))\to(\varsigma-2\theta)\frac{1+2\tau}{4\tau}\), also when
\(\varsigma=\infty\).
\end{proof}

\begin{lemma}[Spread witnesses for \(\cM\)]\label{lem:spreadwit}
Let stable tunings of \(\cM\) with \(L=1\), \(\mu=\varepsilon^2\), \(\beta=1-\theta_\varepsilon\varepsilon\), \(h=y_\varepsilon\) and \(\tau=\tau_\varepsilon\) be given along a sequence
\(\varepsilon\to0\), with \((\theta_\varepsilon,y_\varepsilon,\tau_\varepsilon)\to(\theta,y,\tau)\), \(\theta,y,\tau>0\).  Let \(0<c<\theta/2\), \(r=1-c\varepsilon\), \(t=\log r\), and fix
\(\mathsf s\in(0,\mathsf s_1]\) (\Cref{lem:spreadhalf}).  For each \(\varepsilon\) let \(K\ge3\), \(\omega_0=2\pi/K\), \(B=[(1-\mathsf s)\omega_0,(1+\mathsf s)\omega_0]\),
\(\xi_{\max}=1-\cos((1+\mathsf s)\omega_0)\), and \(\Xi_B=\max_B|\chi_\mu\overline{\chi_L}|\) on \(|z|=r\).  Suppose that, as \(\varepsilon\to0\),
\[
\varepsilon=o(\xi_{\max}^{3/2}),\qquad \max_B|\mathsf M_r|=o(\xi_{\max}^{3/2}),\qquad H_r\le-\lambda_0\varepsilon \Xi_B\ \text{ on }B,
\]
with a constant \(\lambda_0>0\).  Then, for all small \(\varepsilon\), there is a pair with finitely many atoms in \(B\) and
\(|C|<\bar g\) feasible at \(t\).  For some finite \(J=J_\varepsilon\), it is realized from rest by
\(f\in\cF_{\varepsilon^2,1}(\R^{2J+1})\) and \(x_0\ne0\), with
\(\norm{x_N}\ge c_0(1-c\varepsilon)^N\norm{x_0}\) for every \(N\ge0\), where \(c_0>0\)
may depend on \(\varepsilon\); in particular \(\rhoz\ge1-c\varepsilon\).
\end{lemma}

\begin{proof}
Put \(\mathfrak c=4/(y_\varepsilon^2(1-\varepsilon^2)^2)\); by~\eqref{eq:factor}, \(g-2is=-\mathfrak c\chi_\mu\overline{\chi_L}\) at \(t\).  On \(B\), \(\omega\le\frac76\cdot\frac{2\pi}3<\pi\), so
\(\sin\omega>0\) and \(2-\xi\ge\frac15\), and \(\xi\ge\frac12\xi_{\max}\), because \(\sin\frac{(1-\mathsf s)\omega_0}2\ge\frac{1-\mathsf s}{1+\mathsf s}\sin\frac{(1+\mathsf s)\omega_0}2\) and
\((\frac{1-\mathsf s}{1+\mathsf s})^2\ge\frac12\).

\emph{Signs.}  By \Cref{lem:Pscaled}(b), \(\Ima(\chi_\mu\overline{\chi_L})=(y_\varepsilon-y_\mu)r\sin\omega\,\iota_\varepsilon\) with
\(\iota_\varepsilon=\varepsilon(\theta_\varepsilon-2c)+(1+\tau_\varepsilon)c^2\varepsilon^2+2\tau_\varepsilon r\xi\), and \(\iota_\varepsilon>0\) for small \(\varepsilon\), as \(\theta>2c\).  By \Cref{lem:mixed},
\(\xi\mathsf M_r=\xi\Rea(\chi_\mu\overline{\chi_L})+\sin\omega\,\Ima(\chi_\mu\overline{\chi_L})\), and \(\sin^2\omega=\xi(2-\xi)\); hence
\[
\Rea(\chi_\mu\overline{\chi_L})=\mathsf M_r-(2-\xi)(y_\varepsilon-y_\mu)r\,\iota_\varepsilon\le|\mathsf M_r|-2(2-\xi)(y_\varepsilon-y_\mu)\tau_\varepsilon r^2\xi\le-\mathsf a\xi,\qquad \mathsf a=\tfrac{y\tau}{10},
\]
for small \(\varepsilon\), because \(2(2-\xi)(y_\varepsilon-y_\mu)\tau_\varepsilon r^2\to2(2-\xi)y\tau\ge\frac25y\tau\) uniformly on \(B\), and \(|\mathsf M_r|=o(\xi_{\max})\le\frac{y\tau}{10}\xi_{\max}\le\frac{y\tau}5\xi\).
Thus \(g=-\mathfrak c\Rea(\chi_\mu\overline{\chi_L})>0\) and \(2s=\mathfrak c\Ima(\chi_\mu\overline{\chi_L})>0\) on \(B\): \(g-2is=Re^{-i\Theta}\) with \(\Theta\in(0,\frac\pi2)\)
continuous and \(R=\mathfrak c|\chi_\mu\overline{\chi_L}|\ge\mathfrak c \mathsf a\xi\).

\emph{Defect.}  By \Cref{lem:angle,lem:mixed}, \(2R\sin\frac\omega2\sin\frac{\omega-2\Theta}2=\Rea[(1-e^{i\omega})(g-2is)]=-\mathfrak c\xi\mathsf M_r\).  As
\(\sin\frac\omega2\ge\frac\omega\pi\) and \(|\omega-2\Theta|<\pi\),
\(|\omega-2\Theta|\le\pi|\sin\frac{\omega-2\Theta}2|\le\pi^2|\mathsf M_r|/(2\mathsf a\omega)\), so \(\mathsf d\le\pi^2\max_B|\mathsf M_r|/(\mathsf a\omega_0)\).

\emph{Hypotheses of \Cref{lem:spreadhalf}.}  Since \(c_{\mathsf k}\le\frac12\), \(\mathsf s\le\frac16\) and \(\omega_0\le\frac{2\pi}3\),
\(\mathsf m=\min(\frac14,\frac{\omega_0^2}6,c_{\mathsf k}\mathsf s^2)\ge m_{\mathsf s}\omega_0^2\) with \(m_{\mathsf s}=9c_{\mathsf k}\mathsf s^2/(4\pi^2)\).  Moreover
\(\xi_{\max}\le\frac12(1+\mathsf s)^2\omega_0^2\le\omega_0^2\), so \(\mathsf m\ge m_{\mathsf s}\xi_{\max}\) and \(\mathsf w=\mathsf s\omega_0\ge\mathsf s\xi_{\max}^{1/2}\).  Hence
\(\mathsf d\le\pi^2\max_B|\mathsf M_r|/(\mathsf a\xi_{\max}^{1/2})=o(\xi_{\max})\), so \(\mathsf d\le\mathsf m/2\) and \(\mathsf d\to0\); and \(|t|\le2c\varepsilon\), so
\(|t|/(\mathsf m\mathsf w)\le2c\varepsilon/(m_{\mathsf s}\mathsf s\xi_{\max}^{3/2})\to0\).  Thus~\eqref{eq:spreadhyp} holds for every fixed \(\eta\in(0,1)\) and small \(\varepsilon\).

\emph{The lag-one constraint.}  On \(B\), \(\Phi^{(1)}_t=-\mathfrak c H_r\ge\mathfrak c\lambda_0\varepsilon \Xi_B\), so \(\int_B\Phi^{(1)}_t\,d\nu\ge\mathfrak c\lambda_0\varepsilon \Xi_B\) for every probability measure
\(\nu\) on \(B\); while \(\bar g\le\max_Bg\le\mathfrak c \Xi_B\) and \(e^t\sh|t|=\frac12(1-r^2)\le c\varepsilon\).  With \(\eta=\min(\frac12,\frac{\lambda_0}{2c})\),
\(\eta\,e^t\sh|t|\,\bar g\le\frac12\lambda_0\varepsilon\mathfrak c \Xi_B<\mathfrak c\lambda_0\varepsilon \Xi_B\), as \(\Xi_B\ge \mathsf a\xi_{\max}/2>0\); so~\eqref{eq:spreadl1} holds.  \Cref{lem:spreadhalf} gives a
pair with finitely many atoms in \(B\), feasible at \(t\), with \(|C'|<\bar g'\), and \Cref{lem:startup} gives \(f\in\cF_{\varepsilon^2,1}(\C^J\oplus\R)\) and \(x_0\)
with \(\norm{x_N}\ge c_0r^N\norm{x_0}\); so \(\rhoz\ge1-c\varepsilon\).
\end{proof}

\begin{proof}[Proof of \Cref{prop:Ltau}(b)]
\(u_0=N_\tau\) gives \(u_0'(1)=1+\tau\), and \(q_0=-1/y\) is constant, so \(\teff=\tau\).  Since \(\phi_\lambda=\chi_\lambda/h\), the
lag-one symbol here is that of \Cref{sec:closure} divided by \(y_\varepsilon^2\); at \(\varepsilon=0\), \Cref{lem:Pscaled}(a)
gives \(\sfP_0=\sfP_1=0\), \(\sfP_2=2y(2\tau-1)\), \(\sfP_3=4(1-\tau y)\), so \(\Szero=\xi^2(\sfP_2+\sfP_3\xi)/y^2\).  The affine function
\(\sfP_2+\sfP_3\xi\) is positive on \((0,2]\) if and only if \(\sfP_2\ge0\) and \(\sfP_2+2\sfP_3=8-2y(1+2\tau)>0\), and nonnegative on
\([0,2]\) if and only if \(\sfP_2\ge0\) and \(\sfP_2+2\sfP_3\ge0\).  Finally
\(u_0-2v_0=\frac2y\bigl(z^2-(2-\frac y2(1+\tau))z+1-\frac{y\tau}2\bigr)\), whose roots lie in the open unit disk by the
Schur--Cohn conditions \(|1-\frac{y\tau}2|<1\), i.e.\ \(0<y\tau<4\), and \(|2-\frac y2(1+\tau)|<2-\frac{y\tau}2\), i.e.\ \(\frac y2>0\) and
\(4-\frac y2(1+2\tau)>0\), all implied by \(y(1+2\tau)<4\) and \(\tau\ge\frac12\).
\end{proof}

\begin{lemma}[Quadratic constant up to the boundary]\label{lem:quadboundary}
Let a family~\eqref{eq:scaling} have \(\tau\ge\frac12\), \(y(1+2\tau)=4\) and \(\varsigma_\varepsilon\to\varsigma\in(2\theta,\infty]\), with
\(\varsigma_\varepsilon\) and \(\cpi\) as in~\eqref{eq:margin}--\eqref{eq:cpi}.  Then
\(\sqrt\kappa\,(1-\rho_q)\to\min(\cQ(\theta,y),\cpi)\).
\end{lemma}

\begin{proof}
Write \(\ell=h\lambda\in[y_\mu,y_\varepsilon]\) and
\(P_\ell(z)=(z-1)(z-\beta)+\ell N_{\tau_\varepsilon}(z)\), so
\(\chi_\lambda=P_{h\lambda}\).
The proof of \Cref{lem:cq} already bounds the positive real roots by
the slow endpoint contribution
\(1-\varepsilon\cQ(\theta,y)+o(\varepsilon)\), and the complex
roots by \(\sqrt\beta=1-\theta\varepsilon/2+o(\varepsilon)\).
Both estimates hold up to the limiting stability boundary.

For a negative real root \(w\) of \(P_\ell\),
\[
\chi_L(w)=P_\ell(w)+(y_\varepsilon-\ell)N_{\tau_\varepsilon}(w)
\le0,
\]
because \(N_{\tau_\varepsilon}(w)<0\).
By \Cref{lem:altmode}(b), \(\chi_L\) has two real roots for small
\(\varepsilon\), the smaller being
\(z_L=-1+\delta_\varepsilon\).  A monic quadratic is nonpositive
only between its roots, so \(w\ge z_L\).
Thus the largest modulus of a negative real mode is attained at
\(\lambda=L\), and
\(\delta_\varepsilon/\varepsilon\to\cpi\), also when
\(\cpi=\infty\).  Since \(\cQ\le\theta/2\), combining the three
types of roots gives the claimed minimum.
\end{proof}

\begin{proof}[Proof of \Cref{thm:B}]
(a) The claim for some \(\delta\) implies it for every larger \(\delta\), so \(\delta\) may be taken small.  The family is a second-order
family with slow data \((\theta,y)\) (\Cref{sec:beyond}).
\emph{Case \(\theta/3<\cQ\).}  Let \(0<\delta<\cQ-\theta/3\) and \(c=\theta/3+\delta=\cFF+\delta\).  \Cref{lem:flowspirals} gives \(J\), depending on
\((\theta,y,\delta)\) only, and for small \(\varepsilon\) a pair with \(J\) atoms that is feasible at \(t=-c\varepsilon\) and has \(|C|<\bar g\).
\Cref{lem:startup} with \(t<0\) gives \(f\in\cF_{\varepsilon^2,1}(\C^J\oplus\R)=\cF_{\varepsilon^2,1}(\R^{2J+1})\) and
\(x_0\) with \(\norm{x_N}\ge c_0e^{-c\varepsilon N}\norm{x_0}\).

\emph{Case \(\cQ\le\theta/3\).}  Then \(\theta^2>4y\), and the
roots \(a_\varepsilon>b_\varepsilon\) of \(\chi_\mu\) lie in
\((0,1)\) for small \(\varepsilon\), with
\(a_\varepsilon=1-\cQ\varepsilon+o(\varepsilon)\).
On \(f(x)=\mu x^2/2\), normalize the rest initialization to
\(x_{-1}=x_0=1\).  The exact solution is
\[
x_N=Aa_\varepsilon^N+Bb_\varepsilon^N,\qquad
A=\frac{a_\varepsilon(1-b_\varepsilon)}{a_\varepsilon-b_\varepsilon}
\ge1,\qquad B=1-A\le0.
\]
Hence \(x_N/a_\varepsilon^N
=A+B(b_\varepsilon/a_\varepsilon)^N\ge A+B=1\).
For each \(\delta>0\), small \(\varepsilon\) gives
\(a_\varepsilon\ge e^{-(\cQ+\delta)\varepsilon}\), proving the
required bound for every \(N\), with \(c_0=1\).
Embed this scalar example in \(\R^3\) to match the stated dimension.

(b) Let \(0<c<\cFF\) and \(r=1-c\varepsilon\).  Then \(c<\cQ\) and \(\theta-3c>0\).  By \Cref{lem:Pscaled}(a),
\(\sfP_1=2y(\theta-3c)\varepsilon+o(\varepsilon)>0\) for small \(\varepsilon\), \(\sfP_2\to2y(2\tau-1)>0\), because \(\tau>\frac12\), and
\(\sfP_2+2\sfP_3\to2y(2\tau-1)+8(1-\tau y)=2(4-y(1+2\tau))>0\); so the affine function
\(\xi\mapsto\sfP_2+\sfP_3\xi\) is positive at \(\xi=0\) and \(\xi=2\), hence on \([0,2]\), for small \(\varepsilon\), and
\(\sfP_1+\sfP_2\xi+\sfP_3\xi^2>0\) on \([0,2]\).  \Cref{lem:lagscaled} gives \(\rhoF<r\).

(c) By (a), \(\rhoF\ge\rhoz\ge e^{-(\cFF+\delta)\varepsilon}\) for small \(\varepsilon\) (\Cref{lem:envelope}), so
\(\limsup\sqrt\kappa(1-\rhoz)\le\cFF+\delta\) for every \(\delta>0\), and the same for \(\rhoF\).  By (b),
\(\liminf\sqrt\kappa(1-\rhoF)\ge c\) for every \(c<\cFF\), and the same for \(\rhoz\le\rhoF\).
\end{proof}

\begin{proof}[Proof of \Cref{prop:boundary}]
Let \(0<c<\min(\cFF(\theta,y),\cpi)\), put
\(r=1-c\varepsilon\), and let
\(\Psi(\xi)=\sfP_1+\sfP_2\xi+\sfP_3\xi^2\).
By \Cref{lem:Pscaled}(a),
\(\Psi(0)=2y(\theta-3c)\varepsilon+o(\varepsilon)>0\)
and
\(\sfP_3\to-4(2\tau-1)/(1+2\tau)<0\).
Thus \(\Psi\) is concave, and it remains to check \(\Psi(2)>0\).
By \Cref{lem:cubic},
\[
2\Psi(2)=(1+r)\chi_\mu(-r)\chi_L(-r)-\sfP_0.
\]
Here \(\chi_\mu(-r)\to4\) and \(\sfP_0=O(\varepsilon^3)\).
With the roots of \Cref{lem:altmode}(b),
\[
\chi_L(-r)=(\delta_\varepsilon-c\varepsilon)(r+z_2),
\qquad r+z_2\to\frac{4\tau}{1+2\tau}>0.
\]
Since \(\delta_\varepsilon/\varepsilon\to\cpi>c\), this product is
bounded below by a positive multiple of \(\varepsilon\), including
the case \(\cpi=\infty\).  Hence \(\Psi(2)>0\), so concavity gives
\(\Psi>0\) on \([0,2]\).  \Cref{lem:lagscaled} proves
\(\rhoF<1-c\varepsilon\), certified by lag one.

For the limit, \Cref{thm:B}(a) gives
\(\limsup\sqrt\kappa(1-\rhoz)\le\cFF\), and the same for \(\rhoF\).  If \(\cpi<\infty\), \Cref{lem:altmode}(b) gives, from rest on \(f(x)=\frac12x^2\) on \(\R\),
\(\rhoz\ge1-\delta_\varepsilon=1-\cpi\varepsilon+o(\varepsilon)\).  By the first part, \(\liminf\sqrt\kappa(1-\rhoF)\ge c\) for every \(c<\min(\cFF,\cpi)\), and the
same for \(\rhoz\le\rhoF\).
\end{proof}

\begin{proof}[Proof of \Cref{cor:B}]
(i) \(1-\beta=2/(\sqrt\kappa+1)\), so \(\theta_\varepsilon=2/(1+\varepsilon)\to2\), \(y=1\), \(\tau_\varepsilon\to1\);
the tunings are stable since \(1+2\beta<2(1+\beta)\), and \(y(1+2\tau)=3<4\).  As
\(\theta^2=4y\), \(\cQ=1\) and \(\cFF=\frac23\).  Squaring the distance doubles the
constant.
(ii) Substitution in~\eqref{eq:scaling} and~\eqref{eq:margin} gives the stated slow data and margin,
with \(\cQ=\cpi=2/\sqrt3\).  \Cref{prop:boundary} gives the constant
\(\min(4/(3\sqrt3),2/\sqrt3,2/\sqrt3)=4/(3\sqrt3)\).
(iii) By \Cref{rem:twostep}, QHM with \((\alpha,\beta,\nu_{\rm qh})\) has the query representation~\eqref{eq:M} with
\(h=\alpha(1-\beta)\), the same \(\beta\) and \(\tau=\tau_{\rm qh}\); its iterates are the queries of~\eqref{eq:M} from rest, and their rates
satisfy \(\rho^0_{\tilde x}\le\rhoz\le\max(\rho^0_{\tilde x},z_\tau)\), \(z_\tau=\tau_{\rm qh}/(1+\tau_{\rm qh})<1\).  If \(\beta\ne z_\tau\), then on each
\(\frac\lambda2x^2\), \(\lambda\in[\mu,L]\), their asymptotic root factor is \(\rho(\chi_\lambda)\) (\Cref{lem:restquad}).
(a) Let the iterates accelerate on quadratics with constant \(c\).  Then \(\beta\ne z_\tau\) for large \(\kappa\): if \(\beta=z_\tau\), then
\(h_2=\beta h_1\) in the gradient-correction form, so \(e_k=\tilde x_{k+1}-\tilde x_k+h_1\nabla f(\tilde x_k)\) satisfies \(e_k=\beta e_{k-1}\); from rest
\(e_0=0\), hence \(e_k=0\) for all \(k\), and the iterates are those of gradient descent with step \(h_1\), whose
factor on quadratics is at least \((\kappa-1)/(\kappa+1)\).  So \(\rho_q\le1-c\kappa^{-1/2}\), and \Cref{lem:accel} gives \(\tau_{\rm qh}<2/c^2\) and \(1-\beta<\max(4,8/c)\kappa^{-1/2}\), hence
\(1-\nu_{\rm qh}=\tau_{\rm qh}(1-\beta)/\beta=O(\kappa^{-1/2})\).  If \(\limsup\tau_{\rm qh}<\frac12\), \Cref{thm:wall} gives \(\rhoz>1\), hence
\(\rho^0_{\tilde x}>1\), for all large \(\kappa\); so acceleration on \(\cF_{\mu,L}\) forces \(\liminf\tau_{\rm qh}\ge\frac12\).  (b) This is
\Cref{thm:B}(c) with \(hL=\alpha L(1-\beta)\to y\); as \(\tau_{\rm qh}\) is bounded, \(z_\tau\) stays below \(1\) while
\(\rhoz\to1\), so \(\rho^0_{\tilde x}=\rhoz\) for large \(\kappa\).
(iv) The gradient-correction dictionary gives \(h=\gamma\), \(\tau=1\) and the same \(\beta\).
The tuning is stable by \Cref{lem:qstab} and has slow data \((2,1,1)\), so
\Cref{thm:B}(c) gives the constant \(2/3\).

(v) Write \(\alpha\) for the HNAG time step, \(B\) for its gradient coefficient, and \(\gamma\) for the
coefficient of \(y_k-x_{k+1}\) in its first update (\(\gamma=1\) for HNAG\(^{++}\), \(\gamma=2\) for HNAG\(^{+}\)).
With \(\alpha B=1/L\), eliminating the auxiliary state gives the gradient-correction parameters
\[
\begin{gathered}
D=(1+\alpha)(1+\gamma\alpha),\qquad \beta=D^{-1},\\
hL=\frac{\gamma\alpha^2\kappa+\alpha}{D},\qquad
\tau=\frac1{\gamma\alpha^2\kappa+\alpha}.
\end{gathered}
\]
The first update recovers the auxiliary state from two consecutive iterates.
The resulting state maps have norm bounds uniform over the function class, so they preserve the uniform rate.
For HNAG\(^{++}\), \(\alpha=\sqrt2\varepsilon=\eta\), giving the parameters in~(v) and
\[
(\theta,a,b)=\left(2\sqrt2,\frac{\sqrt2}{4},3\sqrt2\right).
\]
These coordinates satisfy the flow-region inequalities of \Cref{rem:regions}, so \(c=\theta/3=2\sqrt2/3\).
The quadratic constant is \(\cQ(2\sqrt2,2)=\sqrt2\); the alternating constant \(\cpi=3\sqrt2\) is larger
(\Cref{lem:quadboundary}).  For HNAG\(^{+}\), \(\alpha=\varepsilon/(1-\varepsilon)\), giving exactly
\(\beta=(1-\varepsilon)^2/(1+\varepsilon)\), \(hL=2-\varepsilon\) and
\(\tau=(1-\varepsilon)^2/((2-\varepsilon)(1+\varepsilon))\), the TMM parameters of \Cref{cor:TMM}.
\end{proof}

\begin{samepage}
\begin{corollary}[Ceilings at fixed extrapolation]\label{cor:frontier}
\leavevmode
\begin{enumerate}[label=\textup{(\roman*)}]
\item Every family in \(\cM\) with \(\tau\to\tau_\infty\) has quadratic constant
\(\limsup\sqrt\kappa(1-\rho_q)\le2/\sqrt{1+2\tau_\infty}\).
\item If a family in \(\cM\) is accelerated on \(\cF_{\mu,L}\), its quadratic
constant \(\limsup\sqrt\kappa(1-\rho_q)\) is at most \(\sqrt2\).
\item Fix \(\tau>\frac12\) and \(y_\tau=4/(1+2\tau)\).  The exact constants of the families of \Cref{thm:B}(c) with this
limit \(\tau\) have supremum \(\sqrt{y_\tau/2}=\sqrt{2/(1+2\tau)}<1\), and none of them attains it.  It is attained on the
stability boundary by the families of \Cref{prop:boundary} with \(\theta=3\sqrt{y_\tau/2}\) and \(\cpi\ge\theta/3\), i.e.\
\(\varsigma\ge2\theta+\frac{4\tau}{1+2\tau}\cdot\frac\theta3\).
\item Fix \(\tau\ge\frac12\) and \(y_\tau=4/(1+2\tau)\).  For every family~\eqref{eq:scaling} with this
limit \(\tau\) (so that \(y\le y_\tau\)), the pair formed by \(\limsup\sqrt\kappa(1-\rho_q)\) and
\(\limsup\sqrt\kappa(1-\rhoz)\) is dominated componentwise by a point of the arc
\(\mathcal C_\tau=\{(\zeta_-(\theta,y_\tau),\theta/3):2\sqrt{y_\tau}\le\theta\le3\sqrt{y_\tau/2}\}\).  For \(\tau>\frac12\), every point of
\(\mathcal C_\tau\) is attained by a family on the stability boundary; for \(\tau=\frac12\), see \Cref{cor:arc}.
\item If a family of stable tunings of \(\cM\) satisfies \(\liminf_{\kappa\to\infty}\sqrt\kappa(1-\rhoF)\ge1\), then, in the
notation of~\eqref{eq:scaling}, \((\theta_\varepsilon,y_\varepsilon,\tau_\varepsilon)\to(3,2,\frac12)\).
\end{enumerate}
\end{corollary}
\end{samepage}

\begin{proof}[Proof of \Cref{cor:nag-asymptotic-optimal}]
Put \(\varepsilon=\kappa^{-1/2}\) and \(\theta=\sqrt6\).
For the proposed tuning, \(y=4/3\), \(\tau_\varepsilon=\beta
=1-\theta\varepsilon\), and
\[
2(1+\beta)-y(1+2\tau)=\frac{2\theta}{3}\varepsilon>0,
\qquad \varsigma_\varepsilon=\frac{8\theta}{3}.
\]
Thus \(c_\pi=\theta/2\) and
\(\cQ(\theta,4/3)=\theta/3=\sqrt{2/3}\).
For \(\kappa\ge6\), \(\beta\in[0,1)\) and \(\tau=\beta\), so
\Cref{thm:tonguefree} also gives linear convergence at each such condition number.
\Cref{prop:boundary} gives the two asymptotic expansions.

For optimality, suppose a family with \(\tau=\beta\) has a larger
positive asymptotic constant.  It is accelerated on quadratics, so
\Cref{lem:accel} gives, along a subsequence, convergent slow data
\((\theta,y,\tau)\) with \(\beta\to1\), hence \(\tau=1\).
Stability forces \(y\le4/3\), while \Cref{thm:B}(a) and
\Cref{cor:damping} bound the constant by
\(\cFF(\theta,y)\le\sqrt{y/2}\le\sqrt{2/3}\), a contradiction.
\end{proof}

\begin{figure}[t]
\centering
\begin{tikzpicture}
\begin{axis}[width=0.65\textwidth,height=0.38\textwidth,
  xmin=0,xmax=3.4,ymin=0,ymax=1.5,xlabel={\(\theta\)},ylabel={\(y\)},
  ylabel style={rotate=-90},
  xtick={0,1,2,3},ytick={0,0.5,1,4/3},yticklabels={\(0\),\(\frac12\),\(1\),\(\frac43\)},
  tick label style={font=\footnotesize},label style={font=\small},
  axis lines=left,clip=true,
  legend style={at={(0.5,-0.20)},anchor=north,font=\scriptsize,draw=none,
    /tikz/every even column/.append style={column sep=1em}},
  legend columns=2,legend cell align=left]
\addplot[draw=none,fill=gray!12,forget plot] coordinates {(0,{4/3}) (3.4,{4/3}) (3.4,1.5) (0,1.5)} --cycle;
\addplot[thin,dashed,domain=0:{4/3},samples=80,variable=t,forget plot] ({2*sqrt(t)},{t});
\addplot[thick,domain=0:{4/3},samples=80,variable=t,forget plot] ({3*sqrt(t/2)},{t});
\addplot[thin,gray!70!black,forget plot] coordinates {({3/4},{4/3}) ({3/4},{1/8})};
\addplot[thin,gray!70!black,domain={3/4}:3.4,forget plot] {x/4-1/16};
\addplot[thin,gray!70!black,forget plot] coordinates {({3/2},{4/3}) ({3/2},{1/2})};
\addplot[thin,gray!70!black,domain={3/2}:{19/6},forget plot] {x/2-1/4};
\addplot[thin,gray!70!black,forget plot] coordinates {(2,{4/3}) (2,{8/9})};
\addplot[thin,gray!70!black,domain=2:{8/3},forget plot] {2*x/3-4/9};
\addplot[thin,gray!70!black,forget plot] coordinates {({9/4},{4/3}) ({9/4},{9/8})};
\addplot[thin,gray!70!black,domain={9/4}:{91/36},forget plot] {3*x/4-9/16};
\addplot[only marks,mark=*,mark size=1.7pt] coordinates {(2,1)};
\addlegendentry{Nesterov, \(h=1/L\)}
\addplot[only marks,mark=square*,mark size=1.6pt] coordinates {({3/sqrt(2)},1)};
\addlegendentry{robust momentum, \(\tau=1\)}
\addplot[only marks,mark=triangle*,mark size=2pt] coordinates {({4/sqrt(3)},{4/3})};
\addlegendentry{Nesterov, quadratic-optimal}
\addplot[only marks,mark=diamond*,mark size=2pt] coordinates {({sqrt(6)},{4/3})};
\addlegendentry{Nesterov, class-optimal}
\node[font=\scriptsize,anchor=south] at (axis cs:0.75,1.35) {\(\frac14\)};
\node[font=\scriptsize,anchor=south] at (axis cs:1.5,1.35) {\(\frac12\)};
\node[font=\scriptsize,anchor=south] at (axis cs:2,1.35) {\(\frac23\)};
\node[font=\scriptsize,anchor=south] at (axis cs:2.25,1.35) {\(\frac34\)};
\node[font=\scriptsize] at (axis cs:3.0,1.42) {unstable};
\end{axis}
\end{tikzpicture}
\caption{Design level sets \(\cFF=c\) at \(\tau=1\), for
\(c=\frac14,\frac12,\frac23,\frac34\) (labels at the top).
Each joins \(\theta=3c\), \(y\ge2c^2\), to \(y=c\theta-c^2\), \(\theta\ge3c\),
on the balance curve \(\theta=3\sqrt{y/2}\) (solid).
The dashed curve is critical damping \(\theta=2\sqrt y\); the shaded region is unstable.
At \(y=4/3\), the alternating mode also constrains the constant (\Cref{prop:boundary}).
The marked tunings are described in \Cref{cor:B,cor:nag-asymptotic-optimal,sec:robust}.}
\label{fig:levels}
\end{figure}
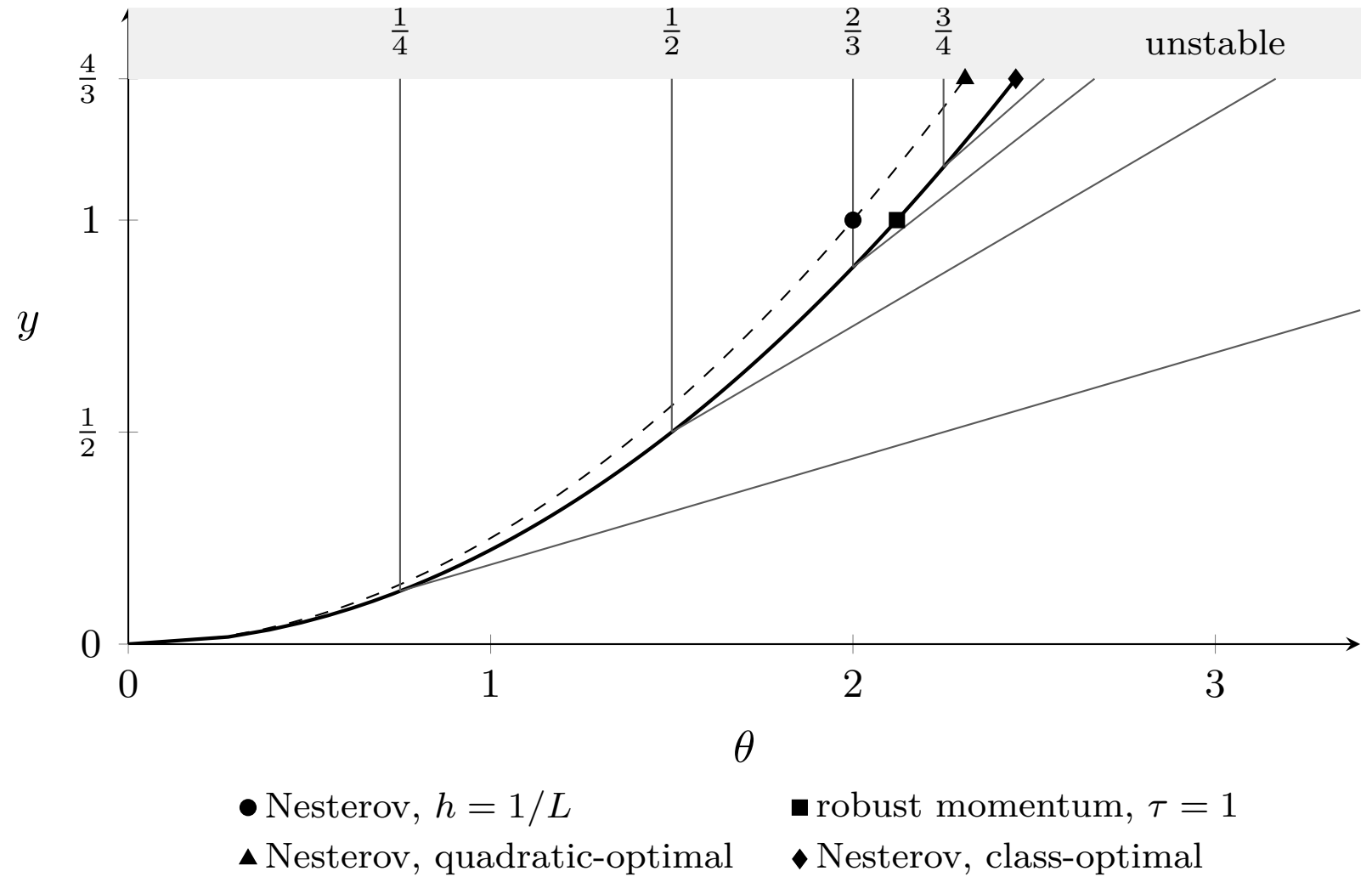

\begin{proof}[Proof of \Cref{cor:frontier}]
(i) By \Cref{thm:price},
\(1-\rho_q\le2/(\sqrt{\kappa_\tau}+1)\), and
\(\sqrt\kappa/\sqrt{\kappa_\tau}=(1+2\tau-2\tau/\kappa)^{-1/2}\to(1+2\tau_\infty)^{-1/2}\).
(ii) By \Cref{thm:wall,lem:accel}, \(\tau_\kappa\) is bounded with
\(\liminf\tau_\kappa\ge\frac12\); apply (i) along subsequences.
(iii) By \Cref{cor:damping} with \(m=y\), \(\cFF(\theta,y)\le\sqrt{y/2}\), with equality only at \(\theta=3\sqrt{y/2}\).  The families
of \Cref{thm:B}(c) with limit \(\tau\) have \(y<y_\tau\), so their constants are at most \(\sqrt{y/2}<\sqrt{y_\tau/2}\), and the families with
constant parameters \(\theta=3\sqrt{y/2}\), \(y\uparrow y_\tau\), show that the supremum is \(\sqrt{y_\tau/2}\); it is below \(1\) since
\(y_\tau<2\) for \(\tau>\frac12\).  On the stability boundary, \Cref{prop:boundary} gives the constant
\(\min(\theta/3,\cQ(\theta,y_\tau),\cpi)\).  At \(\theta=3\sqrt{y_\tau/2}\), \(\theta/3=\cQ(\theta,y_\tau)=\sqrt{y_\tau/2}\) by \Cref{cor:damping}, so the
constant equals \(\sqrt{y_\tau/2}\) if and only if \(\cpi\ge\theta/3\), i.e.\ \((\varsigma-2\theta)(1+2\tau)/(4\tau)\ge\theta/3\).

(iv) \emph{Domination.}  \Cref{lem:cq} and \Cref{thm:B}(a) bound the two
constants componentwise by \((\cQ,\cFF)(\theta,y)\), and stability
gives \(y\le y_\tau\).  By \Cref{lem:cQprop}(b), this pair is
\(\alpha(\cQ,\cFF)(\theta/\alpha,y_\tau)\), where
\(\alpha=\sqrt{y/y_\tau}\le1\).  \Cref{cor:damping} now gives
domination by a point of \(\mathcal C_\tau\).

\emph{Attainment.}  Let \(\tau>\frac12\) and \((\zeta_-(\theta,y_\tau),\theta/3)\in\mathcal C_\tau\); on this arc
\(\theta/3\le\zeta_-(\theta,y_\tau)\).  Take \(\theta_\varepsilon=\theta\), \(\tau_\varepsilon=\tau\) and \(y_\varepsilon=y_\tau-\varepsilon^{1/2}\), so that
\(\varsigma_\varepsilon=(1+2\tau)\varepsilon^{-1/2}\to\infty\) and \(\cpi=\infty\); \Cref{prop:boundary} gives the constant \(\cFF(\theta,y_\tau)=\theta/3\) and
\Cref{lem:quadboundary} the quadratic constant \(\zeta_-(\theta,y_\tau)\).

(v) Let a family of stable tunings (with \(L=1\),
\(\mu=\varepsilon^2\), \(\theta_\varepsilon=(1-\beta)/\varepsilon\), \(y_\varepsilon=h\)) satisfy
\(\liminf\sqrt\kappa(1-\rhoF)\ge1\), and suppose that along a sequence \(\varepsilon\to0\) the triple
\((\theta_\varepsilon,y_\varepsilon,\tau_\varepsilon)\) stays away from \((3,2,\frac12)\).  Fix \(c\in(0,1)\).  For small \(\varepsilon\),
\(\rho_q\le\rhoF\le1-c\varepsilon\), so by \Cref{lem:accel}(b) a subsequence converges to some \((\theta,y,\tau)\) with \(\theta\ge2c>0\) and is a
family~\eqref{eq:scaling}.  By \Cref{thm:wall} (applied along the subsequence, which
is accelerated on quadratics), \(\tau\ge\frac12\); by stability, \(y(1+2\tau)\le4\), so \(y\le2\).
\Cref{thm:B}(a) gives \(\limsup\sqrt\kappa(1-\rhoz)\le\cFF(\theta,y)\le\sqrt{y/2}\le1\).  Since
\(\rhoz\le\rhoF\), the assumption forces \(\sqrt{y/2}=1\), i.e.\ \(y=2\) and then \(\tau=\frac12\), and
\(\cFF(\theta,2)=1\), i.e.\ \(\theta=3\).  This contradicts the choice of the sequence.
\end{proof}

\begin{proof}[Proof of \Cref{prop:onedisc}]
Suppose \(\varepsilon_k\to0\), \(c_k\in(\theta/3,c_\sharp]\) and \(\gamma_k\in[0,\gamma_{\max}]\) with
\((\delta_{\varepsilon_k\gamma_k},C_k)\) feasible at \(-c_k\varepsilon_k\).  Pass to a subsequence with
\(c_k\to c_\infty\), \(\gamma_k\to\gamma_\infty\).  By \Cref{lem:Cforced},
\(C_k=C_{*,k}+O(\varepsilon_k^3)\); in the normalization of~\eqref{eq:gtst} the lag
constraints become
\(\gt\cA_{\gamma_k}(q)\pm\st\cB_{\gamma_k}(q)+O(\varepsilon_k)\sh(c_kq)\ge0\) for \(q\in\varepsilon_k\N\), and
the minimizer constraint becomes \(|\gamma_k\st|\le c_k\gt+O(\varepsilon_k)\).  For fixed
\(q>0\) take \(q_k=\varepsilon_k\lfloor q/\varepsilon_k\rfloor\to q\); by \Cref{lem:reduce} the limits
show that \(\gamma_\infty\) satisfies~\eqref{eq:P12} at \(c_\infty\), and \(|(\theta-2c_\infty)\gamma_\infty^2|\le c_\infty(\gamma_\infty^2-\Delta_{c_\infty})\), which
is~\eqref{eq:P3} (with \(m=y\)).  If \(c_\infty>\theta/3\) this contradicts \Cref{prop:one}; if \(c_\infty=\theta/3\),
\eqref{eq:P3} reads \(0\ge\frac\theta3\Delta_{\theta/3}>0\), a contradiction.
\end{proof}
\begin{proof}[Proof of \Cref{thm:layer}]
Throughout, \(2\tau_\varepsilon-1=-2A_\varepsilon\varepsilon^{1/2}\), \(\tau_\varepsilon\to\frac12\), and
\(4(1-\tau_\varepsilon y_\varepsilon)\to2(2-y)>0\).  We use repeatedly that
\(\varepsilon^{1/2}\xi\le\frac12(\varepsilon+\xi^2)\) and \(\varepsilon\xi\le\frac12\varepsilon^{1/2}(\varepsilon+\xi^2)\).

(a) It suffices to take \(0<c<\clayer(A)\); put \(r=1-c\varepsilon\).  Then \(c<\cQ\) and \(3c<\Lambda(A)\).  By
\Cref{lem:Pscaled}(a), uniformly for \(\xi\in[0,2]\),
\[
\sfP_1+\sfP_2\xi+\sfP_3\xi^2\ \ge\ 2y(\theta-3c)\varepsilon-4yA_+\varepsilon^{1/2}\xi+2(2-y)\xi^2-o(\varepsilon+\xi^2),
\]
because \(\sfP_2\ge-4y_\varepsilon(A_\varepsilon)_+\varepsilon^{1/2}-O(\varepsilon)\) and \(y_\varepsilon(A_\varepsilon)_+\to yA_+\) (also for
\(A=-\infty\)).  The binary form \(2y(\theta-3c)a^2-4yA_+ab+2(2-y)b^2\) is positive definite,
since \(3c<\Lambda(A)\) is \(yA_+^2<(\theta-3c)(2-y)\); evaluated at \((\varepsilon^{1/2},\xi)\) it is at
least \(m_1(\varepsilon+\xi^2)\) with \(m_1>0\).  Hence \(\sfP_1+\sfP_2\xi+\sfP_3\xi^2\ge(m_1-o(1))(\varepsilon+\xi^2)>0\) on
\([0,2]\) for small \(\varepsilon\), and \Cref{lem:lagscaled} gives \(\rhoF<r\).

(b) Let \(A\in(0,\infty)\), and let \(c>\Lambda(A)_+/3\); we may assume \(c<\theta/2\) (if \(c\) works, so does every
larger \(c\)).  Put \(r=1-c\varepsilon\), \(\varphi_*=Ay/(2-y)\) and
\(m_3(\varphi)=2y(\theta-3c)-4Ay\varphi+2(2-y)\varphi^2\).  Then
\(m_3(\varphi_*)=\frac{2y}{2-y}\bigl[(\theta-3c)(2-y)-yA^2\bigr]<0\), because \(3c>\Lambda(A)\).
Put \(m_*=-\frac12m_3(\varphi_*)>0\), and choose
\(\mathsf s\in(0,\mathsf s_1]\) (\Cref{lem:spreadhalf}) so small that \(m_3\le-m_*\) on
\(I_\varphi=[\varphi_*(1-2\mathsf s)^2,\varphi_*(1+2\mathsf s)^2]\).

\emph{The band.}  Let \(\omega_*=\arccos(1-\varepsilon^{1/2}\varphi_*)\); as \(\varepsilon\to0\), \(\omega_*\to0\) and
\(\omega_*^2/(2\varepsilon^{1/2}\varphi_*)\to1\).  For \([\omega_a,\omega_b]=[(1-\frac32\mathsf s)\omega_*,(1+\frac32\mathsf s)\omega_*]\),
\(2\pi(1-\mathsf s)/\omega_a-2\pi(1+\mathsf s)/\omega_b=\frac{2\pi}{\omega_*}\bigl(\frac{1-\mathsf s}{1-3\mathsf s/2}-\frac{1+\mathsf s}{1+3\mathsf s/2}\bigr)\to\infty\), so \Cref{lem:band}(a) gives \(K\ge3\)
with \(B=[(1-\mathsf s)\omega_0,(1+\mathsf s)\omega_0]\subset[\omega_a,\omega_b]\), \(\omega_0=2\pi/K\).  Since \(\xi=1-\cos\omega=\frac{\omega^2}2(1+O(\omega^2))\), every \(\omega\in B\) has
\(\xi=\varepsilon^{1/2}\varphi\) with \(\varphi\in I_\varphi\), for small \(\varepsilon\).

\emph{Symbols on the band.}  By \Cref{lem:Pscaled}(a), \(\sfP_0=O(\varepsilon^3)\),
\(\sfP_1=2y(\theta-3c)\varepsilon+o(\varepsilon)\), \(\sfP_2=-4y_\varepsilon A_\varepsilon\varepsilon^{1/2}+O(\varepsilon)\) and
\(\sfP_3=2(2-y)+o(1)\), so on \(B\)
\[
H_r=\varepsilon^{3/2}\varphi\bigl[m_3(\varphi)+o(1)\bigr]+O(\varepsilon^3)\le-\tfrac14\varepsilon^{3/2}\varphi_*m_* .
\]
By \Cref{lem:Pscaled}(b), \(\mathsf M_r=O(\varepsilon)\) on \(B\) (its coefficients are \(O(\varepsilon)\), \(-4y_\varepsilon A_\varepsilon\varepsilon^{1/2}+O(\varepsilon)\) and \(O(1)\)),
and \(|\chi_\mu\overline{\chi_L}|\le C_2\varepsilon^{1/2}\) on \(B\), as \(|\chi_\mu|=O(\xi)\) and \(|\chi_L|=O(1)\) there.  Since \(\xi\asymp\varepsilon^{1/2}\) on \(B\), the
hypotheses of \Cref{lem:spreadwit} hold, with \(\lambda_0=\varphi_*m_*/(4C_2)\): \(\varepsilon\) and \(\max_B|\mathsf M_r|\) are \(O(\varepsilon)=o(\varepsilon^{3/4})\).  It gives
\(\rhoz\ge1-c\varepsilon\).

(c) By \Cref{prop:Ltau}, \(L_\tau=\mathsf A_0+\mathsf A_1\xi+\mathsf A_2\xi^2\) with
\(\mathsf A_0=2\theta_\varepsilon y_\varepsilon\varepsilon+O(\varepsilon^2)>0\),
\(\mathsf A_1=2y_\varepsilon(2\tau_\varepsilon-1)+O(\varepsilon)=-4y_\varepsilon A_\varepsilon\varepsilon^{1/2}(1+o(1))<0\) and
\(\mathsf A_2=4(\beta-\tau_\varepsilon y_\varepsilon)\to2(2-y)\); the \(o(1)\) is uniform because
\(\liminf A_\varepsilon>0\).  Hence
\(\mathsf A_1^2/(4\mathsf A_0\mathsf A_2)=\frac{y_\varepsilon A_\varepsilon^2}{\theta_\varepsilon(2-y_\varepsilon)}(1+o(1))\ge(A'/A_*)^2(1+o(1))\) for some
\(A'>A_*\), so \(L_\tau\) has two roots \(0<\xi_1<\xi_2\) with
\(\xi_2/\xi_1=(1+\sqrt d)/(1-\sqrt d)\), \(d=1-4\mathsf A_0\mathsf A_2/\mathsf A_1^2\ge d_0>0\), and \(L_\tau<0\) on
\((\xi_1,\xi_2)\).  Moreover \(\xi_2<-\mathsf A_1/\mathsf A_2\to0\).  With \(\omega_i=\arccos(1-\xi_i)\), therefore
\(\omega_2\to0\) and \(\omega_2/\omega_1\ge1+\eta'\) for small \(\varepsilon\), with \(\eta'>0\); so
\(2\pi/\omega_1-2\pi/\omega_2\ge2\pi\eta'/\omega_2\to\infty\), and \((2\pi/\omega_2,2\pi/\omega_1)\) contains an integer \(K\ge3\), at which
\(L_\tau(1-\cos\frac{2\pi}K)<0\), and \Cref{lem:band}(b) gives the orbits and \(\rhoz>1\).  If
\(A_\varepsilon\to A<\infty\), then \(\xi_1\xi_2=\mathsf A_0/\mathsf A_2\asymp\varepsilon\) and \(\xi_2/\xi_1\) is bounded, so
\(\xi_i\asymp\varepsilon^{1/2}\), \(\omega_i\asymp\varepsilon^{1/4}\) and \(K\asymp\varepsilon^{-1/4}=\kappa^{1/8}\).

\emph{Consequences.}  Let \(A_\varepsilon\to A\in[-\infty,A_*)\).  By (a),
\(\liminf\sqrt\kappa(1-\rhoF)\ge\clayer(A)\).  \Cref{thm:B}(a), which holds for every \(\tau\), gives
\(\limsup\sqrt\kappa(1-\rhoz)\le\cFF(\theta,y)\), and for \(A>0\), (b) gives \(\limsup\sqrt\kappa(1-\rhoz)\le\Lambda(A)/3\).
Since \(\Lambda(A)\le\theta\), with equality for \(A\le0\), \(\min(\cFF,\Lambda(A)/3)=\clayer(A)\) in both cases, and
\(\rhoz\le\rhoF\) closes the chain.  If \(A_\varepsilon\to A\in[A_*,\infty)\), (b) applies to every \(c>0=\Lambda(A)_+/3\), and if
\(A=\infty\), (c) applies; in both cases the family is not accelerated.  Finally \(\clayer\) is
nonincreasing in \(A\) and positive on \([-\infty,A_*)\); applying the above along subsequences gives the last
claim.
\end{proof}

\begin{remark}[One-atom spectra in the layer]\label{rem:onelayer}
For a family as in \Cref{thm:layer}, let \(A_\varepsilon\to A\in(0,\infty)\), \(c>\Lambda(A)_+/2\) and \(r=1-c\varepsilon\).  By \Cref{lem:Pscaled}(b), with \(\xi=\varepsilon^{1/2}\varphi\),
\[
\varepsilon^{-1}\mathsf M_r=2y(\theta-2c)-4Ay\varphi+2(2-y)\varphi^2+o(1)
\]
uniformly for \(\varphi\) in compact sets.  The quadratic in \(\varphi\) is minimal at \(\varphi_*=Ay/(2-y)\), with the value
\(\frac{2y}{2-y}\bigl[(\theta-2c)(2-y)-yA^2\bigr]\), which is negative because \(2c>\Lambda(A)\).  The points \(\xi_K=1-\cos\frac{2\pi}K\) with
\(K\asymp\varepsilon^{-1/4}\) have consecutive ratios tending to \(1\), so for small \(\varepsilon\) some \(\xi_K=\varepsilon^{1/2}\varphi\) has \(\varphi\) close to
\(\varphi_*\); there \(H^{\rm mix}_r=\xi_K\mathsf M_r<0\), while \(\Ima(\chi_\mu\overline{\chi_L})>0\) by \Cref{lem:Pscaled}(b).  \Cref{lem:mixed} shows that
\((\delta_{2\pi/K},0)\) is feasible at \(t=\log r\), and \Cref{lem:startup} gives \(\rhoz\ge1-c\varepsilon\).  So the one-atom spectra of
\Cref{lem:mixed} give \(\limsup\sqrt\kappa(1-\rhoz)\le\Lambda(A)_+/2\), and with \Cref{thm:B}(a) the bound \(\min(\cFF,\Lambda(A)_+/2)\)
drawn dashed in \Cref{fig:frontier}.
\end{remark}

\section{Proofs for Section~\ref{sec:TMM}}
\label{app:corner}

\subsection{The exact robust-momentum certificate}
\label{app:robust-certificate}

Robust momentum~\cite{CHVL18} is~\eqref{eq:M} with
\begin{equation}\label{eq:robust}
hL=\kappa(1-\varrho)^2(1+\varrho),\qquad \beta=\frac{\kappa\varrho^3}{\kappa-1},\qquad
\tau=\frac{\varrho^3}{(\kappa-1)(1-\varrho)^2(1+\varrho)},\qquad \varrho\in\bigl[1-\kappa^{-1/2},\,1-\kappa^{-1}\bigr],
\end{equation}
where \(\varrho\) is its design rate; for \(\varrho=1-\kappa^{-1/2}\) it is TMM, and for \(\varrho=1-\kappa^{-1}\) its
queries follow gradient descent with step \(1/L\)~\cite{CHVL18}.  For this family the lag-one certificate can be
computed in closed form at every condition number.

\begin{theorem}[Robust momentum]\label{thm:robust}
Let \(\kappa>1\) and \(\varrho\in[1-\kappa^{-1/2},1-\kappa^{-1}]\), and consider the tuning~\eqref{eq:robust}.  Then
\[
\chi_\mu(z)=(z-\varrho)(z-\varrho^2),\qquad \chi_L(z)=z(z-z_L),\qquad z_L=1-\kappa(1-\varrho)^2(1+\varrho),
\]
and on the circle \(|z|=\varrho\), with \(\xi=1-\cos\omega\),
\[
\Rea\bigl[(1-z)\chi_\mu(z)\overline{\chi_L(z)}\bigr]=\varrho^3(1+\varrho)\bigl(1-\kappa(1-\varrho)^2\bigr)\,\xi\,\bigl((1-\varrho)^2+2\varrho\xi\bigr)\ \ge0 .
\]
Moreover, on every circle \(|z|=r\) with \(\varrho<r<1\),
\[
\Rea\bigl[(1-z)\chi_\mu(z)\overline{\chi_L(z)}\bigr]>0 .
\]
The tuning is stable and has exact rates
\[
\rhoF=\rhoz=\rho_q=\varrho .
\]
At TMM, \(\varrho=1-\kappa^{-1/2}\), and the lag-one symbol vanishes on \(|z|=\varrho\).
\end{theorem}

\begin{proof}[Proof of \Cref{thm:robust}]
\emph{Characteristic polynomials.}  From~\eqref{eq:robust}, \(hL\tau=\kappa\varrho^3/(\kappa-1)=\beta\) and
\(h\mu\tau=\beta/\kappa\).  Hence \(Q_L=\beta-hL\tau=0\) and \(S_L=1+\beta-hL(1+\tau)=1-hL=z_L\), so
\(\chi_L(z)=z(z-z_L)\).  Further \(Q_\mu=\beta(1-\frac1\kappa)=\varrho^3\) and
\(S_\mu=1+\beta-h\mu-\frac\beta\kappa=1+\varrho^3-(1-\varrho)^2(1+\varrho)=\varrho+\varrho^2\), so \(\chi_\mu(z)=(z-\varrho)(z-\varrho^2)\).

\emph{Admissibility and stability.}  Since \(0<\varrho\le1-\kappa^{-1}\), we have \(\varrho^3<\varrho\) and
\(\beta\in[0,1)\); also \(\tau\ge0\), \(h>0\).  Since \(\varrho\ge1-\kappa^{-1/2}\), \(\kappa(1-\varrho)^2\le1\) and
\(hL\le1+\varrho<2\).  As \(hL\tau=\beta\), \(y(1+2\tau)-2(1+\beta)=hL-2<0\), and the tuning is stable by
\Cref{lem:qstab}.  Moreover \(-\varrho\le z_L\le\varrho\): indeed
\[
\varrho+z_L=(1+\varrho)\bigl(1-\kappa(1-\varrho)^2\bigr)\ge0,\qquad
\varrho-z_L=(1-\varrho)\bigl(\kappa(1-\varrho^2)-1\bigr)\ge0,
\]
since \(\kappa(1-\varrho)^2\le1\) and \(1-\varrho^2\ge1-\varrho\ge\kappa^{-1}\).

\emph{The symbol.}  Let \(z=\varrho e^{i\omega}\).  Since \(1-\varrho e^{i\omega}=e^{i\omega}\,\overline{(e^{i\omega}-\varrho)}\),
\[
(1-z)(z-\varrho^2)=\varrho(1-\varrho e^{i\omega})(e^{i\omega}-\varrho)=\varrho\,e^{i\omega}\,|e^{i\omega}-\varrho|^2,
\qquad (z-\varrho)\bar z=\varrho^2(1-e^{-i\omega}).
\]
Hence
\[
(1-z)\chi_\mu(z)\overline{\chi_L(z)}=\varrho^3|e^{i\omega}-\varrho|^2\,(e^{i\omega}-1)(\varrho e^{-i\omega}-z_L),
\]
and \(\Rea[(e^{i\omega}-1)(\varrho e^{-i\omega}-z_L)]=(\varrho+z_L)(1-\cos\omega)=(\varrho+z_L)\xi\).  With
\(|e^{i\omega}-\varrho|^2=(1-\varrho)^2+2\varrho\xi\) and \(\varrho+z_L=(1+\varrho)(1-\kappa(1-\varrho)^2)\) this is the stated formula.

\emph{Lower bounds.}  \(\varrho\) is a root of \(\chi_\mu\), so \(\rhoF\ge\rho_q\ge\varrho\) (\Cref{lem:envelope}), and \(\rhoz\ge\varrho\) by
\Cref{lem:restquad} with \(\lambda=\mu\), since the roots of \(\chi_\mu\) are \(\varrho\) and \(\varrho^2\).

\emph{Positive circles and the upper bound.}  Fix \(\varrho<r<1\) and \(\xi\in[0,2]\).
Let \(H_r(s;\xi)\) be the circle symbol with \(\chi_\mu=(z-\varrho)(z-\varrho^2)\) and \(\chi_L=z(z-s)\).
It is affine in \(s\), so
\[
2\varrho H_r(z_L;\xi)=(\varrho-z_L)H_r(-\varrho;\xi)+(\varrho+z_L)H_r(\varrho;\xi).
\]
At the two endpoints, polynomial expansion gives
\begin{align*}
H_r(-\varrho;\xi)&=r(r^2-\varrho^2)\bigl[(r-\varrho^2)(1-r)+\xi\bigl(r^2+\varrho^2+2r\varrho(2-\xi)\bigr)\bigr],\\
H_r(\varrho;\xi)&=r\bigl[(r-\varrho)^2+2r\varrho\xi\bigr]
\bigl[(r-\varrho^2)(1-r)+(r^2+\varrho^2)\xi\bigr].
\end{align*}
Both are strictly positive: \((r-\varrho^2)(1-r)>0\), and all remaining bracketed terms are nonnegative,
with \((r-\varrho)^2>0\).  Since \(-\varrho\le z_L\le\varrho\), the interpolation weights are nonnegative
and sum to \(2\varrho>0\).  Therefore \(H_r(z_L;\xi)>0\) for every \(\varrho<r<1\).
The roots \(\varrho,\varrho^2\) of \(\chi_\mu\) lie in \(|z|<r\), so \Cref{lem:separate} gives
\(\rho_q<r\) and \(\rhoF<r\).  Letting \(r\downarrow\varrho\) proves \(\rhoF\le\varrho\).
With the lower bounds and \Cref{lem:envelope}, \(\varrho\le\rho_q\le\rhoF\le\varrho\) and \(\varrho\le\rhoz\le\rhoF\).
At \(\varrho=1-\kappa^{-1/2}\), \(\kappa(1-\varrho)^2=1\), so the symbol at radius \(\varrho\) vanishes identically.
\end{proof}

\begin{proof}[Proof of \Cref{cor:TMM}]
TMM~\cite{VFL18} is~\eqref{eq:M} with \(\beta=\varrho^2/(2-\varrho)\), \(\tau=\varrho^2/((1+\varrho)(2-\varrho))\) and
\(hL=1+\varrho\), where \(\varrho=1-\kappa^{-1/2}\).  (It reports the point \((1+\delta_{\rm T})x_k-\delta_{\rm T}x_{k-1}\),
\(\delta_{\rm T}=\varrho^2/(1-\varrho^2)\); as in \Cref{rem:twostep}, this point has the same rate as \(x_k\), because
\(\delta_{\rm T}/(1+\delta_{\rm T})=\varrho^2<\varrho\).)  This is~\eqref{eq:robust} at \(\varrho=1-\kappa^{-1/2}\): then
\(\kappa(1-\varrho)^2=1\) and \(\kappa-1=\varrho(2-\varrho)/(1-\varrho)^2\), so \(\kappa(1-\varrho)^2(1+\varrho)=1+\varrho\),
\(\kappa\varrho^3/(\kappa-1)=\varrho^2/(2-\varrho)\) and \(\varrho^3/((\kappa-1)(1-\varrho)^2(1+\varrho))=\varrho^2/((2-\varrho)(1+\varrho))\).
\Cref{thm:robust} gives \(\rhoF=\rhoz=\rho_q=\varrho\).  With \(\varepsilon=1-\varrho\),
\(1-\beta=\varepsilon(3-\varepsilon)/(1+\varepsilon)\), \(2-hL=\varepsilon\), and
\[
\tfrac12-\tau=\frac{(1-\varrho)(2+3\varrho)}{2(2+\varrho-\varrho^2)}=\tfrac54\varepsilon\,(1+O(\varepsilon));
\]
in the notation of~\eqref{eq:corner}, \((\theta_\varepsilon,a_\varepsilon,b_\varepsilon)\to(3,\frac54,1)\).  In particular
\(\tau<\frac12\) for every \(\varrho\in(0,1)\).  Finally \(\zeta_-(3,2)=1=\theta/3\), and optimality among deterministic first-order methods with value and gradient access
follows from~\cite{DT22}.
\end{proof}
\begin{lemma}[Expansions at the corner]\label{lem:cornerexp}
Under~\eqref{eq:corner}, let \(r=1-c\varepsilon\), \(c\in\R\), and write \(\doteq\) for equality up to
\(o(\varepsilon)\), uniformly for \(\xi\in[0,2]\) and for \(c\) in compact sets.  (With \(\theta_\varepsilon,a_\varepsilon,b_\varepsilon\) in
place of their limits the remainders are \(O(\varepsilon^2)\).)
Then the coefficients of \Cref{lem:cubic} and \Cref{lem:mixed} satisfy
\[
\sfP_1\doteq4(\theta-3c)\varepsilon,\quad \sfP_2\doteq2(\theta+2c-4a)\varepsilon,\quad \sfP_3\doteq2(4a+b-2\theta)\varepsilon,
\]
\[
\sfP_1+\sfP_2\xi+\sfP_3\xi^2\doteq\varepsilon\bigl[\mathsf N(\xi)-4c(3-\xi)\bigr],\qquad
\mathsf M_r(\xi)\doteq\varepsilon\bigl[\mathsf N(\xi)-2c(4-\xi)\bigr].
\]
Moreover \(\sfP_0=c\varepsilon\chi_\mu(r)\chi_L(r)\) with \(\chi_\mu(r)=\varepsilon^2(c^2-\theta_\varepsilon c+2+O(\varepsilon))\) and \(\chi_L(r)=2+O(\varepsilon)\),
\(\chi_L(-1)=\cpi\varepsilon+o(\varepsilon)\), \(S_L\to-1\), and \(\chi_L\) has a real root \(z_L=-1+\cpi\varepsilon+o(\varepsilon)\).
\end{lemma}

\begin{proof}
Substitute \(y_\varepsilon=2-b_\varepsilon\varepsilon\) and
\(\tau_\varepsilon=1/2-a_\varepsilon\varepsilon\) in \Cref{lem:Pscaled}.
The corner coordinates are bounded, so collecting the polynomial coefficients gives
\(O(\varepsilon^2)\) remainders with the member parameters;
replacing them by their limits costs \(o(\varepsilon)\), uniformly on the stated compact sets.
The characteristic formulas at \(r\) also follow from \Cref{lem:Pscaled}(a);
retaining \(\theta_\varepsilon\) gives the asserted \(O(\varepsilon)\) remainder inside
\(\chi_\mu(r)/\varepsilon^2\).
Finally,
\[
\chi_L(-1)=(4a_\varepsilon+2b_\varepsilon-2\theta_\varepsilon)\varepsilon
-2a_\varepsilon b_\varepsilon\varepsilon^2
=\cpi\varepsilon+o(\varepsilon),\qquad
\chi_L'(-1)=-1+O(\varepsilon).
\]
The simple real root near \(-1\) therefore satisfies
\(z_L=-1+\chi_L(-1)(1+O(\varepsilon))=-1+\cpi\varepsilon+o(\varepsilon)\).
The expansion holds for every sign of \(\cpi\).
\end{proof}
\begin{proof}[Proof of \Cref{thm:corner}]
The tunings are stable for small \(\varepsilon\): by \Cref{lem:qstab}, stability is
\(\chi_L(-1)>0\), and \(\chi_L(-1)=\cpi\varepsilon+o(\varepsilon)\) with \(\cpi>0\).

(a) Let \(0<c<\min(\cQ(\theta,2),\clag)\) and \(r=1-c\varepsilon\).  Since \(3-\xi\ge1\) on \([0,2]\),
\[
\mathsf N(\xi)-4c(3-\xi)=4(3-\xi)\bigl(\mathsf q(\xi)-c\bigr)\ge4(\clag-c)=:m_1>0 ,
\]
so by \Cref{lem:cornerexp} \(\sfP_1+\sfP_2\xi+\sfP_3\xi^2\ge\varepsilon m_1-o(\varepsilon)>0\) on \([0,2]\) for small \(\varepsilon\).  The family has
\((\theta_\varepsilon,y_\varepsilon,\tau_\varepsilon)\to(\theta,2,\frac12)\), and \Cref{lem:lagscaled} gives \(\rhoF<r\).

(b) For a requested constant \(\widehat c\), it suffices to construct a witness
at some \(0<c<\widehat c\), since
\(\rhoz\ge1-c\varepsilon>1-\widehat c\varepsilon\).
We treat the terms in the minimum separately and choose \(c\) strictly between
the selected term (or zero, if larger) and \(\widehat c\).
The family~\eqref{eq:corner} is a family~\eqref{eq:scaling} with \(y=2\)
and \(\tau=\frac12\), so \Cref{thm:B}(a) gives
\(\rhoz\ge e^{-(\cFF(\theta,2)+\delta)\varepsilon}\) for every \(\delta>0\) and small \(\varepsilon\), with the witnesses described there; note that
\(\cFF(\theta,2)=\min(\theta/3,\cQ(\theta,2))\).

For \(\cpi\), \Cref{lem:altmode}(a) applies, since \(\chi_L(-1)=\cpi\varepsilon+o(\varepsilon)>0\) and \(S_L\to-1\) (\Cref{lem:cornerexp}): from rest on \(f(x)=\frac12x^2\) on \(\R\) (here \(L=1\)),
\(\rhoz\ge|z_L|=1-\cpi\varepsilon+o(\varepsilon)\), and the dominant mode \(z_L^k\) alternates in sign.

For \(\mathsf q(\xi_K)\), let \(K\ge3\) with \(\mathsf q(\xi_K)<\theta/3\) (otherwise every \(c>\mathsf q(\xi_K)\) exceeds \(\cFF(\theta,2)\), and the first term
applies), \(\omega_0=2\pi/K\), and \(\max(0,\mathsf q(\xi_K))<c<\theta/3\).  Then \(\mathsf N(\xi_K)-4c(3-\xi_K)<0\); fix
\(\mathsf s\in(0,\mathsf s_1]\) such that \(\mathsf N(\xi)-4c(3-\xi)\le-m_*<0\) for \(\omega\in B=[(1-\mathsf s)\omega_0,(1+\mathsf s)\omega_0]\), and put
\(r=1-c\varepsilon\).  On \(B\), \(\xi\ge\xi_{\min}>0\), and by \Cref{lem:cornerexp}
\(H_r=\sfP_0+\xi\varepsilon[\mathsf N(\xi)-4c(3-\xi)+o(1)]\le-\frac12\xi_{\min}m_*\varepsilon\) and \(\mathsf M_r=O(\varepsilon)\) for small \(\varepsilon\), while
\(|\chi_\mu\overline{\chi_L}|\le C_3\) on \(B\).  Since \(B\) is fixed across \(\varepsilon\), \Cref{lem:spreadwit} applies with \(\lambda_0=\xi_{\min}m_*/(2C_3)\) and gives
\(\rhoz\ge1-c\varepsilon\).

(c) The upper bound on \(\limsup\) is (b).  If \(\clag>0\), (a) gives
\(\liminf\sqrt\kappa(1-\rhoF)\ge\min(\cQ(\theta,2),\clag)\), and \(\rhoz\le\rhoF\) gives the middle inequality.  If
\(\mathsf q\) attains its minimum at a point of \(\{0,2\}\cup\{\xi_K\}\), then \(\clag=\cdisc\); if \(\cQ(\theta,2)\le\clag\), both bounds equal
\(\cQ(\theta,2)\).  If \(\mathsf N(\xi)\ge4\ccorner(3-\xi)\) on \([0,2]\), then \(\clag\ge\ccorner\), while
\(\cdisc\le\min(\mathsf q(0),\mathsf q(2))=\min(\theta/3,\cpi)\); so both bounds equal \(\ccorner\).  If \(4a+b\le2\theta\), then \(\mathsf N\) is
concave, so \(\mathsf N(\xi)-4\ccorner(3-\xi)\) is concave on \([0,2]\); at \(\xi=0\) it equals
\(4(\theta-3\ccorner)\ge0\), and at \(\xi=2\) it equals \(8(2a+b-\theta)-4\ccorner=4(\cpi-\ccorner)\ge0\).  Hence it
is nonnegative on \([0,2]\).

\emph{Explicit regions (\Cref{rem:regions}).}  For \(\theta\le3\),
\(\cQ(\theta,2)\ge\theta/3\) by \Cref{cor:damping}, so
\(\ccorner=\min(\theta/3,\cpi)\).  The identities
\[
\begin{aligned}
\mathsf N(\xi)-\tfrac{4\theta}3(3-\xi)&=2\xi\Bigl[\bigl(\tfrac{5\theta}3-4a\bigr)+(4a+b-2\theta)\xi\Bigr],\\
\mathsf N(\xi)-4\cpi(3-\xi)&=2(2-\xi)\Bigl[(7\theta-12a-6b)-(4a+b-2\theta)\xi\Bigr],
\end{aligned}
\]
reduce the criterion of (c) to the signs of affine functions at \(\xi=0\) and \(\xi=2\).  If
\(a\le\frac{5\theta}{12}\) and \(2a+b\ge\frac{7\theta}6\), the first bracket is nonnegative at both ends (its
value at \(\xi=2\) is \(2(2a+b)-\frac{7\theta}3\)), and \(\cpi=2(2a+b-\theta)\ge\theta/3\), so \(\ccorner=\theta/3\) and (c)
applies.  If \(2a+b<\frac{7\theta}6\), then \(\cpi<\theta/3\), so \(\ccorner=\cpi\); the second bracket is positive
at \(\xi=0\) and equals \(11\theta-20a-8b\ge0\) at \(\xi=2\) when \(5a+2b\le\frac{11\theta}4\), and (c) applies.

(d) Let \(\theta=3\) and \(a\le\frac54\), \(2a+b\ge\frac72\).  By \Cref{rem:regions}, the constant is \(\theta/3=1\).  Conversely, if \(\theta\ne3\), then
\(\min(\theta/3,\cQ(\theta,2))=\cFF(\theta,2)<1\) by \Cref{cor:damping} (with \(m=2\), the maximum \(\sqrt{m/2}=1\) of \(\cFF(\cdot,2)\) is
attained only at \(\theta=3\)), and (b) gives \(\limsup\sqrt\kappa(1-\rhoz)<1\).  Let \(\theta=3\).  If \(2a+b<\frac72\),
then \(\cpi<1\).  If \(a>\frac54\), then \(\mathsf q(\xi)-1=2\xi[(5-4a)+(4a+b-6)\xi]/(4(3-\xi))<0\) for small \(\xi>0\), hence
\(\mathsf q(\xi_K)<1\) for all large \(K\), as \(\xi_K\to0\).  In both cases \(\cdisc<1\), and (b) again gives
\(\limsup\sqrt\kappa(1-\rhoF)\le\limsup\sqrt\kappa(1-\rhoz)<1\).
\end{proof}

\begin{proof}[Proof of \Cref{cor:arc}]
(i) These tunings form a family~\eqref{eq:corner} with \(\theta=2\sqrt2\), \(a=0\), \(b=4\), so
\(2a+b=4>\theta\), \(4a+b=4\le2\theta\) and \(\cpi=2(4-2\sqrt2)>\sqrt2\).  Since \(\theta^2=8=4y\),
\(\cQ(\theta,2)=\sqrt2\), and \Cref{thm:corner}(c) gives the constant
\(\ccorner=\min(\frac{2\sqrt2}3,\sqrt2,\cpi)=\frac{2\sqrt2}3\) on \(\cF_{\mu,L}\).  By \Cref{lem:quadboundary}
the quadratic constant is \(\min(\cQ(\theta,2),\cpi)=\sqrt2\).
(ii) The arc \(\mathcal C_{1/2}\) has \(y_{1/2}=2\) and \(2\sqrt2\le\theta\le3\).  For such \(\theta\) take the
family~\eqref{eq:corner} with \(a=0\) and \(b=2\theta\): then \(2a+b>\theta\), \(4a+b=2\theta\), \(\cpi=2\theta>\zeta_-(\theta,2)\), and
\Cref{thm:corner}(c) and \Cref{lem:quadboundary} give the pair \((\zeta_-(\theta,2),\theta/3)\).
\end{proof}

\begin{lemma}[Near the corner]\label{lem:nearcorner}
Let a family~\eqref{eq:scaling} satisfy
\((\theta_\varepsilon,y_\varepsilon,\tau_\varepsilon)\to(\theta,2,\frac12)\), \(\theta>0\),
with \(\alpha_\varepsilon,\varsigma_\varepsilon,\mathsf N_\varepsilon\) as in~\eqref{eq:alphavs},
allowing unbounded \(\alpha_\varepsilon\) and \(\varsigma_\varepsilon\).  Let
\(r=1-c\varepsilon\) with \(c\) in a compact set.  Then, uniformly for \(\xi\in[0,2]\) and for \(c\) in that set,
\[
\sfP_1+\sfP_2\xi+\sfP_3\xi^2=\varepsilon\bigl[\mathsf N_\varepsilon(\xi)-4c(3-\xi)+o(1)\bigr],\qquad
\mathsf M_r(\xi)=\varepsilon\bigl[\mathsf N_\varepsilon(\xi)-2c(4-\xi)+o(1)\bigr],
\]
with \(\sfP_k\) as in \Cref{lem:cubic} and \(\mathsf M_r\) as in \Cref{lem:mixed}.  Moreover \(\sfP_0=c\varepsilon\,\chi_\mu(r)\chi_L(r)\)
with \(\chi_\mu(r)=\varepsilon^2(c^2-\theta c+2+o(1))\) and \(\chi_L(r)=2+o(1)\), and \(\chi_L(-1)=(\varsigma_\varepsilon-2\theta_\varepsilon)\varepsilon\), \(\chi_L'(-1)\to-1\).
\end{lemma}

\begin{proof}
Apply \Cref{lem:Pscaled}(c) to the physical parameters
\((\theta_\varepsilon,y_\varepsilon,\tau_\varepsilon)\), which stay in a compact set.
The terms involving the possibly unbounded corner coordinates cancel exactly:
\[
2y_\varepsilon(2\tau_\varepsilon-1)=-8\alpha_\varepsilon\varepsilon,
\qquad
4(1-\tau_\varepsilon y_\varepsilon)
=(\varsigma_\varepsilon+4\alpha_\varepsilon)\varepsilon.
\]
The first-order coefficient polynomials converge to their values at \((\theta,2,1/2)\),
uniformly for \(c\) in the prescribed compact set.  Combining them with the two exact identities
gives both symbol expansions, with remainders uniform on \(\xi\in[0,2]\), including
unbounded \(\alpha_\varepsilon\) and \(\varsigma_\varepsilon\).
The formulas at \(r\) follow from \Cref{lem:Pscaled}(a), while
\(\chi_L(-1)=2(1+\beta)-y_\varepsilon(1+2\tau_\varepsilon)
=(\varsigma_\varepsilon-2\theta_\varepsilon)\varepsilon\) and
\(\chi_L'(-1)=-2-S_L\to-1\).
\end{proof}

\begin{proof}[Proof of \Cref{thm:optimal}]

\emph{Reformulation.}  Suppose
\((\theta_\varepsilon,y_\varepsilon,\tau_\varepsilon)
\to(3,2,1/2)\), and let
\[
c_{{\rm lag},\varepsilon}=\min_{0\le\xi\le2}
\frac{\mathsf N_\varepsilon(\xi)}{4(3-\xi)}.
\]
Define
\[
\mathsf B_\varepsilon(\xi)
=(5-4\alpha_\varepsilon)+\mathsf B'_\varepsilon\xi,
\qquad
\mathsf B'_\varepsilon=
\tfrac12(\varsigma_\varepsilon+4\alpha_\varepsilon-12).
\]
Direct expansion gives
\[
\mathsf N_\varepsilon(\xi)-4(3-\xi)
=2\xi\mathsf B_\varepsilon(\xi)
+(\theta_\varepsilon-3)(4+2\xi-4\xi^2),
\qquad \mathsf B_\varepsilon(2)=\varsigma_\varepsilon-7.
\]
Since \(4\le4(3-\xi)\le12\) and \(\xi\mathsf B_\varepsilon(\xi)=0\)
at \(\xi=0\), the condition \(c_{{\rm lag},\varepsilon}\to1\)
is equivalent to \(\min_{[0,2]}\xi\mathsf B_\varepsilon(\xi)\to0\).
It requires \(\liminf\varsigma_\varepsilon\ge7\), by evaluation at
\(\xi=2\).  Assume this inequality, and put
\(d_\varepsilon=(4\alpha_\varepsilon-5)_+\).
When \(d_\varepsilon=0\), both endpoint values of the affine
\(\mathsf B_\varepsilon\) are at least \(-o(1)\).
More generally, if \(d_\varepsilon\to0\), then
\(\mathsf B_\varepsilon(0)\ge-d_\varepsilon\) and
\(\mathsf B_\varepsilon(2)\ge-o(1)\), so
\(\min_{[0,2]}\xi\mathsf B_\varepsilon(\xi)
\ge-2d_\varepsilon-o(1)\to0\).
Along any sequence with \(d_\varepsilon\ge d_*>0\), the vertex
\[
\xi_v=\frac{d_\varepsilon}
{\varsigma_\varepsilon+5+d_\varepsilon-12}
\]
lies in \((0,1+o(1)]\), and
\[
\min_{[0,2]}\xi\mathsf B_\varepsilon(\xi)
=-\frac{d_\varepsilon^2}
{2(\varsigma_\varepsilon+5+d_\varepsilon-12)}.
\]
The quantity \(\varsigma_\varepsilon+5+d_\varepsilon-12\) is at least
\(d_*-o(1)\).  Adding \(12\) therefore changes the ratio by a bounded
factor, proving equivalence with
\(d_\varepsilon^2/(\varsigma_\varepsilon+5+d_\varepsilon)\to0\).
If this ratio fails to vanish, a subsequence has
\(d_\varepsilon\ge d_*>0\) and
\(m_\varepsilon:=d_\varepsilon^2/(2\mathsf B'_\varepsilon)\)
bounded below by a positive constant.  We may pass to further
subsequences in the cases below whenever a bounded/unbounded alternative
is used.

\emph{Sufficiency.}  Let the conditions hold and \(0<c<1\).  By \Cref{lem:nearcorner} with \(c_{{\rm lag},\varepsilon}\to1\),
\(\sfP_1+\sfP_2\xi+\sfP_3\xi^2=\varepsilon[\mathsf N_\varepsilon(\xi)-4c(3-\xi)+o(1)]\ge\varepsilon[4(3-\xi)(c_{{\rm lag},\varepsilon}-c)+o(1)]\ge2(1-c)\varepsilon\) for small
\(\varepsilon\).  Since \(c<1=\cQ(3,2)\), \Cref{lem:lagscaled} gives \(\rhoF<1-c\varepsilon\).  Hence \(\liminf\sqrt\kappa(1-\rhoF)\ge1\).  Conversely,
\Cref{thm:B}(a), applied to the family~\eqref{eq:scaling} with limit \((3,2,\frac12)\), gives \(\limsup\sqrt\kappa(1-\rhoz)\le\cFF(3,2)=1\), and \(\rhoz\le\rhoF\).

\emph{Necessity.}  Suppose \(\sqrt\kappa(1-\rhoF)\to1\).
By \Cref{cor:frontier}(v), the physical parameters tend to \((3,2,1/2)\).
First, \(\liminf\varsigma_\varepsilon\ge7\): otherwise a subsequence has
\(\varsigma_\varepsilon\to s<7\), with \(s\ge6\) by stability.
By \Cref{lem:nearcorner},
\(\chi_L(-1)=(\varsigma_\varepsilon-2\theta_\varepsilon)\varepsilon>0\) and
\(\chi_L'(-1)\to-1\).  The alternating quadratic mode of \Cref{lem:altmode}(a)
then gives \(\rhoz\ge1-(s-6+o(1))\varepsilon\), contradicting
\(\rhoz\le\rhoF\) and \(s-6<1\).

If the ratio criterion fails, pass to a sequence on which
\[
D_\varepsilon=\varsigma_\varepsilon+5+d_\varepsilon,
\qquad R_\varepsilon=\frac{d_\varepsilon^2}{D_\varepsilon}\ge q,
\qquad 0<q\le1.
\]
We separate bounded and unbounded deficits, passing to subsequences as needed.

\emph{Bounded deficit.}  If \(d_\varepsilon\le M\), stability and
\(qD_\varepsilon\le d_\varepsilon^2\) bound \(\varsigma_\varepsilon\) above
and keep \(d_\varepsilon\) away from zero.
By compactness, \(d_\varepsilon\to d_*>0\) and
\(\varsigma_\varepsilon\to s\ge7\) along a subsequence.  Put
\(a=(d_*+5)/4\) and \(b=(s-4a)/2\).
The exact identities
\[
\frac{1/2-\tau_\varepsilon}{\varepsilon}
=\frac{2\alpha_\varepsilon}{y_\varepsilon}\to a,
\qquad
\frac{2-y_\varepsilon}{\varepsilon}
=\frac{\varsigma_\varepsilon-4\alpha_\varepsilon}{2}\to b
\]
place this subsequence in \Cref{thm:corner}, since
\(2a+b=s/2\ge7/2>3\).
Part (d) requires \(a\le5/4\), contradicting \(d_*>0\).

It remains to consider \(d_\varepsilon\to\infty\).
Then \(D_\varepsilon\to\infty\), \(\mathsf B'_\varepsilon=(D_\varepsilon-12)/2>0\), and
\[
m_\varepsilon=\frac{d_\varepsilon^2}{D_\varepsilon-12}
=R_\varepsilon\frac{D_\varepsilon}{D_\varepsilon-12},
\qquad \frac{m_\varepsilon}{R_\varepsilon}\to1.
\]
Thus the bounded and divergent depth alternatives are the same for
\(m_\varepsilon\) and \(R_\varepsilon\).
Set \(\xi_{\mathsf B}=d_\varepsilon/\mathsf B'_\varepsilon\);
since \(\varsigma_\varepsilon\ge7-o(1)\), \(\xi_{\mathsf B}/2\le1+o(1)\).
For \(0<\xi\le\xi_{\mathsf B}/2\),
\(\mathsf B_\varepsilon(\xi)\le-d_\varepsilon/2\), so the reformulation gives
\begin{equation}\label{eq:dip}
\mathsf N_\varepsilon(\xi)-4c(3-\xi)
\le-d_\varepsilon\xi+12(1-c)+o(1)
\qquad(0<\xi\le\tfrac{\xi_{\mathsf B}}2,\ 0\le c\le1).
\end{equation}
The error is uniform on these ranges.

\emph{Moderate dip.}  Suppose \(q\le R_\varepsilon\le M\).
Let \(\xi_0=q/d_\varepsilon\) and \(c'=1-q/48\).
Since \(m_\varepsilon\ge R_\varepsilon\ge q\),
\(\xi_0\le\xi_{\mathsf B}/2\).
On \([\xi_0/2,\xi_0]\), \eqref{eq:dip} gives
\[
\mathsf N_\varepsilon(\xi)-4c'(3-\xi)
\le-\frac q2+\frac q4+o(1)\le-\frac q8,
\qquad |\mathsf N_\varepsilon(\xi)|\le15
\]
for small \(\varepsilon\); the second bound follows from
\(-d_\varepsilon\le\mathsf B_\varepsilon(\xi)\le0\) and
\(d_\varepsilon\xi\le q\).
By~\eqref{eq:alphavs}, \(\varepsilon D_\varepsilon\to0\).
Consequently
\(\varepsilon d_\varepsilon^2\le M\varepsilon D_\varepsilon\to0\), and hence
\(\varepsilon\xi_0^{-3/2}\to0\).

Choose \(\mathsf s=\min(\mathsf s_1,1/20)\).
The frequencies corresponding to \([\xi_0/2,\xi_0]\) have ratio tending to
\(\sqrt2\) and tend to zero.
Thus \Cref{lem:band}(a) supplies \(K\ge3\) and a band
\(B=[(1-\mathsf s)2\pi/K,(1+\mathsf s)2\pi/K]\) inside that interval of frequencies.
At \(r=1-c'\varepsilon\), \Cref{lem:nearcorner} gives, on \(B\),
\[
H_r\le-\frac q{64}\xi_0\varepsilon,
\qquad |\mathsf M_r|\le24\varepsilon,
\qquad |\chi_\mu\overline{\chi_L}|\le C\xi_0.
\]
Here \(\sfP_0=O(\varepsilon^3)=o(\varepsilon\xi_0)\); the last estimate follows from
\(\chi_\mu=(z-1)(z-\beta)+y_\varepsilon\varepsilon^2N_{\tau_\varepsilon}(z)
=O(\xi_0+\varepsilon\sqrt{\xi_0}+\varepsilon^2)=O(\xi_0)\)
and uniform boundedness of \(\chi_L\).
Since \(\xi_{\max}\in[\xi_0/2,\xi_0]\) and
\(\varepsilon\xi_0^{-3/2}\to0\), \Cref{lem:spreadwit} applies with
\(\lambda_0=q/(64C)\), giving \(\rhoz\ge1-c'\varepsilon\), a contradiction.

\emph{Deep dip.}  Otherwise pass to a subsequence with \(R_\varepsilon\to\infty\).
Take integers \(K_\varepsilon\sim2\pi\sqrt{d_\varepsilon/24}\) and put
\(\xi_\varepsilon=1-\cos(2\pi/K_\varepsilon)\).
Then \(K_\varepsilon\ge3\), \(\xi_\varepsilon\to0\) and
\(d_\varepsilon\xi_\varepsilon\to12\).  The exact reformulation yields
\[
\mathsf N_\varepsilon(\xi_\varepsilon)
=12-4\xi_\varepsilon-2d_\varepsilon\xi_\varepsilon
+\frac{(d_\varepsilon\xi_\varepsilon)^2}{R_\varepsilon}
-12\xi_\varepsilon^2+o(1)\to-12.
\]
By \Cref{lem:nearcorner} at \(c=0\),
\(L_\tau(\xi_\varepsilon)=\varepsilon[\mathsf N_\varepsilon(\xi_\varepsilon)+o(1)]<0\).
\Cref{lem:band}(b) gives \(K_\varepsilon\)-periodic witnesses and \(\rhoz>1\), again a contradiction.
Thus the ratio criterion is necessary.
\end{proof}

\section{Proofs for Section~\ref{sec:finite}}
\label{app:finite}

\begin{lemma}[Structure of \(L_\tau\)]\label{lem:struct}
\(L_\tau\) is affine in \(y_\mu\), and at \(y_\mu=y\) it equals
\(|\chi_L(e^{i\omega})|^2\).  Hence, for \(y>0\) and \(0\le y_\mu\le y\),
\(L_\tau(\xi;y_\mu)=(1-\frac{y_\mu}y)L_\tau(\xi;0)+\frac{y_\mu}y|\chi_L(e^{i\omega})|^2\), where
\(L_\tau(\xi;0)=2(1-\beta)y+\xi[2(1-\beta)^2+y(2\tau(3-\beta)-2)]+4\xi^2(\beta-\tau y)\).
With \(y_{\max}=2(1+\beta)/(1+2\tau)\),
\[
L_\tau(\xi;0)\big|_{y=0}=2\xi\bigl[(1-\beta)^2+2\beta\xi\bigr],\qquad
L_\tau(\xi;0)\big|_{y=y_{\max}}=\frac{2(2-\xi)}{1+2\tau}\bigl[(1-\beta^2)+2\xi(2\tau-\beta)\bigr].
\]
\end{lemma}

\begin{proof}
The coefficients in~\eqref{eq:Acoef} are affine in \(y_\mu\), and
\(\chi_\mu=\chi_L\) at \(y_\mu=y\), proving the interpolation
formula.  Setting \(y_\mu=0\) gives the displayed polynomial.
At \(y=0\) its value is immediate; substitution of
\(y_{\max}=2(1+\beta)/(1+2\tau)\) factors it as
\(2(2-\xi)[(1-\beta^2)+2\xi(2\tau-\beta)]/(1+2\tau)\).
\end{proof}

\begin{proof}[Proof of \Cref{thm:tonguefree}]
Use the multiplier weight \(\lambda\) specified in \Cref{thm:tonguefree}, and write
\(B(\xi)=(1-\beta^2)+2\xi(2\tau-\beta)\), so that
\(L_\tau(\xi;0)|_{y=y_{\max}}=\frac{2(2-\xi)}{1+2\tau}B(\xi)\) (\Cref{lem:struct}); note
\(B(2)=8(\tau-\tau_*)\) and \(B(\frac32)=6(\tau-\tau_3)\).

(a), \(\tau\ge\tau_*\).  Fix \(\xi\in(0,2]\).  \(L_\tau(\xi;0)\) is affine in \(y\), positive
at \(y=0\), and nonnegative at \(y=y_{\max}\), because \(B\) is affine with
\(B(0)=1-\beta^2>0\) and \(B(2)\ge0\).  So \(L_\tau(\xi;0)>0\) for \(0\le y<y_{\max}\), which by
\Cref{lem:qstab} is the stability range.  For a stable tuning \(\chi_L\) has no
root on the unit circle, and \Cref{lem:struct} gives \(L_\tau(\xi;y_\mu)>0\) for
\(0<y_\mu<y\); \Cref{prop:Ltau} gives \(\rhoF<1\).

(a), \(\tau_3\le\tau<\tau_*\).  Then \(\beta-2\tau>0\), because \(\tau_*\le\beta/2\) (equivalent to \(\beta^2\le1\)), so
\(\lambda>0\).  On \(|z|=1\) write \(\chi_\mu\overline{\chi_L}=\mathcal R_0+iq\sin\omega\) with
\(q(\xi)=(y-y_\mu)[(1-\beta)+2\tau\xi]\) by~\eqref{eq:imab}, and recall \(\mathcal R_0=L_\tau-(2-\xi)q\) from
the proof of \Cref{prop:Ltau}.  Since \(1-\bar z^2=2\sin^2\omega+2i\sin\omega\cos\omega\),
\[
\Rea\bigl[\chi_\mu\overline{\chi_L}(1-\bar z^2)\bigr]=2\sin^2\omega\,(\mathcal R_0-q\cos\omega)
=2\xi(2-\xi)\bigl(L_\tau-(3-2\xi)q\bigr),
\]
and together with \(\Rea[\chi_\mu\overline{\chi_L}(1-z)]=\xi L_\tau\),
\[
\Rea\bigl[\chi_\mu\overline{\chi_L}\,((1-z)+\lambda(1-\bar z^2))\bigr]=\xi F(\xi),\qquad
F=(1+2\lambda(2-\xi))L_\tau+2\lambda(2-\xi)(2\xi-3)q .
\]
We show \(F>0\) on \((0,2]\) for every stable tuning.  For fixed \(\xi\), \(F\) is affine in
\(y_\mu\in[0,y]\) (both \(L_\tau\) and \(q\) are), and at \(y_\mu=0\) it is affine in
\(y\in[0,y_{\max}]\).  At \(y_\mu=y\), \(q=0\) and \(F=(1+2\lambda(2-\xi))|\chi_L|^2>0\), since \(\chi_L\) has no
root on the unit circle.  At \(y_\mu=y=0\), \(F=(1+2\lambda(2-\xi))\,2\xi[(1-\beta)^2+2\beta\xi]>0\).  At
\(y_\mu=0\), \(y=y_{\max}\), using \Cref{lem:struct} and \(y_{\max}(1+2\tau)/2=1+\beta\),
\(F=\frac{2(2-\xi)}{1+2\tau}G(\xi)\) with
\[
G(\xi)=(1+2\lambda(2-\xi))B(\xi)+2\lambda(2\xi-3)(1+\beta)\bigl[(1-\beta)+2\tau\xi\bigr].
\]
Substitution of the chosen \(\lambda\) and completion of the square give the exact identity
\[
G(\xi)=6(1+\lambda)(\tau-\tau_3)
+4\lambda\beta(1+2\tau)\bigl(\xi-\tfrac32\bigr)^2\ \ge0 .
\]
So \(F(\xi;0,y)\) is affine in \(y\in[0,y_{\max}]\), positive at \(y=0\) and nonnegative at
\(y=y_{\max}\), hence positive for \(y<y_{\max}\); and \(F(\xi;y_\mu,y)\), affine in \(y_\mu\) between
the positive values at \(y_\mu=0\) and \(y_\mu=y\), is positive.  Thus
\(\Rea[\chi_\mu\overline{\chi_L}((1-z)+\lambda(1-\bar z^2))]>0\) for \(z\ne1\), and \Cref{prop:lagone}(c) gives
\(\rhoF<1\).

(b) If \(\tau<\tau_3\), then \(L_\tau(\frac32;0)|_{y=y_{\max}}=B(\frac32)/(1+2\tau)<0\); by
continuity in \(y\) and \Cref{lem:struct}, \(L_\tau(\frac32;y_\mu)<0\) for some \(y<y_{\max}\)
and small \(y_\mu>0\).  Take \(h=1\), \(L=y\), \(\mu=y_\mu\); the tuning is stable, and
\Cref{lem:band}(b) with \(K=3\) (\(\xi=\frac32\)) gives \(\Theta(2\pi/3)<\pi/3\), the \(3\)-periodic orbits and \(\rhoz>1\).

(c) If \(\tau<\tau_*\), \(B<0\) on some \((\xi_1,2]\), so \(L_\tau(\xi;0)<0\) at \(y=y_{\max}\)
for \(\xi\in(\xi_1,2)\), and the same perturbation argument applies.

The special cases follow from \(\tau_*-\tau_3=(1-\beta^2)/24\), \(\tau_*(\beta)\le\beta\iff\beta^2-4\beta-1\le0\),
\(\tau_3(\beta)\le0\iff\beta\le(\sqrt{13}-3)/2\) and \(\tau_*(\beta)\le0\iff\beta\le\sqrt5-2\).
\end{proof}

\begin{proof}[Proof of \Cref{cor:QHM}]
Use the query-rate dictionary from the proof of \Cref{cor:B}(iii).
For \(\beta\ne z_\tau\), \Cref{lem:restquad} identifies stability with convergence of the iterates
on every quadratic in \(\cF_{\mu,L}\).
For fixed \(\beta,\nu_{\rm qh}\), the extrapolation \(\tau_{\rm qh}\) is fixed and varying \(\alpha\)
covers all \(h>0\).  Apply \Cref{thm:tonguefree}: linear convergence of \(\cM\) gives
\(\rho^0_{\tilde x}\le\rhoz\le\rhoF<1\), and the divergent orbits from rest of \Cref{thm:tonguefree} (a \(3\)-cycle
tongue) have \(\rhoz>1\), hence \(\rho^0_{\tilde x}>1\) by the second inequality.
\end{proof}

\begin{proof}[Proof of \Cref{prop:mono}]
For any fixed method, inclusion of function classes gives
\[
\mu\le\mu'<L
\quad\Longrightarrow\quad
\Gamma_{\mu',L}(n)\le\Gamma_{\mu,L}(n)
\quad\Longrightarrow\quad
\rho_{\cF_{\mu',L}}\le\rho_{\cF_{\mu,L}}.
\]
Thus \(I=\{\mu\in(0,L):\rhoF<1\}\) is an upper interval.
It is open: \Cref{thm:A}(b) supplies a strict certificate at \(t=0\)
for each \(\mu\in I\), and \Cref{lem:persist} preserves it for nearby
\(\mu\).  Stability, independent of \(\mu\) by \Cref{lem:qstab},
gives \(r_c<1\), so \Cref{thm:A}(c) applies there.

By \Cref{lem:struct}, \(L_\tau(\xi;y_\mu)\) is affine in \(y_\mu\)
and is strictly positive at \(y_\mu=y\), where it equals
\(|\chi_L(e^{i\omega})|^2\).  Therefore
\(\{y_\mu:\min_{[0,2]}L_\tau(\xi;y_\mu)>0\}\) is a nonempty open
upper interval with right endpoint \(y\).  The lag-one criterion
places it inside \(I\).  Since \(y_\mu=h\mu\), this proves both
interval descriptions.  The class-inclusion argument also shows that
the rate itself is nonincreasing in \(\mu\), for any fixed stationary
method.
\end{proof}
\end{document}